\documentclass[authoryear,review,12pt]{elsarticle}
\usepackage{appendix}
\usepackage{comment}
\usepackage{setspace}
\usepackage[margin=0cm]{caption}
\usepackage[font={sf,bf},textfont=md]{caption}
\usepackage[tbtags]{amsmath}
\usepackage{amsthm,amssymb,amsfonts}
\usepackage{natbib}
\usepackage{multirow}
\usepackage{lscape}
\usepackage{graphicx}
\usepackage{subcaption}
\usepackage{makecell}
\usepackage{exscale}
\usepackage{booktabs}
\usepackage{array}
\usepackage{fullpage}
\usepackage{url}
\usepackage{algorithm}
\usepackage{algpseudocode}
\usepackage{bm}
\usepackage{mathtools}
\usepackage{wrapfig}
\usepackage{lipsum}
\usepackage{mathrsfs}
\usepackage{relsize}
\usepackage{dsfont}
\usepackage{multirow}
\usepackage{apalike}
\usepackage[usenames,dvipsnames,svgnames,table]{xcolor}
\usepackage[colorlinks=true,
linkcolor=blue,
urlcolor=blue,
citecolor=blue]{hyperref}

\title{Testing Kronecker Product Covariance Matrices for High-dimensional Matrix-Variate Data}
\date{}

\newtheorem{theorem}{Theorem}
\newtheorem{lemma}{Lemma}
\newtheorem{corollary}{Corollary}
\newtheorem{condition}{Condition}

\newtheorem{remark}{Remark}

\newcommand{\ba}{\bm{a}}

\newcommand{\be}{\bm{e}}

\newcommand{\bh}{\bm{h}}

\newcommand{\bq}{\bm{q}}
\newcommand{\br}{\bm{r}}

\newcommand{\bw}{\bm{w}}
\newcommand{\bx}{\bm{x}}
\newcommand{\by}{\bm{y}}

\newcommand{\Ab}{\mathbf{A}}
\newcommand{\Bb}{\mathbf{B}}
\newcommand{\Cb}{\mathbf{C}}
\newcommand{\Db}{\mathbf{D}}

\newcommand{\Gb}{\mathbf{G}}
\newcommand{\Hb}{\mathbf{H}}
\newcommand{\Ib}{\mathbf{I}}

\newcommand{\Mb}{\mathbf{M}}

\newcommand{\Qb}{\mathbf{Q}}

\newcommand{\Sbb}{\mathbf{S}}

\newcommand{\Ub}{\mathbf{U}}
\newcommand{\Vb}{\mathbf{V}}
\newcommand{\Wb}{\mathbf{W}}
\newcommand{\Xb}{\mathbf{X}}
\newcommand{\Yb}{\mathbf{Y}}
\newcommand{\Zb}{\mathbf{Z}}

\newcommand{\bB}{\bm{B}}

\newcommand{\bU}{\bm{U}}
\newcommand{\bV}{\bm{V}}

\newcommand{\bX}{\bm{X}}
\newcommand{\bY}{\bm{Y}}

\newcommand{\bgamma}{\bm{\gamma}}

\newcommand{\bxi}{\bm{\xi}}

\newcommand{\bGamma}{\bm{\Gamma}}

\newcommand{\bTheta}{\bm{\Theta}}

\newcommand{\bSigma}{\bm{\Sigma}}

\newcommand{\bPhi}{\bm{\Phi}}
\newcommand{\bPsi}{\bm{\Psi}}
\newcommand{\bOmega}{\bm{\Omega}}

\begin{document}
	\begin{frontmatter}

		\author[myfirstaddress]{Yu Long}
		\ead{fduyulong@163.com}
				\author[myfirstaddress]{Xie Jiahui}
		\ead{jiahui.xie@u.nus.edu}
		\author[myfirstaddress]{Zhou Wang}
		\ead{wangzhou@nus.edu.sg}
		\address[myfirstaddress]{Department of Statistics and Data Science, National University of Singapore,
			Singapore.}

		\begin{abstract}
			Kronecker product covariance structure provides an efficient way to modeling the inter-correlations of matrix-variate data. In this paper, we propose testing statistics for Kronecker product covariance matrix based on linear spectral statistics of renormalized sample covariance matrices. Central limit theorem is proved for the linear spectral statistics with explicit formulas for mean and covariance functions, which fills the gap in the literature. We then theoretically justify that  the proposed testing statistics have well-controlled sizes and strong powers. To facilitate practical usefulness, we further propose a bootstrap resampling algorithm to approximate the limiting distributions of associated linear spectral statistics. Consistency of the bootstrap procedure is guaranteed under mild conditions. A more general model which allows the existence of noises will also be discussed. In the simulations,   the empirical sizes of the proposed testing procedure and its bootstrapped version are close to corresponding theoretical values, while the powers converge to one quickly as the dimension and sample size grow. 
		\end{abstract}
		
		
	\end{frontmatter}

\section{Introduction}
Estimation of covariance matrix is one of the most fundamental problems in statistical learning and related applications. Recent decades have seen fruitful research in this field, especially under the high-dimensional settings where the conventional sample covariance matrix is no longer a consistent estimator. To name a few examples, see \cite{bickel2008regularized,cai2010optimal,fan2013large} and the references therein. However, most of these estimation procedures are designed only for data of vector form. Thanks to the rapid advance in data science and information technology, there is growing demand for the analysis of matrix-variate or higher-order tensor-variate data.  In this paper, we will mainly focus on matrix-variate data.

In matrix-variate regime, a simple approach is to first vectorize the data matrices and then apply the conventional vector-based procedures. However, naively stacking the rows or columns of matrix usually leads to a great loss of the information contained in the matrix structure, and also a higher risk of  the ``curse of dimensionality".  Instead, to characterize the inter-connections of  matrix-variate data, the Kronecker product covariance matrix structure attracts more and more attention nowadays.  We say the $p\times q$ data matrix $Y$ has Kronecker product covariance matrix if
\begin{equation}\label{kronecker}
	\text{cov}\{\text{Vec}(Y)\}=\Sigma_{Y}=\Sigma_{V}\otimes \Sigma_{U},
\end{equation}
where $\text{Vec}(Y)$ stands for stacking the columns of $Y$ into a $pq$-dimensional vector, $\Sigma_{U}$ and $\Sigma_{V}$ are $p\times p$ and $q\times q$ {cross-row and cross-column covariance matrices}, respectively.  It's also common in the literature to directly assume
\begin{equation}\label{separable}
	Y=UXV^\prime,\quad UU^\prime=\Sigma_{U},\quad VV^\prime=\Sigma_{V},
\end{equation}
where $U$ and $V$ are $p\times p$ and $q\times q$ deterministic matrices, and $X$ is $p\times q$ random matrix composed of independent and identically distributed  entries. See \cite{zhou2014gemini} and \cite{leng2018covariance}.  If the entries of $X$ are standard normal variables, we say $Y$  follows matrix-variate normal distribution.  Sometimes, the decomposition in  (\ref{separable}) is also referred as separable covariance model. It's easy to see that  the data generating model (\ref{separable}) satisfies the Kronecker product covariance matrix structure in (\ref{kronecker}). The Kronecker product assumption retains the matrix structure of the data, meanwhile effectively reducing the number of unknown parameters in the covariance matrix from $(pq)(pq+1)/2$ to $p(p+1)/2+q(q+1)/2$.

Estimating covariance matrix with Kronecker product assumption has been considered in the literature both in low-dimensional and high-dimensional settings. For example, the flip-flop algorithm in \cite{lu2004likelihood} and its extensions in \cite{srivastava2008models} and \cite{werner2008estimation} are suitable to the low-dimensional cases. On the other hand, estimation in high dimensions usually relies on sparsity assumption on the population covariance or correlation matrix  and penalized optimization algorithms, see \cite{tsiligkaridis2012sparse,tsiligkaridis2013covariance,leng2018covariance}. Some other papers focus on the estimation of the precision matrix $\Sigma_{Y}^{-1}$ with sparsity assumptions on $\Sigma_{U}^{-1}$ and $\Sigma_{V}^{-1}$, see \cite{allen2010transposable}, \cite{leng2012sparse}, and \cite{zhou2014gemini} to list a few. With all these estimation procedures, it's natural to ask which one is more preferred in real applications, especially when they give significantly disparate results. This paper is partially motivated by the portfolio data example in our supplementary material, where various approaches output different {guesses} for the covariance matrix which sometimes leads to contradictory investing strategies. In other words, it's in urgent need to propose some testing procedure which can evaluate { preliminary hypothesis of the covariance matrix}.

Unfortunately, the testing of Kronecker product covariance matrix is more challenging  and  only  few studies have ever tried  to propose some testing procedures. Up to our knowledge, most of them are only for the low-dimensional settings and require normality assumptions so that the likelihood ratio test may work. See \cite{lu2005likelihood}, \cite{srivastava2008models} and \cite{hao2016testing}.  These constraints are not  easily fulfilled in real applications, while the existence and uniqueness of the maximum likelihood estimation should also be considered. See \cite{ros2016existence}.

Motivated by the above arguments,  we aim to propose some distribution-free testing procedure for the Kronecker product covariance matrix in high dimensions, which is the first contribution of this paper. Our method is based on the  column-column or row-row sample covariance matrix, defined by 
$S=(Tq)^{-1}\sum_{t=1}^TY_tY_t^\prime$ where $Y_t$ {$(1\le t\le T)$ are independent observations of $Y$,  with $T$ being the number of observations}.   We renormalize  $S$ and construct testing statistics based on associated  linear spectral statistics. Detailed definitions and procedures are presented in the next section, which can be viewed as an extension of \cite{bai2004clt} to the high-dimensional matrix-variate regime. The method is flexible and general since many conventional  testing statistics can be written as special cases of linear spectral statistics. {A more general model which allows the existence of noises in (\ref{separable}) will also be discussed in the paper.}

Our second contribution is to derive the central limit theorem  for the linear spectral statistics mentioned above.  We observe that  under the data generating model (\ref{separable})$, S$ is close to the separable covariance model ever studied in \cite{bai2019central} and \cite{li2021central}. However, they require that the ratio of dimension over sample size  converges to a constant, which is usually not fulfilled in matrix-variate regime. The definition of $S$ is equivalent to regarding each column of $Y_t$ as individual observation, which increases the effective sample size. Hence, the ratio of dimension {$(p)$} over sample size {$(Tq)$} typically converges to zero.  The derivation of  central limit theorem under this case is more challenging  than that in \cite{bai2019central}, see \cite{bai1988convergence} and \cite{chen2015clt} for intuition. We formally prove the central limit theorem, with explicit formulas for the mean and covariance functions, which fills the gap in the literature of random matrix theory.

Although the theoretical mean and covariance functions of the central limit theorem are presented in our theorems, they involve complex number integration which is hard to calculate and  unknown parameters which need to be estimated. Then, the third contribution of this paper is to propose a bootstrap algorithm to approximate the limiting distribution of the linear spectral statistics. We borrow idea from \cite{lopes2019bootstrapping} and design a bootstrap resampling procedure which can output accurate critical values for the testing statistics. This facilitates to the practical usefulness of  the proposed testing procedure.

The rest of this paper is organized as follows. Section \ref{sec2} illustrates the motivation of our testing statistics and the detailed procedure. The testing statistics can be written as linear spectral statistics of renormalized sample covariance matrices. Hence, in Section \ref{sec3}, we show the theoretical results on the central limit theorem of corresponding linear spectral statistics. With the central limit theorem, we discuss the asymptotic sizes and powers of the testing procedure in Section \ref{sec4}. Section \ref{sec5} provides a bootstrap algorithm for the approximation of associated limiting distributions. {Section \ref{sec:noise} discusses the more general model which allows the existence of noises.} Section \ref{sec6}  verifies the empirical sizes and powers of the proposed testing procedure and its bootstrapped version with simulated data under different settings. Section \ref{sec8} discusses some extensions. Some additional simulation results, a real data example and  all the technical proofs of the theorems, lemmas and corollaries are put into our supplementary material.

\section{Testing procedure}\label{sec2}
Let $\{Y_t\}_{t=1}^T$ be  independent observations of $p\times q$ random matrix $Y$ which satisfies the separable structure (\ref{separable}). That is,
\begin{equation}\label{model}
	Y_t=UX_tV^\prime\quad (t=1,\ldots,T),
\end{equation}
where the entries of $X_t$ are independent and identically distributed with mean 0 and variance 1. We are interested in hypothesis testing of $\Sigma_{Y}=\Sigma_{Y,0}$, where $\Sigma_{Y,0}$ is a preliminary guess of the population covariance matrix. Under the Kronecker product assumption, {it's more informative to test $\Sigma_{U}=\Sigma_{U_0}$  and $\Sigma_{V}=\Sigma_{V_0}$ separately, given some matrices $\Sigma_{U_0}$ and $\Sigma_{V_0}$}.  Considering the exchangeability  of $U$ and $V$, we  mainly  focus on the testing of $\Sigma_{U}$ in this paper, while that of $\Sigma_{V}$ follows a parallel procedure by transposing $Y_t$.

Let's consider the trivial case $q=1$ first,  which reduces to the testing of covariance matrix of conventional vector-valued data.  Many approaches are available in the literature no matter in low dimensions or high dimensions, such as the testing statistics in \cite{john1971some}, \cite{ledoit2002some} and  \cite{chen2010tests}.  It's well known that many testing statistics are spectral statistics of the sample covariance matrix,  defined by
	$S=(Tq)^{-1}\sum_{t=1}^TY_tY_t^\prime$.
To this end, for any Hermitian matrix $S$ of size $p\times p$, its empirical spectral distribution is defined by
$F^{S}(x)=p^{-1}\sum_{j=1}^pI(\lambda_j^{S}\le x)$,
where $\lambda_j^{S}$ is the $j$-th largest eigenvalue of $S$. Moreover,  linear spectral statistics corresponding to $S$ are quantities of the form
$	p^{-1}\sum_{j=1}^pf(\lambda_j^{S})=\int f(x)dF^{S}(x)$,
with some  continuous and bounded real function $f$ on $(-\infty,\infty)$.

In matrix-variate regime where $q>1$,  $S$ is referred {to} as the column-column sample covariance matrix, since it regards each column of $Y_t$ as individual observation. Obviously, the expectation of $S$ is $q^{-1}\text{tr}(\Sigma_{V})\Sigma_{U}$. Due to identifiability, we may let  $q^{-1}\text{tr}(\Sigma_{V})=1$, then $S$ is unbiased estimate of $\Sigma_{U}$.  This motivates us to test $\Sigma_{U}=\Sigma_{U_0}$ still based on $S$ even for the matrix-variate  data.  {Other identification condition such as that in \cite{srivastava2008models} is also available, with minor adjustment on our calculations.} A closer look at $S$ shows that
\begin{equation}\label{S:separable}
	S=(Tq)^{-1}U(X_1,\ldots,X_T)(I_T\otimes V^\prime)(I_T\otimes V^\prime)^\prime(X_1,\ldots,X_T)^\prime U^\prime.
\end{equation}
The representation (\ref{S:separable}) is similar to the separable sample covariance matrix in \cite{bai2019central}, where the central limit theorem of associated linear spectral statistics is proved. However, in \cite{bai2019central}, {the ratio of  dimension $(p)$ over sample size $(Tq)$ } is required to converge to some constant $\gamma\in(0,\infty)$. {For  high-dimensional matrix-variate  data considered in the current paper, the ratio  usually tends to zero.}  It's well known in random matrix theory that the spectral properties of sample covariance matrices are totally different under the two cases.

When $p/(Tq)\rightarrow 0$, motivated by \cite{bai1988convergence} and \cite{chen2015clt}, it's more convenient to first normalize $S$ by defining
\begin{equation}\label{renormalize}
	\bar S=\{(Tq)/p\}^{1/2}\{S-E(S)\},
\end{equation}
where $E(\cdot)$ denotes the expectation. In this paper, we propose to test the null hypothesis $\Sigma_{U}=\Sigma_{U_0}$ using linear spectral statistics of the renormalized $\bar S$.
The limiting distribution of linear spectral statistics associated with $\bar S$ has ever been studied by \cite{chen2015clt}. However, they only consider the special case where $U=V=I$, which does not cover the matrix-variate scenarios considered in this paper. Up to our knowledge, no results are found in the literature for general  $U$ and $V$. To overcome this difficulty, we derive the central limit theorem of linear spectral statistics associated with $\bar S$ in the next section.

Under high dimensional settings, the estimated covariance matrices from penalized optimization procedures are usually positive definite, see \cite{leng2018covariance}.
Formally, in this paper we are testing
\[
	H_0: \Sigma_{U}=\Sigma_{U_0},\text{ v.s. } H_1: \Sigma_{U}\ne\Sigma_{U_0}, \text{ for some } \Sigma_{U_0}>0.
\]
Given $\Sigma_{U_0}>0$, we can rewrite (\ref{model}) as
\begin{equation}\label{transform}
	\tilde Y_t=\Sigma_{U_0}^{-1/2}Y_t=\Sigma_{U_0}^{-1/2}UX_tV^\prime=\tilde UX_tV^\prime \quad (t=1,\ldots, T).
\end{equation}
Then, the null hypothesis $\Sigma_{U}=\Sigma_{U_0}$ is equivalent to $\Sigma_{\tilde U}=\Sigma_{U_0}^{-1/2}\Sigma_{U}\Sigma_{U_0}^{-1/2}=I$. Moreover, the null hypothesis $\Sigma_{Y}=\Sigma_{V_0}\otimes \Sigma_{U_0}$ is equivalent to $\Sigma_{\tilde Y}=\Sigma_{V_0}\otimes I_p$. Consequently, after the transformation (\ref{transform}), it's sufficient to consider testing the identity of $\Sigma_{\tilde U}$. We summarize the above arguments and our testing procedure  in Algorithm \ref{alg1}.

\begin{algorithm}[H]
	\caption{Testing procedure for $H_0: \Sigma_{U}=\Sigma_{U_0}$}\label{alg1}
	\begin{algorithmic}[1]
		\Require data matrices $\{Y_t\}_{t=1}^T$,  function $f$.
		\Ensure P-value of the testing.
		
		\State  transform $\{Y_t\}_{t=1}^T$ to $\{\tilde Y_t\}_{t=1}^T$ based on (\ref{transform}),
		\State calculate the renormalized column-column sample covariance matrix associated with $\{\tilde Y_t\}_{t=1}^n$ by
		$	\tilde S=\{(Tq)/p\}^{1/2}\big\{(Tq)^{-1}\sum_{t=1}^T\tilde Y_t\tilde Y_t^\prime-I\big\}$,
		
		\State calculate the testing statistic
		$
		\mathcal{T}=\{\sum_{j=1}^pf(\lambda_j^{\tilde S})-\mu\}/\sigma$,
		where $\mu$ and $\sigma^2$ are given in (\ref{mu}) and (\ref{sigma2}) in Section 4,
		\State output the P-value $2\{1-F^{G}(|\mathcal{T}|)\}$, where $F^{G}(x)$ is the cumulative probability function of standard normal distribution.
	\end{algorithmic}
\end{algorithm}

\section{Linear spectral statistics}\label{sec3}
\subsection{Preliminary results}
If $U=\Sigma_{U}^{1/2}$, the renormalized sample covariance matrix $\bar S$ in (\ref{renormalize}) can be viewed as a special case of the more general form
\begin{equation}\label{sp}
	\bar S_p=(np)^{-1/2}\Big\{A_p^{1/2}X_pB_nX_p^\prime A_p^{1/2}-(\text{tr}B_n)A_p\Big\},
\end{equation}
{where $A_p$ and $B_n$ are respective $p\times p$ and $n\times n$ deterministic matrices.  $X_p$ is $p\times n$ random matrix with independent and identically distributed entries. We use  $p$ and $n$ to indicate the respective dimension and sample size}. Letting $A_p=\Sigma_{U}$, $n=Tq$,  $B_n=I_T\otimes \Sigma_{V}$, we then get (\ref{renormalize}). This section proves central limit theorem for linear spectral statistics of $\bar S_p$.
In the below, we propose some assumptions directly on $A_p$, $X_p$ and $B_n$.  We will come back to the matrix-variate model (\ref{model}) in the next section.
\begin{condition}\label{c1}
	In (\ref{sp}), suppose that
	\begin{enumerate}
		\item $X_p=(x_{ij})_{p\times n}$ where $\{x_{ij}: i=1,\ldots,p,j=1\ldots,n\}$ are  independent and identically distributed real random variables with $E(x_{11})=0$, $E(x_{11}^2)=1$, $E(x_{11}^4)=\nu_4$, and $E(|x_{11}|^{4+\delta_0})<\infty$ for some $\delta_0>0$.
		\item $p/n\rightarrow 0$ as $p\rightarrow \infty$.
		\item $A_p$ and $B_n$ are non-negative deterministic real symmetric matrices with bounded eigenvalues $(a_1,\ldots,a_p)$ and $(b_1,\ldots,b_n)$, respectively in decreasing order. The empirical spectral densities of $A_p$ and $B_n$ converge to some  probability functions $F^{A}$ and $F^{B}$ which are not degenerate at 0 as $p\rightarrow\infty$, respectively. 
	\end{enumerate}
\end{condition}

The above conditions are standard and common in random matrix theory, see \cite{bai2019central}. The condition of finite $(4+\delta_0)$-th moment is to derive almost surely upper and lower bounds for the eigenvalues of $\bar S_p$, which is not stringent in real applications. The condition $p/n\rightarrow 0$ is from the matrix-variate setting where $n=Tq$. Hence, the fundamental large number in this paper is $p$ rather than $n$.
The central limit theorem for linear spectral statistics of $\bar S_p$ under the above general conditions has not been studied in the literature. Hence, this section fills in this gap, and also provides the mean and variance parameters in Step 3 of  Algorithm \ref{alg1}.

 {Define the Stieltjes transform of any distribution function $F(x)$ as
$m_{F}(z)=\int (x-z)^{-1}dF(x)$, $z\in \mathbb{C}^+$.
Then, for  the empirical spectral distribution $F^{\bar S_p}(x)$, its Stieltjes transform can be written as
$
m_{F^{\bar S_p}}(z)=p^{-1}\text{tr}(\bar S_p-zI)^{-1}=p^{-1}\sum_{j=1}^p(\lambda_j^{\bar S_p}-z)^{-1}$, $z\in \mathbb{C}^+$.
Given $z\in\mathbb{C}^+$, define $m_p(z)$ and $s_p(z)$ as the solution in $\mathbb{C}^+$ to the equations
\[
m_p(z)=-\int\frac{1}{z+x\bar\lambda_{B_n^2}s_{p}(z)}dF^{A_p}(x),\quad s_p(z)=-\int\frac{x}{z+x\bar\lambda_{B_n^2}s_{p}(z)}dF^{A_p}(x),
\]
where $\bar\lambda_{B_n^2}=n^{-1}\text{tr}(B_n^2)$. Indeed, $m_p(z)$ is the Stieltjes transform of some probability function $F_p(x)$, which works as an approximation to $F^{\bar S_p}(x)$.  Let
$\tilde G_p(f)=p\int_{-\infty}^{+\infty}f(x)d(F^{\bar S_p}(x)-F_p(x))$,
where $f\in\mathcal{M}=\{\text{functions which are  analytic in an open domain containing} [-2c,2c]\}$ with $c=\lim\sup_{p}a_1b_1$ being an almost surely upper bound for the spectral norm of $\bar S_p$.  Then, $\tilde G_p(f)$ can be viewed as normalized linear spectral statistics. 

To provide the asymptotic mean of $\tilde G_p(f)$, we define

	\[
	\mathcal{X}_p(z)=-\frac{1}{p}\mathcal{A}_p(z)\times \frac{1}{p}\sum_{k=1}^pa_k^2\tilde\epsilon_k(z)^3-\mathcal{Y}_p(z)\times \frac{1}{p}\sum_{k=1}^p\frac{a_k\tilde\epsilon_k(z)^2}{1-\mathcal{Y}_p(z)a_k\tilde\epsilon_k(z)},
	\]
	where
	\[
	\begin{split}
		\mathcal{A}_p(z)=&(\nu_4-3)n^{-1}\sum_{j=1}^n\mathrm{B}_{n,jj}^2+\bar\lambda_{B_n^2}+\mathcal{B}_p(z),\quad 		\tilde\epsilon_k(z)=\big\{z+a_k\bar\lambda_{B_n^2}s_p(z)\big\}^{-1},\\
		\mathcal{B}_p(z)=&\frac{1}{p}\sum_{k=1}^pa_k^2\bar\lambda^2_{B_n^2}\tilde\epsilon_k(z)^2\bigg\{1-\frac{1}{p}\sum_{k=1}^pa_k^2\bar\lambda_{B_n^2}\tilde\epsilon_k(z)^2\bigg\}^{-1},
	\end{split}
	\]
	and $\mathcal{Y}_p(z)$  is the solution satisfying $\mathcal{Y}_p(z)=o(1)$ as $p\rightarrow\infty$ to the equation
	\[
	\begin{split}
		x=\bar\lambda_{B_n^2}\bigg\{-\frac{1}{p}\mathcal{A}_p(z)\times\frac{1}{p}\sum_{k=1}^pa_k^3\tilde\epsilon_k^3(z)+s_p(z)\bigg\}+\frac{1}{n}\sum_{j=1}^n\frac{(\lambda_j^{B_n})^2\mathcal{D}_{p2}(x,z)}{1-\lambda_j^{B_n}(\frac{p}{n})^{1/2}\mathcal{D}_{p2}(x,z)}
	\end{split}
	\]
	with $
	\mathcal{D}_{p2}(x,z)=p^{-1}\sum_{k=1}^p(a_k\tilde\epsilon_k)/(1-xa_k\tilde\epsilon_k)$.
	Further let
	\begin{equation}\label{gpf}
	G_{p}(f)=\tilde G_p(f)+\frac{p}{2\pi i}\oint_{\mathcal{C}}f(z)\mathcal{X}_p(z)dz,
	\end{equation}
	where $\mathcal{C}$ is the contour formed by the boundary of the rectangle with four vertices $(\pm u_0,\pm iv_0)$. Here $u_0=2\lim\sup_{p}a_1b_1+\epsilon_0$ with sufficiently small $\epsilon_0$, and $v_0$ is any positive number so that $f$ is analytic in a neighborhood of $\mathcal{C}$.
Then, the next theorem shows that $G_p(f)$ converges weakly to a Gaussian limit.
	
\begin{theorem}[Linear spectral statistics]\label{thm2}
	Under Condition \ref{c1}, further assume that either of the following two assumptions holds: 
(1).$A_p$ is diagonal;(2).$E(x_{11}^4)=3$,
	then for any $f_1,\ldots,f_l\in\mathcal{M}$, the finite dimensional random vector $\big(G_{p}(f_1),\ldots G_p(f_l)\big)$ converges weakly to a Gaussian vector $\big(Y(f_1), \ldots, Y(f_l)\big)$ with mean function $E\{ Y(f)\}=0$ and covariance function
	\[
	\text{cov}\big(Y(f_j),Y(f_k)\big)=-\frac{1}{4\pi^2}\oint_{\mathcal{C}_1}\oint_{\mathcal{C}_2}f_j(z_1)f_k(z_2)\Lambda(z_1,z_2)dz_1dz_2,
	\]
	where the contours $\mathcal{C}_1$ and $\mathcal{C}_2$ are non-overlapping, counterclockwise, and enclosing the interval  $[-2\lim\sup a_1b_1,2\lim\sup a_1b_1]$, and $\Lambda(z_1,z_2)$ is defined by
		\[
	\begin{split}
	\Lambda(z_1,z_2)=&\frac{\partial^2}{\partial z_1\partial z_2}\lim_{p\rightarrow\infty}\frac{1}{p}\sum_{k=1}^p \vartheta_k(z_1,z_2)\quad \text{with}\\
	\vartheta_k(z_1,z_2)=&\tilde\epsilon_k(z_1)\tilde\epsilon_k(z_2)\bigg\{\frac{\frac{2}{p}a_k^2\bar\lambda^2_{B_n^2}\sum_{i<k}a_i^2\tilde\epsilon_i(z_1)\tilde\epsilon_i(z_2)}{1-\frac{1}{p}\bar\lambda_{B_n^2}\sum_{i<k}a_i^2\tilde\epsilon_i(z_1)\tilde\epsilon_i(z_2)}+a_k^2\bigg(\frac{\nu_4-3}{n}\sum_{j=1}^n\mathrm{B}_{n,jj}^2+2\bar\lambda_{B_n^2}\bigg)\bigg\}.
	\end{split}
	\]
\end{theorem}

}

\section{Sizes and powers of the testing procedure}\label{sec4}
\subsection{Calculate parameters $\mu$ and $\sigma^2$}
It seems that the mean and covariance functions in Theorem \ref{thm2} are very complicated. Based on the argument above Algorithm \ref{alg1}, it's sufficient to consider the case where $\Sigma_{\tilde U}=I$, or equivalently, $\tilde U$ is an orthogonal matrix. Orthogonal transformation has no effects on the eigenvalues of $\tilde S$. Then, it's sufficient to consider the case with $A_p=I$ in Theorem \ref{thm2}, although the more general result in the last section has its own interest in random matrix theory. The next corollary follows directly.
{
\begin{corollary}[Calculations for  $A_p=I_p$]\label{cor1}
	When $A_p=I_p$, we have
	\begin{equation}\label{simplify}
	\begin{split}
	m_p(z)=&s_p(z)=-\frac{1}{z+\bar\lambda_{B_n^2}s_p(z)}=-\tilde\epsilon_k(z), \quad\mathcal{B}_p(z)=\frac{\bar\lambda^2_{B_n^2}m_p^2(z)}{1-\bar\lambda_{B_n^2}m_p^2(z)}=\bar\lambda^2_{B_n^2}m_p^\prime(z),
	\end{split}
	\end{equation}
while $\mathcal{Y}_p(z)$ is the solution satisfying $\mathcal{Y}_p(z)=o(1)$ as $p\rightarrow\infty$ to equation
	\[
	\begin{split}
		&x=\bigg\{\frac{1}{p}\mathcal{A}_p(z)\times m_p^3(z)+m_p(z)\bigg\}\times \bar\lambda_{B_n^2}+\bigg\{\frac{1}{n}\sum_{j=1}^n\frac{(\lambda_j^{B_n})^2}{1-\lambda_j^{B_n}(\frac{p}{n})^{1/2}\mathcal{D}_{p2}(x,z)}\bigg\}\mathcal{D}_{p2}(x,z)
	\end{split}	\]
	with
$	\mathcal{D}_{p2}(x,z)=-m_p(z)/\{1+xm_p(z)\}$,
and 
	\[
	\begin{split}
		\Lambda(z_1,z_2)
		=&m^\prime(z_1)m^\prime(z_2)\bigg\{(\nu_4-3)\lim_{n\rightarrow\infty}n^{-1}\sum_{j=1}^n\mathrm{B}_{n,jj}^2+\frac{2\bar\lambda_{B_n^2}}{(1-\bar\lambda_{B_n^2}m(z_1)m(z_2))^2}\bigg\}.
	\end{split}
	\]
Further, 
the covariance function in Theorem \ref{thm2} is reduced to
	\[
	\text{cov}\big(Y(f_j),Y(f_k)\big)=\frac{1}{4\pi^2}\int_{-\sqrt{4\bar\lambda_{B_n^2}}}^{\sqrt{4\bar\lambda_{B_n^2}}}\int_{-\sqrt{4\bar\lambda_{B_n^2}}}^{\sqrt{4\bar\lambda_{B_n^2}}}f_j^\prime(t_1)f_k^\prime(t_2)H(t_1,t_2)dt_1dt_2,
	\]
	where
	\[
	\begin{split}
		H(t_1,t_2)=&\frac{1}{\bar\lambda^2_{B_n^2}}\bigg\{(\nu_4-3)\lim_{n\rightarrow \infty}\frac{1}{n}\sum_j\mathrm{B}_{n,jj}^2\bigg\}\sqrt{4\bar\lambda_{B_n^2}-t_1^2}\sqrt{4\bar\lambda_{B_n^2}-t_2^2}\\
		&+2\log\bigg\{\frac{4\bar\lambda_{B_n^2}-t_1t_2+\sqrt{(4\bar\lambda_{B_n^2}-t_1^2)(4\bar\lambda_{B_n^2}-t_2^2)}}{4\bar\lambda_{B_n^2}-t_1t_2-\sqrt{(4\bar\lambda_{B_n^2}-t_1^2)(4\bar\lambda_{B_n^2}-t_2^2)}}\bigg\}.
	\end{split}
	\]
	Moreover, if $f_j(x)=f_k(x)=x^2$, we have $\text{cov}\big(Y(f_j),Y(f_k)\big)=4\bar\lambda^2_{B_n^2}$.
\end{corollary}

}

Now we are ready to provide the parameters $\mu$ and $\sigma^2$ in the testing algorithm.
First, by (\ref{simplify}), we conclude that $F_p(x)$ satisfies a rescaled semi-circle law with density function
\begin{equation}\label{hpx}
	h_p(x)=\frac{1}{2\pi \bar\lambda_{B_n^2}}\sqrt{4\bar\lambda_{B_n^2}-x^2}=\frac{1}{2\pi \bar\lambda_{\Sigma_{V}^2}}\sqrt{4\bar\lambda_{\Sigma_{V}^2}-x^2},
\end{equation}
by noting that $B_n=I\otimes \Sigma_{V}$ in the matrix-variate model.
Then, 
\[
G_{p}(f)=\sum_{j=1}^pf(\lambda_j^{\tilde S})-p\int f(x) h_p(x)dx+\frac{p}{2\pi i}\oint_{\mathcal{C}}f(z)\mathcal{X}_p(z)dz.
\]
This suggests that
\begin{equation}\label{mu}
	\mu=p\int f(x) h_p(x)dx-\frac{p}{2\pi i}\oint_{\mathcal{C}}f(z)\mathcal{X}_p(z)dz,
\end{equation}
where  $\mathcal{X}_p$ is defined {with $A_p=I$, $B_n=I_T\otimes \Sigma_{V}$}, $n=Tq$ and $\mathcal{C}$ is a contour enclosing interval $[-2\lambda_{\max}(\Sigma_{V})-\epsilon_0,2\lambda_{\max}(\Sigma_{V})+\epsilon_0]$.
The variance parameter is
\begin{equation}\label{sigma2}
	\sigma^2=\frac{1}{4\pi^2}\int_{-\sqrt{4\bar\lambda_{\Sigma_{V}^2}}}^{\sqrt{4\bar\lambda_{\Sigma_{V}^2}}}\int_{-\sqrt{4\bar\lambda_{\Sigma_{V}^2}}}^{\sqrt{4\bar\lambda_{\Sigma_{V}^2}}}f^\prime(t_1)f^\prime(t_2)H(t_1,t_2)dt_1dt_2,
\end{equation}
where $H(t_1,t_2)$ is given in Corollary \ref{cor1} by replacing $B_n$ with $I_T\otimes \Sigma_{V}$. Further if $f(x)=x^2$,{
\begin{equation}\label{mu and sigma2}
	\mu=(p+1)\bar\lambda_{\Sigma_{V}^2}+q^{-1}(\nu_4-3)\sum_{j=1}^q\sigma_j^4,\quad \sigma^2=4\bar\lambda^2_{\Sigma_{V}^2},
\end{equation}
where $\text{diag}(\Sigma_{V})=(\sigma_1^2,\ldots,\sigma_q^2)$. The simplified formula for $\mu$ in (\ref{mu and sigma2}) is from direct calculation of expectation. If further $\Sigma_V=I_q$, the results in Corollary \ref{cor1} will be consistent with those in Theorem 1.1 of \cite{chen2015clt}.  See more details in Section I of our supplementary material.}

\subsection{Sizes and powers}
It's time to study the asymptotic sizes and powers of our testing procedure. Firstly, we propose assumptions on the data-generating model (\ref{model}), which are parallel to the previous Condition \ref{c1}.

\begin{condition}\label{c2}
	{Suppose that model (\ref{model}) satisfies Condition \ref{c1} as $\min\{p,T,q\}\rightarrow \infty$ by letting $n=Tq$, $A_p=\Sigma_U$ and $B_n=I_T\otimes\Sigma_V$. 
	 For identifiability, further let $q^{-1}\text{tr}(\Sigma_{V})=1$.}
\end{condition}

We remark that $T$ is allowed to be smaller than $p$ or $q$ in the above condition. We have the next two theorems on the asymptotic behavior of the testing statistics $\mathcal{T}$.

\begin{theorem}[Null hypothesis]\label{null}
	If Condition \ref{c2} holds, under the null hypothesis that $\Sigma_{U}=\Sigma_{U_0}>0$, we have
	$\mathcal{T}\longrightarrow\mathcal{N}(0,1)$ in distribution,
	where $\mathcal{T}$ is the testing statistic defined in Algorithm \ref{alg1} with $\mu$ and $\sigma^2$ given by (\ref{mu}) and (\ref{sigma2}).
\end{theorem}

\begin{theorem}[Alternative hypothesis]\label{alternative}
Under Condition \ref{c2} and the alternative hypothesis, if $f(x)=x^2$ and $
	p^{-1}\text{tr}(\Sigma_{U_0}^{-1/2}\Sigma_{U}\Sigma_{U_0}^{-1/2}-I)^2\ge c$ 
	for some constant $c>0$, then
	$P(|\mathcal{T}|>c_{\alpha})\rightarrow 1$ for any $\alpha\in(0,1)$,
	where $c_{\alpha}$ is the $\alpha$-th quantile of standard normal distribution.
\end{theorem}

Consequently, we reject the null hypothesis as long as the P-value from algorithm \ref{alg1} is smaller than some predetermined significant level $\alpha$.

\subsection{Estimating unknown parameters}
It's seen from (\ref{mu}) and (\ref{sigma2}) that  $\mu$ and $\sigma^2$ are dependent on the covariance matrix $\Sigma_{V}$ and the fourth moment $\nu_4$ of the data. In real applications, these parameters are unknown and need to be estimated. In this subsection, we discuss the estimation of these parameters. For simplicity, we only consider the case where $f(x)=x^2$ so that (\ref{mu and sigma2}) holds.

Indeed, by (\ref{mu and sigma2}), it's sufficient to find consistent estimators for the two parts of $\mu$. {For the first part, it suffices to consider  $p\bar\lambda_{\Sigma_{V}^2}$, while a natural estimator is
$ (p/q)\times \|(Tp)^{-1}\sum_{t=1}^T\tilde Y_t^\prime\tilde Y_t\|_F^2$.}
The consistency is demonstrated in the next lemma.
\begin{lemma}\label{consistcy: mu 1}
	Under Condition \ref{c2} and the null hypothesis $\Sigma_{U}=\Sigma_{U_0}$,  we have
	\[
	\frac{p}{q}\bigg\|\frac{1}{Tp}\sum_{t=1}^T\tilde Y_t^\prime\tilde Y_t\bigg\|_F^2-p\bar\lambda_{\Sigma_{V}^2}=\frac{q}{T}+o_p(1).
	\]
\end{lemma}

{Therefore, the estimator is consistent after a bias correction in the right hand side. 
It also indicates that a consistent estimator of $\sigma^2$ is
$\hat\sigma^2=(4/q)\times\{\|(Tp)^{-1}\sum_{t=1}^T\tilde Y_t^\prime\tilde Y_t\|_F^2-q^2/(Tp)\}$.
For the second part of $\mu$, note the relationship 
\begin{equation}\label{relation}
p^{-1}\text{var}\Big\{\text{tr}(\tilde Y_t^\prime\tilde Y_t)\Big\}=(\nu_4-3)\sum_{j=1}^q\sigma_j^4+2\|\Sigma_{V}\|_F^2,
\end{equation}
where $\sigma_j^2$ is defined in (\ref{mu and sigma2}).
Define $
\zeta=p^{-1}\text{var}\{\text{tr}(\tilde Y_t^\prime\tilde Y_t)\}$, $\omega=\sum_{j=1}^q\sigma_j^4$, $\tau=\|\Sigma_{V}\|_F^2$,
and the finite sample versions
\[
\begin{split}
\hat\zeta=&\frac{1}{Tp}\sum_{t=1}^T\bigg\{\text{tr}(\tilde Y_t^\prime\tilde Y_t)-\frac{1}{T}\sum_{t=1}^T\text{tr}(\tilde Y_t^\prime\tilde Y_t)\bigg\}^2,\quad
\hat\tau=\bigg\|\frac{1}{Tp}\sum_{t=1}^T\tilde Y_t^\prime\tilde Y_t\bigg\|_F^2.
\end{split}
\]
Then, a natural estimator for the second part of $\mu$ is given by $\hat \mu_2=q^{-1}(\hat\zeta-2\hat\tau)$. Moreover, if $\nu_4$ is of interest, a natural estimator is 
\[
\hat\nu_4=\max\{3+\hat\omega^{-1}(\hat\zeta-2\hat\tau),1\},\quad\text{where}\quad \hat\omega=\sum_{j=1}^q \{(Tp)^{-1}\sum_{i=1}^p\sum_{t=1}^T \tilde Y_{t,ij}^2\}^2.
\]
The consistency is guaranteed by the next lemma.
\begin{lemma}\label{consistency:nu}
	Write $\mu_2=q^{-1}(\nu_4-3)\sum_{j=1}^q\sigma_j^4$. Under Condition \ref{c2} and the null hypothesis, we have $
	\hat\mu_2-\mu_2=o_p(1)$.
Further if $E(x_{1,11}^8)<\infty$, we have $\hat \nu_4-\nu_4=o_p(1)$.
\end{lemma}
}

\section{Bootstrapping}\label{sec5}
{By (\ref{mu}) and (\ref{sigma2}), the derivations of $\mu$ and $\sigma^2$ are not easy even if  $\Sigma_{V}$ and $\nu_4$ are given. }The mean correction term $\mathcal{X}_p(z)$ involves integration over
complex number contour, where the term $\mathcal{Y}_p(z)$ has no closed-form. Therefore, numerical approximation with high accuracy is needed. To  overcome this challenge, in this section we propose a bootstrap resampling algorithm to generate critical values for the testing procedure. The algorithm can be regarded as  an extension of the technique in \cite{lopes2019bootstrapping} to the matrix-variate regime.

The algorithm is motivated by the fact that $\mu$ and $\sigma^2$ only depend on $\Sigma_{V}$ and $\nu_4$ under the null hypothesis. Therefore, if $\Sigma_{V}$ and $\nu_4$ are given, we can regenerate some random matrices whose linear spectral statistics have exactly the same limiting distributions as the original ones. Now we assume some estimators of $\Sigma_{V}$ and $\nu_4$ are available,  denoted by $\hat\Sigma_{V}$ and $\hat\nu_4$ satisfying $q^{-1}\text{tr}(\hat\Sigma_{V})=1$.  Then, the bootstrap procedure is presented in  Algorithm \ref{alg2} below.

Define the L{\'e}vy-Prokhorov distance between two probability measure on $\mathbb{R}^d$ by
\[
d_{LP}(P,Q)=\inf\{\epsilon: P(\mathcal{A})\le Q(\mathcal{A}^{\epsilon})+\epsilon, Q(\mathcal{A})\le P(\mathcal{A}^\epsilon)+\epsilon, \forall \mathcal{A}\in \mathcal{B}(\mathbb{R}^d)\},
\]
where $\mathcal{A}^{\epsilon}$ is defined by
\[
\mathcal{A}^{\epsilon}=\{x\in \mathbb{R}^d:\text{there exists } y\in \mathcal{A},\text{ such that } \|x-y\|<\epsilon\}.
\]
If this distance converges to 0, it essentially indicates convergence in distribution. We have the next theorem on the consistency of the bootstrap procedure.
{\begin{theorem}[Bootstrap consistency]\label{bootstrap}
	Let $\mathcal{K}=\sum_{j=1}^pf(\lambda_j^{\tilde S})-p\int f(x)h_p(x)dx$. Under Condition \ref{c2} and the null hypothesis, if $\hat\Sigma_{V}$ and $\hat \nu_4$ satisfy 
		\[
	q^{-1}\text{tr}(\hat\Sigma_V)=1,\quad \|\hat \Sigma_{V}-\Sigma_{V}\|=o_p(\min\{1,\sqrt{T}/p\}),\quad \hat\nu_4-\nu_4=o_p(1),
	\] we have $
	d_{LP}(\mathcal{K},\mathcal{K}_1^*\mid\{X_t\}_{t=1}^T)\longrightarrow 0$
	in probability.
\end{theorem}

Therefore, we reject the null hypothesis as long as
\[
\mathcal{K}<c_{\alpha/2}^*,\quad \text{ or }\quad\mathcal{K}>c_{1-\alpha/2}^*,
\]
where $c_{\alpha}^*$ is the $\alpha$-th sample quantile of $\{\mathcal{K}_b^*\}_{b=1}^B$. It's worth mentioning that the convergence rate in Theorem \ref{bootstrap} can be easily fulfilled with sparsity conditions. For instance, the optimal convergence rate in \cite{cai2011adaptive} is $\{\log q/(Tp)\}^{(1-s)/2}$, where $s$ is related to the sparsity and can be 0.  Then,  the rate for $\hat\Sigma_V$ holds as long as
$(p\log q)^{1/2}/T= o(1).$}

\begin{algorithm}[H]
	\caption{Bootstrap resampling algorithm}\label{alg2}
	\begin{algorithmic}[1]
		\Require data matrices $\{\tilde Y_t\}_{t=1}^T$,  function $f$, estimators $\hat\Sigma_{V}$ and $\hat\nu_4$, bootstrap sample size $m$.
		\Ensure  A series of bootstrapped testing statistics $\{\mathcal{K}_b^*\}_{b=1}^m$.
		
		\State for $b=1$, generate a series of $p\times q$ matrices $\{
		Z_t\}_{t=1}^T$ whose entries are from independent and identically distributed  Pearson distribution system with parameters $(0,1,0,\hat\nu_4)$,
		\State calculate the renormalized sample covariance matrix
		\[
		S^*_1=\{(Tq)/p\}^{1/2}\{(Tq)^{-1}\sum_{t=1}^TZ_t\hat\Sigma_{V}Z_t^\prime-I\},
		\]
		\State calculate the normalized linear spectral statistic
		\[
		\mathcal{K}_1^*=\sum_{j=1}^pf(\lambda_j^{S_1^*})-p\int f(x)\hat h_p(x) d(x),
		\]
		where $\hat h_p(x)$ is defined in (\ref{hpx}) by replacing $\Sigma_{V}$ with $\hat\Sigma_{V}$,
		\State repeat the above procedure $m$ times, and output the series $\{\mathcal{K}_b^*\}_{b=1}^m$.
	\end{algorithmic}
\end{algorithm}
\begin{remark}
	If $\hat\nu_4=1$, we generate the entries of $Z_t$ from independent and identically distributed  Bernoulli distribution with $\Pr (z_{t,ij}=\pm 1)=0.5$.
\end{remark}

\section{Noised model}\label{sec:noise}
{The standard model (\ref{separable}) is sometimes too idealistic since it does not allow any noise. In this section, we  generalize the model by letting
\begin{equation}\label{noise model}
\check y_{ij}=y_{ij}+\sigma_\alpha\varphi+\sigma_\beta\phi_{ij}, \quad \sigma_\alpha,\sigma_\beta\ge 0,
\end{equation}
where $Y=(y_{ij})$ has the separable structure  in (\ref{separable}). $\varphi$ is a common random noise independent of $Y$ with $E(\varphi)=0$, $E(\varphi^2)=1$ and $E(\varphi^4)<\infty$. $\phi_{ij}$'s are independent individual noises with $E(\phi_{ij})=0$, $E(\phi_{ij}^2)=1$ and $E(\phi_{ij}^4)=\tilde\nu_4<\infty$. $\phi_{ij}$'s are independent of $\varphi$ and $Y$.  $\sigma_{\alpha}$ and $\sigma_\beta$ are two constants.  The proposal of common noise $\varphi$ is partially motivated by our real data analysis where a spiked eigenvalue is found in both the column-column and row-row sample covariance matrices. The target of this section is to test the hypothesis $\Sigma_{U}=\Sigma_{U_0}$ for some $\Sigma_{U_0}>0$ under the noised model (\ref{noise model}). As a byproduct, we will also discuss how to estimate  $\sigma_\alpha$, $\sigma_\beta$. 

Given independent observations, model (\ref{noise model}) can also be written in matrix form as 
\begin{equation}\label{noise matrix}
\check Y_t=Y_t+\sigma_\alpha\varphi_t\mathbf{1}_p\mathbf{1}_q^\prime+\sigma_\beta \Phi_t\quad (1\le t\le T),\text{ with }\quad \Phi_t=(\phi_{t,ij}),
\end{equation}
where $\mathbf{1}_p$ is a p-dimensional vector with all entries being 1.
When $\sigma_\alpha>0$, we are actually considering a spiked model where the common noise dominates in the population covariance matrix, and may further dominate in the  linear spectral statistics. Therefore, in order to test $\Sigma_{U}$, the first step is to remove the common noise from the system. 

A direct approximation to $(\sigma_\alpha\varphi_t)$ is $\hat\sigma_{\alpha,t}=(pq)^{-1}\mathbf{1}_p^\prime\check Y_t\mathbf{1}_q$. Therefore, we subtract the common noise by defining \[
\hat Y_t=\check Y_t-\hat\sigma_{\alpha,t}\mathbf{1}_p\mathbf{1}_q^\prime=Y_t+\sigma_{\beta}\Phi_t-(pq)^{-1}\mathbf{1}_p^\prime(Y_t+\sigma_\beta\Phi_t)\mathbf{1}_q\times\mathbf{1}_p\mathbf{1}_q^\prime.
\]
	Consequently, we can assume $\sigma_\alpha=0$ without loss of generality if only $\Sigma_U$ is of concern. Otherwise, one can always have a good estimator for $\sigma_{\alpha}^2$ simply by $\hat\sigma_\alpha^2=T^{-1}\sum_{t=1}^T\hat\sigma_{\alpha,t}^2$. 
	However, it's more challenging to handle the individual noises $\phi_{ij}$'s. Firstly, we estimate $\sigma_\beta^2$. Write $\Omega=	E\{\text{Vec}(\check Y_1)\text{Vec}(\check Y_1)^\prime\}$. When $\sigma_{\alpha}=0$, elementary calculations lead to
	\[
\text{diag}(\Omega)=\text{diag}(\Sigma_V)\otimes\text{diag}(\Sigma_U)+\sigma_\beta^2I.
	\]
We reshape $\text{diag}(\Omega)$ into a $p\times q$ matrix column by column so that
\[
\text{Re}(\Omega)=\vec{u}\vec{v}^{\,\prime}+\sigma_\beta^2\mathbf{1}_p\mathbf{1}_q^\prime,\text{ where } \vec{u}=(\Sigma_{U,11},\ldots,\Sigma_{U,pp}),\vec{v}=(\Sigma_{V,11},\ldots,\Sigma_{V,qq}).
\]
Denote the leading left and right singular vectors of $(I-p^{-1}\mathbf{1}_p\mathbf{1}_p^\prime)\text{Re}(\Omega)(I-q^{-1}\mathbf{1}_q\mathbf{1}_q^\prime)$ as $\tilde u$  and $\tilde v$, respectively. Then, there exist coefficients $c_1,c_2,k_1,k_2$ such that
\[ 
p^{-1/2}\vec{u}=c_1\tilde u+p^{-1/2}k_1\mathbf{1}_p,\quad q^{-1/2}\vec{v}=c_2\tilde v+q^{-1/2}k_2\mathbf{1}_q.
\]
Moreover, $\tilde u$ is orthogonal of $\mathbf{1}_p$, and $\tilde v$ is orthogonal of $\mathbf{1}_q$. Therefore,
\[
\begin{split}
&(pq)^{-1/2}\tilde u^\prime \text{Re}(\Omega)\tilde v=c_1c_2,\quad (p^2q)^{-1/2}\mathbf{1}_p^\prime \text{Re}(\Omega)\tilde v=k_1c_2,\\
&(pq^2)^{-1/2}\tilde u^\prime \text{Re}(\Omega)\mathbf{1}_q=c_1k_2,\quad (pq)^{-1}\mathbf{1}_p^\prime \text{Re}(\Omega)\mathbf{1}_q=k_1k_2+\sigma_{\beta}^2.
\end{split}
\]
Consequently, 
\[
\sigma_{\beta}^2=(pq)^{-1}\mathbf{1}_p^\prime \text{Re}(\Omega)\mathbf{1}_q-\{pq\tilde u^\prime \text{Re}(\Omega)\tilde v\}^{-1}\{\mathbf{1}_p^\prime \text{Re}(\Omega)\tilde v\times \tilde u^\prime \text{Re}(\Omega)\mathbf{1}_q\}.
\]
A natural estimator for $\sigma_{\beta}^2$ is then obtained by replacing $\Omega$ with $
\hat\Omega=T^{-1}\sum_{t=1}^T\text{Vec}(\hat Y_t)\text{Vec}(\hat Y_t)^\prime$
in the above process, denoted by $\hat\sigma_{\beta}^2$. The next lemma demonstrates the consistency.
\begin{lemma}\label{sigma beta}
	Under condition \ref{c2}, if $\vec{u}^\prime(I-p^{-1}\mathbf{1}_p\mathbf{1}_p^\prime)\vec{u}/p\ge \epsilon_1$ and $\vec{v}^\prime(I-p^{-1}\mathbf{1}_q\mathbf{1}_q^\prime)\vec{v}/p\ge \epsilon_2$ for some constants $\epsilon_1,\epsilon _2>0$, we have
	\[
	\hat\sigma_{\beta}^2-\sigma_{\beta}^2=O_p[(pq)^{-1}+T^{-1}\times( \min\{T,p,q\})^{-1/2}+(Tpq)^{-1/2}].
	\]
	\end{lemma}
\begin{remark}
	The conditions for $\vec{u}$ and $\vec{v}$ in the above lemma are mainly for identifiability between the Kronecker product part and the individual noises.  For instance, if $U$ and $V$ are both identity matrices, we will have $\epsilon_1=\epsilon_2=0$, but the model is not identifiable. Under such cases, it's more meaningful to test $\Sigma_U$ and $\sigma_\beta^2$ jointly based on Theorem \ref{noise clt} below.  
	\end{remark}

	Motivated by the normalization in (\ref{renormalize}) and (\ref{transform}), we define
	\begin{equation}\label{noise normalization}
	\mathcal{Y}_t=\Sigma_{U_0}^{-1/2}\hat Y_t,\quad \bar{\mathcal{S}}=\{(Tq)/p\}^{1/2}\{(Tq)^{-1}\sum_{t=1}^T\mathcal{Y}_t\mathcal{Y}_t^\prime-\mathcal{E}_0\},
	\end{equation}
	where $\mathcal{E}_0=q^{-1}E(\mathcal{Y}_1\mathcal{Y}_1^\prime)$ under the null $\Sigma_U=\Sigma_{U_0}$ is given by
	\[
\begin{split}
		\mathcal{E}_0=&I+\sigma_\beta^2\Sigma_{U_0}^{-1}+(pq)^{-1}(p^{-1}\mathbf{1}_p^ \prime\Sigma_{U_0}\mathbf{1}_p\times q^{-1}\mathbf{1}_q^\prime\Sigma_V\mathbf{1}_q-\sigma_\beta^2)\Sigma_{U_0}^{-1/2}\mathbf{1}_p\mathbf{1}_p^\prime\Sigma_{U_0}^{-1/2}\\
		&-(pq^2)^{-1}\mathbf{1}_q^\prime\Sigma_V\mathbf{1}_q(\Sigma_{U_0}^{1/2}\mathbf{1}_p\mathbf{1}_p^\prime\Sigma_{U_0}^{-1/2}+\Sigma_{U_0}^{-1/2}\mathbf{1}_p\mathbf{1}_p^\prime\Sigma_{U_0}^{1/2}).
\end{split}
	\]
	The normalization heavily relies on $\sigma_{\beta}^2$. This is the reason why we need to estimate $\sigma_\beta^2$ first. In the following,  we will focus on linear spectral statistics of $\bar{\mathcal{S}}$ in a special case where $f(x)=x^2$, equivalently, $\text{tr}\bar{\mathcal{S}}^2$. The next theorem gives the asymptotic distribution under the null hypothesis.
	\begin{theorem}\label{noise clt}
		Assume that Condition \ref{c2} and the null hypothesis hold, while $T\le (\max\{p,q\})^C$ for some constant $C>0$. Then, $\tilde\sigma^{-1}\{\text{tr}(\bar{\mathcal{S}}^2-\tilde\mu)\}\rightarrow \mathcal{N}(0,1)$ in distribution, with
		\[
		\begin{split}
		\tilde\mu=&q^{-1}(\nu_4-3)\sum_{j=1}^q(\Sigma_V)_{jj}^2+(p+1)\bar\lambda_{\Sigma_V^2}+2(p+1)\sigma_\beta^2\bar\lambda_{\Sigma_{U_0}^{-1}}\\&+\sigma_{\beta}^4\{p^{-1}(\tilde\nu_4-3)\sum_{j=1}^p(\Sigma_{U_0}^{-1})_{jj}^2+\bar\lambda_{\Sigma_{U_0}^{-2}}+p\bar\lambda^2_{\Sigma_{U_0}^{-1}}\},\\
		\tilde\sigma^2=&4\bar\lambda^2_{\Sigma_{V}^2}+4\sigma_{\beta}^8\bar \lambda^2_{\Sigma_{U_0}^{-2}}+8\sigma_{\beta}^4(\bar\lambda^2_{\Sigma_{U_0}^{-1}}+\bar\lambda_{\Sigma_{U_0}^{-2}}+\bar\lambda^2_{\Sigma_{U_0}^{-1}}\bar\lambda_{\Sigma_{V}^2})+16\sigma_{\beta}^2(\bar\lambda_{\Sigma_{U_0}^{-1}}\bar\lambda_{\Sigma_{V}^2}+\sigma_{\beta}^4\bar\lambda^2_{\Sigma_{U_0}^{-1}}).
		\end{split}
		\]
	\end{theorem}
The asymptotic power is verified by the next  theorem.
	\begin{theorem}\label{noise power}
	Assume that Condition \ref{c2} holds and $T\le (\max\{p,q\})^C$ for some constant $C>0$. If $
	p^{-1}\text{tr}(\Sigma_{U_0}^{-1/2}\Sigma_{U}\Sigma_{U_0}^{-1/2}-I)^2\ge c$ 
	for some constant $c>0$, we have $|\tilde\sigma^{-1}\{\text{tr}(\bar{\mathcal{S}}^2-\tilde\mu)\}|\ge C$ with probability tending to 1 for any constant $C>0.$
\end{theorem}
With the above two theorems, it's easy to conduct the test based on $\tilde\sigma^{-1}\{\text{tr}(\bar{\mathcal{S}}^2-\tilde\mu)$ if all the parameters
are given. The unknown parameters in the above procedure are $\mathcal{E}_0$, $\tilde\mu$ and $\tilde\sigma$. Estimation of these parameters is standard and organized in Section A of our supplementary material. Since the asymptotic mean and variance have closed-form expressions, it's not necessary to use the bootstrap technique.}

\section{Numerical studies}\label{sec6}
In this section, we investigate  empirical performance of the proposed testing procedure and bootstrap technique. {Here we only show the results for the model without noise, while those for the noised model are in Section B.3 of our supplementary material. }According to data-generating model (\ref{model}), we sample the entries of $X_t$ from independent and identically distributed standard normal distribution $\mathcal{N}(0,1)$ or Bernoulli distribution with $\Pr(x_{t,ij}=\pm 1)=0.5$ so that both continuous and discrete distributions are considered. For the null hypothesis, we generate covariance matrix $\Sigma_{U}$ by spectral decomposition $\Sigma_{U}=\Sigma_{U_0}=\Gamma\Lambda\Gamma^\prime$ where $\Gamma$ is from Haar distribution and the diagonal entries of $\Lambda$ are sampled independently from uniform distribution $\mathcal{U}[1,2]$. 

For the alternative hypothesis, we consider two common cases.  The first one is related to a low-rank perturbation such that
\[
H_{A1}:\Sigma_{U}=\Sigma_{U_1}=\Sigma_{U_0}+p^{-1/2}\beta_1\gamma\gamma^\prime,
\]
where $\gamma$ is the $p$-dimensional vector whose entries are sampled independently from  $\mathcal{N}(0,1)$. For the second case, we consider small but dense perturbations on eigenvalues and let
\[
H_{A2}:\Sigma_{U}=\Sigma_{U_2}=\Gamma(\Lambda+\beta_2I_p)\Gamma^\prime.
\]
It's easy to see under both alternative hypothesis, the condition in Theorem \ref{alternative} is satisfied as long as $\beta_1$ and $\beta_2$ are some positive constants.
In the simulations, $\beta_1=\beta_2=0.1$.  To ease the computation of $\mathcal{X}_p(z)$, we set a sparse block diagonal structure for  $\Sigma_{V}$  as 
\[
(\Sigma_{V})_{ij}=\left\{
\begin{aligned}
1,&\quad i=j,&\\
0.5,&\quad (i,j)\in\{(2k,2k-1),(2k-1,2k)\}_{k=1,\ldots,q/2},&\\
0,& \quad \text{otherwise,}&\\
\end{aligned}\right.
\]
so that $q^{-1}\text{tr}(\Sigma_{V})=1$ and the empirical spectral distribution of $\Sigma_{V}$ is
$0.5\delta_{1.5}+0.5\delta_{0.5}$,
where $\delta_x$ is the Dirac measure at $x$. We test the null hypothesis $\Sigma_{U}=\Sigma_{U_0}$ under all scenarios with various combinations of $(T,p,q)$ and under the two distribution families.

The testing procedure can be implemented using the formulas for $\mu$ and $\sigma^2$ in (\ref{mu and sigma2}) with $f(x)=x^2$. We denote the results by $\mathcal{T}_{FO}$. It can also be implemented with the bootstrap algorithm. For the bootstrap procedure, we will further compare the results if taking $\Sigma_{V}$ and $\nu_4$ as given (denoted by $\mathcal{T}_{BG}$) with those if the parameters are estimated (denoted by $\mathcal{T}_{BE}$). The adaptive thresholding technique in \cite{cai2011adaptive} is applied to the estimation of $\Sigma_{V}$ with data matrix $(\tilde Y_1^\prime,\ldots,\tilde Y_T^\prime)$, while $\nu_4$ is estimated according to Lemma \ref{consistency:nu}. We rescale $\hat\Sigma_{V}$ so that $q^{-1}\text{tr}(\hat\Sigma_{V})=1$. Critical values from bootstrap algorithm are based on $B=200$ replications. As a competitor, another naive but natural approach is to first estimate $\Sigma_{V}$ and then right-multiply $\hat\Sigma_{V}^{-1/2}$ to each $Y_t$ so that the problem is transformed to testing conventional vector-variate covariance matrix. Then, various approaches in the literature  are available.  We denote this naive method as $\mathcal{T}_{NA}$. The approach in \cite{chen2015clt} is selected for the testing after transformation, while $\hat\Sigma_{V}$ and $\hat\nu_4$ are estimated similarly  in our bootstrap algorithm. All the simulation results hereafter are based on 1000 replications.

\begin{table*}[hbtp]
	\begin{center}
		\small
		\addtolength{\tabcolsep}{4pt}
		\caption{Empirical sizes of different testing procedures under the null hypothesis. }\label{tab1}
		\renewcommand{\arraystretch}{1}
		\scalebox{0.9}{ 	
			\begin{tabular*}{16cm}{lllllllllllll}
				&&&\multicolumn{4}{l}{$\alpha=0.05$}&\multicolumn{4}{l}{$\alpha=0.10$}\\
				family&$p(q)$&$T$&$\mathcal{T}_{FO}$&$\mathcal{T}_{BG}$&$\mathcal{T}_{BE}$&$\mathcal{T}_{NA}$&$\mathcal{T}_{FO}$&$\mathcal{T}_{BG}$&$\mathcal{T}_{BE}$&$\mathcal{T}_{NA}$\\
				Normal&20&20&0.063&0.048&0.037&0.691&0.120&0.094&0.079&0.768
				\\
				&20&100&0.052&0.064&0.052&0.676&0.103&0.118&0.103&0.753
				\\
				&100&20&0.069&0.052&0.036&1&0.126&0.116&0.080&1
				\\
				&100&100&0.047&0.060&0.049&1&0.106&0.110&0.095&1\\
				Bernoulli&20&20&0.040&0.068&0.050&0.666&0.093&0.120&0.101&0.755
				\\
				&20&100&0.043&0.053&0.050&0.701&0.089&0.111&0.098&0.784
				\\
				&100&20&0.044&0.055&0.043&1&0.091&0.108&0.096&1
				\\
				&100&100&0.047&0.051&0.054&1&0.094&0.109&0.104&1\\
		\end{tabular*}}		
	\end{center}
\end{table*}

\begin{table*}[hbtp]
	\begin{center}
		\small
		\addtolength{\tabcolsep}{8pt}
		\caption{Empirical powers of different testing procedures under the alternative hypothesis $H_{A1}$.}\label{tab2}
		\renewcommand{\arraystretch}{1}
		\scalebox{0.9}{ 	
			\begin{tabular*}{16cm}{lllllllllllll}
				&&&\multicolumn{3}{l}{$\alpha=0.05$}&\multicolumn{4}{l}{$\alpha=0.10$}\\
				family&$p(q)$&$T$&$\mathcal{T}_{FO}$&$\mathcal{T}_{BG}$&$\mathcal{T}_{BE}$&$\mathcal{T}_{FO}$&$\mathcal{T}_{BG}$&$\mathcal{T}_{BE}$\\
				Normal&20&20&0.075&0.043&0.038&0.134&0.079&0.067
				\\
				&20&60&0.392&0.231&0.321&0.499&0.353&0.426
				\\
				&20&100&0.883&0.765&0.833&0.926&0.843&0.902
				\\
				&100&20&0.283&0.149&0.213&0.392&0.225&0.314
				\\
				&100&60&0.978&0.944&0.974&0.991&0.965&0.989
				\\
				&100&100&1&1&1&1&1&1\\
				Bernoulli&20&20&0.050&0.048&0.048&0.091&0.101&0.089
				\\
				&20&60&0.361&0.253&0.339&0.485&0.365&0.479
				\\
				&20&100&0.897&0.817&0.882&0.945&0.900&0.936
				\\
				&100&20&0.289&0.168&0.286&0.400&0.260&0.382
				\\
				&100&60&0.978&0.947&0.980&0.992&0.971&0.992\\
				&100&100&1&1&1&1&1&1\\
		\end{tabular*}}		
	\end{center}
\end{table*}

Table \ref{tab1} displays the empirical sizes of the mentioned testing procedures under the null hypothesis. It's seen that when $T,p,q$ are sufficiently large, the empirical sizes of the proposed testing procedure and its bootstrapped versions are well-controlled. The bootstrap method can  be even more accurate sometimes when the sample size $T$ is small. On the contrary, the naive testing approach $\mathcal{T}_{NA}$ fails to control the sizes under all settings. It indicates that the estimation error of $\Sigma_{V}$ has non-negligible effects  on this approach.

The performance of the proposed testing procedure and its bootstrapped versions under the alternative hypothesis $H_{A1}$  is reported in Table \ref{tab2}.  The empirical powers can be weak when $T,p,q$ are small, but they converge quickly to one as $T,p,q$ grow, as claimed in Theorem \ref{alternative}. Moreover, there is no significant difference between the powers by bootstrapped  approaches  and the formula-based one. Similar results are observed for the alternative hypothesis $H_{A2}$. To save space, we leave the details into our supplementary material. In conclusion, the proposed  procedure can effectively test the two alternative hypothesis, while the bootstrap algorithm provides accurate approximation for the testing statistics. The supplementary material also contains a real-data analysis section, where the proposed  procedure is applied to testing covariance matrices of well-structured portfolio returns.

\section{Discussion}\label{sec8}
In current settings, we assume the observations are independent across time. This condition can be further relaxed to  separable temporal correlations, i.e.,
\[
\{\text{Vec}(Y_1),\ldots,\text{Vec}(Y_T)\}=(V\otimes U)\{\text{Vec}(X_1),\ldots,\text{Vec}(X_T)\}C,
\]
where $C$ is $T\times T$ deterministic matrix indicating temporal dependence. {The proposed procedure is also applicable to tensor (array) data, such as 
\[
Y_t=X_t\times_1 U_1\times_2\cdots\times_d U_d \quad (t=1,\ldots,T),
\]
where $X_t\in\mathbb{R}^{p_1\times \cdots\times p_d}$ is a tensor of order $d$ with independent and identically distributed entries, $U_1,\ldots,U_d$ are corresponding coefficient matrices,  and $\times_k$ stands for $k$-mode multiplication.  See \cite{kolda2009tensor} for more details. Then, we can calculate sample covariance matrix using tensor matricization, which will satisfy the form in (\ref{sp}) thus Theorem \ref{thm2} is still available.    Another closely related problem is to test whether the separable model (\ref{model}) holds or not. Motivated by the estimation of $\sigma_\beta^2$, this is possible according to the edge eigenvalues of $\text{Re}(\hat\Omega)$. When model (\ref{model}) holds, $\text{Re}(\hat\Omega)$ has only one spiked eigenvalue, while more will be observed if noises exist. A similar method has been studied by \cite{guggenberger2022test} in low dimensions. The extension to high dimensions requires to know the asymptotic distribution of edge eigenvalues of $\text{Re}(\hat\Omega)$, which is not available. } On the other hand, inspired by our real data analysis and the noised model, the Kronecker product covariance model may contain a low rank plus noise structure. Up to our knowledge, the existing literature has not considered this case yet. We are also interested in estimating the Kronecker product covariance matrices under such a structure in high dimensions.  This structure can also potentially model the volatility of high-frequency data by generalizing the factor model in \cite{kong2017on} and \cite{kong2018systematic} to matrix-variate regime. 

\section{Acknowledgments}
We thank the editor, associate editor and three anonymous referees for careful reading and valuable comments.

\section{Supplementary material}
	This supplementary material provides additional simulation results, a real data example and  technical proofs of  all the theorems, corollaries  and lemmas  in the paper ``Testing Kronecker Product Covariance Matrices for High-dimensional Matrix-Variate Data''. There are {ten sections. Section \ref{sec: estimate} presents how to estimate unknown parameters for the test of noised model, corresponding to Section \ref{sec:noise} of the main paper. }Section \ref{sec: addsim} provides more simulation results not shown in the main paper. Section \ref{real data} is for real data analysis.   Sections \ref{seca} to \ref{sece} are for the proof of central limit theorem for linear spectral statistics, corresponding to Theorem \ref{thm2} of the main paper. {For the convenience of readers, we give an outline of the proof in Section \ref{seca}. } The proof for theorems and lemmas in Sections \ref{sec4} and \ref{sec5} of the main paper is devoted to Section \ref{secf}. {Section \ref{sec: proof noise} proves results related to the noised model.

	We first introduce some notation used throughout the proof. We use boldface for matrices and vectors. Let $\|\cdot\|$ denote the spectral norm of a matrix or Euclidean norm of a vector. $\Ab^\prime$ (or $\Ab^*$) denotes the transpose (or conjugate transpose) of a  real (complex) matrix $\Ab$.  $\mathrm{A}_{ij}$ is the $(i,j)$-th entry of $\Ab$. $\lambda_j(\Ab)$ and $\lambda_j^{\Ab}$ stand for the $j$-th largest eigenvalue of $\Ab$.  For complex number $z$, $\Im z$ and $\Re z$ denote its imaginary and real parts, respectively. $i$ is the complex number satisfying $i=\sqrt{-1}$.  $\bar\lambda_{\Bb^h}=n^{-1}\text{tr}(\Bb^h)$ for symmetric matrix $\Bb$ and any integer $h>0$, where $\text{tr}(\cdot)$ denotes the trace.  $\overset{d}{\longrightarrow}$,  $\overset{P}{\longrightarrow}$ ( $\overset{i.p.}{\longrightarrow}$),  $\overset{a.s.}{\longrightarrow}$ denote convergence in distribution, in probability and almost surely, respectively. $\otimes$ denotes Kronecker product. $\mathbb{E}(\cdot)$ is for expectation.

\begin{appendices}
	\section{Estimation of unknown parameters  for noised model}\label{sec: estimate}
	{Under the noised model, in order to construct testing statistics,  we need to first estimate the normalization matrix $\mathcal{E}_0$ in (\ref{noise normalization}), and the asymptotic mean $\tilde \mu$, variance $\tilde\sigma^2$ in Theorem \ref{noise clt}. 
		
		We start with $\mathcal{E}_0$. By (\ref{noise normalization}), it suffices to provide plug-in estimators for $\sigma_\beta^2$ and $q^{-1}\mathbf{1}_q^\prime\bSigma_{\Vb}\mathbf{1}_q$. $\hat\sigma_\beta^2$ is given by Lemma \ref{sigma beta}. Therefore, we focus on the latter. By elementary calculation,
		\[
		\begin{split}
		\mathbb{E}\frac{1}{q}\mathbf{1}_q^\prime\bigg(\frac{1}{Tp}\sum_{t=1}^T\mathcal{Y}_t^\prime\mathcal{Y}_t\bigg)\mathbf{1}_q=&\frac{1}{q}\mathbf{1}_q^\prime\bSigma_{\Vb}\mathbf{1}_q+\frac{\sigma_\beta^2}{p}\text{tr}\bSigma_{\Ub_0}^{-1}+\frac{\mathbf{1}_p^\prime\bSigma_{\Ub_0}^{-1}\mathbf{1}_p}{p^2}\bigg(\frac{\mathbf{1}_p^\prime\bSigma_{\Ub_0}\mathbf{1}_p}{p}\frac{\mathbf{1}_q^\prime\bSigma_{\Vb}\mathbf{1}_q}{q}-\sigma_\beta^2\bigg)\\
		&-\frac{2}{p}\frac{\mathbf{1}_q^\prime\bSigma_{\Vb}\mathbf{1}_q}{q}.
		\end{split}
		\]
		Therefore, a natural estimator for $q^{-1}\mathbf{1}_q^\prime\bSigma_{\Vb}\mathbf{1}_q$ is
		\begin{equation}\label{plug-in 1}
		\bigg\{\frac{1}{q}\mathbf{1}_q^\prime\bigg(\frac{1}{Tp}\sum_{t=1}^T\mathcal{Y}_t^\prime\mathcal{Y}_t\bigg)\mathbf{1}_q-\frac{\hat\sigma_\beta^2}{p}\text{tr}\bSigma_{\Ub_0}^{-1}+\frac{\hat\sigma_\beta^2}{p^2}\mathbf{1}_p^\prime\bSigma_{\Ub_0}^{-1}\mathbf{1}_p\bigg\}\times \bigg(1+\frac{\mathbf{1}_p^\prime\bSigma_{\Ub_0}^{-1}\mathbf{1}_p}{p^2}\frac{\mathbf{1}_p^\prime\bSigma_{\Ub_0}\mathbf{1}_p}{p}-\frac{2}{p}\bigg)^{-1}.
		\end{equation}
		We use the plug-in estimators $\hat\sigma_\beta^2$ and (\ref{plug-in 1}) to calculate $\mathcal{E}_0$, denoted as $\hat{\mathcal{E}}_0$.  The next lemma holds.
		\begin{lemma}\label{lem: plug-in 1}
			Under the  conditions  in Theorem \ref{noise clt}, further assume $T=o(p^2q)$, $q=o(\min\{T^2,Tp\})$. Then, we have $\text{tr}\hat{\mathcal{S}}^2-\text{tr}\bar{\mathcal{S}}^2=o_p(1)$, where $\hat{\mathcal{S}}$ is calculated with $\hat{\mathcal{E}}_0$ instead of $\mathcal{E}_0$ in the definition of $\bar{\mathcal{S}}$. 
		\end{lemma}
		The proof is given in Section \ref{sec: proof estimate}. 
		Next, we estimate $\tilde\mu$. This is more challenging because the order of $\tilde\mu$ is  $O(p)$. Therefore, it's necessary  to do similar bias correction  as we did in Lemma \ref{consistcy: mu 1}. Motivated by Lemma \ref{consistcy: mu 1}, we first calculate the limit of $\|(Tp)^{-1}\sum_{t=1}^T\mathcal{Y}_t^\prime\mathcal{Y}_t\|_F^2$, shown in the next lemma whose proof is given in Section \ref{sec: proof estimate}.
		\begin{lemma}\label{tilde mu}
			Under the conditions  in Theorem \ref{noise clt}, we have
			\[
			\frac{p}{q}\bigg\|\frac{1}{Tp}\sum_{t=1}^T\mathcal{Y}_t^\prime\mathcal{Y}_t\bigg\|_F^2-(p\bar\lambda_{\bSigma_{\Vb}^2}+p\sigma_\beta^4\bar\lambda^2_{\bSigma_{\Ub_0}^{-1}}+2p\sigma_\beta^2\bar\lambda_{\bSigma_{\Ub_0}^{-1}})=\mathbb{B}_0+o_p(1),
			\]
			where $\mathbb{B}_0$ is given by
			\[
			\begin{split}
			\mathbb{B}_0=
			\frac{q}{T}\bigg(\sigma_\beta^4\bar\lambda_{\bSigma_{\Ub_0}^{-2}}+1+2\sigma_\beta^2\bar\lambda_{\bSigma_{\Ub_0}^{-1}}\bigg).
			\end{split}
			\]
		\end{lemma}
		The bias correction term $\mathbb{B}_0$ contains unknown parameters $\sigma_\beta^2$, which can be handled similarly as in Lemma \ref{plug-in 1}, resulting in a plug-in estimator $\hat{\mathbb{B}}_0$.  Then, compared with the definition of  $\tilde\mu$, it remains to estimate terms associated with fourth moments $\nu_4$ and $\tilde\nu_4$. Motivated by (\ref{relation}), under the noised model and null hypothesis, we have
		\[
		\begin{split}
		&\frac{1}{pq}\text{var}\bigg\{\text{tr}(\mathcal{Y}_t^\prime\mathcal{Y}_t)\bigg\}=\frac{1}{pq}\text{var}\bigg\{\text{tr}(\Vb\Xb_t^\prime\Xb_t\Vb^\prime+\sigma_\beta^2\bPhi_t^\prime\bSigma_{\Ub_0}^{-1}\bPhi_t+2\sigma_\beta\Vb\Xb_t^\prime\bSigma_{\Ub_0}^{-1/2}\bPhi_t)\bigg\}+o(1)\\
		=&\frac{\nu_4-3}{q}\sum_{j=1}^q\bSigma_{\Vb,jj}^2+2\bar\lambda_{\bSigma_{\Vb}^2}+\frac{\sigma_\beta^4(\tilde\nu_4-3)}{p}\sum_{j=1}^p(\bSigma_{\Ub_0}^{-1})_{jj}^2+2\sigma_\beta^4\bar\lambda_{\bSigma_{\Ub_0}^{-2}}+4\sigma_\beta^2\bar\lambda_{\bSigma_{\Ub_0}^{-1}}+o(1).
		\end{split}
		\]
		On the other hand, we use observations across $t$ to calculate the sample variance. By elementary calculations,
		\begin{small}
			\[
			\begin{split}
			\frac{1}{Tpq}\sum_{t=1}^T\bigg[\text{tr}(\mathcal{Y}_t^\prime\mathcal{Y}_t)-\frac{1}{T}\sum_{s=1}^T\text{tr}(\mathcal{Y}_s^\prime\mathcal{Y}_s)\bigg]^2
			=&\frac{1}{pq}\text{var}\bigg\{\text{tr}(\mathcal{Y}_t^\prime\mathcal{Y}_t)\bigg\}-\frac{1}{pq}\bigg[\frac{1}{T}\sum_{s=1}^T\big\{\text{tr}(\mathcal{Y}_s^\prime\mathcal{Y}_s)-\mathbb{E}\text{tr}(\mathcal{Y}_s^\prime\mathcal{Y}_s)\big\}\bigg]^2\\
			=&\frac{1}{pq}\text{var}\bigg\{\text{tr}(\mathcal{Y}_t^\prime\mathcal{Y}_t)\bigg\}+o_p(1).
			\end{split}
			\]
		\end{small}
		Therefore, we conclude that
		\[
		\begin{split}
		&\frac{1}{Tpq}\sum_{t=1}^T\bigg[\text{tr}(\mathcal{Y}_t^\prime\mathcal{Y}_t)-\frac{1}{T}\sum_{s=1}^T\text{tr}(\mathcal{Y}_s^\prime\mathcal{Y}_s)\bigg]^2\\
		=&\frac{\nu_4-3}{q}\sum_{j=1}^q\bSigma_{\Vb,jj}^2+2\bar\lambda_{\bSigma_{\Vb}^2}+\frac{\sigma_\beta^4(\tilde\nu_4-3)}{p}\sum_{j=1}^p(\bSigma_{\Ub_0}^{-1})_{jj}^2+2\sigma_\beta^4\bar\lambda_{\bSigma_{\Ub_0}^{-2}}+4\sigma_\beta^2\bar\lambda_{\bSigma_{\Ub_0}^{-1}}+o_p(1).
		\end{split}
		\]
		The last step is to estimate $\bar\lambda_{\bSigma_{\Vb}^2}$. Directly by Lemma \ref{tilde mu},
		\begin{equation}\label{bar lambda sigma V2}
		\frac{1}{q}\bigg\|\frac{1}{Tp}\sum_{t=1}^T\mathcal{Y}_t^\prime\mathcal{Y}_t\bigg\|_F^2-\frac{1}{p}\mathbb{B}_0-(\sigma_\beta^4\bar\lambda^2_{\bSigma_{\Ub_0}^{-1}}+2\sigma_\beta^2\bar\lambda_{\bSigma_{\Ub_0}^{-1}})=\bar\lambda_{\bSigma_{\Vb}^2}+o_p(1).
		\end{equation}
		Combining all the above results, we can obtain the next lemma, whose proof is omitted.
		\begin{lemma}
			Under the same conditions as in Theorem \ref{noise clt}, further assume that $q=o(T^2)\times \min\{\sqrt{T},\sqrt{p},\sqrt{q}\}$. Then, we have $\hat \mu-\tilde\mu=o_p(1)$ where $\hat \mu$ is given by
			\[
			\begin{split}
			\hat \mu=&\frac{p-1}{p}\bigg(\frac{p}{q}\bigg\|\frac{1}{Tp}\sum_{t=1}^T\mathcal{Y}_t^\prime\mathcal{Y}_t\bigg\|_F^2-\hat{\mathbb{B}}_0\bigg)+\frac{1}{Tpq}\sum_{t=1}^T\bigg[\text{tr}(\mathcal{Y}_t^\prime\mathcal{Y}_t)-\frac{1}{T}\sum_{s=1}^T\text{tr}(\mathcal{Y}_s^\prime\mathcal{Y}_s)\bigg]^2\\
			&+\hat\sigma_\beta^4\bar\lambda^2_{\bSigma_{\Ub_0}^{-1}}-\hat\sigma_\beta^4\bar\lambda_{\bSigma_{\Ub_0}^{-2}}.
			\end{split}
			\]
		\end{lemma}
		
		To construct testing statistic, it remains to estimate the asymptotic variance $\tilde\sigma^2$. Under the null hypothesis, the only unknown parameter is $\bar\lambda_{\bSigma_{\Vb}^2}$, which has already been considered in (\ref{bar lambda sigma V2}), by replacing the unknown parameters with plug-in estimators.}

	\section{Additional simulation results}\label{sec: addsim}
	{\subsection{Model without noise}
		For the model without noise, under the null hypothesis,} we also investigate the asymptotic normality of the proposed testing statistics, as shown in Figure \ref{fig:density}. The empirical powers under alternative hypothesis $H_{A2}$ are displayed in Table \ref{tab3}. The results are similar to those under $H_{A1}$.
	
	\begin{figure}[hbpt]
		\begin{subfigure}{.45\textwidth}
			\centering
			\includegraphics[width=7.5cm,height=7.5cm]{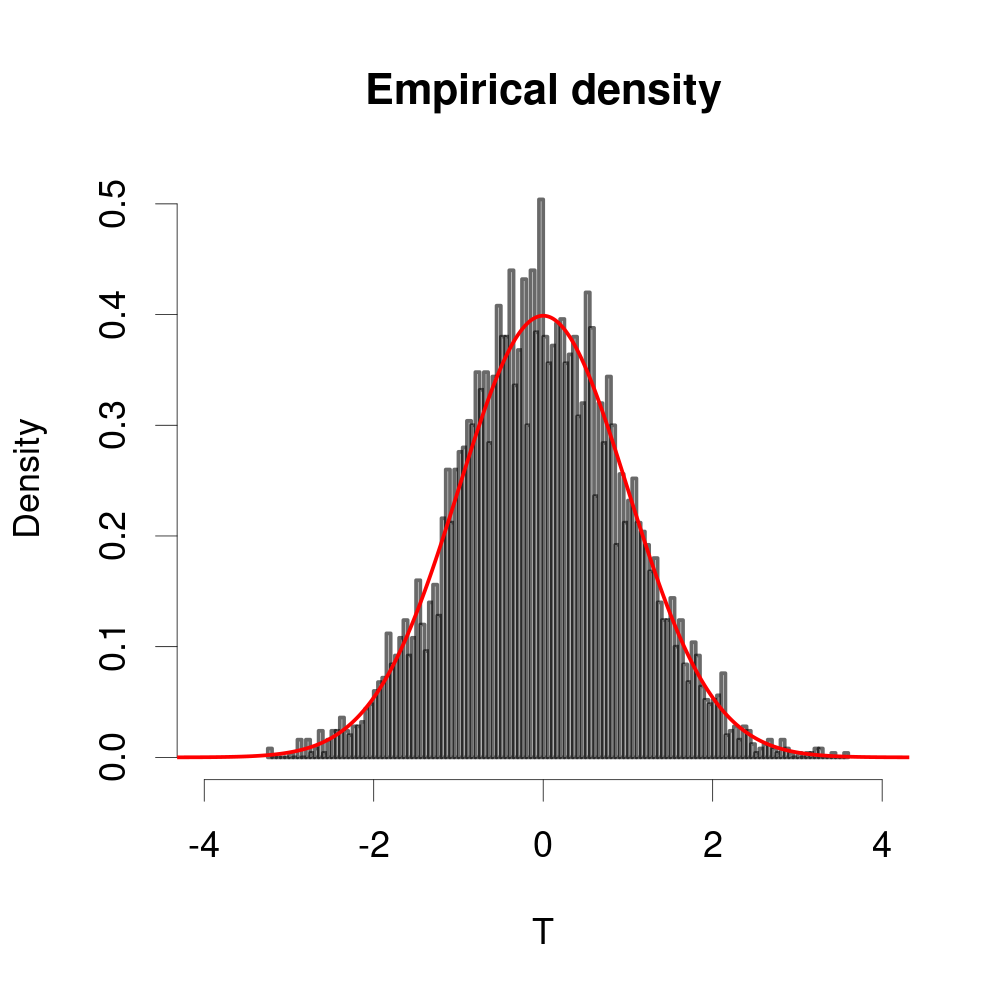}
		\end{subfigure}
		\begin{subfigure}{.45\textwidth}
			\centering
			\includegraphics[width=7.5cm,height=7.5cm]{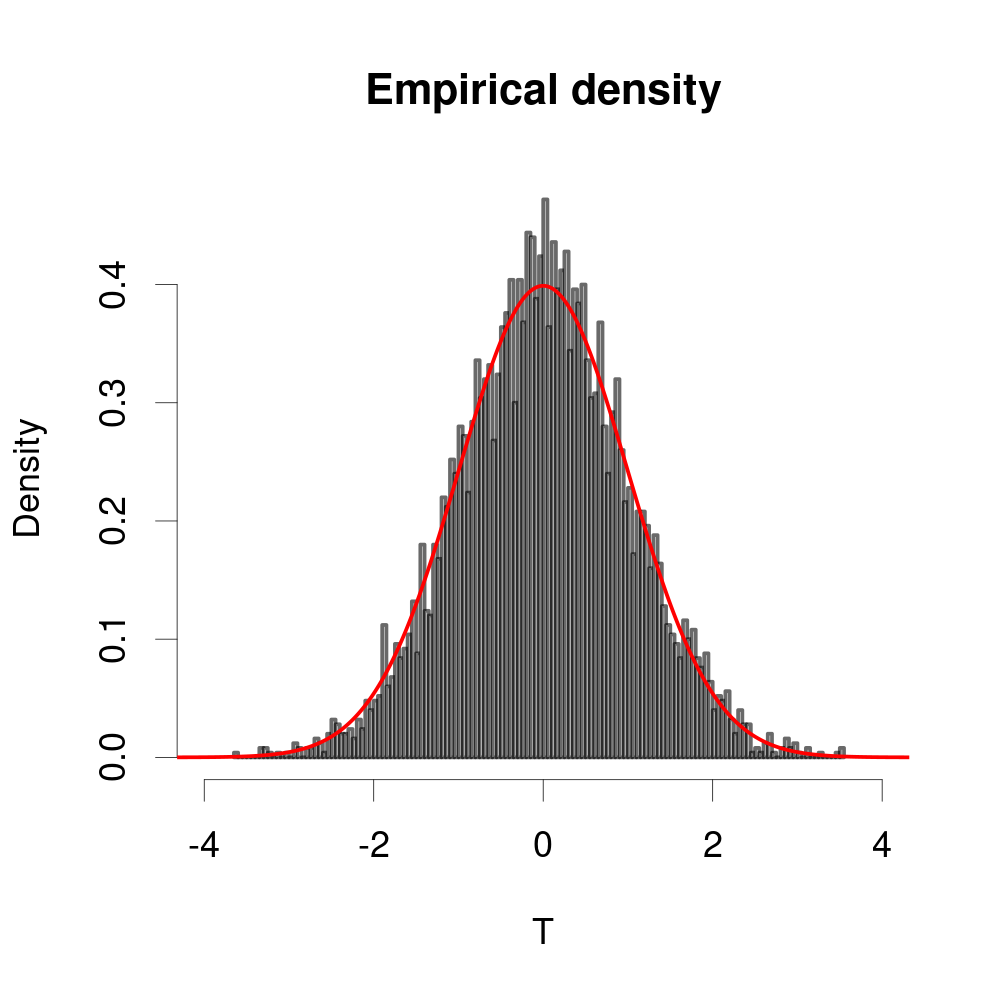}
		\end{subfigure}
		\caption{Empirical densities for the proposed testing statistics under the null hypothesis when $T=p=q=100$. The red real lines correspond to the standard normal density. Left: Normal distribution; Right: Bernoulli distribution.}\label{fig:density}
	\end{figure}
	
	\begin{table*}[hbtp]
		\begin{center}
			\small
			\addtolength{\tabcolsep}{8pt}
			\caption{Empirical powers of different testing procedures under the alternative hypothesis $H_{A2}$ with various $p$ and $T$.}\label{tab3}
			\renewcommand{\arraystretch}{1}
			\scalebox{0.9}{ 	
				\begin{tabular*}{18cm}{lllllllllllll}
					\toprule[1.2pt]
					&&&\multicolumn{3}{l}{$\alpha=0.05$}&\multicolumn{4}{l}{$\alpha=0.10$}\\\cmidrule(lr){4-6}\cmidrule(lr){7-9}
					family&$p(q)$&$T$&$\mathcal{T}_{FO}$&$\mathcal{T}_{BG}$&$\mathcal{T}_{BE}$&$\mathcal{T}_{FO}$&$\mathcal{T}_{BG}$&$\mathcal{T}_{BE}$\\\midrule[1.2pt]
					\multirow{9}{*}{Normal}&20&20&0.179&0.061&0.063&0.267&0.115&0.115
					\\
					&20&60&0.594&0.077&0.432&0.695&0.151&0.547
					\\
					&20&100&0.929&0.487&0.845&0.952&0.615&0.919
					\\\cmidrule(lr){2-9}
					&100&20&0.979&0.807&0.948&0.993&0.882&0.971
					\\
					&100&60&1&0.997&1&1&0.999&1
					\\
					&100&100&1&1&1&1&1&1\\\midrule[1.2pt]
					\multirow{9}{*}{Bernoulli}&20&20&0.137&0.064&0.086&0.223&0.119&0.161
					\\
					&20&60&0.586&0.071&0.491&0.704&0.153&0.633
					\\
					&20&100&0.973&0.587&0.947&0.983&0.729&0.976
					\\\cmidrule(lr){2-9}
					&100&20&0.980&0.818&0.962&0.992&0.895&0.984
					\\
					&100&60&1&0.999&1&1&1&1
					\\
					&100&100&1&1&1&1&1&1\\
					\bottomrule[1.2pt]
			\end{tabular*}}		
		\end{center}
	\end{table*}
	
	{
		\subsection{Dependence on the separable assumption}  
		The proposed testing procedure depends heavily on the separable model (\ref{model}). The separability may not hold even if  $\text{cov}\{\text{Vec}(\Yb)\}$ has a Kronecker product form, i.e., $\bSigma_{\Vb}\otimes \bSigma_{\Ub}$. Therefore, in this subsection. we investigate the limiting distribution of the proposed testing statistic without the separable assumption. For the data, we let 
		\[
		\text{Vec}(\Yb_t)\overset{i.i.d.}{\sim}\mathcal{MV}({\bf 0},\bSigma),
		\]
		where $\mathcal{MV}({\bf 0},\bSigma)$ stands for some multivariate distribution with mean ${\bf 0}$ and covariance matrix $\bSigma=\bSigma_{\Vb}\otimes\bSigma_{\Ub}$. $\bSigma_{\Ub}$ and $\bSigma_{\Vb}$ are the same as those in Table \ref{tab1} of the main paper. We reshape $\text{Vec}(\Yb_t)$ to $p\times q$ matrix column by column, and use the proposed procedure to calculate the testing statistic with $f(x)=x^2$. Figure \ref{fig: sensitivity} plots the empirical densities when $\Yb_t$ are from multivariate Gaussian and multivariate $t$ distribution with degree of freedom equal to 8. For better comparison, we further present the empirical Cumulative Distribution Function (CDF) of the testing statistics in Table \ref{tab: sensitivity}, together with the quantiles and CDF of standard normal distribution.  The results are based on 1000 replications with $T=p=q=60$.
		
		\begin{figure}[hbpt]
			\begin{subfigure}{.45\textwidth}
				\centering
				\includegraphics[width=7cm,height=7cm]{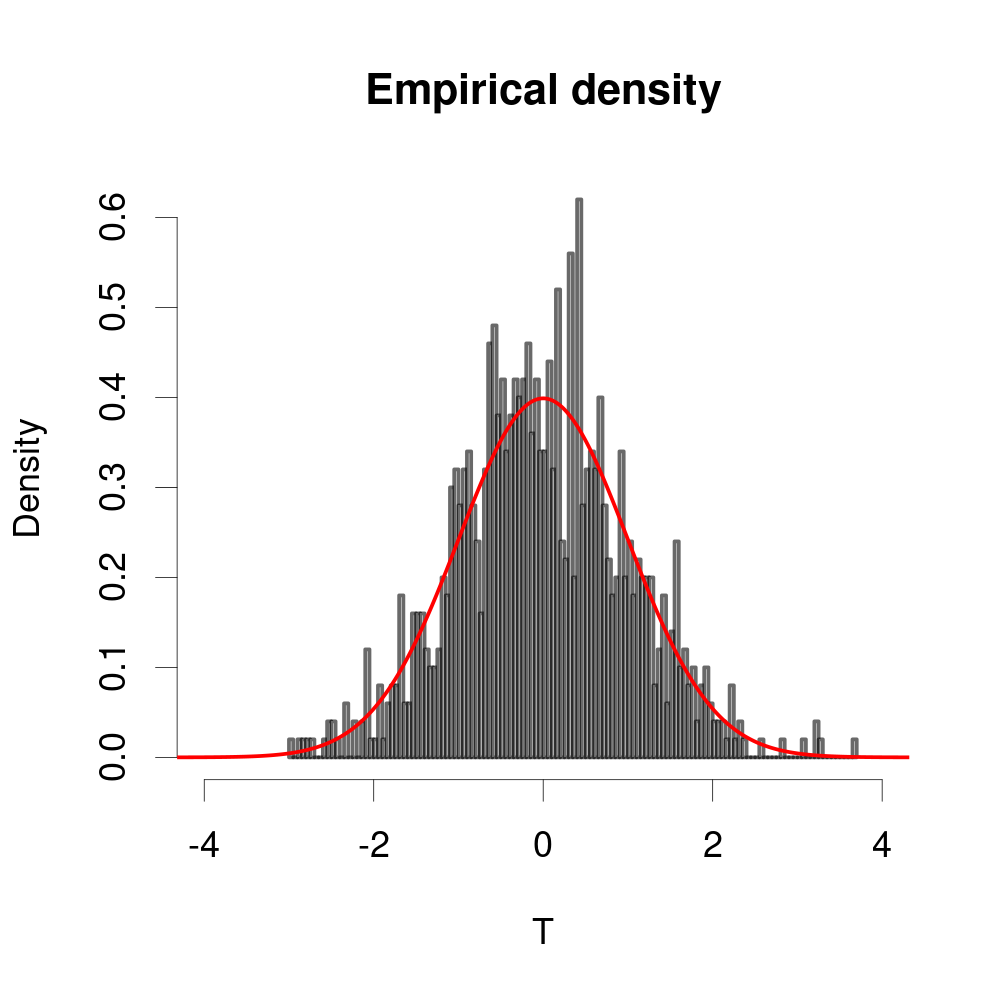}
			\end{subfigure}
			\begin{subfigure}{.45\textwidth}
				\centering
				\includegraphics[width=7cm,height=7cm]{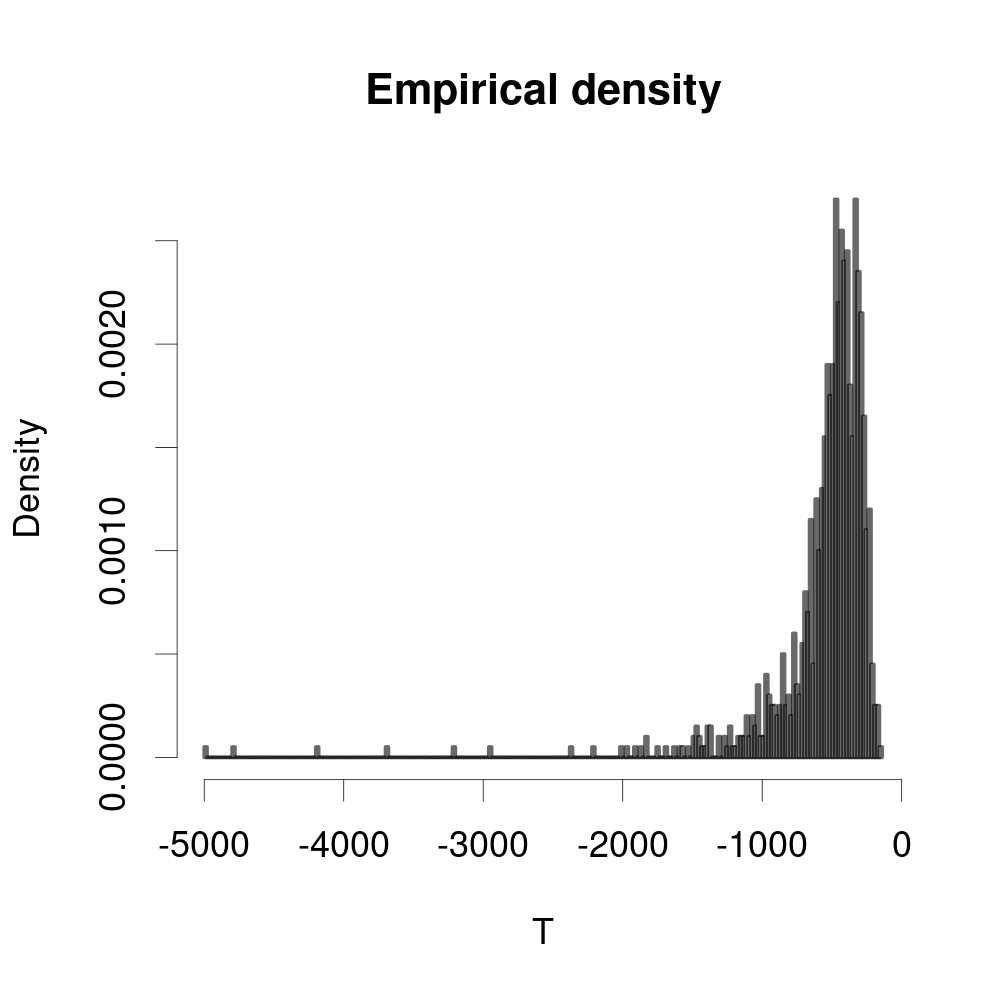}
			\end{subfigure}
			\caption{Empirical densities for the proposed testing statistics without separable assumption but $\text{cov}\{\text{Vec}(\Yb_t)\}=\bSigma_{\Vb}\otimes\bSigma_{\Ub}$. Left: multivariate Gaussian; Right: multivariate $t$. The red real line is the density of standard normal distribution.}\label{fig: sensitivity}
		\end{figure}
		
		\begin{table*}[hbtp]
			\begin{center}
				\small
				\addtolength{\tabcolsep}{0pt}
				\caption{Standard normal quantiles ($\Psi(\alpha)$) and the corresponding empirical cumulative functions of the testing statistics. }\label{tab: sensitivity}
				\renewcommand{\arraystretch}{1}
				\scalebox{1}{ 	
					\begin{tabular*}{16.5cm}{lllllllllllll}
						\toprule[1.2pt]
						$\Psi(\alpha)$&-1.645&-1.282&-1.036&-0.842&-0.675&-0.524&-0.385&-0.253&-0.126&0.000
						\\\midrule[1.2pt]
						$\alpha$&0.050&0.100&0.150&0.200&0.250&0.300&0.350&0.400&0.450&0.500
						\\
						Gaussian&0.051&	0.096&0.143&0.201&0.242&0.309&0.361&0.412&0.467&0.514
						\\
						$t$&1&1&1&1&1&1&1&1&1&1
						\\\midrule[1.2pt]
						$\Psi(\alpha)$&0.126&0.253&0.385&0.524&0.675&0.842&1.036&1.282&1.645&
						\\\midrule[1.2pt]
						$\alpha$&0.550&0.600&0.650&0.700&0.750&0.800&0.850&0.900&0.950&
						\\
						Gaussian&0.563&0.607&0.654&0.707&0.759&0.802&0.852&0.898&0.948&
						\\
						$t$&1&1&1&1&1&1&1&1&1&\\
						\bottomrule[1.2pt]
				\end{tabular*}}		
			\end{center}
		\end{table*}
		
		By Figure \ref{fig: sensitivity} and Table \ref{tab:noise power1}, the testing procedure still works for Gaussian settings but losses effectiveness under $t$ distribution. Actually, under  Gaussian settings, the assumptions of separability and Kronecker product covariance structure are equivalent. After  the transformation $\tilde\Yb_t=\bSigma_{\Ub}^{-1/2}\Yb_t$, each column of $\tilde\Yb_t$ will have covariance matrix equal to $\Ib_p$ while each row of $\tilde\Yb_t$ will have covariance matrix equal to $\bSigma_{\Vb}$. When the entries are Gaussian, we can write $\tilde\Yb_t\overset{d}{=}\Zb_t\bSigma_{\Vb}^{-1/2}$, where the entries of $\Zb_t$ are from i.i.d. $\mathcal{N}(0,1)$. When the entries are not from Gaussian, the equivalence will not hold anymore. Under such cases,  limiting distributions of linear spectral statistics will depend on not only  second-order information of the data, but also higher-order correlations between the entries of $\tilde\Yb_t$. Therefore, there can be great shifts both on the asymptotic mean and variance, which is justified by the simulation results with $t$ distribution. 
		
		As a conclusion, the proposed procedure can be used to directly test Kronecker product covariance matrix if the data are from Gaussian. Otherwise, we need to assume separability or specify higher-order correlations in order to remove the Gaussian assumption.

		\subsection{Model with noise}
		This part investigates the empirical performance of our testing procedure for the noised model, corresponding to Section \ref{sec:noise} of the main paper. We generate independent samples according to model (\ref{noise model}) with $\sigma_\alpha=\sigma_\beta=1$, while the entries of $\Xb_t$, $\varphi_t$ and $\bPhi_t$ are all from standard normal or Bernoulli distribution. Considering the identification condition in Lemma \ref{sigma beta}, we let $\Ub=\Db_1\bGamma_1$ and $\Vb=\Db_2\bGamma_2$, where $\Db_1$, $\Db_2$ are diagonal matrices and $\bGamma_1$, $\bGamma_2$ are orthogonal matrices. The spectral distributions of $\Db_1^2$ and  $\Db_2^2$ are $0.5\delta_{1-c_1}+0.5\delta_{1+c_1}$ and $0.5\delta_{1-c_2}+0.5\delta_{1+c_2}$, respectively for some constants $c_1,c_2\in(0,1)$. When $c_1,c_2$ are small, $\bSigma_{\Ub}$ and $\bSigma_{\Vb}$  will be close to identity matrices, and it's harder to identify the individual noises from the system. The simulations results below are all based on 1000 replications if not specified separately.

		\begin{figure}[hbpt]
			\begin{subfigure}{.45\textwidth}
				\centering
				\includegraphics[width=7cm,height=7cm]{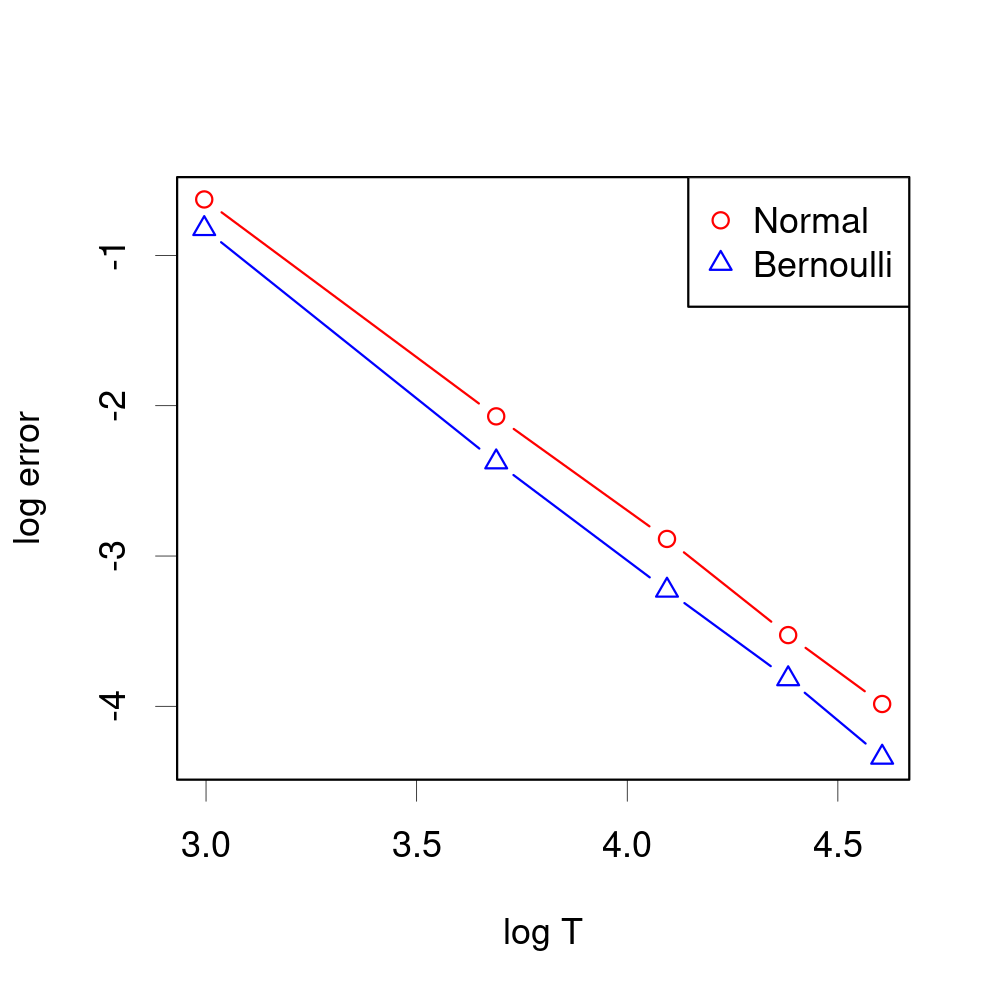}
			\end{subfigure}
			\begin{subfigure}{.45\textwidth}
				\centering
				\includegraphics[width=7cm,height=7cm]{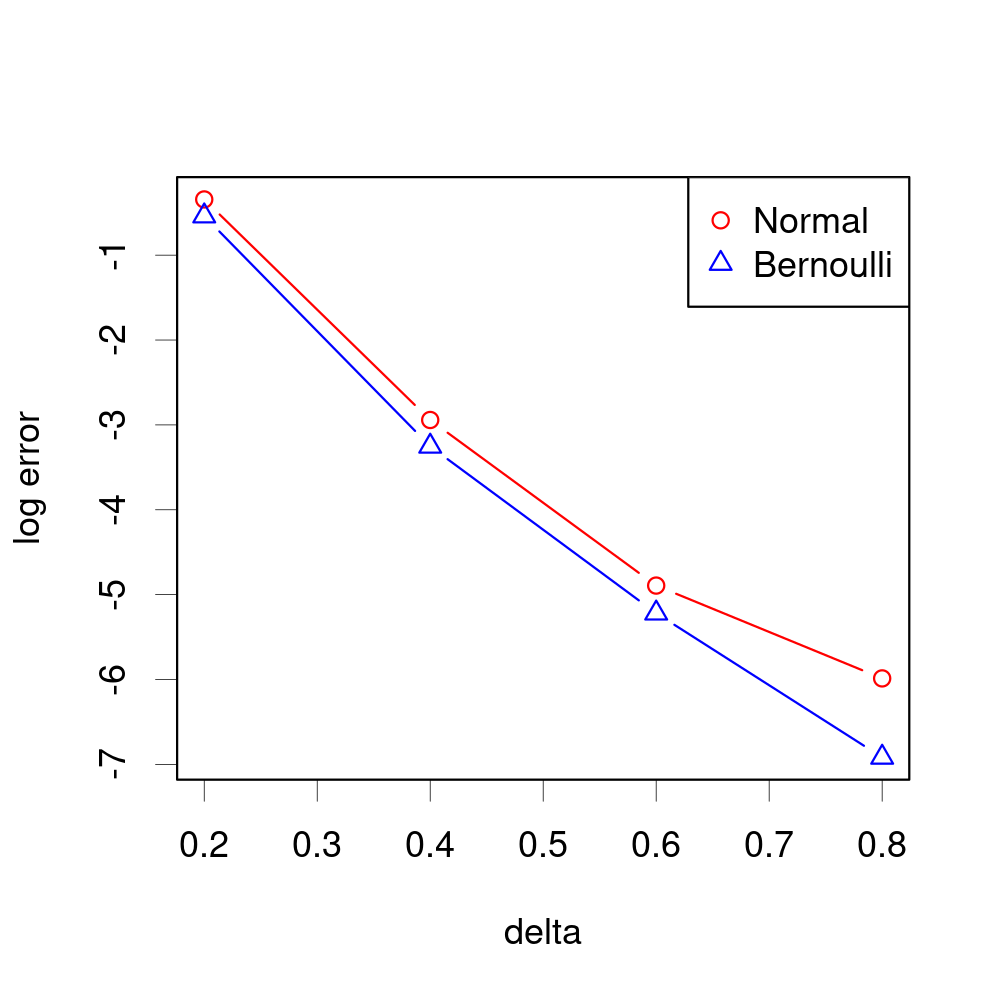}
			\end{subfigure}
			\caption{Average estimation error $\ln|\hat\sigma_\beta^2-\sigma_\beta^2|$. Left: $c_1=c_2=0.5$ while $T=p=q$. Right: $T=p=q=100$ while $c_1=c_2=\delta$. }\label{fig:sigma}
		\end{figure}
		
		\begin{figure}[hbpt]
			\begin{subfigure}{.45\textwidth}
				\centering
				\includegraphics[width=7cm,height=7cm]{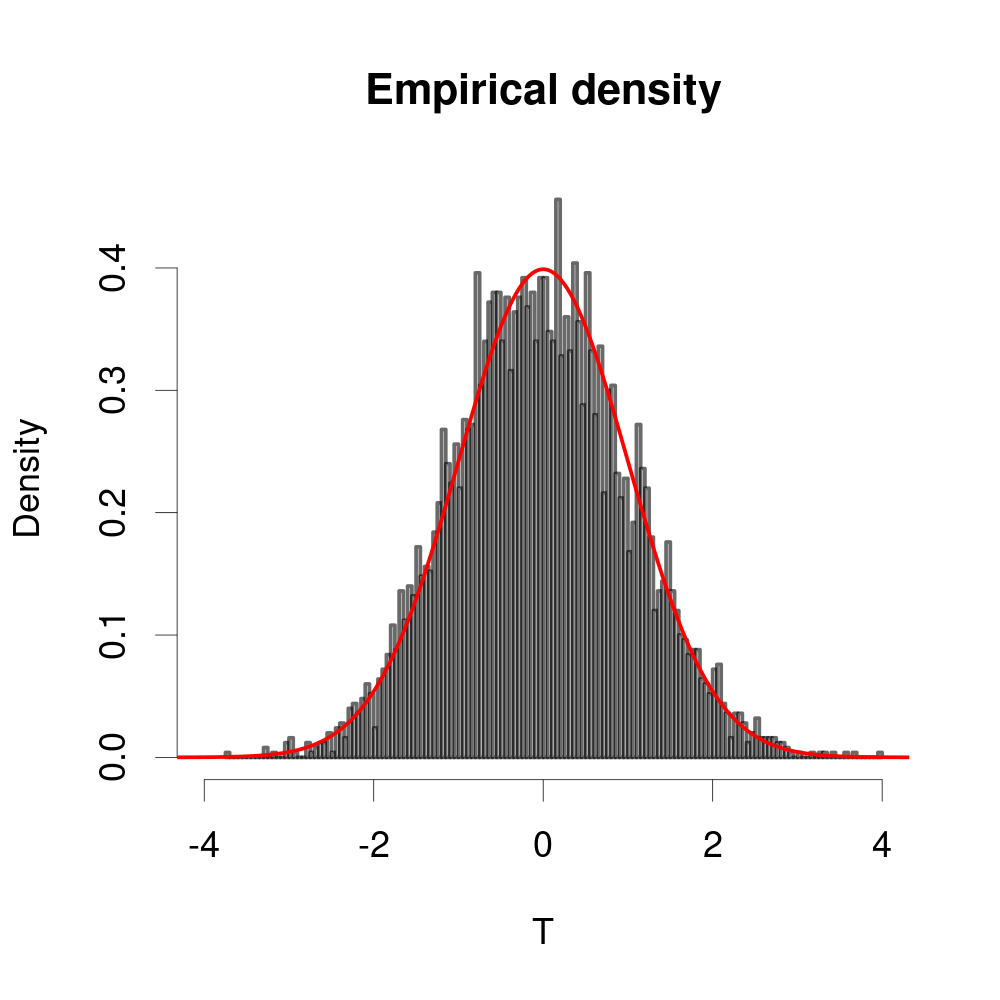}
			\end{subfigure}
			\begin{subfigure}{.45\textwidth}
				\centering
				\includegraphics[width=7cm,height=7cm]{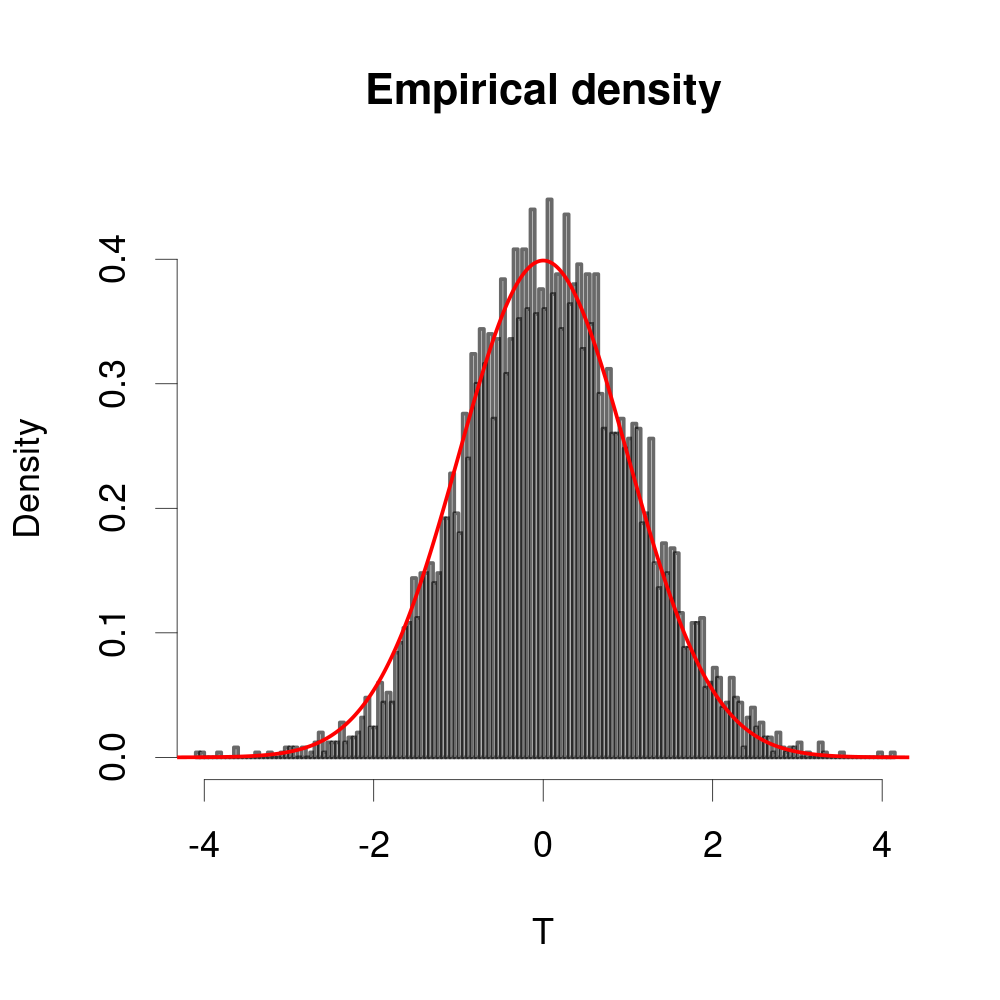}
			\end{subfigure}
			\caption{Empirical densities for the proposed testing statistics under the null hypothesis and noised model. The red real lines correspond to the standard normal density. Left: Normal distribution; Right: Bernoulli distribution. }\label{fig:noise density}
		\end{figure}
		
		Our first experiment is to study the estimation of $\sigma_{\beta}^2$, which plays an important role in the normalization step. Figure \ref{fig:sigma} shows the average absolute estimation error after log transformation ($\ln|\hat\sigma_\beta^2-\sigma_\beta^2|$)  under different settings. As shown in the left panel, the estimation error decreases as the dimension and sample size grow. The two lines have scope smaller than $-1.5$, which implies that the convergence rate of $\hat\sigma_\beta^2$ can be even faster than that in Lemma \ref{sigma beta}. The right panel indicates that $\sigma_\beta^2$ can be estimated more accurately when it's easier to identify the individual noises, as expected. 
		
		Our second experiment is to verify the asymptotic distribution in Theorem \ref{noise clt}. We fix $T=p=q=100$ and $c_1=0.5$, $c_2=0.8$. Figure \ref{fig:noise density} shows the histograms of the normalized testing statistics, $\tilde\sigma^{-1}\{\text{tr}(\bar{\mathcal{S}}^2-\tilde\mu)\}$, over 5000 replications.  It's clear that the empirical densities are close to normal density under the two distribution families. 
		
		Our third experiment is to verify the empirical sizes of the proposed testing statistics under the null. We will also investigate the effects of the estimation error for unknown parameters. When all the unknown parameters ($\sigma_\beta$, $\bSigma_{\Vb}$, $\nu_4$ and $\tilde\nu_4$) are given, we denote the results by $\mathcal{T}_{FG}$, meaning ``Fully Given''. On the contrary, if all the unknown parameters are estimated by the arguments in Section \ref{sec: estimate}, we denote the results by $\mathcal{T}_{FE}$, meaning ``Fully Estimated''. Moreover, note that the estimation of $\hat\sigma_\beta^2$ mainly affects the normalization matrix $\mathcal{E}_0$. Therefore,  we further construct a testing statistic where $\mathcal{E}_0$ is estimated with $\sigma_\beta^2$ given, denoted by $\mathcal{T}_{PG}$, meaning ``Partially Given''. The empirical sizes of the three testing statistics are presented in Table \ref{tab:noise size}.  It's seen that when $T,p,q$ are sufficiently large, our testing statistics have well-controlled sizes no matter the unknown parameters are given or estimated. Comparing $\mathcal{T}_{FE}$ and $\mathcal{T}_{PG}$, we conclude that when $T,p,q$ are small, the testing statistic may lose effectiveness  mainly because the estimated $\sigma_\beta^2$ is not accurate, which is understandable by Lemma \ref{sigma beta}.  
		
		\begin{table*}[hbtp]
			\begin{center}
				\small
				\addtolength{\tabcolsep}{8pt}
				\caption{Empirical sizes of different testing statistics under the null hypothesis and noised model,  with various $p$ and $T$.}\label{tab:noise size}
				\renewcommand{\arraystretch}{1}
				\scalebox{0.9}{ 	
					\begin{tabular*}{18cm}{lllllllllllll}
						\toprule[1.2pt]
						&&&\multicolumn{3}{l}{$\alpha=0.05$}&\multicolumn{4}{l}{$\alpha=0.10$}\\\cmidrule(lr){4-6}\cmidrule(lr){7-9}
						family&$p(q)$&$T$&$\mathcal{T}_{FG}$&$\mathcal{T}_{FE}$&$\mathcal{T}_{PG}$&$\mathcal{T}_{FG}$&$\mathcal{T}_{FE}$&$\mathcal{T}_{PG}$\\\midrule[1.2pt]
						\multirow{9}{*}{Normal}&20&20&0.076&0.632&0.059&0.140&0.699&0.103
						\\
						&20&60&0.063&0.218&0.065&0.121&0.304&0.117
						\\
						&20&100&0.049&0.124&0.047&0.109&0.194&0.101
						\\\cmidrule(lr){2-9}
						&60&20&0.061&0.196&0.055&0.110&0.275&0.118
						\\
						&60&60&0.048&0.071&0.048&0.091&0.139&0.106
						\\
						&60&100&0.041&0.054&0.044&0.082&0.113&0.096
						\\\cmidrule(lr){2-9}
						&100&20&0.078&0.115&0.070&0.146&0.188&0.101
						\\
						&100&60&0.065&0.060&0.055&0.101&0.122&0.101
						\\
						&100&100&0.041&0.047&0.050&0.094&0.097&0.085\\\midrule[1.2pt]
						\multirow{9}{*}{Bernoulli}&20&20&0.062&0.387&0.044&0.104&0.469&0.091
						\\
						&20&60&0.050&0.099&0.052&0.094&0.156&0.096
						\\
						&20&100&0.039&0.061&0.035&0.081&0.104&0.079\\
						\cmidrule(lr){2-9}
						&60&20&0.068&0.089&0.052&0.112&0.154&0.106
						\\
						&60&60&0.050&0.044&0.042&0.099&0.105&0.092
						\\
						&60&100&0.051&0.055&0.053&0.104&0.101&0.100
						\\\cmidrule(lr){2-9}
						&100&20&0.074&0.066&0.059&0.124&0.124&0.110
						\\
						&100&60&0.055&0.051&0.048&0.099&0.101&0.101
						\\
						&100&100&0.053&0.051&0.051&0.101&0.090&0.093\\
						\bottomrule[1.2pt]
				\end{tabular*}}		
			\end{center}
		\end{table*}
		
		Our last experiment is to investigate the empirical powers under alternative hypotheses. Similarly to the model without noise, we consider two types of alternative hypotheses. For $H_{A1}$, we let $\Ub=\Db_1\bGamma_1+\beta_1\bgamma\bgamma^\prime/\sqrt{p}$, where $\bgamma$ is the $p$-dimensional vector whose entries are from i.i.d. $\mathcal{N}(0,1)$. For $H_{A2}$, we let $\Ub=(\Db_1+\beta_2\Ib_p)\bGamma_1$. In the simulations, $\beta_1=\beta_2=0.1$. The empirical sizes of $\mathcal{T}_{FG}$, $\mathcal{T}_{FE}$ and $\mathcal{T}_{PG}$ under $H_{A1}$ and $H_{A2}$ are shown in Table \ref{tab:noise power1} and Table \ref{tab:noise power2}, respectively. We conclude that the empirical powers are pretty strong as long as $T,p,q$ are large.
		
		\begin{table*}[hbtp]
			\begin{center}
				\small
				\addtolength{\tabcolsep}{8pt}
				\caption{Empirical powers of different testing statistics under the noised model and alternative hypothesis $H_{A1}$,  with various $p$ and $T$.}\label{tab:noise power1}
				\renewcommand{\arraystretch}{1}
				\scalebox{0.9}{ 	
					\begin{tabular*}{18cm}{lllllllllllll}
						\toprule[1.2pt]
						&&&\multicolumn{3}{l}{$\alpha=0.05$}&\multicolumn{4}{l}{$\alpha=0.10$}\\\cmidrule(lr){4-6}\cmidrule(lr){7-9}
						family&$p(q)$&$T$&$\mathcal{T}_{FG}$&$\mathcal{T}_{FE}$&$\mathcal{T}_{PG}$&$\mathcal{T}_{FG}$&$\mathcal{T}_{FE}$&$\mathcal{T}_{PG}$\\\midrule[1.2pt]
						\multirow{4}{*}{Normal}&20&20&0.705&0.921&0.623&0.770&0.948&0.713
						\\
						&20&100&0.995&0.995&0.995&0.997&0.997&0.996
						\\
						&100&20&0.911&0.917&0.869&0.941&0.955&0.914
						\\
						&100&100&1&1&1&1&1&1\\\midrule[1.2pt]
						\multirow{4}{*}{Bernoulli}&20&20&0.160&0.635&0.157&0.227&0.703&0.236
						\\
						&20&100&0.918&0.930&0.907&0.947&0.962&0.947
						\\
						&100&20&1&1&1&1&1&1
						\\
						&100&100&1&1&1&1&1&1\\
						\bottomrule[1.2pt]
				\end{tabular*}}		
			\end{center}
		\end{table*}
		
		\begin{table*}[hbtp]
			\begin{center}
				\small
				\addtolength{\tabcolsep}{8pt}
				\caption{Empirical powers of different testing statistics under the noised model and alternative hypothesis $H_{A2}$,  with various $p$ and $T$.}\label{tab:noise power2}
				\renewcommand{\arraystretch}{1}
				\scalebox{0.9}{ 	
					\begin{tabular*}{18cm}{lllllllllllll}
						\toprule[1.2pt]
						&&&\multicolumn{3}{l}{$\alpha=0.05$}&\multicolumn{4}{l}{$\alpha=0.10$}\\\cmidrule(lr){4-6}\cmidrule(lr){7-9}
						family&$p(q)$&$T$&$\mathcal{T}_{FG}$&$\mathcal{T}_{FE}$&$\mathcal{T}_{PG}$&$\mathcal{T}_{FG}$&$\mathcal{T}_{FE}$&$\mathcal{T}_{PG}$\\\midrule[1.2pt]
						\multirow{4}{*}{Normal}&20&20&0.919&0.212&0.351&0.944&0.285&0.446
						\\
						&20&100&1&0.863&1&1&0.893&1
						\\
						&100&20&1&0.917&1&1&0.944&1
						\\
						&100&100&1&1&1&1&1&1\\\midrule[1.2pt]
						\multirow{4}{*}{Bernoulli}&20&20&0.942&0.121&0.310&0.962&0.173&0.424
						\\
						&20&100&1&0.975&1&1&0.985&1
						\\
						&100&20&1&0.980&1&1&0.989&1
						\\
						&100&100&1&1&1&1&1&1\\
						\bottomrule[1.2pt]
				\end{tabular*}}		
			\end{center}
		\end{table*}
		
	}
	
	\section{Real data analysis}\label{real data}
	In this section, we analyze a real data set. The data set consists of monthly returns of 100 portfolios, which can be freely downloaded from \url{http://mba.tuck.dartmouth.edu/pages/faculty/ken.french/data_library.html}. The portfolios can be further categorized into 10 levels of capital sizes and 10 levels of book-to-equity  ratios. Hence, it's naturally structured as $10\times 10$ matrix-variate series, also known as the Fama-French $10\times 10$ portfolios.  Considering the missing rate, we only use the data from January-1964 to  December-2020, with a total of 684 months. We impute the missing values (missing rate is $0.25\%$) by linear interpolation.  The data set was ever analyzed in \cite{wang2019factor} and \cite{yu2020projected}. Following their preprocessing procedures, we first subtract the monthly market excess returns and then standardize the series one by one. Figure \ref{fig1} is an illustration of the standardized series.
	
	\begin{figure}[hbpt]
		\centering
		\includegraphics[width=15cm,height=15cm]{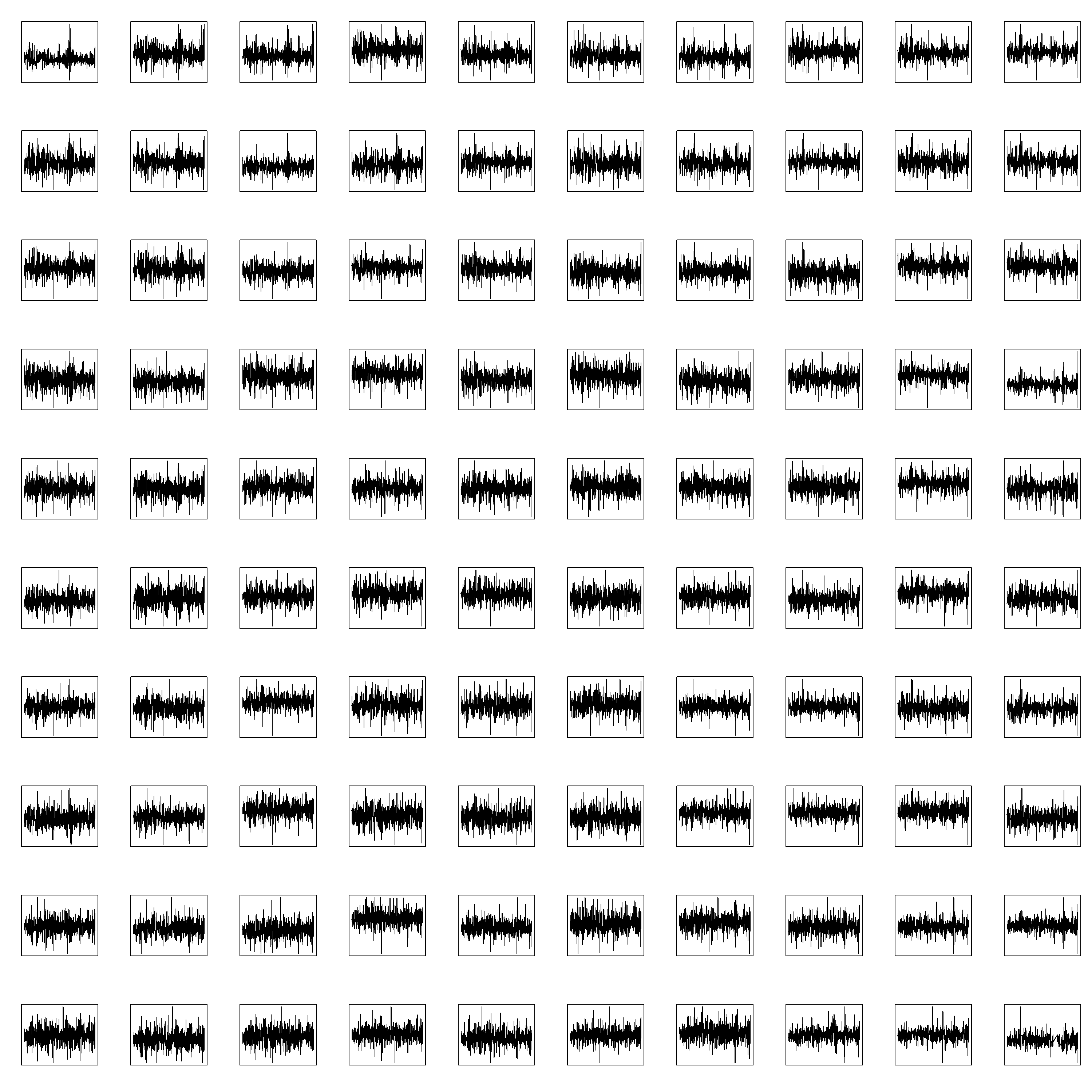}
		
		\caption{Standardized Fama-French $10\times 10$ portfolio series. Each series is associated with  a level of capital size and a level of book-to-equity ratio.  Top to bottom: lower book-to-equity  ratio to higher. Left to right: small capital size to larger.}\label{fig1}
	\end{figure}
	
	{We start with the model (\ref{model}) without noise, which assumes a Kronecker product covariance structure for the $10\times 10$ matrices. We first provide some initial guesses of the row-row and column-column covariance matrices $\bSigma_{\Ub}$ and $\bSigma_{\Vb}$. The target of this study is to select the most convincing one.
		Since the sample size is much larger than dimensions, it's natural to consider the sample versions
		\[
		\bSigma_{\Ub_1}=\frac{1}{6840}\sum_{t=1}^{684}\Yb_t\Yb_t^\prime,\quad  \bSigma_{\Vb_1}=\frac{1}{6840}\sum_{t=1}^{684}\Yb_t^\prime\Yb_t,
		\]
		which are actually unbiased estimators respectively under the identification conditions $\text{tr}(\bSigma_{\Vb})/10=1$ and $\text{tr}(\bSigma_{\Ub})/10=1$.} Although the two identification conditions may not hold simultaneously, it has no effects on this study since we will test $\bSigma_{\Ub}$ and $\bSigma_{\Vb}$ separately. We plot  $\bSigma_{\Ub_1}$, $\bSigma_{\Vb_1}$, their inverse matrices and eigenvalues in Figure \ref{fig2}.
	
	\begin{figure}[hbpt]
		\begin{subfigure}{.3\textwidth}
			\centering
			\includegraphics[width=4cm,height=4cm]{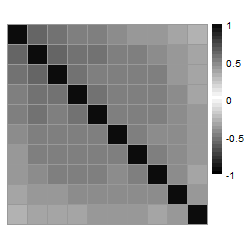}
			\subcaption{heatmap of $\hat\bSigma_{\Ub_1}$}
		\end{subfigure}
		\begin{subfigure}{.3\textwidth}
			\centering
			\includegraphics[width=4cm,height=4cm]{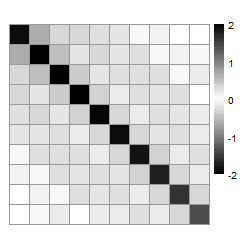}
			\subcaption{heatmap of $\hat\bSigma_{\Ub_1}^{-1}$}
		\end{subfigure}
		\begin{subfigure}{.3\textwidth}
			\centering
			\includegraphics[width=4cm,height=4cm]{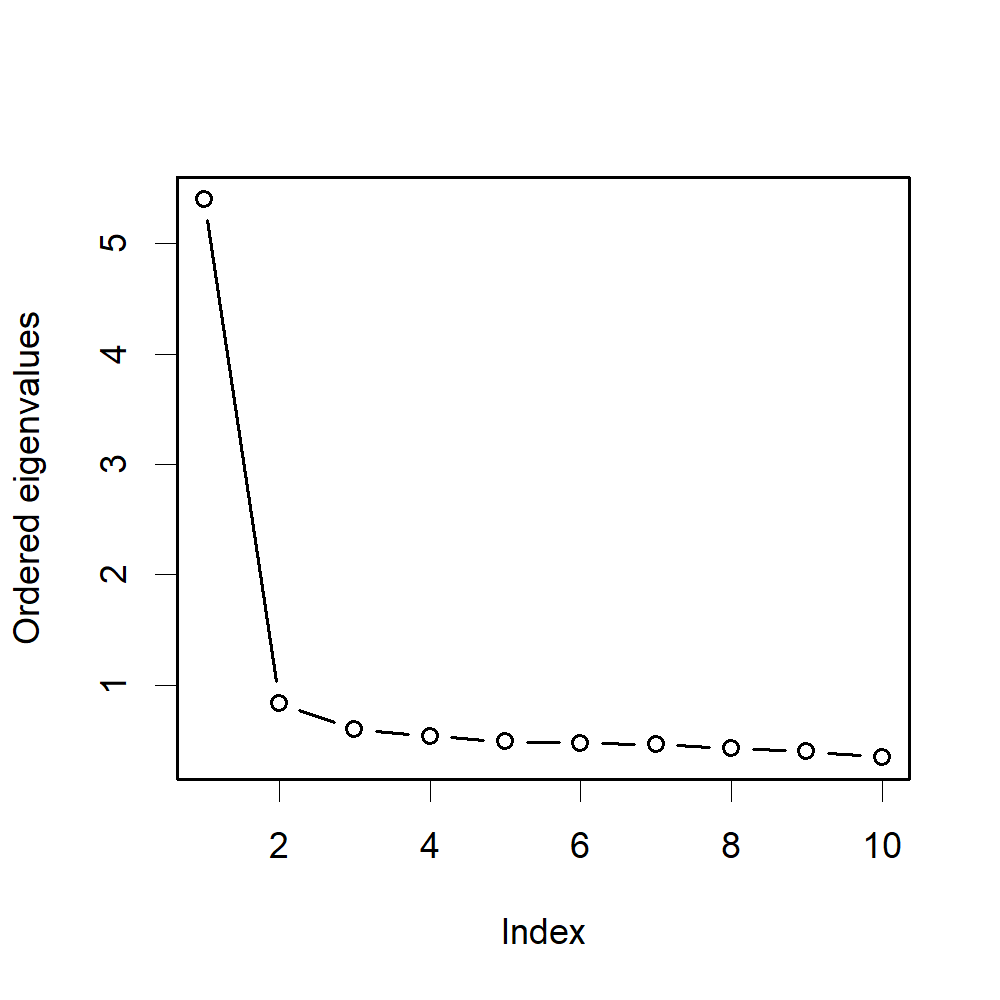}
			\subcaption{Eigenvalues of $\hat\bSigma_{\Ub_1}$}
		\end{subfigure}
		
		\begin{subfigure}{.3\textwidth}
			\centering
			\includegraphics[width=4cm,height=4cm]{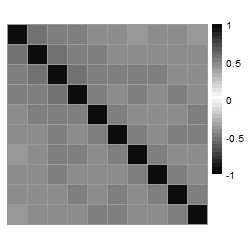}
			\subcaption{heatmap of $\hat\bSigma_{\Vb_1}$}
		\end{subfigure}
		\begin{subfigure}{.3\textwidth}
			\centering
			\includegraphics[width=4cm,height=4cm]{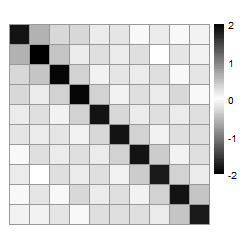}
			\subcaption{heatmap of $\hat\bSigma_{\Vb_1}^{-1}$}
		\end{subfigure}
		\begin{subfigure}{.3\textwidth}
			\centering
			\includegraphics[width=4cm,height=4cm]{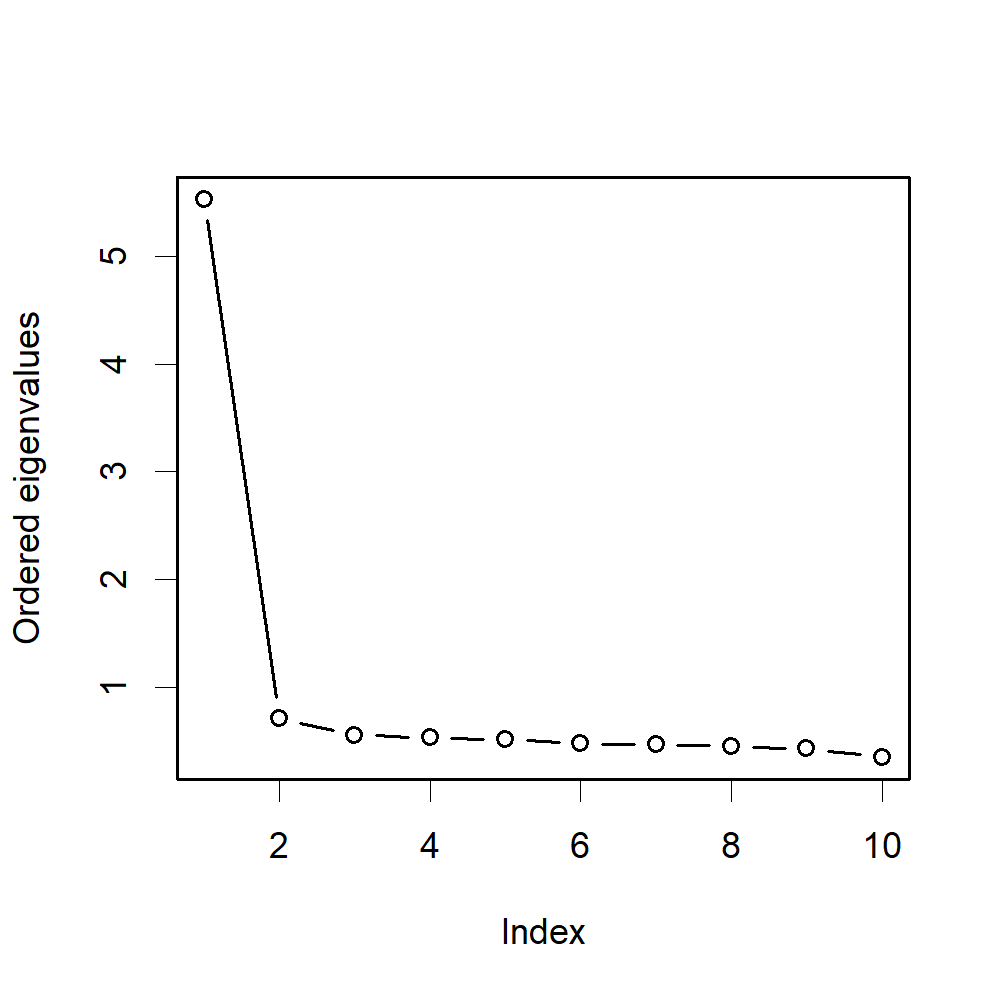}
			\subcaption{Eigenvalues of $\hat\bSigma_{\Vb_1}$}
		\end{subfigure}
		\caption{Sample covariances  $\hat\bSigma_{\Ub_1}$, $\hat\bSigma_{\Vb_1}$, their inverses and corresponding eigenvalues. }\label{fig2}
	\end{figure}
	
	By Figure \ref{fig2}, the largest eigenvalues of $\bSigma_{\Ub_1}$ and $\bSigma_{\Vb_1}$ are well separated from the others. This indicates a potential low-rank plus noise structure, which is motivated by the matrix-valued factor models ever studied in \cite{wang2019factor} and \cite{yu2020projected}. Then, we construct factor-based {guess} of $\bSigma_{\Ub_2}$ by
	\[
	\bSigma_{\Ub_2}=\sum_{j=1}^k\lambda_j(\bSigma_{\Ub_1})\bxi_j\bxi_j^\prime+\text{diag}(\bPsi),
	\]
	where $\bxi_j$ is the $j$-th eigenvector of $\bSigma_{\Ub_1}$ and $\bPsi=\bSigma_{\Ub_1}-\sum_{j=1}^k\lambda_j(\bSigma_{\Ub_1})\bxi_j\bxi_j^\prime$. Motivated by the ``POET'' estimator in \cite{fan2013large}, $\text{diag}(\bPsi)$ can be replaced by penalized versions, e.g., using the hard-thresholding technique in \cite{bickel2008covariance}. In this example, we use $\text{diag}(\bPsi)$ for simplicity. We let $k=5$ so that the leading factors can explain over $75\%$ variance of $\bSigma_{\Ub_1}$. Similar procedure leads to a modified version of $\bSigma_{\Vb}$, denoted by $\bSigma_{\Vb_2}$. Figure \ref{fig2} also shows that  $\bSigma_{\Ub_1}^{-1}$ and $\bSigma_{\Vb_1}^{-1}$ seem to be sparse. Hence, motivated by \cite{zhou2014gemini}, it's reasonable to consider the graphical-lasso estimator, defined by
	\[
	\bOmega_{gl}=\arg\min_{\bOmega> 0}\bigg(\text{tr}(\bSigma_{\Ub_1}\bOmega)-\log\det(\bOmega)+\lambda|\bOmega|_{1,\text{off}}\bigg),
	\]
	where $\lambda$ is a tuning parameter. We let $\lambda=0.05$ in this example for comparison. Then, the graphical-lasso based guess of $\bSigma_{\Ub}$ is $\bSigma_{\Ub_3}:=\bOmega_{gl}^{-1}$. Similar technique can be applied to $\bSigma_{\Vb_1}$ to obtain $\bSigma_{\Vb_3}$.
	
	Currently we have three preliminary guesses for $\bSigma_{\Ub}$ and $\bSigma_{\Vb}$. In portfolio theory, the covariance matrix is critical to finding the risk-minimization investing strategy, which is equivalent to solving the optimization problem
	\[
	\bw_{op}=\arg\min_{\|\bw\|=1}\text{var}\bigg(\bw^\prime\text{Vec}(\Yb)\bigg)=\frac{\bSigma_{\Yb}^{-1}{\bf 1}}{{\bf 1}^\prime\bSigma_{\Yb}^{-1}{\bf 1}},
	\]
	where $\bSigma_{\Yb}=\text{Cov}\{\text{Vec}(\Yb)\}$ and ${\bf 1}$ is a vector with entries all being 1. $\hat\bw$ determines the best weights assigned to $\text{Vec}(\Yb)$.  Then, based on the preliminary guesses, taking $\bSigma_{\Yb}=\bSigma_{\Vb}\otimes \bSigma_{\Ub}$, we calculate the risk-minimization portfolio weights and report them in Figure \ref{fig3}. $\bw_{op}$ is restructured into $10\times 10$ matrices for better representation.
	
	\begin{figure}[hbpt]
		\begin{subfigure}{.42\textwidth}
			\centering
			\scalebox{0.7}{ \begin{tabular*}{9cm}{|c|c|c|c|c|c|c|c|c|c|}
					\cline{1-10}
					1.9&1.4&1.2&1.1&1.5&1.5&1.6&1.8&1.6&1.9\\\cline{1-10}
					0.8&0.5&0.5&0.4&0.6&0.6&0.6&0.7&0.6&0.7\\\cline{1-10}
					0.5&0.3&0.3&0.3&0.4&0.4&0.4&0.4&0.4&0.5\\\cline{1-10}
					0.9&0.6&0.5&0.5&0.7&0.7&0.7&0.8&0.7&0.9\\\cline{1-10}
					0.7&0.5&0.5&0.4&0.6&0.6&0.6&0.7&0.6&0.7\\\cline{1-10}
					0.9&0.6&0.5&0.5&0.6&0.7&0.7&0.8&0.7&0.8\\\cline{1-10}
					1.0&0.7&0.6&0.6&0.7&0.7&0.8&0.9&0.8&0.9\\\cline{1-10}
					1.3&0.9&0.8&0.8&1.0&1.0&1.1&1.3&1.1&1.3\\\cline{1-10}
					1.8&1.3&1.1&1.1&1.4&1.4&1.5&1.7&1.5&1.8\\\cline{1-10}
					2.7&1.9&1.7&1.6&2.1&2.1&2.3&2.6&2.3&2.7\\\cline{1-10}
			\end{tabular*}}
			\subcaption{$\bw_{op}$ from $\bSigma_{\Vb_1}\otimes\bSigma_{\Ub_1}$}
		\end{subfigure}
		\begin{subfigure}{.56\textwidth}
			\centering
			\scalebox{0.7}{ \begin{tabular*}{9cm}{|c|c|c|c|c|c|c|c|c|c|}
					\cline{1-10}
					3.5&1.9&0.8&0.9&3.1&1.8&2.5&3.3&1.9&3.0\\\cline{1-10}
					0.5&0.3&0.1&0.1&0.5&0.3&0.4&0.5&0.3&0.5\\\cline{1-10}
					-0.2&-0.1&0.0&0.0&-0.2&-0.1&-0.1&-0.2&-0.1&-0.2\\\cline{1-10}
					0.5&0.3&0.1&0.1&0.4&0.3&0.4&0.5&0.3&0.4\\\cline{1-10}
					-0.1&-0.1&0.0&0.0&-0.1&-0.1&-0.1&-0.1&-0.1&-0.1\\\cline{1-10}
					1.8&1.0&0.4&0.5&1.6&0.9&1.3&1.7&1.0&1.5\\\cline{1-10}
					1.0&0.5&0.2&0.3&0.9&0.5&0.7&1.0&0.6&0.9\\\cline{1-10}
					2.4&1.3&0.6&0.6&2.1&1.2&1.7&2.3&1.3&2.0\\\cline{1-10}
					2.3&1.2&0.5&0.6&2.0&1.2&1.6&2.2&1.3&1.9\\\cline{1-10}
					3.8&2.0&0.9&0.9&3.3&1.9&2.6&3.5&2.0&3.2\\\cline{1-10}
			\end{tabular*}}
			\subcaption{$ \bw_{op}$ from $\bSigma_{\Vb_2}\otimes\bSigma_{\Ub_2}$}
		\end{subfigure}
		
		\begin{subfigure}{.42\textwidth}
			\centering
			\scalebox{0.7}{ \begin{tabular*}{9cm}{|c|c|c|c|c|c|c|c|c|c|}
					\cline{1-10}
					1.6&1.2&1.1&1.1&1.3&1.3&1.4&1.6&1.4&1.6\\\cline{1-10}
					0.8&0.6&0.6&0.5&0.7&0.7&0.7&0.8&0.7&0.8\\\cline{1-10}
					0.6&0.4&0.4&0.4&0.5&0.5&0.5&0.6&0.5&0.6\\\cline{1-10}
					0.9&0.7&0.6&0.6&0.7&0.7&0.8&0.8&0.8&0.9\\\cline{1-10}
					0.8&0.6&0.5&0.5&0.6&0.6&0.7&0.8&0.7&0.8\\\cline{1-10}
					0.9&0.7&0.6&0.6&0.7&0.7&0.8&0.9&0.8&0.9\\\cline{1-10}
					1.0&0.7&0.7&0.6&0.8&0.8&0.8&0.9&0.8&1.0\\\cline{1-10}
					1.3&0.9&0.9&0.8&1.0&1.0&1.1&1.2&1.1&1.3\\\cline{1-10}
					1.7&1.2&1.1&1.1&1.3&1.3&1.4&1.6&1.4&1.7\\\cline{1-10}
					2.4&1.8&1.6&1.6&1.9&2.0&2.1&2.3&2.1&2.4\\\cline{1-10}
			\end{tabular*}}
			
			\subcaption{$ \bw_{op}$ from $\bSigma_{\Vb_3}\otimes\bSigma_{\Ub_3}$}
		\end{subfigure}
		\begin{subfigure}{.56\textwidth}
			\centering
			\includegraphics[width=7.2cm,height=4.5cm]{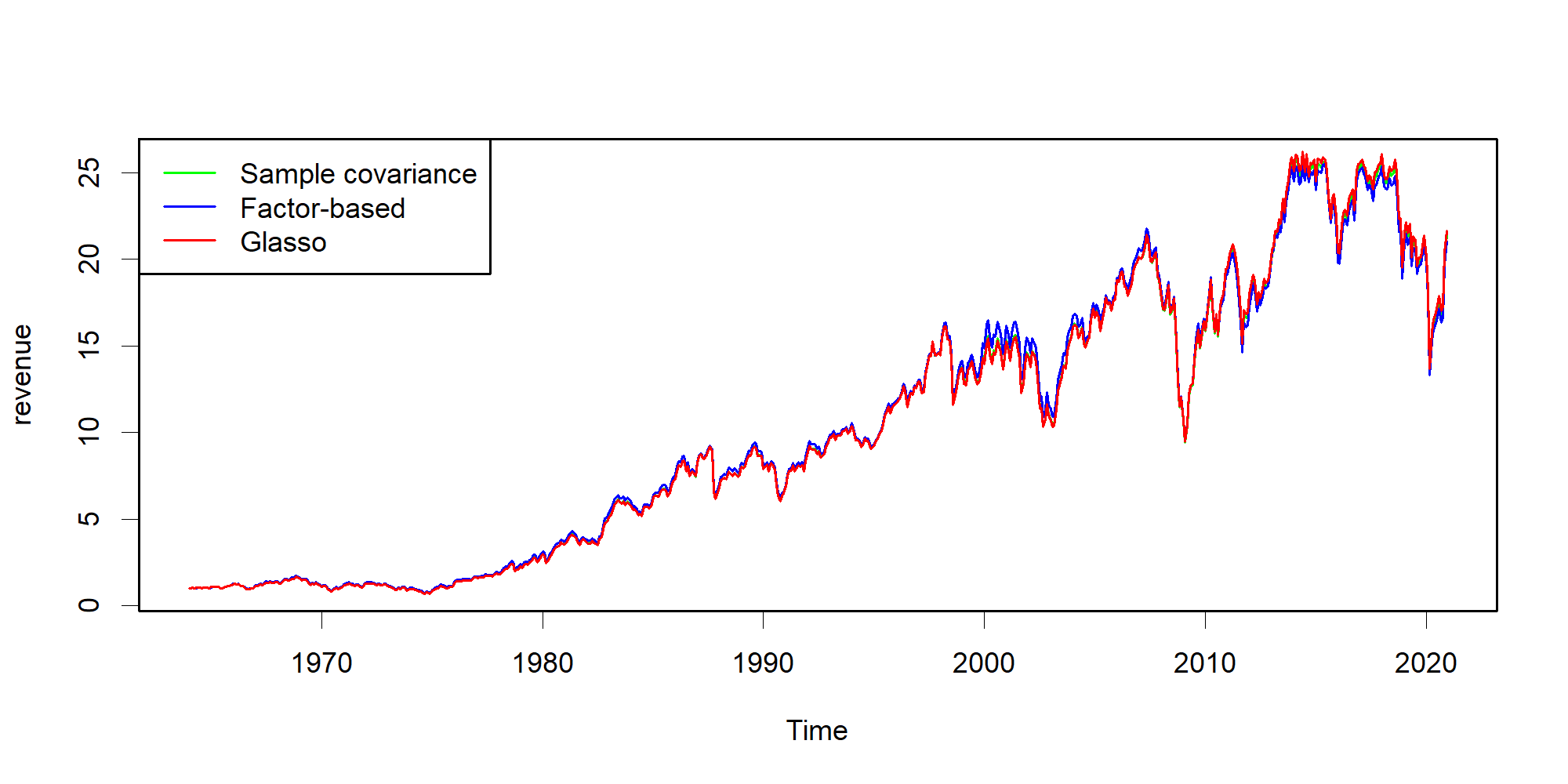}
			\subcaption{historical revenue}
		\end{subfigure}
		
		\caption{(a) to (c): optimal portfolio weights with $\bSigma_{\Yb}$ from three estimators. All the values have been multiplied by 100 for better representation. (d): historical  revenue of corresponding portfolio weights regardless of liquidity and transaction fees.}\label{fig3}
	\end{figure}
	
	It's seen from Figure \ref{fig3} that the optimal portfolio weights can vary a lot when $\bSigma_{\Yb}$ is estimated by different methods. For the factor-based $\hat\bw_2$, it assigns negative weights to some portfolios, referring to ``short'' strategy which are contradictory to the other two methods. It's then of importance to determine which strategy is more convincing. Empirically, the historical  revenues of different strategies can provide a criterion to select the best one. However, in this example, the historical revenue curve of the three strategies are so close to each other, as displayed in Figure \ref{fig3} (d). In the below, we apply the proposed testing procedure to select the most convincing strategy, i.e., to select the covariance matrix nearest to the truth.
	
	{We use the testing algorithm \ref{alg1} to test $\bSigma_{\Ub}=\bSigma_{\Ub_0}$, where $\bSigma_{\Ub_0}$ is any of the above preliminary guess. We clarify here that when saying ``testing'', we regard  $\bSigma_{\Ub_0}$ as a given constant matrix, but not an estimator. The unknown parameters in the testing statistics are  estimated by the procedure illustrated in Section 4.3.  This leads to the testing statistics $\mathcal{T}_1=30.15$, $\mathcal{T}_2=184.23$ and $\mathcal{T}_3=129.31$, associated with $\bSigma_{\Ub_1}$, $\bSigma_{\Ub_2}$ and $\bSigma_{\Ub_3}$ respectively.} That is, all the preliminary guesses of $\bSigma_{\Ub}$ are rejected significantly, which is common in real examples because the model can be misspecified. However, a larger deviation from the interval $[-2,2]$ usually indicates more significant rejection. Hence, in this example, it's more recommended to believe in the sample covariance matrix $\bSigma_{\Ub_1}$. This is partially rationalized by the large sample size of this data set. Actually, if we select the tuning parameter $\lambda$ of the graphical lasso estimator by cross validation, it always leads to $\lambda=0$, which is exactly associated with the sample covariance matrix. {It's also standard in the literature to test where the covariance matrix is an scaled identity matrix, i.e., $\bSigma_{\Ub}=\Ib\times\text{tr}(\bSigma_{\Ub_1})/10$. This leads to a testing statistic $\mathcal{T}_4=212.69$, which deviates most significantly  from the interval $[-2,2]$. In other words, the preliminary guesses  $\bSigma_{\Ub_1}$, $\bSigma_{\Ub_2}$ and $\bSigma_{\Ub_3}$ work at least better than a simple guess of identity matrix. 
		For $\bSigma_{\Vb}$, a parallel procedure is applied to the series $\{\Yb_t^\prime\}_{t=1}^{684}$, which leads to testing statistics $\mathcal{T}_1=35.67$, $\mathcal{T}_2=280.70$ and $\mathcal{T}_3=137.96$. The hypothesis $\bSigma_{\Vb}=\Ib\times\text{tr}(\bSigma_{\Vb_1})/10$ will lead to $\mathcal{T}_4=223.64$. Therefore, the sample covariance matrix estimator $\bSigma_{\Vb_1}$ is still the most convincing one in this example.}

	{Next, we consider the noised model (\ref{noise model}) for this real example. Using the procedure in Section \ref{sec:noise}, we have $\hat\sigma_{\beta}^2=0.47$. Then, we remove the common noise from the model and obtain $\{\hat\Yb_t\}$.  Left and middle panels of Figure \ref{fig: eig noise} plot the eigenvalues of the sample covariance matrices $\bSigma_{\Ub_1}$ and $\bSigma_{\Vb_1}$ from $\{\Yb_t\}$. Compared with Figure \ref{fig2}, the eigenvalues now decay more smoothly, and no significantly spiked eigenvalues are found.

		\begin{figure}[hbpt]
			\begin{subfigure}{.3\textwidth}
				\centering
				\includegraphics[width=4.8cm,height=5cm]{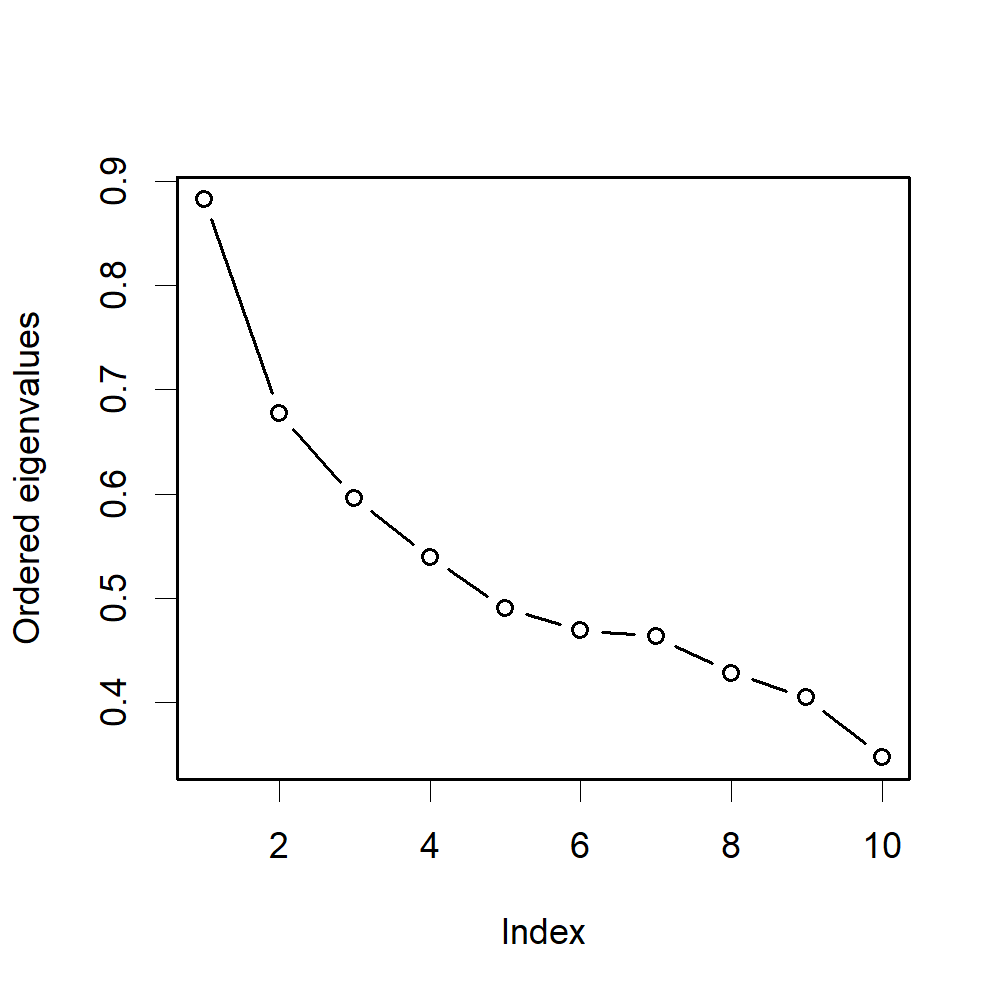}
			\end{subfigure}
			\begin{subfigure}{.3\textwidth}
				\centering
				\includegraphics[width=4.8cm,height=5cm]{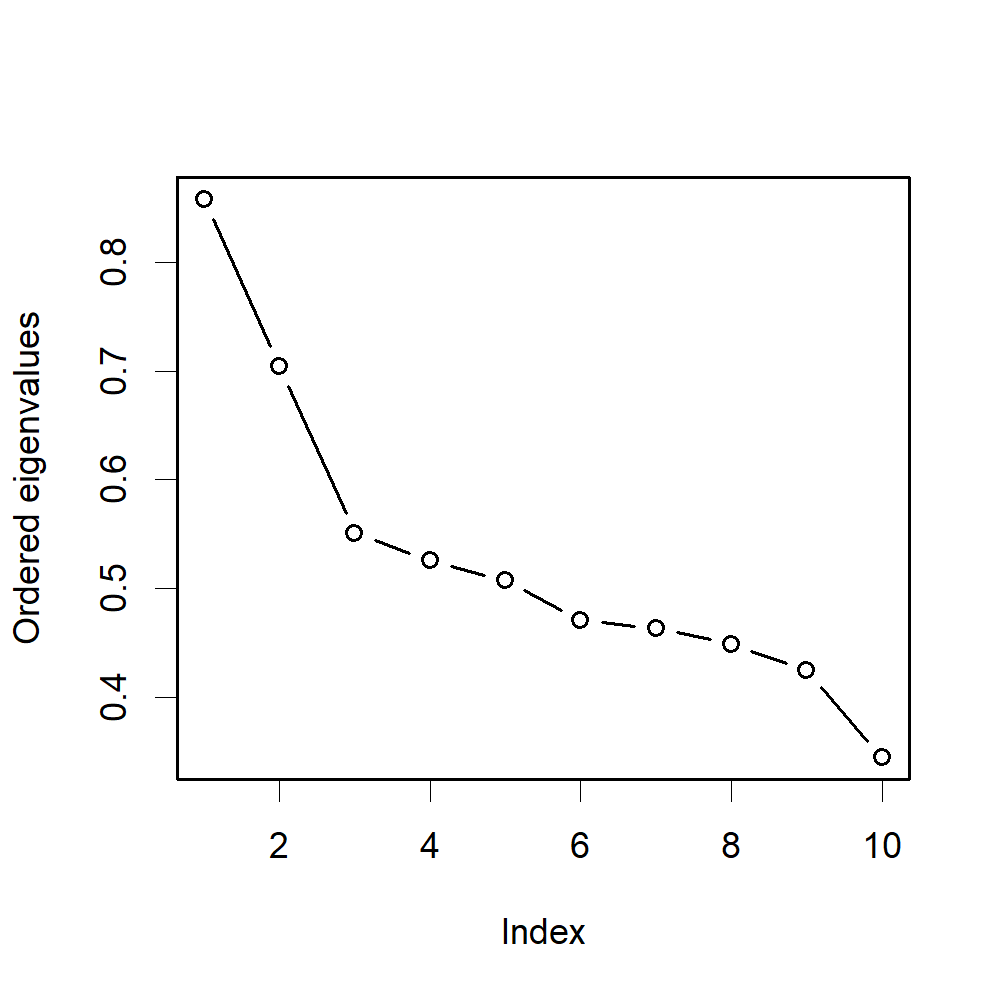}
			\end{subfigure}
			\begin{subfigure}{.3\textwidth}
				\centering
				\includegraphics[width=4.8cm,height=5cm]{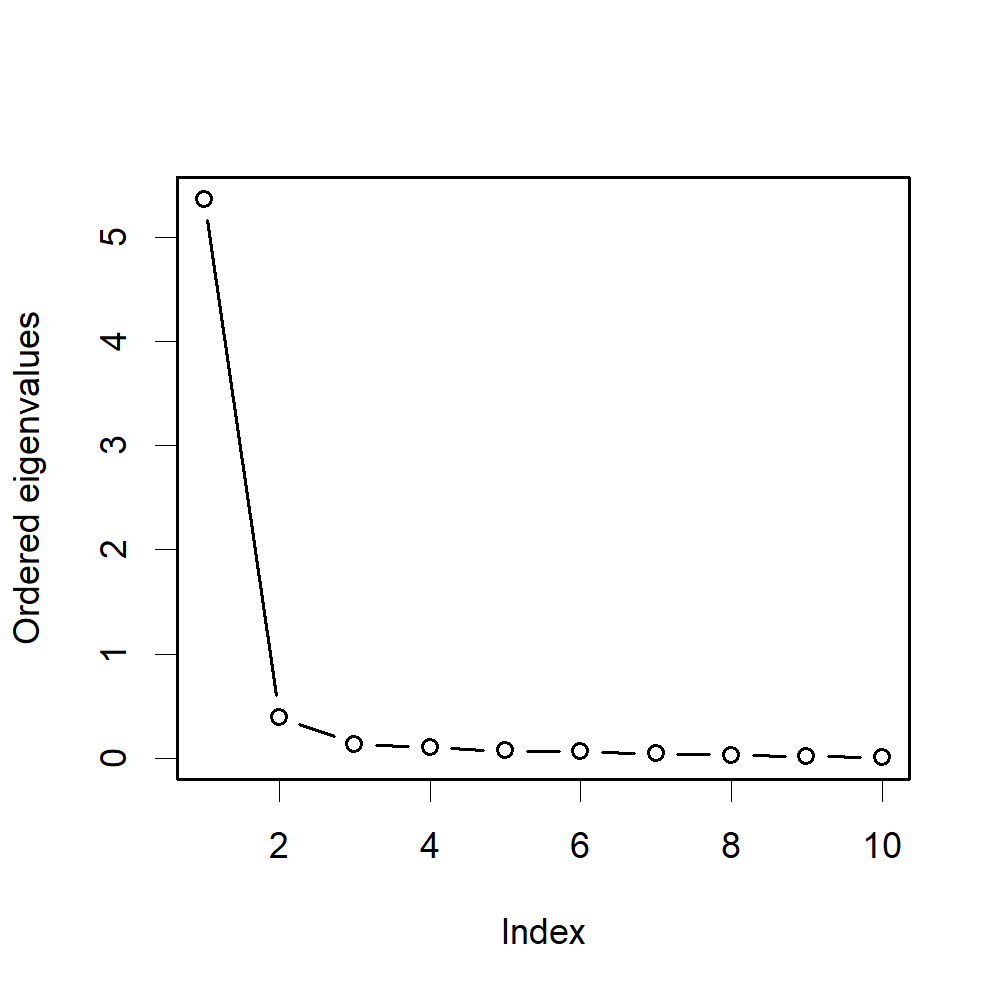}
			\end{subfigure}
			\caption{Eigenvalues of the sample covariance matrices $\bSigma_{\Ub_1}$, $\bSigma_{\Vb_1}$ from $\{\Yb_t\}$, and $\text{Re}(\Omega)$, under noised model.}\label{fig: eig noise}
		\end{figure}
		
		In the next step, we calculate the reshaped matrix $\text{Re}(\hat\Omega)$ according to the definition in the main paper. Right panel of Figure \ref{fig: eig noise} shows the eigenvalues of $\text{Re}(\hat\Omega)$. According to our proof of Lemma \ref{sigma beta} and the arguments in Section \ref{sec8}, $\text{Re}(\hat\Omega)$ should have two spiked eigenvalues if $\sigma_{\beta}^2\ge c$ for some constant $c>0$, which is not observed in the  figure. Therefore, we conclude that either $\sigma_{\beta}^2=0$ or the identification condition in Lemma \ref{sigma beta}  does not hold. Under such cases, it's more meaningful to test $\bSigma_{\Ub}$ and $\sigma_{\beta}^2$ jointly. In this example, we test the null hypothesis that $\bSigma_{\Ub}=\bSigma_{\Ub_1}$ and $\sigma_{\beta}^2=0$ hold simultaneously.  This leads to testing statistic equal to $-23.58$. Compared with the model without noise, the new testing statistic deviates less from the interval $[-2,2]$, which benefits mainly from the existence of common noise.  A parallel procedure for $\bSigma_{\Vb}$ leads to testing statistic equal to $-34.55$ and similar conclusion holds.

	}

	\section{Outline of the proof} \label{seca}
	
	{
		
		The major target of Sections \ref{seca} to \ref{sece} is to prove the CLT of LSS in Theorem \ref{thm2} of the main paper.  We first give two lemmas concerning the spectral properties of $\bar{S}_p$.
		\begin{lemma}[\cite{wang2014limiting}]\label{lem:wang}
			Suppose that Condition \ref{c1} holds. Then the empirical spectral distribution  of $\bar S_p$  almost surely  converges to a non-random probability function $F$ whose Stieltjes transform $m(z)$ is the unique solution satisfying $\Im m(z)>0$ and
			\begin{equation}\label{eq:system}
			m(z)=-\int\frac{1}{z+x\bar\lambda_{B^2}s(z)}dF^{A}(x),\quad s(z)=-\int\frac{x}{z+x\bar\lambda_{B^2}s(z)}dF^{A}(x),\quad z\in \mathbb{C}^+.
			\end{equation}
			where $\bar\lambda_{B^2}=\lim_{n\rightarrow \infty}n^{-1}\text{tr}(B_n^2)=\int x^2dF^{B}(x)$ is a positive constant.
		\end{lemma}
		\begin{lemma}[Bound on $\|\bar S_p\|$]\label{lem1}
			Suppose that Condition \ref{c1} holds. Then the spectral norm of $\bar S_p$ satisfies
			$	\|\bar S_p\|\le 2\lim\sup_{p}a_1b_1$  almost surely.
		\end{lemma}
		Lemma \ref{lem:wang} is borrowed from \cite{wang2014limiting}, while the proof of Lemma \ref{lem1} is in the next section.  Indeed, $m_p(z)$ and $s_p(z)$ in the main paper are finite sample versions of $m(z)$ and $s(z)$, respectively. Recall the definition
		\[
		\tilde G_p(f)=p\int_{-\infty}^{+\infty}f(x)d(F^{\bar S_p}(x)-F_p(x)),
		\]
		where $f\in\mathcal{M}=\{\text{functions which are  analytic in an open domain containing} [-2c,2c]\}$ with $c=\lim\sup_{p}a_1b_1$.  Define an event $U_p=\{\|\bar S_p\|\le 2\lim\sup a_1b_1\}$. Then by Lemma \ref{lem1} and the Cauchy integral formula,  $U_p$ happens with probability tending to one and
		\[
		\tilde G_{p}(f)=-\frac{1}{2\pi i}\oint_{\mathcal{C}}f(z)p\{m_{F^{\bar S_p}}(z)-m_p(z)\}dz,
		\]
		where $m_{F^{\bar S_p}}(z)$ and $m_p(z)$ are Stieltjes transforms associated with $F^{\bar S_p}(x)$ and $F_p(x)$ respectively, $\mathcal{C}$ is the contour formed by the boundary of the rectangle with four vertices $(\pm u_0,\pm iv_0)$ where $u_0=2\lim\sup_{p}a_1b_1+\epsilon_0$ with sufficiently small $\epsilon_0$, and $v_0$ is any positive number so that $f$ is analytic in a neighborhood of $\mathcal{C}$. Hence, it's sufficient to study the Stieltjes transform on the contour $\mathcal{C}$. Our major target is then to prove that with the mean correction term  $\mathcal{X}_p(z)$ defined in the main paper,  the process
		\begin{equation}\label{form}
		M_p(z)=p\{m_{F^{\bar S_p}}(z)-m_p(z)-\mathcal{X}_p(z)\}
		\end{equation}
		converges weakly to some limit which needs to be specified. 
		
		We will start with a simple case where $A_p$ is diagonal.  To avoid dealing with small imaginary part of $z$, we  denote
		\[
		\begin{split}
		&\mathcal{C}_u=\{u+iv_0:u\in[-u_0,u_0]\},\quad \mathcal{C}_l=\{-u_0+iv:v\in[\eta_p/p,v_0]\},\\
		&\mathcal{C}_r=\{u_0+iv:v\in[\eta_p/p,v_0]\},\quad
		\mathcal{C}_0=\{-u_0+iv:v\in[0,\eta_p/p]\}\cup \{u_0+iv:v\in[0,\eta_p/p]\},
		\end{split}
		\]
		where $\{\eta_p\}$ is a sequence decreasing to $0$ satisfying $\eta_p\ge p^{-\alpha}$ for some $\alpha\in(0,1)$, and $v_0$ is a constant.
		The following theorem specifies the limit of the process $M_p(z)$ for $z\in\mathcal{C}_u$.
		
		\begin{theorem}[Stieltjes transform]\label{thm:m}
			With Condition \ref{c1}, if $A_p$ is diagonal, with $\mathcal{X}_p(z)$ defined in the main paper,
			the process $\{M_{p}(z):z\in \mathcal{C}_u\}$ converges weakly to a two-dimensional Gaussian process $M(\cdot)$ with mean zero and the covariance function $\Lambda(z_1,z_2)$,  where $\Lambda(z_1,z_2)$ is defined in Theorem \ref{thm2}. 
		\end{theorem}
		
		Therefore, the proof of Theorem \ref{thm2} is divided into two parts. The first part proves Theorem \ref{thm:m}, organized in Sections \ref{secb} to \ref{non-random}. The second part is in Section \ref{sece}, which considers contours $\mathcal{C}_{l,r,0}$ and general $\Ab_p$.  Here we only present some key ideas.
		
		To prove Theorem \ref{thm:m}, we decompose $M_p(z)$ into random part and non-random part as follows,
		\[
		M_p(z)=p[m_{F^{\bar\Sbb_p}}(z)-\mathbb{E}m_{F^{\bar\Sbb_p}}(z)]+p[\mathbb{E}m_{F^{\bar\Sbb_p}}(z)-m_p(z)-\mathcal{X}_p(z)]=M_{p1}(z)+M_{p2}(z).
		\]
		The random part $M_{p1}(z)$ mainly contributes to the covariance function while  the  non-random part $M_{p2}(z)$ mainly contributes to the mean correction term $\mathcal{X}_p(z)$.  
		
		To deal with $M_{p1}(z)$, we expand the Stieltjes transform and observe that the randomness mainly comes from a sum of martingale  in equation (\ref{eq:main}) below. With the decomposition in this equation, we easily verify the finite-dimensional distribution of $M_{p1}(z)$ and specify the covariance function, following a similar expansion strategy  in \cite{chen2015clt}.  Finally, we finish the proof for $M_{p1}(z)$ by verifying the tightness.  The details are put in Section \ref{random}. 
		
		For $M_{p2}(z)$, we first calculate the diagonal entries of Green function with Schur's complement formula. It's seen that the diagonal entries will converge to the terms $\tilde\epsilon_k(z)$ defined in the main paper. However, the convergence rates are slow, so we need more detailed calculations of the errors to specify those which are not  smaller than $O(p^{-1})$. The proof procedure is more complicated than that in \cite{chen2015clt}, where $\mathbb{E}m_{F^{\bar\Sbb_P}}(z)$ can be explicitly solved by a quadratic equation of one variable. With general $\Ab_p$ and $\Bb_n$, we find that the equation therein is not correct anymore. Instead, we expand the diagonal entries of the Stieltjes transform  using the inverse matrix formula recursively. After each expansion step, we remove some ''big'' errors by rigorous approximation. Finally,  we sum up all these non-negligible errors,  which lead to the mean correction term $\mathcal{X}_p(z)$.
		
		To complete the proof of Theorem \ref{thm2}, we first prove that the integrals corresponding to $\mathcal{C}_{l,r,0}$ are negligible.  Since the event $\Ub_p(\epsilon_0):=\{\|\bar\Sbb_p\|\le 2\lim\sup a_1b_1+\epsilon_0/2\}$ happens with probability tending to 1 for large $p$, we aim  to prove that under this event
		\[
		\lim_{v_0\downarrow 0}\lim_{p\rightarrow\infty}\mathbb{E}\bigg|\int_{\mathcal{C}_{0,l,r}}f(z)M_p(z)dz\bigg|\rightarrow 0,
		\]
		which  is verified in Section \ref{sece}.  Moreover, for non-diagonal $\Ab_p$, note that the theorem still holds if $x_{ij}$ are i.i.d. standard normal variables, due to the orthogonal invariance property. For non-Gaussian variables, we use the Lindeberg replacement technique and compare the characteristic functions of  LSS with Gaussian and non-Gaussian distributions. Under the additional constraint $\nu_4=3$, we can replace  standard normal variables with general $x_{ij}$, which has negligible effects on the limiting distribution of the LSS. Then,  Theorem \ref{thm2} holds. Now we move to the details.

	}

	\section{Truncation and  the proof of Lemma \ref{lem1}}\label{secb}
	
	\subsection{Truncation}
	In this section, we first truncate the random variables $x_{ij}$ by $\delta_p\sqrt[4]{np}$, which hepls in controlling some higher-order moments of $x_{ij}$ in the proof. Then, we  prove  Lemma \ref{lem1}, which provides preliminary upper and lower bounds for the eigenvalues of $\bar\Sbb_p$. In the proof, we may suppress the dependence on $p$ and $z$ if it doesn't cause any confusion.
	
	We follow the truncation technique in \cite{chen2012convergence} and \cite{chen2015clt}.
	By Condition \ref{c1}, there is a series of  $\delta_p$ satisfying
	\[
	\lim_{p\rightarrow \infty}\delta_p^{-4}\mathbb{E}|x_{11}|^4I(|x_{11}|>\delta_p\sqrt[4]{np})=0,\quad \delta_p\downarrow 0,\delta_p\sqrt[4]{np}\uparrow\infty.
	\]
	Define the truncated variables as $\hat x_{ij}=x_{ij}I(x_{ij}\le \delta_p\sqrt[4]{np})$, and the standardized truncated variables as $\tilde x_{ij}=(\hat x_{ij}-\mathbb{E}\hat x_{ij})/\sigma$, where $\sigma^2$ is the variance of $\hat x_{ij}$. Let $\hat\Xb=(\hat x_{ij})$, $\tilde\Xb=(\tilde x_{ij})$,  $\hat\Sbb_p$ and $\tilde\Sbb_p$ be the normalized separable sample covariance matrix by replacing $ \Xb$ with $\hat \Xb$ and $\tilde\Xb$, respectively.  Then, similarly to the proof of Theorem 1 in \cite{chen2012convergence}, $\mathbb{P}(\bar\Sbb_p\ne \hat\Sbb_p,i.o.)=0$  by choosing a proper series of $\delta_p$. Hence, below we aim to control the effects by replacing $\hat\Xb$ with $\tilde\Xb$.
	
	Some elementary calculations yield that
	\begin{eqnarray}\label{eq:ex11}
	|\mathbb{E}\hat x_{ij}|&=&\big|\mathbb{E}x_{ij}-\mathbb{E}x_{ij}I(|x_{ij}|> \delta_p\sqrt[4]{np})\big|\le\mathbb{E}|x_{ij}|I(|x_{ij}|> \delta_p\sqrt[4]{np})|\\
	&\le& (\delta_p\sqrt[4]{np})^{-3}\mathbb{E}|x_{11}|^4I(|x_{11}|>\delta_p\sqrt[4]{np})=o\big((np)^{-3/4}\big),\\
	|1-\sigma^2|&\le &\mathbb{E}x_{11}^2I(|x_{11}|>\delta_p\sqrt[4]{np})|\le o\big((np)^{-1/2}\big).\label{eq:sx11}
	\end{eqnarray}
	Then, by  Wely's theorem,
	\begin{equation}\label{eq:G minus}
	\begin{split}
	&|G_p(f,\tilde\Xb)- G_p(f,\hat\Xb)|\le \sum_{j=1}^p|f(\lambda_j(\hat\Sbb_p))-f(\lambda_j(\tilde\Sbb_p))|\le C_f\sum_{j=1}^p|\lambda_j(\hat\Sbb_p)-\lambda_j(\tilde\Sbb_p)|\\
	\le&C_f\times p\|\hat\Sbb_p-\tilde\Sbb_p\|\le C_f\times \sqrt{\frac{p}{n}}\|\Ab\|\bigg(\|(\hat\Xb-\tilde\Xb)\Bb\hat\Xb^\prime\|+\|(\hat\Xb-\tilde\Xb)\Bb\tilde\Xb^\prime\|\bigg).
	\end{split}
	\end{equation}
	Note that by the i.i.d. assumption,
	\begin{equation}\label{eq: X tilde minus}
	\begin{split}
	\hat\Xb-\tilde\Xb=&\frac{\sigma-1}{\sigma}\hat\Xb+\frac{1}{\sigma}\mathbb{E}\hat x_{11}{\bf 1}_p{\bf 1}_n^\prime.
	\end{split}
	\end{equation}
	Therefore, 
	\[
	\begin{split}
	\sqrt{\frac{p}{n}}\|(\hat\Xb-\tilde\Xb)\Bb\hat\Xb^\prime\|\le &\frac{\sigma-1}{\sigma}\|\Bb\|\times \sqrt{np}\times \bigg\|\frac{1}{n}\hat\Xb\hat\Xb^\prime\bigg\|+\sqrt{\frac{p}{n}}\times \frac{1}{\sigma}|\mathbb{E}\hat x_{11}|\times \sqrt{p}\times O_p(\sqrt{np})\\
	\le &o_p(1),
	\end{split}
	\]
	where the $o_p(1)$ is by (\ref{eq:ex11}), (\ref{eq:sx11}), $p/n\rightarrow 0$ and the fact $\|n^{-1}\hat\Xb\hat\Xb^\prime\|=O_{a.s.}(1)$. Similarly, we can prove 
	\[
	\sqrt{\frac{p}{n}}\|(\hat\Xb-\tilde\Xb)\Bb\tilde\Xb^\prime\|=o_p(1).
	\]
	Consequently, the asymptotic distributions of $G_{p}(f,\tilde\Xb)$ and $ G_p(f,\hat\Xb)$ are the same. Similarly, the asymptotic distributions of $M_p(z,\tilde\Xb)$ and $M_p(z,\hat\Xb)$ are also the same. On the other hand, by (\ref{eq:G minus})
	\[
	\begin{split}
	\|\hat\Sbb_p-\tilde\Sbb_p\|\le  \frac{C}{\sqrt{np}}\|\Ab\|\bigg(\|(\hat\Xb-\tilde\Xb)\Bb\hat\Xb^\prime\|+\|(\hat\Xb-\tilde\Xb)\Bb\tilde\Xb^\prime\|\bigg).
	\end{split}
	\]
	By (\ref{eq: X tilde minus}), 
	\[
	\begin{split}
	\frac{1}{\sqrt{np}}\|(\hat\Xb-\tilde\Xb)\Bb\hat\Xb^\prime\|\le&\frac{\sigma-1}{\sigma}\|\Bb\|\times \sqrt{\frac{n}{p}}\times \bigg\|\frac{1}{n}\hat\Xb\hat\Xb^\prime\bigg\|+\frac{1}{\sqrt{np}}\times \frac{1}{\sigma}|\mathbb{E}\hat x_{11}|\times \sqrt{np}\times O_{a.s.}\big((np)^{3/4}\big)\\
	=&o_{a.s}(1),
	\end{split}
	\]
	where we use $|\hat x_{ij}|\le \delta_p\sqrt[4]{np}$ for the term $\mathbf{1}_n^\prime\hat\Xb^\prime$. Similarly, 
	\[
	\frac{1}{\sqrt{np}}\|(\hat\Xb-\tilde\Xb)\Bb\tilde\Xb^\prime\|=o_{a.s.}(1).
	\]
	In conclusion, replacing $\Xb$ with $\tilde\Xb$ will have negligible effects on the LSS and eigenvalues.
	Hence, we will focus on  $\tilde\Xb$ rather than $\Xb$ in the proof. For simplicity, we still write $\Xb$ but assume that
	\begin{equation}\label{eq:truncation}
	|x_{ij}|\le \delta_p\sqrt[4]{np},\quad\mathbb{E}x_{ij}=0,\quad\mathbb{E}x_{ij}^2=1,\quad \mathbb{E}x_{ij}^4=\nu_4+o(1).
	\end{equation}

	\subsection{The proof of Lemma \ref{lem1}}
	We prove Lemma \ref{lem1} here, which gives a rough bound for the support of $F^{\bar\Sbb_p}$.  By definition of $\bar\Sbb_p$,
	\begin{equation}\label{eq:sp}
	\begin{split}
	\|\bar\Sbb_p\|\le& \bigg\|\frac{1}{\sqrt{np}}(\Ab^{1/2}\Xb\Bb\Xb^\prime\Ab^{1/2}-n\bar\lambda_{\Bb}\Ab)\bigg\|\le \|\Ab\|\bigg\|\frac{1}{\sqrt{np}}(\Xb\Bb\Xb^\prime-n\bar\lambda_{\Bb})\bigg\|\\
	\le &a_1\bigg(\max_{i}\frac{1}{\sqrt{np}}|\bx_i^\prime\Bb\bx_i-n\bar\lambda_{\Bb}|+\bigg\|\frac{1}{\sqrt{np}}\Bb^0_x\bigg\|\bigg),
	\end{split}
	\end{equation}
	where $\Bb_x^0$ is defined by 
	\[
	(\Bb_x^0)_{ij}=\left\{\begin{aligned}
	&0,&i=j,\\
	&\bx_i^\prime\Bb\bx_j,&i\ne j.
	\end{aligned}\right.
	\]
	We first deal with $\Bb_x^0$. For simplicity, we suppress the index $x$. Let $\be_i$ be the $p$-dimensional vector with the $i$-th entry being 1 and the others being 0. Then, 
	\[
	\begin{split}
	&\bigg\|\frac{1}{\sqrt{np}}\Bb^0\bigg\|^2=\sup_{\|\bxi\|=1}\frac{1}{np}\|\Bb^0\bxi\|^2=\sup_{\|\bxi\|=1}\sum_{i=1}^p\frac{1}{np}\bigg(\bx_i^\prime\Bb(\Xb^\prime-\bx_i\be_i^\prime)\bxi\bigg)^2\\
	=&\sup_{\|\bxi\|=1}\sum_{i=1}^p\frac{1}{np}\bigg\|\Bb(\Xb^\prime-\bx_i\be_i^\prime)\bxi\bx_i^\prime\bigg\|^2\le \sup_{\|\bxi\|=1}\sum_{i=1}^p\frac{1}{np}\|\Bb\|^2\bigg\|(\Xb^\prime-\bx_i\be_i^\prime)\bxi\bx_i^\prime\bigg\|^2\\
	\le & \sup_{\|\bxi\|=1}\|\Bb\|^2\sum_{i=1}^p\frac{1}{np}\bigg(\bx_i^\prime(\Xb^\prime-\bx_i\be_i^\prime)\bxi\bigg)^2=\|\Bb\|^2\times \bigg\|\frac{1}{\sqrt{np}}\tilde\Bb_x^0\bigg\|^2,
	\end{split}
	\]
	where $\tilde\Bb_x^0$ is defined by
	\[
	(\tilde\Bb_x^0)_{ij}=\left\{\begin{aligned}
	&0,&i=j,\\
	&\bx_i^\prime\bx_j,&i\ne j,
	\end{aligned}\right.
	\]
	and we use the fact  $\bx_i^\prime(\Xb^\prime-\bx_i\be_i^\prime)\bxi\ne 0$ with  probability one.
	By the proof of Theorem 2 in \cite{chen2012convergence}, we conclude that almost surely
	\begin{equation}\label{eq:as}
	\|(np)^{-1/2}\tilde\Bb_x^0\|\le 2.
	\end{equation} 
	
	Now we come back to the first term in (\ref{eq:sp}) and show that it's $o_{a.s.}(1)$. Write
	\begin{equation}\label{eq:diag}
	\frac{1}{\sqrt{np}}(\bx_i^\prime\Bb\bx_i-n\bar\lambda_{\Bb})=\frac{1}{\sqrt{np}}\sum_{j=1}^n\mathrm{B}_{n,jj}(x_{ij}^2-1)+\frac{1}{\sqrt{np}}\sum_{j\ne l}^n\mathrm{B}_{jl}x_{ij}x_{il}.
	\end{equation}
	Since $\|\Bb\|$ is bounded and $x_{ij}$ are i.i.d., a similar technique to proving equation (9) in \cite{chen2012convergence} leads to
	\begin{equation}\label{eq:as1}
	\max_{i}\frac{1}{\sqrt{np}}\bigg|\sum_{j=1}^n\mathrm{B}_{n,jj}(x_{ij}^2-1)\bigg|=o_{a.s.}(1).
	\end{equation}
	Hence, we only focus on the second term of (\ref{eq:diag}). Using Lemma 5 in \cite{pan2011central}, for  any $j\le n$ and  $k\ge 2$, 
	\[
	\begin{split}
	&\mathbb{P}\bigg(\max_{i}\frac{1}{\sqrt{np}}\bigg|\sum_{l\ne j}^n\mathrm{B}_{jl}x_{ij}x_{il}\bigg|>\epsilon\bigg)\le p\mathbb{P}\bigg(\frac{1}{\sqrt{np}}\bigg|\sum_{l\ne j}^n\mathrm{B}_{jl}x_{1j}x_{1l}\bigg|>\epsilon\bigg)\\
	\le &\frac{p}{\epsilon^{2k}(np)^{k/2}}\mathbb{E}\bigg|\sum_{l\ne j}^n\mathrm{B}_{jl}x_{1j}x_{1l}\bigg|^k\le \frac{C_{k,\epsilon}p}{(np)^{k/2}}(\mathbb{E}|x_{11}|^k)^2\|\Bb\|_F^k.
	\end{split}
	\]
	Take $k=4+\delta_0$ to get
	\[
	\mathbb{P}\bigg(\max_{i}\frac{1}{\sqrt{np}}\bigg|\sum_{l\ne j}^n\mathrm{B}_{jl}x_{ij}x_{il}\bigg|>\epsilon\bigg)\le\frac{C_{\epsilon}}{p^{1+\delta_0/2}}.
	\]
	Since $\delta_0>0$, by the Borel-Cantelli theorem and the fact $\sum_{p=1}^\infty p^{-1-\delta_0/2}<\infty$, we conclude that
	\begin{equation}\label{eq:as2}
	\max_{i}\frac{1}{\sqrt{np}}\bigg|\sum_{l<j}^n\mathrm{B}_{jl}x_{ij}x_{il}\bigg|=o_{a.s.}(1).
	\end{equation}
	The lemma then follows from (\ref{eq:as}), (\ref{eq:diag}), (\ref{eq:as1}) and (\ref{eq:as2}). The next two sections are devoted to the proof of Theorem \ref{thm:m}.

	\section{The random part with diagonal $\Ab_p$}\label{random}
	
	\subsection{Expansion of the Stieltjes transform}
	We rewrite $M_p(z)$ as the sum of random part $M_{p1}(z)$ and  non-random part $M_{p2}(z)$ by
	\[
	M_p(z)=p[m_{F^{\bar\Sbb_p}}(z)-\mathbb{E}m_{F^{\bar\Sbb_p}}(z)]+p[m_{F^{\bar\Sbb_p}}(z)-m_p(z)-\mathcal{X}_p(z)]:=M_{p1}(z)+M_{p2}(z).
	\]
	Then, Theorem \ref{thm:m} can be concluded with detailed analysis of the random part and non-random part. In this section, we focus on the random part and assume $\Ab$ is diagonal. The proof in this part is adapted from \cite{chen2015clt} where both $\Ab$ and $\Bb$ are identity matrices. The differences lie in the following two aspects. First, after expanding the Stieltjes transform by deleting each row of $\Xb$, there will be a coefficient corresponding to  $a_k$ since the diagonal entries of $\Ab$ are no more identical.  Second, we will have more error terms at some steps which are from the off-diagonal entries of $\Bb$.
	
	We first introduce some notation. Write $\Yb=\Ab^{1/2}\Xb\Bb^{1/2}$. Let $\Xb_k$ be the $(p-1)\times n$ matrix after removing  $\bx_k^\prime$ from $\Xb$, where $\bx_k^\prime$ is the $k$-th row vector of $\Xb$.  Similarly, let $\by_k^\prime$ be the $k$-th row of $\Yb$, $\Yb_k$ be the $(p-1)\times n$ matrix after removing  the $k$-th row  from $\Yb$, and $\bar\Sbb_k=\frac{1}{\sqrt{np}}(\Yb_k\Yb_k^\prime-n\bar\lambda_{\Bb}\Ab_k)$, where $\Ab_k$ is obtained by deleting the $k$-th row and column in $\Ab$. Let $\Db:=\bar\Sbb-z\Ib_p$ and $\Db_k=\bar\Sbb_k-z\Ib_{p-1}$, where $z=u+iv\in\mathcal{C}_u$. Note that the $k$-th diagonal entry of $\Db$ is  $d_k:=\frac{1}{\sqrt{np}}(\by_k^\prime\by_k-n\bar\lambda_{\Bb}a_k)-z$ and the $k$-th row of $\Db$ with the $k$-th element deleted is $\bq_k^\prime:=\frac{1}{\sqrt{np}}\by_k^\prime\Yb_k^\prime$.
	
	Define the $\sigma$-field generated by $(\bx_1,\ldots,\bx_k)$ as $\mathcal{F}_k$ and the conditional expectation $\mathbb{E}_k(\cdot):=\mathbb{E}(\cdot|\mathcal{F}_k)$. By Shur's complement formula and the inverse formula,
	\begin{eqnarray}\label{eq:D inverse1}
	\text{tr}(\Db^{-1})&=&\frac{1}{d_k-\bq_k^\prime\Db_k^{-1}\bq_k}+\text{tr}(\Db_k-d_k^{-1}\bq_k\bq_k^\prime)^{-1},\\
	\text{tr}(\Db_k-d_k^{-1}\bq_k\bq_k^\prime)^{-1}&=&\text{tr}(\Db_k^{-1})+d_k^{-1}\bq_k^\prime\Db_k^{-1}(\Db_k-d_k^{-1}\bq_k\bq_k^\prime)^{-1}\bq_k,\\
	\quad\quad d_k^{-1}\bq_k^\prime\Db_k^{-1}(\Db_k-d_k^{-1}\bq_k\bq_k^\prime)^{-1}\bq_k&=&\frac{d_k^{-1}\bq_k^\prime\Db_k^{-2}\bq_k}{1-d_k^{-1}\bq_k^\prime\Db_k^{-1}\bq_k}=\frac{\bq_k^\prime\Db_k^{-2}\bq_k}{d_k-\bq_k^\prime\Db_k^{-1}\bq_k}.\label{eq:D inverse3}
	\end{eqnarray}
	Therefore, we have
	\begin{equation}
	\text{tr}(\Db^{-1}-\Db_k^{-1})=-\frac{1+\bq_k^\prime\Db_k^{-2}\bq_k}{-d_k+\bq_k^\prime\Db_k^{-1}\bq_k}.
	\end{equation}
	As a result,
	\begin{equation}\label{eq:mF expansion}
	\begin{split}
	p[m_{F^{\bar\Sbb}}-\mathbb{E}(m_{F^{\bar\Sbb}})]=&\sum_{k=1}^p(\mathbb{E}_k-\mathbb{E}_{k-1})\text{tr}(\Db^{-1}-\Db_k^{-1}):=\sum_{k=1}^p\rho_k\\
	:=&(\mathbb{E}_k-\mathbb{E}_{k-1})\zeta_k-\mathbb{E}_k\kappa_k,
	\end{split}
	\end{equation}
	where
	\[
	\begin{split}
	\rho_k:=&-(\mathbb{E}_k-\mathbb{E}_{k-1})\beta_k(1+\bq_k^\prime\Db_k^{-2}\bq_k),\quad\beta_k:=\frac{1}{-d_k+\bq_k^\prime\Db_k^{-1}\bq_k},\\
	\zeta_k:=&-\tilde\beta_k\beta_k\eta_k(1+\bq_k^\prime\Db_k^{-2}\bq_k),\\
	\tilde\beta_k:=&\frac{1}{z+\frac{1}{np}\text{tr}\Mb_k^{(1)}},\quad \Mb_k^{(s)}:=a_k\Bb^{1/2}\Yb_k^\prime\Db_k^{-s}\Yb_k\Bb^{1/2},s=1,2,\\
	\eta_k:=&\frac{1}{\sqrt{np}}(\by_k^\prime\by_k-n\bar\lambda_{\Bb}a_{k})-\gamma_{k1},\\
	\gamma_{ks}:=&\bq_k^\prime\Db_k^{-s}\bq_k-(np)^{-1}\text{tr}\Mb_k^{(s)},s=1,2,\quad \kappa_k:=\tilde\beta_k\gamma_{k2}.
	\end{split}
	\]
	Note that in the above derivation we use 
	\begin{equation}\label{eq:beta}
	\beta_k=\tilde\beta_k+\beta_k\eta_k\tilde\beta_k.
	\end{equation}
	
	We now provide some useful bounds for further use. First, since the eigenvalues of $\Db$ have the form of $1/(\lambda_j(\bar\Sbb)-z)$, considering the imaginary part we have $\|\Db^{-1}\|\le v^{-1}$. Similarly $\|\Db_k^{-1}\|\le v^{-1}$. Besides, note that $\beta_k$ is the $k$-th diagonal entry of $\Db$, hence $|\beta_k|\le v^{-1}$. By symmetry and the eigenvalue decomposition of $\Db_k^{-1}$, it's also easy to see that all the  diagonal entries of $\Mb_k^{(1)}$ have positive imaginary part, then $|\tilde\beta_k|\le v^{-1}$. Moreover, let $\bGamma=\{\bxi_1,\ldots,\bxi_{p-1}\}$ and $\bTheta=\text{diag}(\theta_1,\ldots,\theta_{p-1})$ be the  eigenvectors and eigenvalues of $\bar\Sbb_k$.  Then,
	\begin{equation}\label{eq:trace mk}
	\begin{split}
	&\bigg|\frac{1}{np}\text{tr}\Mb_k^{(1)}\bigg|=\bigg|\frac{1}{np}\text{tr}\Bb^{1/2}\Yb_k^\prime\bGamma(\bTheta-z\Ib_{p-1})^{-1}\bGamma^\prime\Yb_k\Bb^{1/2}\bigg|\\
	=&\bigg|\frac{1}{np}\text{tr}\bigg(\Bb^{1/2}\Yb_k^\prime\bGamma(\bTheta-z)^{-1}(\bTheta-\bar z)^{-1}(\bTheta-u)\bGamma^\prime\Yb_k\Bb^{1/2}\\
	&+iv\times \Bb^{1/2}\Yb_k^\prime\bGamma(\bTheta-z)^{-1}(\bTheta-\bar z)^{-1}\bGamma^\prime\Yb_k\Bb^{1/2}\bigg)\bigg|\\
	\le &\frac{1}{np}\text{tr}\bigg(\Bb^{1/2}\Yb_k^\prime\big(\Db_k^{-1}(z)\Db_{k}^{-1}(\bar z)\big)^{1/2}\Yb_k\Bb^{1/2}\bigg)\\
	\le &\frac{b_1}{np}\text{tr}\bigg(\Yb_k^\prime\big(\Db_k^{-1}(z)\Db_{k}^{-1}(\bar z)\big)^{1/2}\Yb_k\bigg).
	\end{split}
	\end{equation}
	Note that for sufficiently large $p$,
	\[
	\begin{split}
	&\frac{1}{np}\text{tr}\bigg(\Yb_k\Yb_k^\prime\big(\Db_k^{-1}(z)\Db_{k}^{-1}(\bar z)\big)^{1/2}\bigg)\\
	=&\frac{1}{np}\text{tr}\bigg(\big(\sqrt{np}(\Db_k(z)+z\Ib_{p-1})+n\bar\lambda_{\Bb}\Ab_k\big)\big(\Db_k^{-1}(z)\Db_{k}^{-1}(\bar z)\big)^{1/2}\bigg)\\
	=&\frac{1}{np}\text{tr}\bigg(\sqrt{np}\big(\Db_k(z)\Db_{k}^{-1}(\bar z)\big)^{1/2}+z\sqrt{np}\big(\Db_k^{-1}(z)\Db_{k}^{-1}(\bar z)\big)^{1/2}+n\bar\lambda_{\Bb}\Ab_k\big(\Db_k^{-1}(z)\Db_{k}^{-1}(\bar z)\big)^{1/2}\bigg)\\
	\le& C\bigg(\sqrt{\frac{p}{n}}+1\bigg)\le C.
	\end{split}
	\]
	With a similar technique applied to $\Mb_{k}^{(2)}$, we conclude
	\[
	\bigg|\frac{1}{np}\text{tr}\Mb_k^{(1)}\bigg|\le C,\quad \bigg|\frac{1}{np}\text{tr}\Mb_k^{(2)}\bigg|\le C.
	\]
	
	On the other hand, 
	\[
	\begin{split}
	|(1+\bq_k^\prime\Db_k^{-2}\bq_k)\beta_k|\le &(1+|\bq_k^\prime\Db_k^{-2}\bq_k|)\times |\beta_k|\\
	\le &\bigg(1+\sum_j\frac{<\bq_k,\bxi_j>^2}{|\theta_j-z|^2}\bigg)\times \frac{1}{\Im(-d_k+\bq_k^\prime\Db_k^{-1}\bq_k)}\\
	=&\bigg(1+\sum_j\frac{<\bq_k,\bxi_j>^2}{|\theta_j-z|^2}\bigg)\times \frac{1}{v\bigg(1+\sum_j\frac{<\bq_k,\bxi_j>^2}{|\theta_j-z|^2}\bigg)}\\
	\le& v^{-1}.
	\end{split}
	\]
	
	Using (\ref{eq:beta}), we further decompose $\zeta_k$ into
	\[
	\begin{split}
	\zeta_k=&-\tilde\beta^2\eta_k[1+(np)^{-1}\text{tr}\Mb_k^{(2)}]-\tilde\beta^2\eta_k\gamma_{k2}-\tilde\beta_k^2\beta_k\eta_k^2(1+\bq_k^\prime\Db_k^{-2}\bq_k)\\
	:=&\zeta_{k1}+\zeta_{k2}+\zeta_{k3}.
	\end{split}
	\]
	To this end, we should first show the following lemma on large deviation bounds.
	\begin{lemma}\label{lemma:gamma and eta}
		For $z=u+iv$ with $v>0$, under Condition \ref{c1} and (\ref{eq:truncation}), we have
		\[
		\begin{split}
		\mathbb{E}|\gamma_{ks}|^2\le &Cp^{-1},\quad \mathbb{E}|\gamma_{ks}|^4\le C\bigg(\frac{1}{p^2}+\frac{p}{n^2}+\frac{1}{np}\bigg),\\
		\mathbb{E}|\eta_k|^2\le &Cp^{-1},\quad \mathbb{E}|\eta_k|^4\le C\bigg(\frac{\delta_p^4}{p}+ \frac{1}{p^2}+\frac{p}{n^2}+\frac{1}{np}\bigg),	
		\end{split}
		\]
		where $\delta_p\rightarrow 0$ slowly is the truncation parameter.
	\end{lemma}
	The proof of  Lemma \ref{lemma:gamma and eta} is given in the next subsection. With this lemma, the Burkholder inequality and the Cauchy-Schwartz inequality, we have
	\[
	\mathbb{E}\bigg|\sum_{k=1}^p(\mathbb{E}_k-\mathbb{E}_{k-1})\zeta_{k3}\bigg|^2\le \sum_{k=1}^p\mathbb{E}\bigg|\tilde\beta_k^2\beta_k\eta_k^2(1+\bq_k^\prime\Db_k^{-2}\bq_k)\bigg|^2\le C\delta_p^4.
	\]
	Similarly, 
	\[
	\mathbb{E}\bigg|\sum_{k=1}^p(\mathbb{E}_k-\mathbb{E}_{k-1})\zeta_{k2}\bigg|^2\le \sum_{k=1}^p\mathbb{E}\bigg|\tilde\beta_k^2\eta_k\gamma_{k2}\bigg|^2\le C\bigg(\delta_p^2+\frac{p}{n}\bigg).
	\]
	Consequently, we can write
	\begin{equation}\label{eq:main}
	\begin{split}
	p[m_{F^{\bar\Sbb}}(z)-\mathbb{E}(m_{F^{\bar\Sbb}}(z))]=&\sum_{k=1}^p\mathbb{E}_k\bigg[-\bigg(1+\frac{1}{np}\text{tr}(\Mb_k^{(2)})\bigg)\tilde\beta_k^2\eta_k-\kappa_k\bigg]+o_{L_2}(1)\\
	:=&\sum_k\mathbb{E}_k(\alpha_k(z))+o_{L_2}(1).
	\end{split}
	\end{equation}
	It's easy to see that $\mathbb{E}_{k-1}\alpha_k(z)=0$. 
	Then, for the random part, it's sufficient to prove the finite-dimensional convergence and the tightness of $M_{p1}(z)$.  For the finite-dimensional convergence,
	we only need to consider the sum 
	\[
	\sum_{j=1}^l\xi_j\sum_{k=1}^p\mathbb{E}_k(\alpha_k(z_j))=\sum_{k=1}^p\bigg(\sum_{j=1}^l\xi_j\mathbb{E}_k(\alpha_k(z_j))\bigg),
	\]
	where $\xi_1,\ldots,\xi_j$ are complex numbers and $l$ is a positive integer. Before moving forward, we first give the proof of Lemma \ref{lemma:gamma and eta} in the next subsection.

	\subsection{Proof of Lemma \ref{lemma:gamma and eta}}
	\begin{proof}
		Note that $\by_k^\prime\by_k-n\bar\lambda_{\Bb}a_k=a_k[\sum_{j=1}^n\mathrm{B}_{n,jj}(x_{kj}^2-1)+\sum_{j\ne l}\mathrm{B}_{jl}x_{kj}x_{kl}]$, hence by  independence,
		\[
		\begin{split}
		\mathbb{E}\bigg|\frac{1}{\sqrt{np}}(\by_k^\prime\by_k-n\bar\lambda_{\Bb}a_k)\bigg|^2\le &\frac{C}{p},\\
		\mathbb{E}\bigg|\frac{1}{\sqrt{np}}(\by_k^\prime\by_k-n\bar\lambda_{\Bb}a_k)\bigg|^4\le& \mathbb{E}\bigg|\frac{C}{\sqrt{np}}\sum_{j=1}^nb_j(x_{kj}^2-1)\bigg|^4+C\mathbb{E}\bigg|\frac{1}{\sqrt{np}}\sum_{i\ne j}\mathrm{B}_{ij}x_{ki}x_{kj}\bigg|^4\\
		\le&\frac{C}{n^2p^2}\mathbb{E}\bigg(\sum_{j=1}^n(x_{kj}^2-1)^2\bigg)^2+\frac{C}{n^2p^2}(\mathbb{E}x_{11}^4)^2\|\Bb\|_F^4
		\\ \le&\frac{C}{n^2p^2}\bigg(\sum_{j=1}^n\mathbb{E}x_{kj}^4\bigg)^2+\frac{C}{n^2p^2}\sum_{j=1}^n\mathbb{E}x_{jk}^8+\frac{C}{p^2}\le C\bigg(\frac{1}{p^2}+\frac{\delta_p^4}{p}\bigg),
		\end{split}
		\]
		where for the third inequality we use Burkholder's inequality and  Lemma 5 in \cite{pan2011central}.
		On the other hand, note that
		\[
		\begin{split}
		&\bq_k^\prime\Db_k^{-s}\bq_k-\frac{1}{np}\text{tr}(\Mb_k^{(s)})=\frac{1}{np}\bigg(\bx_k^\prime\Mb_k^{(s)}\bx_k-\text{tr}(\Mb_k^{(s)})\bigg)\\
		:=&\frac{1}{np}\bx_k^\prime\Hb\bx_k+\frac{1}{np}\sum_{j=1}^n\mathrm{M}_{jj}^{(s)}(x_{kj}^2-1),
		\end{split}
		\]
		where $\Hb=\Mb_k^{(s)}-\text{diag}(\mathrm{M}_{11}^{(s)},\ldots,\mathrm{M}_{nn}^{(s)})$ and $\mathrm{M}_{jj}^{(s)}$ is the $j$-th diagonal entry of $\Mb_k^{(s)}$.  Let $\bB_j^\prime$ be the $j$-th row of $\Bb$, then it's not hard to verify
		\[
		\begin{split}
		\mathbb{E}|\mathrm{M}_{jj}^{(s)}|^4\le C\mathbb{E}\|\bB_j^\prime\Xb_k^\prime\Ab_k\Db_k^{-s}\Ab_k\Xb_k\bB_j\|^4\le C\mathbb{E}\big(\sum_l\|\mathrm{B}_{jl}\bx_{(k),l}\|^2\big)^4=C\mathbb{E}\big(\sum_l\mathrm{B}_{jl}^2\|\bx_{(k),l}\|^2\big)^4,
		\end{split}
		\]
		where $\bx_{(k),l}$ is the $l$-th column of $\Xb_k$. Then, by Burkholder's inequality,
		\[
		\begin{split}
		\mathbb{E}|\mathrm{M}_{jj}^{(s)}|^4\le & C\mathbb{E}\bigg(\sum_l\big(\mathrm{B}_{jl}^2\|\bx_{(k),l}\|^2\big)^2\bigg)^2
		=C\bigg(\sum_l\mathbb{E}\big(\mathrm{B}_{jl}^2\|\bx_{(k),l}\|^2\big)^2\bigg)^2+C\sum_l\mathbb{E}\big(\mathrm{B}_{jl}^2\|\bx_{(k),l}\|^2\big)^4\\
		\le & Cp^4(\sum_l\mathrm{B}_{jl}^4)^2+C(\sum_l\mathrm{B}_{jl}^8)\mathbb{E}\|\bx_{(k),1}\|^8\le C(p^4+np^2),
		\end{split}
		\]
		where we use the facts that
		\[
		\sum_l\mathrm{B}_{jl}^{2h}\le \sum_l\mathrm{B}_{jl}^2\times \|\Bb\|^{2h-2}\le \|\Bb\|^{2h},
		\]
		and  
		\[
		\mathbb{E}\|\bx_{(k),1}\|^8\le (\sum_{h}\mathbb{E}|x_{(k)_,lh}|^4)^2+\sum_h\mathbb{E}|x_{(k),lh}|^8\le C(p^2+\delta_p^4 np^2).
		\]
		Let $\mathbb{E}_{j_k}(\cdot)=\mathbb{E}(\cdot|x_{k1},\ldots,x_{kj})$, then by Burkholder's inequality and the independence between $\Mb_k^{(s)}$ and $\bx_k$, we have
		\[
		\begin{split}
		\mathbb{E}\bigg|\frac{1}{np}\sum_{j=1}^n\mathrm{M}_{jj}^{(s)}(x_{kj}^2-1)\bigg|^2\le& \frac{C}{n^2p^2}\sum_{j=1}^n\mathbb{E}\bigg|\mathrm{M}_{jj}^{(s)}(x_{kj}^2-1)\bigg|^2\le \frac{C}{p},\\
		\mathbb{E}\bigg|\frac{1}{np}\sum_{j=1}^n\mathrm{M}_{jj}^{(s)}(x_{kj}^2-1)\bigg|^4\le&\frac{C}{n^4p^4}\bigg(\sum_{j=1}^n\mathbb{E}\big|\mathrm{M}_{jj}^{(s)}(x_{kj}^2-1)\big|^2\bigg)^2+\frac{C}{n^4p^4}\sum_{j=1}^n\mathbb{E}|\mathrm{M}_{jj}^{(s)}|^4(x_{kj}^2-1)^4\\
		\le &C\bigg(\frac{1}{np}+\frac{p}{n^2}\bigg).
		\end{split}
		\]
		It remains to deal with $\bx_k^\prime\Hb\bx_k$.  By Lemma 5 in \cite{pan2011central}, for any $h\ge 2$,
		\[
		\mathbb{E}|\bx_k^\prime\Hb\bx_k|^h\le C\big(\mathbb{E}|x_{11}|^h\big)^2\mathbb{E}(\text{tr}\Hb\Hb^*)^{h/2}.
		\]
		Let $h=2$, then
		\[
		\mathbb{E}|\bx_k^\prime\Hb\bx_k|^2\le \mathbb{E}\|\Hb\|_F^2\le \mathbb{E}\|\Mb_k^{(s)}\|_F^2.
		\]
		By the definition of $\Mb_k^{(s)}$, 
		\[
		\begin{split}
		\|\Mb_k^{(s)}\|_F^2\le& \|\Db_k^{-s}\Yb_k\Yb_k^\prime\|_F^2=\bigg\|\Db_k^{-s}\big(\sqrt{np}(\Db_k+z\Ib_{p-1})+n\bar\lambda_{\Bb}\Ab_k\big)\bigg\|_F^2\\
		\le&C\bigg(np\|\Db_k^{1-s}\|_F^2+\|\Db_k^{-s}\|_F^2+n^2\|\Db_k^{-s}\|_F^2\bigg)\le Cn^2p.
		\end{split}
		\]
		Therefore, 
		\[
		\mathbb{E}\bigg|\frac{1}{np}\bx_k^\prime\Hb\bx_k\bigg|^2\le \frac{C}{p},\quad 		\mathbb{E}\bigg|\frac{1}{np}\bx_k^\prime\Hb\bx_k\bigg|^4\le \frac{C}{p^2}.
		\]
		Combining the above results, we then conclude the lemma. 
	\end{proof}

	\subsection{Finite-dimensional distribution}
	Recall the expansion in (\ref{eq:main}). By  the central limit theorem for martingale, for the random part it's sufficient to verify the two  conditions in Lemma 9.12 in \cite{bai2010spectral}. The condition (9.9.2) therein is easily verified by
	\[
	\sum_{k=1}^p\mathbb{E}\bigg|\sum_{j=1}^l\xi_j\mathbb{E}_k(\alpha_k(z_j))\bigg|^4\le C\sum_{k=1}^p\mathbb{E}\bigg(|\eta_k|^4+|\gamma_{k2}|^4\bigg)\rightarrow 0.
	\]
	Hence, in the following we aim to check another condition, which is equivalent to finding the limit in probability of the covariance 
	\[
	\Lambda_p(z_1,z_2):=\sum_{k=1}^p\mathbb{E}_{k-1}[\mathbb{E}_k(\alpha_k(z_1))\cdot\mathbb{E}_k(\alpha_k(z_2))].
	\]
	Recall the expression of $\alpha_k(z)$ in (\ref{eq:main}), 
	\[
	\alpha_k(z)=-\bigg(1+\frac{1}{np}\text{tr}\Mb_k^{(2)}\bigg)\tilde\beta_k^2\eta_k-\gamma_{k2}\tilde\beta_k=\frac{\partial}{\partial z}(\tilde\beta_k\eta_k).
	\] By the dominated convergence theorem, we then focus on 
	\[
	\frac{\partial^2}{\partial z_2\partial z_1}\sum_{k=1}^p\mathbb{E}_{k-1}[\mathbb{E}_k(\tilde\beta_k(z_1)\eta_k(z_1))\cdot\mathbb{E}_k(\tilde\beta_k(z_2)\eta_k(z_2))].
	\]
	We first aim to find the limit of $\tilde\beta_k$.  It's easy to get
	\[
	\tilde\beta_k-\epsilon_k=-\epsilon_k\tilde\beta_k\frac{1}{np}\bigg(\text{tr}\Mb_k^{(1)}-\mathbb{E}\text{tr}\Mb_k^{(1)}\bigg),\quad \epsilon_k=\frac{1}{z+\frac{1}{np}\mathbb{E}\text{tr}\Mb_k^{(1)}}.
	\]
	To this end, we need more notation as follows. We remove $a_k$ from $\Mb_k^{(1)}$ and write $\Mb_k$. Abuse the notation a little to let $\{\be_i,i=1,\ldots,k-1,k+1,\ldots,p\}$ be the $(p-1)$-dimensional vector with the $i$-th ( or $(i-1)$-th if $k<i$) entry being 1 and the others being 0.  Then, we have $\Yb_k=\Yb_{kl}+\be_l\by_{l}^\prime$, where $\Yb_{kl}$ is obtained by replacing the $l$-th (or ($l-1$)-th) row of $\Yb_k$ with 0. Further write
	\[
	\begin{split}
	\bh_l^\prime=&\frac{1}{\sqrt{np}}\by_l^\prime\Yb_{kl}^\prime+\frac{1}{\sqrt{np}}(\by_l^\prime\by_l-na_l\bar\lambda_{\Bb})\be_l^\prime, \quad \br_l=\frac{1}{\sqrt{np}}\Yb_{kl}\by_l,\\
	\Db_{kl,r}=&\Db_k-\be_l\bh_l^\prime=\frac{1}{\sqrt{np}}(\Yb_{kl}\Yb_{k}^\prime-n\bar\lambda_{\Bb}\Ab_{(l)})-z\Ib_{p-1}, \quad \zeta_l=\frac{1}{1+\vartheta_l},\quad \vartheta_l=\bh_l^\prime\Db_{kl,r}^{-1}\be_l,\\
	\Db_{kl}=& \Db_{kl,r}-\br_l\be_l^\prime=\frac{1}{\sqrt{np}}(\Yb_{kl}\Yb_{kl}^\prime-n\bar\lambda_{\Bb}\Ab_{(l)})-z\Ib_{p-1}, \quad \Mb_{kl}=\Bb^{1/2}\Yb_{kl}^\prime\Db_{kl}^{-1}\Yb_{kl}\Bb^{1/2},
	\end{split}
	\]
	where $\Ab_{(l)}$ is obtained by replacing the $l$-th diagonal entry of $\Ab$ with 0.  Then, by the fact $|\tilde\beta|\le C$, $|\epsilon_k|\le C$,
	\[
	\begin{split}
	\mathbb{E}|\tilde\beta_k-\epsilon_k|^2\le&C\mathbb{E}\bigg[\frac{1}{np}\bigg(\text{tr}\Mb_k-\mathbb{E}\text{tr}\Mb_k\bigg)\bigg]^2\le C\mathbb{E}\bigg[\frac{1}{np}\sum_{l=1}^p(\mathbb{E}_{l}-\mathbb{E}_{l-1})\big(\text{tr}\Mb_k-\text{tr}\Mb_{kl}\big)\bigg]^2.
	\end{split}
	\]
	By some elementary but tedious calculations (see section 5.2 in \cite{chen2015clt}), we have
	\begin{equation}\label{eq:mk1}
	\begin{split}
	\Mb_k=&\Mb_{kl}-\frac{a_l\zeta_l}{znp}\Mb_{kl}\bx_l\bx_l^\prime\Mb_{kl}+\frac{a_l\zeta_l}{z\sqrt{np}}\Mb_{kl}\bx_l\bx_l^\prime\Bb+\Bb\bx_l\bx_l^\prime\frac{a_l\zeta_l}{z\sqrt{np}}\Mb_{kl}-\frac{a_l\zeta_l}{z}\Bb\bx_l\bx_l\Bb\\
	:=&B_1(z)+B_2(z)+B_3(z)+B_4(z)+B_5(z).
	\end{split}
	\end{equation}
	Moreover, $\vartheta_l$ can be further simplified as 
	\[
	\vartheta_l=\frac{a_l}{znp}\bx_l^\prime\Mb_{kl}\bx_l-\frac{a_l}{z\sqrt{np}}(\bx_l^\prime\Bb\bx_l-n\bar\lambda_{\Bb}).
	\]
	Similarly to the proof of Lemma \ref{lemma:gamma and eta}, we claim
	\begin{equation}\label{eq:zeta theta}
	\vartheta_l-\frac{a_l}{znp}\text{tr}\Mb_{kl}\overset{L_4}{\longrightarrow}0,\quad \zeta_l-\frac{1}{1+\frac{a_l}{znp}\text{tr}\Mb_{kl}}\overset{L_4}{\longrightarrow}0.
	\end{equation}
	The result on $\zeta_l$ holds because the imaginary part of $z\zeta_l^{-1}$ is $\Im z+\Im a_l/(np)\bx_l^\prime\Mb_{kl}\bx_l$, which is always larger than $v_0$. Hence,
	\[
	|\zeta_l|\le |z|v_0^{-1},
	\]
	and similarly 
	\[
	\bigg|\frac{1}{1+\frac{a_l}{znp}\text{tr}\Mb_{kl}}\bigg|\le |z|v_0^{-1}.
	\]

	By the above argument, with Burkholder's inequality it's easy to get
	\begin{equation}\label{eq:beta epsilon}
	\mathbb{E}|\tilde\beta_k-\epsilon_k|^2\le \frac{C}{p}.
	\end{equation}
	That is, $\tilde\beta_k-\epsilon_k=o_p(1)$. Note that $|\beta_k|\le C$ and $|\epsilon_k|\le C$. Then  by dominated convergence theorem, for any integer $t$, 
	\[
	\mathbb{E}|\tilde\beta_k-\epsilon_k|^t\rightarrow 0.
	\]
	Therefore,
	\[
	\begin{split}
	&\mathbb{E}\bigg|\sum_{k=1}^p\mathbb{E}_{k-1}\big[\mathbb{E}_k(\tilde\beta_k(z_1)\eta_k(z_1))\cdot\mathbb{E}_k(\tilde\beta_k(z_2)\eta_k(z_2))-\mathbb{E}_k(\epsilon_k(z_1)\eta_k(z_1))\cdot\mathbb{E}_k(\epsilon_k(z_2)\eta_k(z_2))\big]\bigg|\\
	\le& \sum_{k=1}^p\mathbb{E}\bigg|\mathbb{E}_k\big[(\tilde\beta_k(z_1)-\epsilon_k(z_1))\eta_k(z_1)\big]\cdot\mathbb{E}_k(\tilde\beta_k(z_2)\eta_k(z_2))\bigg|\\
	&+\sum_{k=1}^p\mathbb{E}\bigg|\mathbb{E}_k(\tilde\beta_k(z_1)\eta_k(z_1))\cdot\mathbb{E}_k\big[(\tilde\beta_k(z_2)-\epsilon_k(z_2))\eta_k(z_2)\big]\bigg|\\
	\rightarrow &0. 
	\end{split}
	\]
	As a result, it suffices to consider 
	\[
	\frac{\partial^2}{\partial z_2\partial z_1}\sum_{k=1}^p\epsilon_k(z_1)\epsilon_k(z_2)\mathbb{E}_{k-1}[\mathbb{E}_k\eta_k(z_1)\cdot\mathbb{E}_k\eta_k(z_2)].
	\]

	Actually, in the next section, we prove the convergence of the ESD of $\bar\Sbb_p$, and the results (\ref{eq:z Mk expansion}) to (\ref{sz}) therein show that
	\[
	\epsilon_k(z)\rightarrow \big[z+a_k\bar\lambda_{\Bb_n^2}s_p(z)\big]^{-1}:=\tilde\epsilon_k\rightarrow \big[z+a_k\bar\lambda_{\Bb^2}s(z)\big]^{-1},
	\]
	which implies that it suffices to consider
	\begin{equation}\label{eq:cov}
	\frac{\partial^2}{\partial z_2\partial z_1}\sum_{k=1}^p\tilde\epsilon_k(z_1)\tilde\epsilon_k(z_2)\mathbb{E}_{k-1}[\mathbb{E}_k\eta_k(z_1)\cdot\mathbb{E}_k\eta_k(z_2)].
	\end{equation}
	In the following, we further simplify the expression  (\ref{eq:cov}). We write
	\[
	\begin{split}
	\mathbb{E}_k\eta_k(z)=&\frac{a_k}{\sqrt{np}}\sum_{j=1}^n\mathrm{B}_{n,jj}(x_{kj}^2-1)+\frac{a_k}{\sqrt{np}}\sum_{i\ne j}^n\mathrm{B}_{ij}x_{ki}x_{kj}\\
	&-\frac{1}{np}\bigg(\sum_{i\ne j}x_{ki}x_{kj}\mathbb{E}_k\mathrm{M}_{k,ij}^{(1)}(z)+\sum_{i}(x_{ki}^2-1)\mathbb{E}_k\mathrm{M}_{k,ii}^{(1)}(z)\bigg).
	\end{split}
	\]
	Hence, after some calculations
	\[
	\mathbb{E}_{k-1}[\mathbb{E}_k\eta_k(z_1)\cdot \mathbb{E}_k\eta_k(z_2)]=\frac{a_k^2}{np}\bigg(\mathbb{E}(x_{11}^2-1)^2\sum_{j}\mathrm{B}_{n,jj}^2+2\sum_{i\ne j}\mathrm{B}_{ij}^2\bigg)+A_1+A_2+A_3+A_4+A_5+A_6,
	\]
	where 
	\[
	\begin{split}
	A_1=&-\frac{a_k}{np\sqrt{np}}\mathbb{E}(x_{11}^2-1)^2\sum_{j}\mathrm{B}_{n,jj}\mathbb{E}_k\mathrm{M}_{k,jj}^{(1)}(z_1), \quad A_2=-\frac{a_k}{np\sqrt{np}}\mathbb{E}(x_{11}^2-1)^2\sum_{j}\mathrm{B}_{n,jj}\mathbb{E}_k\mathrm{M}_{k,jj}^{(1)}(z_2), \\
	A_3=&\frac{2}{n^2p^2}\sum_{i\ne j}\mathbb{E}_k\mathrm{M}_{k,ij}^{(1)}(z_1)\mathbb{E}_k\mathrm{M}_{k,ij}^{(1)}(z_2),\quad A_4=\frac{1}{n^2p^2}\mathbb{E}(x_{11}^2-1)^2\sum_j\mathbb{E}_k\mathrm{M}_{k,jj}^{(1)}(z_1)\mathbb{E}_k\mathrm{M}_{k,jj}^{(1)}(z_2),\\
	A_5=&-\frac{2a_k}{np\sqrt{np}}\sum_{i\ne j}\mathrm{B}_{ij}\mathbb{E}_k\mathrm{M}_{k,ij}^{(1)}(z_1),\quad A_6=-\frac{2a_k}{np\sqrt{np}}\sum_{i\ne j}\mathrm{B}_{ij}\mathbb{E}_k\mathrm{M}_{k,ij}^{(1)}(z_2).
	\end{split}
	\]
	By the proof of Lemma \ref{lemma:gamma and eta}, we already know that $\mathbb{E}|\mathrm{M}_{k,jj}^{(1)}|^4\le C(p^4+np^2)$, then it's easy to conclude 
	\[
	\mathbb{E}\bigg|\sum_{k=1}^pA_j\bigg|\rightarrow 0, \text{ for } j=1,2,4.
	\]
	Moreover, 
	\[
	\mathbb{E}\bigg|\sum_{k=1}^pA_5\bigg|=\mathbb{E}\bigg|\sum_{k=1}^p\frac{2a_k}{np\sqrt{np}}\text{tr}(\Bb\Mb_{k}^{(1)})\bigg|-o(1)\rightarrow 0.
	\]
	Similar results hold for $A_6$. On the other hand,
	\[
	\sum_{k=1}^pA_3=\frac{2}{p}\sum_{k=1}^p\mathbb{Z}_k-\frac{2}{n^2p^2}\sum_{k=1}^p\sum_{j=1}^n\mathbb{E}_k\mathrm{M}_{k,jj}^{(1)}(z_1)\mathbb{E}_k\mathrm{M}_{k,jj}^{(1)}(z_2)=\frac{2}{p}\sum_{k=1}^p\mathbb{Z}_k+o_{L_1}(1),
	\]
	where $\mathbb{Z}_k$ is defined by
	\[
	\mathbb{Z}_k=\frac{1}{n^2p}\text{tr}\bigg(\mathbb{E}_k\Mb_k^{(1)}(z_1)\cdot \mathbb{E}_k\Mb_k^{(1)}(z_2)\bigg).
	\]
	That is, we only need to find the limit in probability of 
	\[
	\tilde \Lambda_p(z_1,z_2):=\frac{1}{p}\sum_{k=1}^p\tilde\epsilon_k(z_1)\tilde\epsilon_k(z_2)\bigg(2\mathbb{Z}_k+a_k^2\big[(\nu_4-3)\frac{1}{n}\sum_j\mathrm{B}_{n,jj}^2+2\bar\lambda_{\Bb^2}\big]\bigg).
	\]
	
	\subsection{Decomposition for $\mathbb{Z}_k$} 
	Note that
	\[
	\Db_k=\sum_{i\ne k}^p\be_i\bh_i^\prime-z\Ib_{p-1}.
	\]
	Multiplying $\Db_k^{-1}$ to both sides leads to
	\[
	z\Db_k^{-1}=-\Ib_{p-1}+\sum_{i\ne k}^p\be_i\bh_i^\prime\Db_k^{-1}.
	\] 
	Therefore,
	\begin{equation}\label{eq:z Mk}
	\begin{split}
	z\Mb_{k}=&-\Bb^{1/2}\Yb_k^\prime\Yb_k\Bb^{1/2}+\sum_{i\ne k}^p\bigg(\Bb^{1/2}\Yb_k^\prime\be_i\bh_i^\prime\Db_k^{-1}\Yb_k\Bb^{1/2}\bigg)\\=&-\Bb^{1/2}\Yb_k^\prime\Yb_k\Bb^{1/2}+\sum_{i\ne k}^p\bigg(\zeta_i\Bb^{1/2}\by_i\bh_i^\prime\Db_{ki,r}^{-1}(\Yb_{ki}+\be_i\by_i^\prime)\Bb^{1/2}\bigg)\\
	=&-\Bb^{1/2}\Yb_k^\prime\Yb_k\Bb^{1/2}+\sum_{i\ne k}^p\bigg(\zeta_i\Bb^{1/2}\by_i\frac{1}{\sqrt{np}}\by_i^\prime\Yb_{ki}^\prime\Db_{ki}^{-1}\Yb_{ki}\Bb^{1/2}\bigg)\\
	&+\sum_{i\ne k}^p\bigg(\zeta_i\vartheta_i\Bb^{1/2}\by_i\by_i^\prime\Bb^{1/2}\bigg)\\
	=&-\sum_{i\ne k}^p\bigg(\zeta_i\Bb^{1/2}\by_i\by_i^\prime\Bb^{1/2}\bigg)+\sum_{i\ne k}^p\bigg(\zeta_i\Bb^{1/2}\by_i\frac{1}{\sqrt{np}}\by_i^\prime\Yb_{ki}^\prime\Db_{ki }^{-1}\Yb_{ki}\Bb^{1/2}\bigg)\\
	=&-\sum_{i\ne k}^p\bigg(\zeta_i\Bb^{1/2}\by_i\by_i^\prime\Bb^{1/2}\bigg)+\sum_{i\ne k}^p\bigg(\zeta_i\Bb\bx_i\frac{a_i}{\sqrt{np}}\bx_i^\prime\Mb_{ki}\bigg).
	\end{split}
	\end{equation}
	On the other hand, by taking $l=i$ in (\ref{eq:mk1}), we write
	\[
	\begin{split}
	z_1\mathbb{Z}_k=\frac{a_k^2}{n^2p}\text{tr}\bigg(\mathbb{E}_kz_1\Mb_k(z_1)\cdot \mathbb{E}_k\Mb_k(z_2)\bigg):=C_1(z_1,z_2)+C_2(z_1,z_2),
	\end{split}
	\]
	with
	\[
	\begin{split}
	C_1(z_1,z_2)=&-\frac{a_k^2}{n^2p}\sum_{i<k}\mathbb{E}_k\bigg(\zeta_i(z_1)\by_i^\prime\Bb^{1/2}\mathbb{E}_k\big[\sum_{j=1}^5B_j(z_2)\big]\Bb^{1/2}\by_i\bigg)\\
	&-\frac{a_k^2}{n^2p}\sum_{i>k}\mathbb{E}_k\bigg(\zeta_i(z_1)\by_i^\prime\Bb^{1/2}\mathbb{E}_k\Mb_k(z_2)\Bb^{1/2}\by_i\bigg):=\sum_{j=1}^6C_{1j},\\
	C_2(z_1,z_2)=&\frac{a_k^2}{n^2p}\sum_{i<k}\mathbb{E}_k\bigg(\frac{a_i\zeta_i(z_1)}{\sqrt{np}}\bx_i^\prime\Mb_{ki}(z_1)\mathbb{E}_k\big[\sum_{j=1}^5B_j(z_2)\big]\Bb\bx_i\bigg)\\
	&+\frac{a_k^2}{n^2p}\sum_{i>k}\mathbb{E}_k\bigg(\frac{a_i\zeta_i(z_1)}{\sqrt{np}}\bx_i^\prime\Mb_{ki}(z_1)\mathbb{E}_k\Mb_k(z_2)\Bb\bx_i\bigg):=\sum_{j=1}^6C_{2j},
	\end{split}
	\]
	where $C_{1j}$, $C_{2j}$  correspond to $B_j$ for $j=1,\ldots,5$, while $C_{16}$ and $C_{26}$ correspond to $i>k$.  Now we aim to control these terms one by one.
	
	We first introduce some useful bounds. For any $n\times n$ matrix $\Qb$ independent of $\bx_i$, we claim
	\begin{eqnarray}
	\label{mm1}\mathbb{E}(\bx_i^\prime\Qb\bx_i-\text{tr}\Qb)^2&\le &C\mathbb{E}\|\Qb\|_F^2,\\
	\mathbb{E}(\bx_i^\prime\Qb\bx_i)^2&\le& C\mathbb{E}(|\text{tr}\Qb|^2+\|\Qb\|_F^2)\le Cn\mathbb{E}\|\Qb\|_F^2,\\
	\quad\label{mm3}\mathbb{E}(\bx_i^\prime\Mb_{ki}\Qb\bx_i-\text{tr}\Mb_{ki}\Qb)^2&\le &C\mathbb{E}\|\Mb_{ki}\Qb\|_F^2\le Cn^2p\mathbb{E}\|\Qb\|^2,\\
	\label{mm4}\mathbb{E}(\bx_i^\prime\Mb_{ki}\Qb\bx_i)^2&\le &C\mathbb{E}(|\text{tr}\Mb_{ki}\Qb|^2+\|\Mb_{ki}\Qb\|_F^2)\le Cn^2p^2\mathbb{E}\|\Qb\|^2.
	\end{eqnarray}
	The proof of (\ref{mm1}) to (\ref{mm3}) is similar to Lemma \ref{lemma:gamma and eta},  while (\ref{mm4}) is concluded from $|(np)^{-1}\text{tr}\Mb_{ki}\Qb|\le C\|\Qb\|$ whose proof is similar to (\ref{eq:trace mk}). Note that $\mathbb{E}_k\Mb_k(z_2)$ is also independent of $\bx_i$ for $i>k$. 
	Therefore, by the Cauchy-Schwartz inequality,
	\[
	\begin{split}
	\mathbb{E}|C_{1j}|\le &C\sqrt{\frac{p}{n}}, \quad j=1,2,3,4,6,\\
	\mathbb{E}|C_{2j}|\le &C\sqrt{\frac{p}{n}}, \quad j=1,2,4,5,6.
	\end{split} 
	\]
	Hence, we only need to consider $C_{15}$ and $C_{23}$.  Specifically,
	\[
	\begin{split}
	C_{15}=&\frac{a_k^2}{n^2p}\sum_{i<k}\mathbb{E}_k\bigg(\zeta_i(z_1)\by_i^\prime\Bb^{1/2}\times \frac{a_i\mathbb{E}_k\zeta_i(z_2)}{z_2}\Bb\bx_i\bx_i^\prime\Bb\Bb^{1/2}\by_i\bigg)\\
	\overset{i.p.}{\longrightarrow}&z_1a_k^2\bar\lambda^2_{\Bb^2}\times \frac{1}{p}\sum_{i<k}a_i^2\tilde\epsilon_i(z_1)\tilde\epsilon_i(z_2).
	\end{split}
	\]
	
	On the other hand,
	\[
	\begin{split}
	C_{23}=&\frac{a_k^2}{n^2p}\sum_{i<k}\mathbb{E}_k\bigg(\frac{a_i\zeta_i(z_1)}{\sqrt{np}}\bx_i^\prime\Mb_{ki}(z_1)\mathbb{E}_k\bigg[\frac{a_i\zeta_i(z_2)}{z_2\sqrt{np}}\Mb_{ki}\bx_i\bx_i^\prime\Bb\bigg]\Bb\bx_i\bigg)\\
	\overset{i.p.}{\longrightarrow}&z_1\mathbb{Z}_k\bar\lambda_{\Bb^2}\times \frac{1}{p}\sum_{i<k}a_i^2\tilde\epsilon_i(z_1)\tilde\epsilon_i(z_2).
	\end{split}
	\]
	That is,
	\[
	\mathbb{Z}_k\overset{i.p.}{\longrightarrow}a_k^2\bar\lambda^2_{\Bb^2}\times \frac{1}{p}\sum_{i<k}a_i^2\tilde\epsilon_i(z_1)\tilde\epsilon_i(z_2)+\mathbb{Z}_k\bar\lambda_{\Bb^2}\times \frac{1}{p}\sum_{i<k}a_i^2\tilde\epsilon_i(z_1)\tilde\epsilon_i(z_2),
	\]
	which implies
	\[
	\mathbb{Z}_k\overset{i.p.}{\longrightarrow}\frac{a_k^2\bar\lambda^2_{\Bb^2}\times \frac{1}{p}\sum_{i<k}a_i^2\tilde\epsilon_i(z_1)\tilde\epsilon_i(z_2)}{1-\bar\lambda_{\Bb^2}\times \frac{1}{p}\sum_{i<k}a_i^2\tilde\epsilon_i(z_1)\tilde\epsilon_i(z_2)},
	\]
	and 
	\[
	\begin{split}
	&\Lambda(z_1,z_2)\overset{i.p.}{\longrightarrow}\\
	&\lim_{p\rightarrow \infty}\frac{\partial^2}{\partial z_1\partial z_2}\frac{1}{p}\sum_{k=1}^p\tilde\epsilon_k(z_1)\tilde\epsilon_k(z_2)\bigg(\frac{2a_k^2\bar\lambda^2_{\Bb^2}\times \frac{1}{p}\sum_{i<k}a_i^2\tilde\epsilon_i(z_1)\tilde\epsilon_i(z_2)}{1-\bar\lambda_{\Bb^2}\times \frac{1}{p}\sum_{i<k}a_i^2\tilde\epsilon_i(z_1)\tilde\epsilon_i(z_2)}+a_k^2\big[(\nu_4-3)\frac{1}{n}\sum_j\mathrm{B}_{n,jj}^2+2\bar\lambda_{\Bb^2}\big]\bigg).
	\end{split}
	\]

	\subsection{Tightness of $M_{p1}(z)$}
	We end this section with a proof of the tightness of $M_{p1}(z)$. By Burkholder's inequality, 
	\[
	\mathbb{E}\bigg|\sum_{k=1}^p\sum_{j=1}^l\xi_j\mathbb{E}_{k-1}(\alpha_k(z_j))\bigg|^2\le C,
	\]
	which ensures the first condition in Theorem 12.3 of \cite{billingsley2013convergence}. For the second condition,  similarly to \cite{bai2010spectral} and \cite{chen2015clt}, we aim to verify
	\[
	\frac{\mathbb{E}|M_{p1}(z_1)-M_{p1}(z_2)|^2}{|z_1-z_2|^2}\le C,\quad z_1,z_2\in \mathcal{C}_u.
	\] 
	
	By (\ref{eq:mF expansion}), 
	\[
	\begin{split}
	&M_{p1}(z_1)-M_{p1}(z_2)=\sum_k\big(\rho_k(z_1)-\rho_k(z_2)\big)\\
	=&\sum_k-(\mathbb{E}_k-\mathbb{E}_{k-1})\bigg(\beta_k(z_1)(1+\bq_k^\prime\Db_k^{-2}(z_1)\bq_k)-\beta_k(z_2)(1+\bq_k^\prime\Db_k^{-2}(z_2)\bq_k)\bigg).
	\end{split}
	\]
	Note that
	\[
	\begin{split}
	\beta_k(z_1)-\beta_k(z_2)=&\beta_k(z_1)\beta_k(z_2)\bigg(z_2-z_1+\bq_k\Db_k(z_2)^{-1}\bq_k-\bq_k\Db_k(z_1)^{-1}\bq_k\bigg),\\
	\bq_k\Db_k(z_2)^{-1}\bq_k-\bq_k\Db_k(z_1)^{-1}\bq_k=&\bq_k^\prime\Db_k(z_2)^{-1}(z_2-z_1)\Db_k(z_1)^{-1}\bq_k.
	\end{split}
	\]
	Moreover,
	\[
	\begin{split}
	\beta_k(z_1)\beta_k(z_2)=&\tilde\beta_k(z_1)\beta_k(z_2)+\tilde\beta_k(z_1)\beta_k(z_1)\eta_k(z_1)\beta_k(z_2)\\
	=&\tilde\beta_k(z_1)\tilde\beta_k(z_2)+\tilde\beta_k(z_1)\tilde\beta_k(z_2)\beta_k(z_2)\eta_k(z_2)+\tilde\beta_k(z_1)\beta_k(z_1)\eta_k(z_1)\beta_k(z_2).
	\end{split}
	\]
	Then, by Burkholder's inequality and the facts that $(\mathbb{E}_k-\mathbb{E}_{k-1})\tilde\beta_k(z_1)\tilde\beta_k(z_2)=0$, $|\beta_k|\le C$, $\mathbb{E}|\eta_k|^2\le Cp^{-1}$, we conclude
	\[
	\mathbb{E}\bigg|\sum_k(\mathbb{E}_k-\mathbb{E}_{k-1})\beta_k(z_1)\beta_k(z_2)(z_2-z_1)\bigg|^2\le C|z_1-z_2|^2.
	\]
	For the other terms, actually it's very similar though more tedious calculations are necessary. For example, we write
	\[
	\begin{split}
	&(z_2-z_1)\beta_k(z_1)\beta_k(z_2)\bq_k^\prime\Db_k(z_2)^{-1}\Db_k(z_1)^{-1}\bq_k\\
	=&(z_2-z_1)\beta_k(z_1)\beta_k(z_2)\bigg(\bq_k^\prime\Db_k(z_2)^{-1}\Db_k(z_1)^{-1}\bq_k-\frac{a_k}{np}\text{tr}\Bb^{1/2}\Yb_k^\prime\Db_k(z_2)^{-1}\Db_k(z_1)^{-1}\Yb_k\Bb^{1/2}\bigg)\\
	&+(z_2-z_1)\beta_k(z_1)\beta_k(z_2)\frac{a_k}{np}\text{tr}\Bb^{1/2}\Yb_k^\prime\Db_k(z_2)^{-1}\Db_k(z_1)^{-1}\Yb_k\Bb^{1/2}.
	\end{split}
	\]
	The first term can be bounded with Burkholder's inequality. For the second term, use the expansion for $\beta_k(z_1)\beta_k(z_2)$ and note that
	\[
	\begin{split}
	&(\mathbb{E}_k-\mathbb{E}_{k-1})\tilde\beta_k(z_1)\tilde\beta_k(z_2)\frac{a_k}{np}\text{tr}\Bb^{1/2}\Yb_k^\prime\Db_k(z_2)^{-1}\Db_k(z_1)^{-1}\Yb_k\Bb^{1/2}=0.
	\end{split}
	\]
	For the remaining terms, we omit the details.

	\section{The non-random part with diagonal $\Ab_p$}\label{non-random}
	
	\subsection{The limit of $\mathbb{E}m_{F^{\bar{\mathbf{S}}_p}}(z)$}
	In this section, we focus on the non-random part $M_{p2}(z)$ under the special case where $\Ab$ is diagonal.  This is more challenging than the trivial case in \cite{chen2015clt} where $\Ab_p$ and $\Bb_n$ are both identity matrices. First, we show how to find the limit of  $\mathbb{E}(m_{F^{\bar\Sbb_p}}(z))$ using a method which is different from \cite{wang2014limiting}. 
	
	We already know that
	\[
	\begin{split}
	\mathbb{E}\text{tr}(\Db^{-1})
	=-\sum_{k=1}^p\mathbb{E}(\beta_k)=&-\sum_{k=1}^p\mathbb{E}\bigg(\frac{1}{z+\frac{1}{np}\by_k^\prime\Yb_k^\prime\Db_k^{-1}\Yb_k\by_k-\frac{1}{\sqrt{np}}(\by_k^\prime\by_k-n\bar\lambda_{\Bb}a_k)}\bigg).
	\end{split}
	\]
	Let
	\begin{equation}\label{eq:epsilon_k mu_k}
	\begin{split}
	\epsilon_k=&\frac{1}{z+\frac{1}{np}\mathbb{E}\text{tr}\Mb_{k}^{(1)}},\\ \mu_k=&\frac{1}{\sqrt{np}}(\by_k^\prime\by_k-n\bar\lambda_{\Bb}a_k)-\frac{1}{np}\by_k^\prime\Yb_k^\prime\Db_k^{-1}\Yb_k\by_k+\frac{1}{np}\mathbb{E}\text{tr}\Mb_{k}^{(1)}.
	\end{split}
	\end{equation}

	Then, by the identity $\beta_k=\epsilon_k+\epsilon_k\beta_k\mu_k$ and the fact that
	\[
	\frac{1}{p}\sum_{k=1}^p\mathbb{E}|\epsilon_k\beta_k\mu_k|\le \frac{C}{p}\sum_k\big(\mathbb{E}|\mu_k|^2\big)^{1/2}\rightarrow 0,
	\]
	we have
	\begin{equation}\label{eq:beta minus epsion}
	\mathbb{E}(m_{F^{\bar\Sbb_p}}(z))=-\frac{1}{p}\sum_k\epsilon_k+o(1).
	\end{equation}
	Moreover, by (\ref{eq:z Mk}), we have
	\[
	\frac{z}{np}\mathbb{E}\text{tr}\Mb_{k}^{(1)}=\frac{a_k}{np}\mathbb{E}\text{tr}\bigg[-\sum_{i\ne k}^p\bigg(\zeta_i\Bb^{1/2}\by_i\by_i^\prime\Bb^{1/2}\bigg)+\sum_{i\ne k}^p\bigg(\zeta_i\Bb\bx_i\frac{a_i}{\sqrt{np}}\bx_i^\prime\Mb_{ki}\bigg)\bigg].
	\]
	Recall that
	\[
	\zeta_i-\frac{1}{1+\frac{a_i}{znp}\text{tr}\Mb_{ki}}\overset{L_4}{\longrightarrow}0.
	\]
	Hence, after some calculations, we have
	\begin{equation}\label{eq:z Mk expansion}
	\begin{split}
	\frac{z}{np}\mathbb{E}\text{tr}\Mb_{k}^{(1)}=&-\frac{a_k}{np}\sum_{i\ne k}^p\mathbb{E}\frac{a_i\text{tr}\Bb^2}{1+\frac{a_i}{znp}\text{tr}\Mb_{ki}}+o(1)
	=za_k\bar\lambda_{\Bb^2}\frac{1}{p}\mathbb{E}\text{tr}\Db_k^{-1}\Ab_k+o(1).
	\end{split}
	\end{equation}
	On the other hand, by Burkholder's inequality and the  expansion in (\ref{eq:mk1}), it's not difficult to verify 
	\[
	\mathbb{E}\bigg|\frac{1}{np}\text{tr}\Mb_{ki}-\frac{1}{np}\mathbb{E}\text{tr}\Mb_{ki}\bigg|^2\rightarrow 0,\quad \bigg|\frac{1}{np}\mathbb{E}\text{tr}\Mb_{k}-\frac{1}{np}\mathbb{E}\text{tr}\Mb_{ki}\bigg|^2\rightarrow 0.
	\]
	Therefore, (\ref{eq:z Mk expansion}) also implies 
	\begin{equation}\label{eq: Dk Ak}
	\begin{split}
	\frac{1}{p}\mathbb{E}\text{tr}\Db_k^{-1}\Ab_k=
	=&-\frac{1}{p}\sum_{j\ne k}\frac{a_j}{z+a_j\bar\lambda_{\Bb^2}\frac{1}{p}\mathbb{E}\text{tr}\Db_{k}^{-1}\Ab_{k}}+o(1).
	\end{split}
	\end{equation}

	Let $s_p(z)$ and $m_{p}(z)$ be the respective solutions in $\mathbb{C}^+$ to 
	\begin{equation}\label{eq:s and m}
	s_p(z)=-\int \frac{x}{z+x\bar\lambda_{\Bb_n^2}s_p(z)}dF^{\Ab_p}, \quad m_p(z)=-\int \frac{1}{z+x\bar\lambda_{\Bb_n^2}s_p(z)}dF^{\Ab_p}\quad z\in\mathbb{C}^+.
	\end{equation}
	Then, as $p\rightarrow\infty$, for $s(z)$ defined in Lemma \ref{lem:wang},
	\[
	\begin{split}
	s_p(z)-s(z)=&\int\bigg(\frac{x}{z+x\bar\lambda_{\Bb^2}s(z)}-\frac{x}{z+x\bar\lambda_{\Bb_n^2}s_p(z)}\bigg)dF^{\Ab_p}+\int \frac{x}{z+x\bar\lambda_{\Bb^2}s(z)}d(F^{\Ab}-F^{\Ab_p})\\
	=&\int\bigg(\frac{x^2\lambda_{\Bb^2}[s_p(z)-s(z)]}{[z+x\bar\lambda_{\Bb^2}s(z)][z+x\bar\lambda_{\Bb^2}s_p(z)]}\bigg)dF^{\Ab_p}+o(1),
	\end{split}
	\]
	where we use the fact that $F^{\Ab_p}\rightarrow F^{\Ab}$ and for $x\in \text{supp}(F^{\Ab})\cup \text{supp}(F^{\Ab_p})$,
	\[
	\bigg|\frac{x}{z+x\bar\lambda_{\Bb^2}s(z)}\bigg|\le \frac{x}{\Im z}\le C.
	\]
	Moreover, by the definition in (\ref{eq:s and m}),   there exists some positive constant $C_z$ satisfying $|s_p(z)|\le C_z$. Similar conclusion holds for $s(z)$. Then, considering the imaginary part, for sufficiently large $p$ we have
	\[
	\begin{split}
	&\bigg|1-\int\bigg(\frac{x^2\lambda_{\Bb^2}}{[z+x\bar\lambda_{\Bb^2}s(z)][z+x\bar\lambda_{\Bb^2}s_p(z)]}\bigg)dF^{\Ab_p}\bigg|
	\ge \int\frac{\lambda_{\Bb^2}x^2}{4[\|z\|^2+(x\lambda_{\Bb^2}C_z)^2]}dF^{\Ab_p}>0.
	\end{split}
	\]
	Then, as $p\rightarrow \infty$
	\[
	s_p(z)-s(z)\rightarrow 0.
	\]
	Similarly, we have
	\begin{equation}\label{sz}
	\quad s(z)-\frac{1}{p}\mathbb{E}\text{tr}\Db_k^{-1}\Ab_k\rightarrow0, \quad m_p(z)\rightarrow m(z),\quad \frac{1}{p}\mathbb{E}\text{tr}\Db^{-1}-m(z)\rightarrow 0.
	\end{equation}
	where $m(z)$ is defined in Lemma \ref{lem:wang}.
	
	\subsection{Convergence of $\mathbb{E}M_{p2}(z)$}
	Now we aim to find the limit of $\mathbb{E}M_{p2}(z)$, which is more challenging due to the multiplication with $p$. That is, we need to study the $o(1)$ terms in (\ref{eq:beta minus epsion}), (\ref{eq:z Mk expansion}) and (\ref{eq: Dk Ak}). For simplicity, we write $m_{\bar F_p}$ for $m_{F^{\bar\Sbb_p}}(z)$ and define 
	\[
	\begin{split}
	\bar\epsilon_k=&\frac{1}{z+\frac{a_k}{np}\mathbb{E}\text{tr}\Mb}\rightarrow\frac{1}{z+a_k\bar\lambda_{\Bb^2}s_p(z)}:=\tilde\epsilon_k,\quad \Mb:=\Bb^{1/2}\Yb^\prime\Db^{-1}\Yb\Bb^{1/2},\\ \bar\mu_k=&\frac{1}{\sqrt{np}}(\by_k^\prime\by_k-n\bar\lambda_{\Bb}a_k)-\frac{1}{np}\by_k^\prime\Yb_k^\prime\Db_k^{-1}\Yb_k\by_k+\frac{a_k}{np}\mathbb{E}\text{tr}\Mb.
	\end{split}
	\]
	Note the difference among $\epsilon_k,\bar\epsilon_k$ and $\tilde\epsilon_k$. Then, 
	\begin{equation}\label{eq:mF-mp}
	\begin{split}
	\mathcal{W}_1:=\mathbb{E}m_{\bar F_p}-m_p=&-\frac{1}{p}\sum_{k=1}^p\mathbb{E}\beta_k+\frac{1}{p}\sum_k\frac{1}{z+a_k\bar\lambda_{\Bb^2}s_p(z)}\\
	=&-\mathbb{E}\frac{1}{p}\sum_{k=1}^p\beta_k\tilde\epsilon_k\bigg[\bar\mu_k+a_k\big(\bar\lambda_{\Bb^2}s_p(z)-\frac{1}{np}\mathbb{E}\text{tr}\Mb\big)\bigg].
	\end{split}
	\end{equation}
	
	Before moving forward, we first introduce several bounds in the next lemma which are useful in later proofs.  The proof of this lemma is postponed to the next  subsection.
	\begin{lemma}\label{lem: bounds}
		Under Condition \ref{c1} and (\ref{eq:truncation}), for $z\in \mathcal{C}_u$, 
		\[
		\begin{split}
		\frac{1}{np}|\text{tr}\Mb\Bb^{h}|\le &C,\text{ for any fixed } h\ge 0,\\
		\bigg|\frac{1}{np}(\mathbb{E}\text{tr}\Mb-\mathbb{E}\text{tr}\Mb_k)\bigg|\le &\frac{C}{p},\text{ for any }k,\\
		\mathbb{E}|\bar\mu_k|^3=&o_(p^{-1}),\quad \mathbb{E}|\bar\mu_k|^4=o(p^{-1}).
		\end{split}
		\]
		The results also hold if we replace $\Mb$ and $\Mb_k$ with $\Mb_k$ and $\Mb_{ki}$, respectively.
	\end{lemma}
	
	We first focus on $p^{-1}\mathbb{E}\sum_k\beta_k\tilde\epsilon_k\bar\mu_k$. By the decomposition $\beta_k=\bar\epsilon_k+\bar\epsilon_k^2\bar\mu_k+\bar\epsilon_k^3\bar\mu_k^2+\beta_k\bar\epsilon_k^3\bar\mu_k^3$,
	\[
	\mathbb{E}\sum_{k=1}^p\beta_k\tilde\epsilon_k\bar\mu_k=\sum_{k=1}^p\tilde\epsilon_k\bigg(\bar\epsilon_k\mathbb{E}\bar\mu_k+\bar\epsilon_k^2\mathbb{E}\bar\mu_k^2+\bar\epsilon_k^3\mathbb{E}\bar\mu_k^3+\bar\epsilon_k^3\mathbb{E}\beta_k\mu_k^4\bigg):=H_1+H_2+H_3+H_4.
	\]
	For $H_4$, Lemma \ref{lem: bounds} implies $\mathbb{E}|\beta_k||\bar\mu_k|^4=o(p^{-1})$ , where the $o(p^{-1})$ is actually uniform over $k$.  A similar conclusion holds for $H_3$. For $H_1$, 
	\[
	\mathbb{E}\bar\mu_k=\frac{a_k}{np}\bigg(\mathbb{E}\text{tr}\Mb-\mathbb{E}\text{tr}\Mb_k\bigg).
	\]
	By a similar expansion  in (\ref{eq:mk1}), we have
	\[
	\begin{split}
	\mathbb{E}\bar\mu_k=&\frac{a_k}{np}\mathbb{E}\bigg(-\frac{a_k\zeta_k}{znp}\bx_k^\prime\Mb_k^2\bx_k+\frac{a_k\zeta_k}{z\sqrt{np}}\bx_k^\prime\Bb\Mb_k\bx_k+\frac{a_k\zeta_k}{z\sqrt{np}}\bx_k^\prime\Mb_k\Bb\bx_k-\frac{a_k\zeta_k}{z}\bx_k^\prime\Bb^2\bx_k\bigg).
	\end{split}
	\]
	Therefore, after some calculations
	\[
	\begin{split}
	\sum_k\tilde\epsilon_k\bar\epsilon_k\mathbb{E}\bar\mu_k=&\sum_k\frac{a_k\tilde\epsilon_k\bar\epsilon_k}{np}\mathbb{E}\bigg(-\frac{a_k\zeta_k}{znp}\text{tr}\Mb_k^2-\frac{a_k\zeta_k}{z}\text{tr}\Bb^2\bigg)+O\bigg(\sqrt{\frac{p}{n}}\bigg)\\
	\rightarrow &\sum_k\frac{a_k^2\tilde\epsilon_k^3}{np}\bigg(-\frac{1}{np}\mathbb{E}\text{tr}\Mb_k^2-\text{tr}\Bb^2\bigg)+o(1).
	\end{split}
	\]
	Moreover, by (\ref{eq:mk1}) and (\ref{eq:z Mk}), 
	\[
	\begin{split}
	\frac{z}{n^2p}\mathbb{E}\text{tr}\Mb_k^2
	=&\frac{z}{n^2p}\mathbb{E}\bigg(\sum_{i}a_i^2\tilde\epsilon_i^2(\bx_i^\prime\Bb^2\bx_i)^2\bigg)+\frac{z}{n^2p}\mathbb{E}\bigg(\sum_{i}\frac{a_i^2\tilde\epsilon_i^2}{np}\bx_i^\prime\Mb_{ki}^2\bx_i\bx_i^\prime\Bb^2\bx_i\bigg)+o(1).
	\end{split}
	\]
	Then, we conclude that
	\[
	\frac{1}{n^2p}\mathbb{E}\text{tr}\Mb_k^2=\frac{1}{p}\sum_{i}\frac{a_i^2\bar\lambda^2_{\Bb^2}}{(z+a_i\bar\lambda_{\Bb^2}s_p(z))^2}+\frac{1}{p}\sum_{i}\frac{a_i^2\bar\lambda_{\Bb^2}}{(z+a_i\bar\lambda_{\Bb^2}s_p(z))^2}\times \frac{1}{n^2p}\mathbb{E}\text{tr}\Mb_k^2+o(1).
	\]
	On the other hand,  $(n^2p)^{-1}\mathbb{E}\text{tr}\Mb_k^2=(n^2p)^{-1}\mathbb{E}\text{tr}\Mb^2+o(1)$.
	
	Now we calculate $H_2$. By the definition of $\bar\mu_k$,
	\[
	\begin{split}
	\mathbb{E}(\bar\mu_k)^2=&\mathbb{E}(\bar\mu_k-\mathbb{E}\bar\mu_k)^2+O(p^{-2}),
	\mathbb{E}(\bar\mu_k-\mathbb{E}\bar\mu_k)^2=S_1+S_2,
	\end{split}
	\]
	where 
	\[
	\begin{split}
	S_1=&\frac{1}{np}\mathbb{E}(\by_k^\prime\by_k-n\bar\lambda_{\Bb}a_k)^2+\mathbb{E}\gamma_{k1}^2,\quad S_2=S_{21}+S_{22},\\
	S_{21}=&\frac{a_k^2}{n^2p^2}\mathbb{E}(\text{tr}\Mb_k-\mathbb{E}\text{tr}\Mb_k)^2,\\
	S_{22}=&-\frac{2a_k}{np\sqrt{np}}\mathbb{E}(\by_k^\prime\by_k-n\bar\lambda_{\Bb}a_k)(\bx_k^\prime\Mb_k\bx_k-\mathbb{E}\text{tr}\Mb_k).
	\end{split}
	\]
	We start with $S_{21}$. By Burkholder's inequality and (\ref{eq:mk1}),
	\[
	|S_{21}|\le \frac{1}{n^2p^2}\sum_{i\ne k}\mathbb{E}\bigg(\mathbb{E}_i-\mathbb{E}_{i-1})(\text{tr}\Mb_k-\text{tr}\Mb_{ki})\bigg)^2=o(p^{-1}).
	\]
	On the other hand,
	\[
	\begin{split}
	|S_{22}|\asymp&\frac{1}{np\sqrt{np}}\mathbb{E}\bigg(\sum_{i}\mathrm{B}_{ii}(x_{ki}^2-1)+\sum_{i\ne j}\mathrm{B}_{ij}x_{ki}x_{kj}\bigg)\bigg(\sum_{i}\mathrm{M}_{k,ii}(x_{ki}^2-1)+\sum_{i\ne j}\mathrm{M}_{k,ij}x_{ki}x_{kj}\bigg)\\
	\le &\frac{1}{np\sqrt{np}}\mathbb{E}\bigg((\nu_4-2)\sum_i\mathrm{B}_{ii}\mathrm{M}_{k,ii}+2\text{tr}\Mb_k\Bb\bigg)\le O(\frac{1}{\sqrt{np}})=o(p^{-1}).
	\end{split}
	\]
	Hence, we only need to consider $S_1$. Note that
	\[
	\frac{1}{np}\mathbb{E}(\by_k^\prime\by_k-n\bar\lambda_{\Bb}a_k)^2=\frac{a_k^2}{np}\mathbb{E}\bigg(\sum_{i}\mathrm{B}_{ii}(x_{ki}^2-1)+\sum_{i\ne j}\mathrm{B}_{ij}x_{ki}x_{kj}\bigg)^2=\frac{a_k^2}{p}\bigg((\nu_4-3)\frac{1}{n}\sum_j\mathrm{B}_{n,jj}^2+2\bar\lambda_{\Bb^2}\bigg). 
	\]
	Similarly,
	\[
	\sum_{k}\mathbb{E}\gamma_{k1}^2=\frac{1}{p}\sum_k\frac{2a_k^2}{n^2p}\mathbb{E}\text{tr}\Mb_k^2+o(1).
	\]
	As a conclusion,
	\[
	\mathbb{E}\sum_{k=1}^p\beta_k\tilde\epsilon_k\bar\mu_k=\frac{1}{p}\sum_ka_k^2\tilde\epsilon_k^3\bigg(\frac{\nu_4-3}{n}\sum_j\mathrm{B}_{n,jj}^2+\bar\lambda_{\Bb^2}+\frac{1}{n^2p}\mathbb{E}\text{tr}\Mb^2\bigg)+o(1):=\mathcal{A}\times \frac{1}{p}\sum_ka_k^2\tilde\epsilon_k^3+o(1).
	\]

	Hence, 
	\begin{equation}\label{w1}
	\mathcal{W}_1=-\frac{1}{p}\mathcal{A}\times \frac{1}{p}\sum_ka_k^2\tilde\epsilon_k^3-\mathcal{W}_2\times \frac{1}{p}\sum_{k=1}^p\mathbb{E}a_k\beta_k\tilde\epsilon_k+o(p^{-1}),\quad \mathcal{W}_2:=\big(\bar\lambda_{\Bb^2}s_p(z)-\frac{1}{np}\mathbb{E}\text{tr}\Mb\big).
	\end{equation}

	Now we aim to find the limit of $p\mathcal{W}_2$.  By a similar expansion  in (\ref{eq:z Mk}),
	\[
	\begin{split}
	zp\mathcal{W}_2=&\frac{1}{n}\sum_{i}\mathbb{E}\text{tr}\bigg(\zeta_i\Bb^{1/2}\by_i\by_i^\prime\Bb^{1/2}\bigg)-\frac{1}{n}\mathbb{E}\sum_{i}\text{tr}\bigg(\zeta_i\Bb\bx_i\frac{a_i}{\sqrt{np}}\bx_i^\prime\Mb_i\bigg)+p\bar\lambda_{\Bb^2}s_p(z).
	\end{split}
	\]
	Actually, here $z^{-1}\zeta_i$ is exactly $\beta_i$, then
	\[
	\begin{split}
	p\mathcal{W}_2=&\sum_{i}a_i\mathbb{E}(\beta_i-\tilde\epsilon_i)\frac{1}{n}\bx_i^\prime\Bb^2\bx_i-\frac{1}{n}\mathbb{E}\sum_{i}\text{tr}\bigg(\beta_i\Bb\bx_i\frac{a_i}{\sqrt{np}}\bx_i^\prime\Mb_i\bigg)+o(1).
	\end{split}
	\]
	Note that
	\[
	\begin{split}
	&\sum_{i}a_i\mathbb{E}(\beta_i-\tilde\epsilon_i)\bigg(\frac{1}{n}\bx_i^\prime\Bb^2\bx_i-\bar\lambda_{\Bb^2}\bigg)=\sum_{i}a_i\mathbb{E}(\beta_i-\bar\epsilon_i)\bigg(\frac{1}{n}\bx_i^\prime\Bb^2\bx_i-\bar\lambda_{\Bb^2}\bigg)\\
	=&\sum_{i}a_i\mathbb{E}\beta_i\bar\epsilon_i\bar\mu_i\bigg(\frac{1}{n}\bx_i^\prime\Bb^2\bx_i-\bar\lambda_{\Bb^2}\bigg)\le O(\sqrt{p/n})\rightarrow 0.
	\end{split}
	\]
	Therefore, 
	\begin{equation}\label{pw2}
	\begin{split}
	p\mathcal{W}_2=&\bar\lambda_{\Bb^2}\sum_{i}a_i\mathbb{E}(\beta_i-\tilde\epsilon_i)-\frac{1}{n}\mathbb{E}\sum_{i}\text{tr}\bigg(\beta_i\Bb\bx_i\frac{a_i}{\sqrt{np}}\bx_i^\prime\Mb_i\bigg)+o(1).
	\end{split}
	\end{equation}
	
	We start with the first term. Similarly to the previous proof, we have
	\[
	\begin{split}
	&\sum_{i}a_i\mathbb{E}(\beta_i-\tilde\epsilon_i)
	=\mathcal{A}\times \frac{1}{p}\sum_ia_i^3\tilde\epsilon_i^3+\mathcal{W}_2\times \sum_ia_i^2\tilde\epsilon_i\mathbb{E}\beta_i+o(1).
	\end{split}
	\]
	On the other hand, by (\ref{eq:z Mk expansion}) to (\ref{sz}), 
	$\mathcal{W}_2=o(1)$. Then, we can  further expand $\sum_ia_i^2\tilde\epsilon_i\mathbb{E}\beta_i$ and repeat the procedure iteratively. In this process, the  $o(1)$ terms are summable and will still be $o(1)$. Then,
	\[
	p\mathcal{W}_2=\mathcal{A}\times \frac{1}{p}\sum_ia_i^3\tilde\epsilon_i^3+ \sum_{j=1}^\infty\mathcal{B}_j-\frac{1}{n}\mathbb{E}\sum_{i}\text{tr}\bigg(\beta_i\Bb\bx_i\frac{a_i}{\sqrt{np}}\bx_i^\prime\Mb_i\bigg)+o(1),
	\]
	where 
	\[
	\mathcal{B}_j=\mathcal{W}_2^{j}\sum_i(a_i\tilde\epsilon_i)^{j+1}.
	\]
	That is, in equation (\ref{pw2}),
	\begin{equation}\label{abj}
	\sum_{i}a_i\mathbb{E}(\beta_i-\tilde\epsilon_i)=\mathcal{A}\times \frac{1}{p}\sum_ia_i^3\tilde\epsilon_i^3+ps_{p}(z)+\sum_i\frac{a_i\tilde\epsilon_i}{1-\mathcal{W}_2a_i\tilde\epsilon_i}+o(1):=\mathcal{D}(\mathcal{W}_2)+o(1).
	\end{equation}

	Now we move to the second term in (\ref{pw2}),
	\[
	\frac{1}{n}\sum_{i=1}^p\mathbb{E}\bigg(\beta_i\frac{a_i}{\sqrt{np}}\bx_i^\prime\Mb_i\Bb\bx_i\bigg)=\frac{1}{n}\sum_{i=1}^p\mathbb{E}\bigg((\tilde\epsilon_i+\beta_i\tilde\epsilon_i(\bar\mu_i+a_i\mathcal{W}_2)\frac{a_i}{\sqrt{np}}\bx_i^\prime\Mb_i\Bb\bx_i\bigg):=J_1+J_2+J_3.
	\]
	For $J_2$, with the bounds on $|\beta_i|$ and $|\tilde\epsilon_i|$ we observe that
	\[
	\bigg|J_2-\frac{1}{n}\sum_{i=1}^p\mathbb{E}\beta_i\tilde\epsilon_i\bar\mu_i\frac{a_i}{\sqrt{np}}\text{tr}\Mb_i\Bb\bigg|\le \frac{1}{n}\sum_{i=1}^p\mathbb{E}|\beta_i||\tilde\epsilon_i||\bar\mu_i|\frac{a_i}{\sqrt{np}}\bigg|\bx_i^\prime\Mb_i\Bb\bx_i-\text{tr}\Mb_i\Bb\bigg|=o(1).
	\]
	Moreover,
	\[
	\frac{1}{n}\sum_{i=1}^p\mathbb{E}\beta_i\tilde\epsilon_i\bar\mu_i\frac{a_i}{\sqrt{np}}\text{tr}\Mb_i\Bb=\frac{1}{n}\sum_{i=1}^p\mathbb{E}(\bar\epsilon_i+\bar\epsilon_i\beta_i\bar\mu_i)\tilde\epsilon_i\bar\mu_i\frac{a_i}{\sqrt{np}}\text{tr}\Mb_i\Bb:=J_{21}+J_{22}.
	\]
	For $J_{22}$, 
	\[
	|J_{22}|\le \frac{1}{n}\sum_{i=1}^p\mathbb{E}|\bar\epsilon_i||\beta_i||\tilde\epsilon_i||\bar\mu_i|^2\bigg|\frac{a_i}{\sqrt{np}}\text{tr}\Mb_i\Bb\bigg|=o(1).
	\]
	On the other hand, for $J_{21}$, note that
	\[
	\begin{split}
	\frac{1}{n\sqrt{np}}\mathbb{E}\bar\mu_i\text{tr}\Mb_i\Bb=&\frac{a_i}{n\sqrt{np}}\mathbb{E}\bigg(\frac{1}{np}\text{tr}\Mb_i-\frac{1}{np}\mathbb{E}\text{tr}\Mb\bigg)\text{tr}\Mb_i\Bb\\
	=&\frac{a_i}{n\sqrt{np}}\mathbb{E}\bigg(\frac{1}{np}\text{tr}\Mb_i-\frac{1}{np}\mathbb{E}\text{tr}\Mb_i\bigg)\text{tr}\Mb_i\Bb+\frac{a_i}{n\sqrt{np}}\mathbb{E}\bigg(\frac{1}{np}\mathbb{E}\text{tr}\Mb_i-\frac{1}{np}\mathbb{E}\text{tr}\Mb\bigg)\text{tr}\Mb_i\Bb\\
	:=&J_{211}+J_{212}.
	\end{split}
	\]
	For $J_{211}$, note that
	\[
	\begin{split}
	J_{211}=&\frac{a_i}{n\sqrt{np}}\mathbb{E}\bigg(\frac{1}{np}\text{tr}\Mb_i-\frac{1}{np}\mathbb{E}\text{tr}\Mb_i\bigg)\bigg(\text{tr}\Mb_i\Bb-\mathbb{E}\text{tr}\Mb_i\Bb\bigg)\\
	\le &C\sqrt{\frac{p}{n}}\sqrt{\mathbb{E}\bigg(\frac{1}{np}\text{tr}\Mb_i-\frac{1}{np}\mathbb{E}\text{tr}\Mb_i\bigg)^2\times \mathbb{E}\bigg(\frac{1}{np}\text{tr}\Mb_i\Bb-\frac{1}{np}\mathbb{E}\text{tr}\Mb_i\Bb\bigg)^2}=o(p^{-1}),
	\end{split}
	\]
	where we use the expansion in (\ref{eq:mk1}), Burkholde's inequality and the martingale decomposition
	\[
	\frac{1}{np}\text{tr}\Mb_i-\frac{1}{np}\mathbb{E}\text{tr}\Mb_i=\sum_{j}(\mathbb{E}_j-\mathbb{E}_{j-1})\bigg(\frac{1}{np}\text{tr}\Mb_i-\frac{1}{np}\text{tr}\Mb_{ij}\bigg).
	\]
	For $J_{212}$, similarly by (\ref{eq:mk1}),
	\[
	|J_{212}|\le C\sqrt{\frac{p}{n}}\bigg(\frac{1}{np}\mathbb{E}\text{tr}\Mb_i-\frac{1}{np}\mathbb{E}\text{tr}\Mb\bigg)\frac{1}{np}\mathbb{E}\text{tr}\Mb_i\Bb=o(p^{-1}).
	\]
	As a result, we conclude that $|J_2|=o(1)$. For $J_3$, we can expand $\beta_i$ once again and write
	\[
	J_3=\mathcal{W}_2\frac{1}{n}\sum_{i=1}^p\mathbb{E}a_i\tilde\epsilon_i\bigg(\tilde\epsilon_i+\tilde\epsilon_i\beta_i(\bar\mu_i+\mathcal{W}_2)\bigg)\frac{a_i}{\sqrt{np}}\bx_i^\prime\Mb_i\Bb\bx_i:=J_{31}+J_{32}+J_{33}.
	\]
	Similarly to $J_2$, we can show that $J_{32}=o(W_2)$. Furthermore, we can keep expanding $\beta_i$ in $J_{33}$ and repeat the procedure iteratively. At each step $k\ge 1$, we get a new error term
	\[
	\mathcal{W}_2^{k-1}\frac{1}{n}\sum_{i=1}^p(a_i\tilde\epsilon_i)^k\frac{1}{\sqrt{np}}\mathbb{E}\text{tr}\Mb_i\Bb,
	\]
	and a negligible term $o(\mathcal{W}_2^{k-1})$. Since $\mathcal{W}_2\rightarrow 0$, such negligible terms are summable and we conclude
	\begin{equation}\label{j3}
	\begin{split}
	\frac{1}{n}\sum_{i=1}^p\mathbb{E}\bigg(\beta_i\frac{a_i}{\sqrt{np}}\bx_i^\prime\Mb_i\Bb\bx_i\bigg)=&\frac{1}{\sqrt{np}}\mathbb{E}\text{tr}\Mb_i\Bb\times \bigg(\sum_{k=1}^\infty\mathcal{W}_2^{k-1}\frac{1}{n}\sum_{i=1}^p(a_i\tilde\epsilon_i)^k\bigg)+o(1)\\
	=&\frac{1}{\sqrt{np}}\mathbb{E}\text{tr}\Mb_i\Bb\times \frac{1}{n}\sum_{i=1}^p\frac{a_i\tilde\epsilon_i}{1-\mathcal{W}_2a_i\tilde\epsilon_i}+o(1)\\
	=&\sqrt{\frac{p}{n}}\times\frac{1}{n}\mathbb{E}\text{tr}\Mb\Bb\times \frac{1}{p}\sum_{i=1}^p\frac{a_i\tilde\epsilon_i}{1-\mathcal{W}_2a_i\tilde\epsilon_i}+o(1).
	\end{split}
	\end{equation}

	Now we focus on the term $\mathbb{E}\text{tr}\Mb\Bb$. Based on (\ref{eq:z Mk}),
	\[
	\begin{split}
	\sqrt{\frac{p}{n}}\times\frac{1}{n}\mathbb{E}\text{tr}\Mb\Bb=&-\sqrt{\frac{p}{n}}\bar\lambda_{\Bb^3}\mathbb{E}\sum_ia_i\beta_i+\sqrt{\frac{p}{n}}\frac{1}{n}\sum_{i=1}^p\mathbb{E}\bigg(\beta_i\frac{a_i}{\sqrt{np}}\bx_i^\prime\Mb_i\Bb^2\bx_i\bigg)+o(1)\\
	:=&P_1+P_2+o(1).
	\end{split}
	\]
	For $P_1$, writing $\beta_i=\beta_i-\tilde\epsilon_i+\tilde\epsilon_i$ and using (\ref{abj}),
	\[
	P_1=-\sqrt{\frac{p}{n}}\bar\lambda_{\Bb^3}\mathcal{D}(\mathcal{W}_2)+\sqrt{\frac{p}{n}}\bar\lambda_{\Bb^3}ps_p(z)+o(1)=-\sqrt{\frac{p}{n}}\bar\lambda_{\Bb^3}\times p\mathcal{D}_2(\mathcal{W}_2)+o(1),
	\]
	where 
	\[
	\mathcal{D}_2(\mathcal{W}_2):= \frac{1}{p}\sum_{i=1}^p\frac{a_i\tilde\epsilon_i}{1-\mathcal{W}_2a_i\tilde\epsilon_i}.
	\]
	For $P_2$, similarly to (\ref{j3}), we conclude 
	\[
	P_2=\frac{p}{n}\frac{1}{n}\mathbb{E}\text{tr}\Mb\Bb^2\times \frac{1}{p}\sum_{i=1}^p\frac{a_i\tilde\epsilon_i}{1-\mathcal{W}_2a_i\tilde\epsilon_i}+o(1).
	\]
	Hence, we can keep expanding $\text{tr}\Mb\Bb^2$ iteratively, and at each step $k\ge 1$ we get a term
	\[
	-p\times\bigg(\sqrt{\frac{p}{n}}\bigg)^k\bar\lambda_{\Bb^{k+2}}\times\bigg( \frac{1}{p}\sum_{i=1}^p\frac{a_i\tilde\epsilon_i}{1-\mathcal{W}_2a_i\tilde\epsilon_i}\bigg)^{k+1},
	\]
	while the other error terms are summable to $o(1)$. Hence, we conclude
	\[
	\begin{split}
	-\frac{1}{n}\sum_{i=1}^p\mathbb{E}\bigg(\beta_i\frac{a_i}{\sqrt{np}}\bx_i^\prime\Mb_i\Bb\bx_i\bigg)=&\bigg(\frac{1}{n}\sum_{j=1}^n\frac{(\lambda_j^{\Bb})^2}{1-\lambda_j^{\Bb}\sqrt{\frac{p}{n}}\mathcal{D}_2(\mathcal{W}_2)}-\bar\lambda_{\Bb^2}\bigg)\times p\mathcal{D}_2(\mathcal{W}_2)+o(1).
	\end{split}
	\]

	Now let $\mathcal{Y}_p$ be  the solution to 
	\[
	x=\bar\lambda_{\Bb^2}\times\frac{1}{p}\mathcal{D}(x)+\bigg(\frac{1}{n}\sum_{j=1}^n\frac{(\lambda_j^{\Bb})^2}{1-\lambda_j^{\Bb}\sqrt{\frac{p}{n}}\mathcal{D}_2(x)}-\bar\lambda_{\Bb^2}\bigg)\mathcal{D}_2(x),
	\]
	and satisfy $\mathcal{Y}_p=o(1)$. Then,
	\[
	p(\mathcal{W}_2-\mathcal{Y}_p)=o(1),
	\]
	which further implies
	\[
	p(\mathcal{W}_1+\frac{1}{p}\mathcal{A}\times\frac{1}{p}\sum_ka_k^2\tilde\epsilon_k^3+\frac{1}{p}\sum_ka_k\tilde\epsilon_k\beta_k\times \mathcal{Y}_p)=o(1).
	\]
	Moreover, note that $\mathcal{Y}_p=o(1)$, then similarly to (\ref{abj}), we can prove that
	\[
	\begin{split}
	\mathcal{Y}_p\sum_{k=1}^pa_k\tilde\epsilon_k\beta_k=&\mathcal{Y}_p\bigg(\sum_{j=1}^\infty \mathcal{W}_2^{j-1}\sum_{k=1}^pa_k^j\tilde\epsilon_k^{j+1}\bigg)+o(1)=\mathcal{Y}_p\sum_k\frac{a_k\tilde\epsilon_k^2}{1-\mathcal{W}_2a_k\tilde\epsilon_k}+o(1)\\
	=&\mathcal{Y}_p\sum_k\frac{a_k\tilde\epsilon_k^2}{1-\mathcal{Y}_pa_k\tilde\epsilon_k}+o(1).
	\end{split}
	\]
	Therefore,
	\[
	p\bigg(\mathcal{W}_1+\frac{1}{p}\mathcal{A}\times\frac{1}{p}\sum_ka_k^2\tilde\epsilon_k^3+\mathcal{Y}_p\times \sum_k\frac{a_k\tilde\epsilon_k^2}{1-\mathcal{Y}_pa_k\tilde\epsilon_k}\bigg)=o(1).
	\]
	We then find the limit of the non-random part and complete the proof of Theorem \ref{thm:m}. It remains to prove the bounds in Lemma \ref{lem: bounds}.
	
	\subsection{Proof of Lemma \ref{lem: bounds}}
	\begin{proof}
		The result for $(np)^{-1}|\text{tr}\Mb\Bb^h|$ is easy by using the similar technique  in deriving (\ref{eq:trace mk}). For the second result, it follows directly  from the expansion in (\ref{eq:mk1}), and the bounds for $\zeta_k$ and $(np)^{-1}\text{tr}\Mb\Bb^h$. Hence, we only prove the bounds for the moments of $\bar\mu_k$.
		
		By definition,
		\[
		\begin{split}
		\bar\mu_k=&\frac{1}{\sqrt{np}}(\by_k^\prime\by_k-n\bar\lambda_{\Bb}a_k)-\bigg(\frac{1}{np}\by_k^\prime\Yb_k^\prime\Db_k^{-1}\Yb_k\by_k-\frac{a_k}{np}\text{tr}\Mb_k\bigg)\\
		&+\bigg(\frac{a_k}{np}\mathbb{E}\text{tr}\Mb_k-\frac{a_k}{np}\text{tr}\Mb_k\bigg)+\bigg(\frac{a_k}{np}\mathbb{E}\text{tr}\Mb-\frac{a_k}{np}\mathbb{E}\text{tr}\Mb_k\bigg)\\
		:=&L_1+L_2+L_3+L_4.
		\end{split}
		\]
		For $L_1$ and $L_2$, by the proof of Lemma \ref{lemma:gamma and eta},
		\[
		\mathbb{E}|L_1|^4=o(p^{-1}),\quad \mathbb{E}|L_2|^4=o(p^{-1}).
		\]
		For $L_3$, we write 
		\[
		\frac{1}{np}\mathbb{E}\text{tr}\Mb_k-\frac{1}{np}\text{tr}\Mb_k=-\frac{1}{np}\sum_{i=1}^p(\mathbb{E}_i-\mathbb{E}_{i-1})\text{tr}(\Mb_k-\Mb_{ki}).
		\]
		Therefore, by Burkholder's inequality and the expansion in (\ref{eq:mk1}), together with the bounds in the proof of Lemma \ref{lemma:gamma and eta}, we conclude
		\[
		\mathbb{E}|L_3|^4=o(p^{-1}).
		\]
		The last term $L_4$ can be handled similarly by Burkholder's inequality, the expansion in (\ref{eq:mk1}) and the bounds in the proof of Lemma \ref{lemma:gamma and eta}. Consequently, we have
		\[
		\mathbb{E}|\bar\mu_k|^4=o(p^{-1}).
		\]
		The result for $\mathbb{E}|\bar\mu_k|^3$ then follows directly from the Cauchy-Schwartz inequality and the fact that $\mathbb{E}|\bar\mu_k|^2=O(p^{-1})$.
	\end{proof}

	\section{Proof of Theorem \ref{thm2}}\label{sece}
	\subsection{Proof of Theorem \ref{thm2} for diagonal $\Ab_p$}
	In the last two sections, we have proved the weak convergence of the process $M_p(z)$ on $\mathcal{C}_u$. By the Cauchy integral formula, to  complete the proof for diagonal $\Ab_p$,  we still need  to show that the integral on $\mathcal{C}_{0,l,r}$ is negligible. Let $u_0\ge 2\lim\sup a_1b_1+\epsilon_0$ for some small positive constant $\epsilon_0$. Since the event $U_p(\epsilon_0):=\{\|\bar\Sbb_p\|\le 2\lim\sup a_1b_1+\epsilon_0/2\}$ happens with probability 1 for large $p$, below we aim  to prove that under this event
	\[
	\lim_{v_0\downarrow 0}\lim_{p\rightarrow\infty}\mathbb{E}\bigg|\int_{\mathcal{C}_{0,l,r}}f(z)M_p(z)dz\bigg|\rightarrow 0.
	\]

	We start with $\mathcal{C}_0$. Since $\|\bar\Sbb_p\|\le 2\lim\sup a_1b_1+\epsilon_0/2$, we have $|m_{F^{\bar\Sbb_p}}(z)|\le 2/\epsilon_0$. Therefore,
	\[
	\mathbb{E}|M_{p1}(z)|\le p\times \sqrt{\mathbb{E}|m_{F^{\bar\Sbb_p}}(z)|^2}\le pC.
	\]
	On the other hand,  $|M_{p2}(z)|\le pC$. Hence,
	\[
	\lim_{v_0\downarrow 0}\lim_{p\rightarrow\infty}\mathbb{E}\bigg|\int_{\mathcal{C}_{0}}f(z)M_{p}(z)dz\bigg|\le \lim_{v_0\downarrow 0}\lim_{p\rightarrow\infty}\int_{\mathcal{C}_{0,}}f(z)\mathbb{E}|M_{p_1}(z)|^2+|M_{p2}(z)|^2dz\le C\times p\epsilon_p\rightarrow 0.
	\]

	Next, for $z\in \mathcal{C}_{l,r}$, it's sufficient to prove that under $U_p(\epsilon_0)$, for sufficiently large $p$,
	\[
	\mathbb{E}|M_{p1}(z)|\le C, \quad |M_{p2}(z)|\le C,\text{ for } z\in \mathcal{C}_{l,r}.
	\]
	We start with $M_{p2}(z)$. Actually, by careful investigation of the previous proof in Section \ref{non-random}, we observe that the condition $\Im z\ge v_0$ ($v_0$ is some constant) only contributes to the bounds  like
	\begin{equation}\label{eq: bound}
	|\beta_k|\le C,\quad \bigg|\frac{1}{np}\text{tr}\Mb\bigg| \le C,\quad |\bar\epsilon_k(z)|\le C,\quad |\tilde\epsilon_k|\le C,\quad \mathcal{W}_2\rightarrow 0.
	\end{equation}
	Hence, we only need to check these bounds for $z\in \mathcal{C}_{l,r}$ and verify the results in Lemma \ref{lem: bounds}. For $\beta_k$, note that $-\beta_k$ is still the $k$-th diagonal element of $\Db^{-1}$. Under the event $U_p(\epsilon_0)$, we have $\|\bar\Sbb_p\|\le 2\lim\sup a_1b_1+\epsilon_0/2$, which implies $|\lambda_j^{\bar\Sbb_p}-z|\ge \epsilon_0/2$ and $\|\Db\|^{-1}\le 2/\epsilon_0$. Therefore, $|\beta_k|$ is still bounded. Moreover, it still holds that $(np)^{-1}|\text{tr}\Mb|\le C$ by similar technique  in (\ref{eq:trace mk}).  
	
	Next, we focus on $\tilde\epsilon_k$. Since $\Im s(z)>0$, the support of $F^{\Ab}$ is bounded and the function $\frac{x}{z+x\bar\lambda_{\Bb^2}s(z)}$ is continuous, from the equations satisfied by $s(z)$, we conclude that there exists a positive constant $\delta_1$ such that for any $t$ in the support of $F^{\Ab}$, 
	\[
	\inf_{z\in\mathcal{C}_{l,r}}|z+t\bar\lambda_{\Bb^2}s(z)|\ge \delta_1.
	\]
	{On the other hand, $s_p(z)\rightarrow s(z)$, then $|\tilde\epsilon_k|\le C$ as long as $a_k$ is in the support of $F^{\Ab}$. For those outside of the support, note that $\Re z=\pm u_0$ for $z\in \mathcal{C}_{l,r}$, while $|u_0|>\|\bar S_p\|$. Therefore, as $\Im z\rightarrow 0$, $\Im s(z)\rightarrow 0$, and it suffices to consider the real part. However, $\Re(z)$ and $\Re s(z)$ are totally determined by $u_0$ and $F^{\Ab}$. Therefore, we can always choose some $u_0$ such that   $\inf_{z\in\mathcal{C}_{l,r}}|z+t\bar\lambda_{\Bb^2}s(z)|\ge \delta_1$ for $t$ equal to the spikes. }
	
	Lastly, we check $\mathcal{W}_2=o(1)$. By the expansion in (\ref{eq:z Mk}),
	\[
	\mathcal{W}_2=\bar\lambda_{\Bb^2}\frac{1}{p}\sum_ia_i(\mathbb{E}\beta_i-\tilde\epsilon_i)+o(1)=-\bar\lambda_{\Bb^2}\big[\frac{1}{p}\mathbb{E}\text{tr}\Db^{-1}\Ab-s_p(z)\big]+o(1).
	\]
	Actually, by the definition of $\beta_i$, (\ref{eq:z Mk}) and the continuous mapping theorem, we know that
	\[
	\beta_i\overset{i.p.}{\longrightarrow} \frac{1}{z+\frac{1}{np}\text{tr}\Mb_i^{(1)}}\overset{i.p.}{\longrightarrow} \frac{1}{z+\frac{a_i}{np}\mathbb{E}\text{tr}\Mb}\overset{i.p.}{\longrightarrow} \frac{1}{z+a_i\bar\lambda_{\Bb^2}\frac{1}{p}\mathbb{E}\text{tr}\Db^{-1}\Ab}.
	\]
	Then, $|\bar\epsilon_k(z)|\le C$ for large $p$. Note that $\beta_i$ and $p^{-1}\text{tr}\Db^{-1}\Ab$ are bounded and $\text{tr}\Db^{-1}\Ab=-\sum_ia_i\beta_i$. Hence, by the dominated convergence theorem,
	\[
	\frac{1}{p}\mathbb{E}\text{tr}\Db^{-1}\Ab=\frac{1}{p}\sum_i\frac{1}{z+a_i\bar\lambda_{\Bb^2}\frac{1}{p}\mathbb{E}\text{tr}\Db^{-1}\Ab}+o(1),
	\]
	which further implies $p^{-1}\mathbb{E}\text{tr}\Db^{-1}\Ab-s_p(z)=o(1)$ and $\mathcal{W}_2=o(1)$.
	Furthermore, with these preliminary  bounds, we observe that the results in Lemma \ref{lem: bounds} still hold. Then, $|M_{p2}(z)|\le C$ for $z\in \mathcal{C}_{l,r}$. 
	
	Now we move to the calculation of $\mathbb{E}|M_{p1}(z)|$ for $z\in \mathcal{C}_{l,r}$. We can not use the decomposition in (\ref{eq:main}) because the bound for $\tilde\beta_k$ is not guaranteed to hold anymore. The strategy is to replace $\tilde\beta_k$ with $\epsilon_k$.  Recall the definition of $\epsilon_k$ and $\mu_k$ in (\ref{eq:epsilon_k mu_k}). Based on (\ref{eq:mF expansion}) and the relationship $\beta_k=\epsilon_k+\beta_k\epsilon_k\mu_k$, 
	\[
	\begin{split}
	M_{p1}(z)=&-\sum_{k=1}^p(\mathbb{E}_k-\mathbb{E}_{k-1})\beta_k(1+\bq_k^\prime\Db_k^{-2}\bq_k)\\
	=&-\sum_{k=1}^p(\mathbb{E}_k-\mathbb{E}_{k-1})\beta_k(1+\frac{1}{np}\text{tr}\Mb_k^{(2)}+\gamma_{k2})\\
	=&-\sum_{k=1}^p(\mathbb{E}_k-\mathbb{E}_{k-1})(\beta_k\gamma_{k2})-\sum_{k=1}^p(\mathbb{E}_k-\mathbb{E}_{k-1})\beta_k\epsilon_k\mu_k\bigg(1+\frac{1}{np}\text{tr}\Mb_k^{(2)}\bigg).
	\end{split}
	\]
	Under the event $U_p(\epsilon_0)$, with the bounds in (\ref{eq: bound}), it's easy to see the second moments of $\gamma_{k2}$ and $\mu_k$ are still $O(p^{-1})$. Hence, by Burkholder's inequality, it reduces  to proving $|\epsilon_k|\le C$. Using the expansions in (\ref{eq:mk1}) and (\ref{eq:z Mk}), we conclude that
	\[
	\epsilon_k\rightarrow \frac{1}{z+a_k\bar\lambda_{\Bb^2}\frac{1}{p}\mathbb{E}\text{tr}\Db^{-1}\Ab}\rightarrow \mathbb{E}\beta_k\le C.
	\]
	As a result, $|M_{p1}(z)|\le C$ for large $p$ and the proof for diagonal $\Ab_p$ has been completed.

	\subsection{Proof of Theorem \ref{thm2} for general $\Ab_p$}
	Up to now we only proved the results for diagonal $\Ab_p$. In this subsection, we extend the results to general non-negative definite matrix $\Ab_p$. Note that if the entries of $\Xb_p$ are i.i.d. standard Gaussian variables, we can regard $\Ab_p$ as diagonal because standard  Gaussian vectors are orthogonally invariant. Therefore, the results hold for Gaussian case. We then follow the interpolation strategy in \cite{bai2019central} to compare  the characteristic functions of linear spectral statistics under the Gaussian case and general case.
	
	The proof is essentially adapted from \cite{bai2019central} and we use similar notation therein. Let $\Xb_p$ be the $p\times n$ random matrix whose entries are i.i.d. from some general distributions. Let $\Zb_p$ be $p\times n$ random matrix whose entries are i.i.d. standard Gaussian variables. Define
	\[
	\begin{split}
	\Wb_p(\theta)=&(w_{jk})=\Xb_p\sin\theta+\Zb_p\cos\theta,\quad \Yb_p(\theta)=\Ab_p^{1/2}\Wb_p(\theta)\Bb_p^{1/2},\\
	\Gb_p(\theta)=&\frac{1}{\sqrt{np}}\Yb_p(\theta)\Yb_p(\theta)^\prime,\quad \bar\Sbb_p(\theta)=\frac{1}{\sqrt{np}}(\Yb_p(\theta)\Yb_p(\theta)^\prime-n\bar\lambda_{\Bb_p}\Ab_p).
	\end{split}
	\]
	Then, $\bar\Sbb_p(0)$ is the matrix of interest. The proofs are very similar to those in \cite{bai2019central} by the observation that $\partial \Gb/\partial \theta=\partial \bar\Sbb/\partial\theta$. Hence, we only show the necessary steps for self-completeness.
	
	Furthermore, let
	\[
	\begin{split}
	\Hb_{p}(t,\theta)=&e^{it\bar\Sbb_p(\theta)},  \quad S(\theta)=\text{tr}\big(\bar\Sbb_p(\theta)\big),\\
	S^0(\theta)=&S(\theta)-p\int f(x)dF_p(x),\quad Z_p(x,\theta)=\mathbb{E}^{ixS^0(\theta)}.
	\end{split}
	\]
	We may suppress the dependence on $p$ and $\theta$ for simplicity.  Therefore, it's sufficient to prove that
	\[
	\frac{\partial Z_p(x,\theta)}{\partial \theta}\rightarrow 0
	\]
	uniformly in $\theta$ over the interval $[0,\pi/2]$, because 
	\[
	Z_p(x,\pi/2)-Z_p(x,0)=\int_{0}^{\pi/2}\frac{\partial Z_p(x,\theta)}{\partial \theta}d\theta.
	\]
	
	Let $f(\lambda)$ be a smooth function with the Fourier transform
	\[
	\hat f(t)=\frac{1}{2\pi}\int_{-\infty}^{+\infty}f(\lambda)e^{-it\lambda}d\lambda.
	\]
	Firstly, we calculate the derivative of $S(\theta)$. By the inverse Fourier transform formula, 
	\[
	\frac{\partial S(\theta)}{\partial w_{jk}}=\int_{-\infty}^{+\infty}\hat f(u)\text{tr}\frac{\partial \Hb_p(u)}{\partial w_{jk}} du,
	\]
	while 
	\[
	\frac{\partial [\Hb_p(u)]_{dl}}{\partial w_{jk}} =\sum_{a,b=1}^p\frac{\partial [\Hb_p(u)]_{dl}}{\partial \bar s_{ab}}\times \frac{\partial \bar s_{ab}}{\partial w_{jk}},
	\]
	where $\bar s_{ab}$ is the $(a,b)$-th entry of $\bar\Sbb_p$.  By Lemma 0.16 in \cite{bai2019central},
	\[
	\begin{split}
	\frac{\partial \Hb_p(t)}{\partial \bar s_{ab}}=i\int_0^te^{is\bar\Sbb_p}\be_a\be_b^\prime e^{i(t-s)\bar\Sbb_p}ds=i\int_0^t\Hb_p(s)\be_a\be_b^\prime \Hb_p(t-s)ds.
	\end{split}
	\]
	On the other hand,
	\[
	\frac{\partial \bar s_{ab}}{\partial w_{jk}}= \frac{\partial g_{ab}}{\partial w_{jk}}= \bigg[\frac{\Gb_p}{\partial w_{jk}}\bigg]_{ab}=\frac{1}{\sqrt{np}}[\Ab^{1/2}]_{aj}[\Bb\Wb^\prime\Ab^{1/2}]_{kb}+\frac{1}{\sqrt{np}}[\Ab^{1/2}\Wb\Bb]_{ak}[\Ab^{1/2}]_{jb}.
	\]
	Let $h_{ij}$ be the $(i,j)$-th entry of $\Hb$, and $f*g(t)=\int_0^tf(s)g(t-s)ds$. Then
	\[
	\begin{split}
	\frac{\partial h_{dl}(t)}{\partial w_{jk}} =&\frac{i}{\sqrt{np}}\sum_{a,b}h_{da}*h_{bl}(t)\bigg\{[\Ab^{1/2}]_{aj}[\Bb\Wb^\prime\Ab^{1/2}]_{kb}+\frac{1}{\sqrt{np}}[\Ab^{1/2}\Wb\Bb]_{ak}[\Ab^{1/2}]_{jb}\bigg\}\\
	=&\frac{i}{\sqrt{np}}[\Hb\Ab^{1/2}]_{dj}*[\Bb\Wb^\prime\Ab^{1/2}\Hb]_{kl}(t)+\frac{i}{\sqrt{np}}[\Hb\Ab^{1/2}\Wb\Bb]_{dk}*[\Ab^{1/2}\Hb]_{jl}(t).
	\end{split}
	\]
	Therefore,
	\[
	\text{tr}\frac{\partial \Hb_p(u)}{\partial w_{jk}}=\frac{2i}{\sqrt{np}}\sum_{d=1}^p[\Hb\Ab^{1/2}]_{dj}*[\Bb\Wb^\prime\Ab^{1/2}\Hb]_{kd}(u)=\frac{2iu}{\sqrt{np}}[\Ab^{1/2}\Hb(u)\Ab^{1/2}\Wb\Bb]_{jk}.
	\]
	We then conclude
	\[
	\frac{\partial S(\theta)}{\partial w_{jk}}=\frac{2i}{\sqrt{np}}\int_{-\infty}^{+\infty}u\hat f(u)[\Ab^{1/2}\Hb(u)\Ab^{1/2}\Wb\Bb]_{jk} du=\frac{2}{\sqrt{np}}\big[\Ab^{1/2}\tilde f(\bar\Sbb_p)\Ab^{1/2}\Wb\Bb\big]_{jk},
	\]
	where
	\[
	\tilde f(\bar\Sbb_p)=i\int_{-\infty}^{+\infty}u\hat f(u)\Hb(u)du.
	\]
	Furthermore,
	\[
	\frac{\partial Z_p(x,\theta)}{\partial \theta}=\frac{2xi}{\sqrt{np}}\sum_{j=1}^p\sum_{k=1}^n\mathbb{E}w_{jk}^\prime\big[\Ab^{1/2}\tilde f(\bar\Sbb_p)\Ab^{1/2}\Wb\Bb\big]_{jk}e^{ixS^0(\theta)},
	\]
	where 
	\[
	w_{jk}^\prime =\frac{d w_{jk}}{d\theta}=x_{jk}\cos\theta-y_{jk}\sin\theta. 
	\]
	Let $\Wb_{p,jk}(w,\theta)$ denote the corresponding matrix $\Wb_p$ by replacing the $(j,k)$-th entry ($w_{jk}$) with $w$. Let
	\[
	\Yb_{p,jk}(w,\theta)=\Ab^{1/2}\Wb_{p,jk}(w,\theta)\Bb_p^{1/2},
	\]
	and define $\Gb_{p,jk}(w,\theta)$, $\bar\Sbb_{p,jk}(w,\theta)$, $\Hb_{p,jk}(w,\theta)$, $S_{jk}(w,\theta)$, $S^0_{jk}(w,\theta)$ similarly. Further let
	\[
	\psi_{jk}(w)=\big[\Ab^{1/2}\tilde f(\bar\Sbb_{p,jk}(w,\theta))\Ab^{1/2}\Wb_{p,jk}(w,\theta)\Bb\big]_{jk}e^{ixS^0_{jk}(w,\theta)}.
	\]
	Then, by Taylor's expansion,
	\[
	\psi_{jk}(w_{jk})=\sum_{l=0}^3\frac{1}{l!}w_{jk}^l\psi_{jk}^{(l)}(0)+\frac{1}{4!}w_{jk}^4\psi_{jk}^{(4)}(\rho w_{jk}),\quad \rho\in(0,1),
	\]
	which implies
	\[
	\frac{\partial Z_p(x,\theta)}{\partial \theta}=\frac{2xi}{\sqrt{np}}\sum_{j=1}^p\sum_{k=1}^n\mathbb{E}w_{jk}^\prime\sum_{l=0}^3\frac{1}{l!}w_{jk}^l\psi_{jk}^{(l)}(0)+\frac{1}{4!}w_{jk}^4\psi_{jk}^{(4)}(\rho w_{jk}).
	\]
	It's easy to see 
	\[
	\begin{split}
	\mathbb{E}w_{jk}^\prime w_{jk}^0=&0,\quad \mathbb{E}w_{jk}^\prime w_{jk}^1=0,\quad\mathbb{E}w_{jk}^\prime w_{jk}^2=\mathbb{E}w_{jk}^3\sin^2\theta\cos\theta, \quad\mathbb{E}w_{jk}^\prime w_{jk}^3=o(1),
	\end{split}
	\]
	under the condition $\nu_4=3$.  Therefore, it's sufficient to prove that
	\[
	\bigg|\frac{1}{\sqrt{np}}\sum_j\sum_k\mathbb{E}w_{jk}^\prime w_{jk}^l\psi_{jk}^l(0)\bigg|\rightarrow 0, l=2,4.
	\]
	
	To this end, we claim that the expansions for $\psi_{jk}^l(0)$ are almost the same as those in \cite{bai2019central}, except that the scaling coefficient in the denominator is $\sqrt{np}$ here rather than $n$. This has no effects on the results because we can slightly modify  Lemma 0.8 in \cite{bai2019central} to a rate of $n^{1/2}p^{3/4}$ when $p/n\rightarrow 0$. Below we write down the modified lemma and its proof, to conclude this subsection. For further detailed proofs, one can refer to \cite{bai2019central}.
	
	\begin{lemma}
		Suppose $\Ab$, $\Bb$ and $\Cb$ are $p\times p$, $p\times n$ and $n\times n$ random matrices respectively, $\tilde\Ab$, $\tilde\Bb$ and $\tilde\Cb$ are $p\times n$ random matrices. Moreover, the moments of their spectral norms are bounded and $p/n\rightarrow 0$. Then we get as $p\rightarrow\infty$,
		\[
		\bigg|\sum_{jk}\mathbb{E}\mathrm{A}_{jj}\Bb_{jk}\Cb_{kk}\bigg|\le C n^{1/2}p^{3/4},\quad 	\bigg|\sum_{jk}\mathbb{E}\Ab_{jk}\Bb_{jk}\Cb_{jk}\bigg|\le C n^{1/2}p^{3/4}.
		\]
	\end{lemma}
	\begin{proof}
		The results follow from the Cauchy-Schwartz inequality directly. Specifically,
		\[
		\begin{split}
		\bigg|\sum_{jk}\mathbb{E}\Ab_{jj}\Bb_{jk}\Cb_{kk}\bigg|\le&\bigg( \sum_{k}\mathbb{E}|\Cb_{kk}|^2\bigg)^{1/2}\times \bigg(\sum_k\mathbb{E}\big|\sum_j \Ab_{jj}\Bb_{jk}\big|^2\bigg)^{1/2}\\
		\le &C\sqrt{n}\times \bigg(\mathbb{E}\sum_{j_1,j_2}\Ab_{j_1j_1}\Ab_{j_2,j_2}\sum_k\Bb_{j_1k}\Bb_{j_2k}\bigg)^{1/2}\\
		\le &C\sqrt{n}\times \bigg(\sum_{j_1,j_2}\mathbb{E}\Ab_{j_1j_1}^2\Ab_{j_2,j_2}^2\bigg)^{1/4}\bigg(\sum_{j_1,j_2}\mathbb{E}(\Bb\Bb^\prime)^2_{j_1,j_2}\bigg)^{1/4}\\
		\le &Cn^{1/2}p^{3/4}.
		\end{split}
		\]
		On the other hand,
		\[
		\begin{split}
		\bigg|\sum_{jk}\mathbb{E}\Ab_{jk}\Bb_{jk}\Cb_{jk}\bigg|\le&\bigg( \sum_{j,k}\mathbb{E}|\Ab_{jk}||\Bb_{jk}|\bigg)^{1/2}\times \bigg(\sum_{j,k}\mathbb{E} \Ab_{jk}||\Bb_{jk}||\Cb_{jk}|^2\bigg)^{1/2}\\
		\le &\bigg( \sum_{j,k}\mathbb{E}|\Ab_{jk}||\Bb_{jk}|\bigg)^{3/4}\times \bigg(\sum_{j,k}\mathbb{E} \Ab_{jk}||\Bb_{jk}||\Cb_{jk}|^4\bigg)^{1/4}\\
		\le &C\bigg( \sum_{j,k}\mathbb{E}|\Ab_{jk}|^2\bigg)^{3/8}\times \bigg( \sum_{j,k}\mathbb{E}|\Bb_{jk}|^2\bigg)^{3/8}\times (np)^{1/4}\\
		\le &Cn^{1/2}p^{3/4},
		\end{split}
		\]
		which concludes the lemma. 
	\end{proof}

	\section{Proof of results in Sections 4 and 5}\label{secf}
	\subsection{Proof of Corollaries \ref{cor1} and $\mu$ in (\ref{mu and sigma2})}
	\begin{proof}
		We first prove Corollary \ref{cor1}. The mean correction term is straightforward from the expressions in Lemma \ref{thm:m}. For the covariance function, note that
		\[
		\begin{split}
		\Lambda(z_1,z_2)=&\frac{\partial^2}{\partial z_1\partial z_2}m(z_1)m(z_2)\bigg(\frac{2}{p}\lim_{p\rightarrow\infty}\sum_{k=1}^p\frac{\frac{k}{p}\bar\lambda^2_{\Bb^2}m(z_1)m(z_2)}{1-\frac{k}{p}\bar\lambda_{\Bb^2}m(z_1)m(z_2)}+\frac{\nu_4-3}{n}\sum_{j=1}^n\mathrm{B}_{n,jj}^2+2\bar\lambda_{\Bb^2}\bigg)\\
		=&\frac{\partial^2}{\partial z_1\partial z_2}m(z_1)m(z_2)\bigg(-\frac{2\log (1-\bar\lambda_{\Bb^2}m(z_1)m(z_2))}{m(z_1)m(z_2)}+(\nu_4-3)n^{-1}\sum_{j=1}^n\mathrm{B}_{n,jj}^2\bigg)\\
		=&m^\prime(z_1)m^\prime(z_2)\bigg((\nu_4-3)\lim_{n\rightarrow}n^{-1}\sum_{j=1}^n\mathrm{B}_{n,jj}^2+\frac{2\bar\lambda_{\Bb^2}}{(1-\bar\lambda_{\Bb^2}m(z_1)m(z_2))^2}\bigg),
		\end{split}
		\]
		where in the second step we use the fact
		\[
		\begin{split}
		\lim_{p\rightarrow\infty}\frac{1}{p}\sum_{k=1}^p\frac{\frac{k}{p}\bar\lambda^2_{\Bb^2}m(z_1)m(z_2)}{1-\frac{k}{p}\bar\lambda_{\Bb^2}m(z_1)m(z_2)}=&\int_0^1\frac{t\bar\lambda_{\Bb^2}^2m(z_1)m(z_2)}{1-t\bar\lambda_{\Bb^2}m(z_1)m(z_2)}dt\\
		=&\bar\lambda_{\Bb^2}\bigg(-1-\big(\bar\lambda_{\Bb^2}m(z_1)m(z_2)\big)^{-1}\log\big(1-\bar\lambda_{\Bb^2}m(z_1)m(z_2)\big)\bigg).
		\end{split}
		\]
		{The simplification of $\Lambda(z_1,z_2)$ is then verified.  Now we move to the covariance function. }Denote
		\[
		\tilde\Lambda(z_1,z_2):=-2\log (1-\bar\lambda_{\Bb^2}m(z_1)m(z_2))+(\nu_4-3)n^{-1}\sum_{j=1}^n\mathrm{B}_{n,jj}^2m(z_1)m(z_2),
		\]
		and consider the covariance function in Theorem \ref{thm2}  with two functions $f,g$. It can be written as 
		\begin{equation}\label{gamma z1 z2}
		\Gamma(z_1,z_2)=-\frac{1}{4\pi^2}\oint_{\mathcal{C}_1}\oint_{\mathcal{C}_2}f(z_1)g(z_2)\frac{\partial^2}{\partial z_1\partial z_2}\tilde\Lambda(z_1,z_2)dz_1dz_2,
		\end{equation}
		where $\mathcal{C}_i$ are disjoint contours formed by vertex $(\pm (2\lim\sup b_1+\epsilon_i)\pm i v_i)$ with some small $\epsilon_i$ and $v_i$. 
		Integrating by parts,  we have
		\[
		\Gamma(z_1,z_2)=-\frac{1}{4\pi^2}\oint_{\mathcal{C}_1}\oint_{\mathcal{C}_2}f^\prime(z_1)g^\prime(z_2)\tilde\Lambda(z_1,z_2)dz_1dz_2.
		\]
		Denote $A(z_1,z_2):=f^\prime(z_1)g^\prime(z_2)\tilde\Lambda(z_1,z_2)$ and $c=2\lim\sup b_1$. Let $v_j\rightarrow 0$ and $\epsilon_j\rightarrow 0$. We have
		\begin{equation}\label{gamma z1 z2 expand}
		\Gamma(z_1,z_2)=-\frac{1}{4\pi^2}\int_{-c}^c\int_{-c}^c\big[A(t_1^-,t_2^-)-A(t_1^-,t_2^+)-A(t_1^+,t_2^-)+A(t_1^+,t_2^+)\big]d{t_1}d{t_2},
		\end{equation}
		where $t_j^{\pm}:=t_j\pm i0$. We first consider $f$ and $g$ to be real-valued functions. Note that
		\[
		m(t\pm i0)=\left\{
		\begin{aligned}
		&\frac{-t+\text{sign}(t)\sqrt{t^2-4\bar\lambda_{\Bb^2}}}{2\bar\lambda_{\Bb^2}},&t^2\ge 4\bar\lambda_{\Bb^2},\\
		&\frac{-t\pm i\sqrt{4\bar\lambda_{\Bb^2}-t^2}}{2\bar\lambda_{\Bb^2}},&t^2<4\bar\lambda_{\Bb^2}.
		\end{aligned}
		\right.
		\]
		Let $\tilde c=\sqrt{4\bar\lambda_{\Bb^2}}$. Then we have
		\[
		\begin{split}
		&\int_{-c}^c\int_{-c}^cf^\prime(t_1)g^\prime(t_2)\big[m(t_1^-)m(t_2^-)-m(t_1^-)m(t_2^+)-m(t_1^+)m(t_2^-)+m(t_1^+)m(t_2^+)\big]dt_1dt_2\\
		=&-\int_{-\tilde c}^{\tilde c}\int_{-\tilde c}^{\tilde c}f^\prime(t_1)g^\prime(t_2)\frac{1}{\bar\lambda_{\Bb^2}^2}\sqrt{4\bar\lambda_{\Bb^2}-t_1^2}\sqrt{4\bar\lambda_{\Bb^2}-t_2^2}dt_1dt_2,\\
		&\int_{-c}^c\int_{-c}^cf^\prime(t_1)g^\prime(t_2)\bigg[\log(1-\bar\lambda_{\Bb^2}m(t_1^-)m(t_2^-))-\log(1-\bar\lambda_{\Bb^2}m(t_1^-)m(t_2^+))\\
		&-\log(1-\bar\lambda_{\Bb^2}m(t_1^+)m(t_2^-))+\log(1-\bar\lambda_{\Bb^2}m(t_1^+)m(t_2^+))\bigg]dt_1dt_2\\
		=&\int_{-\tilde c}^{\tilde c}\int_{-\tilde c}^{\tilde c}f^\prime(t_1)g^\prime(t_2)\log\bigg|\frac{1-\bar\lambda_{\Bb^2}m(t_1^-)m(t_2^-)}{1-\bar\lambda_{\Bb^2}m(t_1^-)m(t_2^+)}\bigg|^2dt_1dt_2\\
		=&\int_{-\tilde c}^{\tilde c}\int_{-\tilde c}^{\tilde c}f^\prime(t_1)g^\prime(t_2)\log\bigg(\frac{4\bar\lambda_{\Bb^2}-t_1t_2+\sqrt{(4\bar\lambda_{\Bb^2}-t_1^2)(4\bar\lambda_{\Bb^2}-t_2^2)}}{4\bar\lambda_{\Bb^2}-t_1t_2-\sqrt{(4\bar\lambda_{\Bb^2}-t_1^2)(4\bar\lambda_{\Bb^2}-t_2^2)}}\bigg)dt_1dt_2.
		\end{split}
		\]
		Hence, the expression in Corollary \ref{cor2} follows. Moreover, the above argument still holds for complex-valued functions $f$ and $g$.
		
		When $f(x)=g(x)=x^2$, it's easy to see that the first term of $H(t_1,t_2)$ is an odd function of $t_1$ or $t_2$, then the corresponding integral is zero. For the second term, take $t_1=\sqrt{4\bar\lambda_{\Bb^2}}\cos\theta$ and $t_2=\sqrt{4\bar\lambda_{\Bb^2}}\sin\theta$. Then, the result follows from  Remark 1.3 in \cite{chen2015clt}. { Meanwhile, by directly calculating expectation, we have
			\[
			\begin{split}
			\mu=&\mathbb{E}\text{tr}\tilde{\Sbb}^2=\frac{1}{Tpq}\sum_{t=1}^T\bigg(\sum_{i=1}^p\mathbb{E}(\|\bx_{t,i\cdot}^\prime\bSigma_{\Vb}\bx_{t,i\cdot}-\text{tr}\bSigma_{\Vb})^2+\sum_{j\ne k}\mathbb{E}(\bx_{t,j\cdot}^\prime\bSigma_{\Vb}\bx_{t,k\cdot})^2\bigg)\\
			=&\frac{\nu_4-3}{q}\sum_{j=1}^q\bSigma_{\Vb,jj}^2+\frac{2}{q}\|\bSigma_{\Vb}\|_F^2+\frac{p-1}{q}\|\bSigma_{\Vb}\|_F^2=\frac{\nu_4-3}{q}\sum_{j=1}^q\sigma_j^4+(p+1)\bar\lambda_{\bSigma_{\Vb}^2},
			\end{split}
			\]
			where $\bx_{t,i\cdot}$ is the $i$-th row vector of $\Xb_t$. This concludes (\ref{mu and sigma2}). }
	\end{proof}
	
	\subsection{Equivalence to \cite{chen2015clt} when $\Ab=\Ib$ and $\Bb=\Ib$}
		As claimed in the main paper, when both $\Ab_p$ and $\Bb_n$ are identity matrices, our results will be consistent with those in \cite{chen2015clt}. Formally, we claim the next lemma. 
		
		\begin{lemma}[Comparison with \cite{chen2015clt}]\label{cor2}
			Under the conditions in Corollary \ref{cor1}, if further $\Bb_n=\Ib_n$, the mean correction term will satisfy
			\begin{equation}\label{check}
			\tilde{\mathcal{A}}\mathcal{X}_p^2+\tilde{\mathcal{B}}\mathcal{X}_p+\tilde{\mathcal{C}}=o(p^{-1}),
			\end{equation}
			where $\tilde{\mathcal{A}}$, $\tilde{\mathcal{B}}$ and $\tilde{\mathcal{C}}$ are from \cite{chen2015clt} and defined by
			\[
			\tilde{\mathcal{A}}=m-\sqrt{\frac{p}{n}}(1+m^2),\quad \tilde{\mathcal{B}}=m^2-1-\sqrt{\frac{p}{n}}m(1+2m^2),\quad\tilde{\mathcal{C}}=\frac{m^3}{p}\mathcal{A}-\sqrt{\frac{p}{n}}m^4,
			\]
			and $m, \mathcal{A}$  are the abbreviations for $m(z)$, $\mathcal{A}_p(z)$, respectively.
		\end{lemma}
		We remark that in \cite{chen2015clt}, they assume $p/n\rightarrow \infty$ rather than $p/n\rightarrow 0$, and consider the normalized sample covariance matrix 
		\[
		\bar S_n= \sqrt{\frac{p}{n}}\bigg(\frac{1}{p}\Xb^\prime\Xb-\Ib_n\bigg),
		\]
		where $\Xb$ is the ${p\times n}$ observed data matrix. This is equivalent to the settings considered in the current paper by simply transposing $\Xb$. Using Lemma \ref{cor2}, one can construct a joint test of $\bSigma_{\Ub}=\bSigma_{\Ub_0}$ and $\bSigma_{\Vb}=\bSigma_{\Vb_0}$ by first defining 
		\[
		\tilde\Yb_t=\bSigma_{\Ub_0}^{-1/2}\Yb_t\bSigma_{\Vb_0}^{-1/2}.
		\]
		However, it's worth mentioning that the joint test is less informative than our individual tests of $\bSigma_{\Ub}$ and $\bSigma_{\Vb}$. The latter provides more details when the joint test rejects the null hypothesis, and also involves more technical innovations. In the following, we prove Lemma \ref{cor2}.
		\begin{proof}
			When $\Bb_n=\Ib_n$, 
			\[
			m_p(z)=s_p(z)=s(z)=m(z)=-\frac{1}{z+m(z)},
			\]
			then the mean correction term can be further simplified as 
			\[
			\mathcal{X}_p=\frac{1}{p}\mathcal{A}\times m^3-\frac{\mathcal{Y}\times m^2}{1+m\mathcal{Y}},
			\]
			where we suppress $z$ for simplicity and 
			\[
			\mathcal{A}=\nu_4-2+m^\prime.
			\]
			Moreover, $\mathcal{D}_2(x)$ is simplified to 
			\[
			\mathcal{D}_2(x)=-\frac{m}{1+mx}.
			\]
			Therefore, $\mathcal{Y}$ is the solution to
			\[
			\begin{split}
			x=&-\frac{1}{p}\mathcal{A}m^3+m-\frac{m}{1+mx}-\bigg(\frac{1}{1+\sqrt{\frac{p}{n}}\frac{m}{1+mx}}-1\bigg)\frac{m}{1+mx}\\
			=&-\frac{1}{p}\mathcal{A}m^3+m-\frac{\frac{m}{1+mx}}{1+\sqrt{\frac{p}{n}}\frac{m}{1+mx}},
			\end{split}
			\]
			and  $\mathcal{Y}=o(1)$.  The above equation is equivalent to
			\[
			\begin{split}
			(x-m+\frac{1}{p}\mathcal{A}m^3)\bigg(1+\sqrt{\frac{p}{n}}\frac{m}{1+mx}\bigg)+\frac{m}{1+mx}=&0.\\
			\Leftrightarrow (x-m+\frac{1}{p}\mathcal{A}m^3)\bigg(1+mx+\sqrt{\frac{p}{n}}m\bigg)+m=&0.
			\end{split}
			\]
			We remove the $o(p^{-1})$ terms and let $\tilde{\mathcal{Y}}$ be the solution to 
			\begin{equation}\label{y}
			mx^2+\bigg(1+\sqrt{\frac{p}{n}}m-m^2\bigg)x+\frac{1}{p}\mathcal{A}m^3-\sqrt{\frac{p}{n}}m^2=0,
			\end{equation}
			and define $\tilde{\mathcal{X}}_p$ by replacing $\mathcal{Y}$ with $\tilde{\mathcal{Y}}$. Then, it's sufficient to verify
			\[
			\mathcal{H}:=\tilde{\mathcal{A}}\tilde{\mathcal{X}}_p^2+\tilde{\mathcal{B}}\tilde{\mathcal{X}}_p+\tilde{\mathcal{C}}
			=o(p^{-1}).
			\]
			
			Note that
			\[
			\begin{split}
			\mathcal{H}=&\bigg(m-\sqrt{\frac{p}{n}}(1+m^2)\bigg)\bigg(\frac{1}{p}\mathcal{A}m^3
			-\frac{\tilde{\mathcal{Y}}m^2}{1+m\tilde{\mathcal{Y}}}\bigg)^2\\
			&+\bigg(m^2-1-\sqrt{\frac{p}{n}}m(1+2m^2)\bigg)\bigg(\frac{1}{p}\mathcal{A}m^3-\frac{\tilde{\mathcal{Y}}m^2}{1+m\tilde{\mathcal{Y}}}\bigg)+\frac{m^3}{p}\mathcal{A}-\sqrt{\frac{p}{n}}m^4.
			\end{split}
			\]
			We first remove all the $o(p^{-1})$ terms and write
			\[
			\begin{split}
			\mathcal{H}=&\bigg(m-\sqrt{\frac{p}{n}}(1+m^2)\bigg)\bigg(\frac{\tilde{\mathcal{Y}}m^2}{1+m\tilde{\mathcal{Y}}}\bigg)^2\\
			&-\bigg(m^2-1-\sqrt{\frac{p}{n}}m(1+2m^2)\bigg)\frac{\tilde{\mathcal{Y}}m^2}{1+m\tilde{\mathcal{Y}}}+\frac{1}{p}\mathcal{A}m^5-\sqrt{\frac{p}{n}}m^4+o(p^{-1}).
			\end{split}
			\]
			Therefore, we only need to prove
			\[
			\begin{split}
			\mathcal{H}_1:=&m^2\bigg(m-\sqrt{\frac{p}{n}}(1+m^2)\bigg)\tilde{\mathcal{Y}}^2-\bigg(m^2-1-\sqrt{\frac{p}{n}}m(1+2m^2)\bigg)\tilde{\mathcal{Y}}(1+m\tilde{\mathcal{Y}})\\
			&+\bigg(\frac{1}{p}\mathcal{A}m^3-\sqrt{\frac{p}{n}}m^2\bigg)(1+m\tilde{\mathcal{Y}})^2=o(p^{-1}).
			\end{split}
			\]
			Once again, we remove all the $o(p^{-1})$ terms and after some calculations,
			\begin{equation}\label{y1}
			\begin{split}
			\mathcal{H}_1:=&m\tilde{\mathcal{Y}}^2-\bigg(m^2-1-\sqrt{\frac{p}{n}}m\bigg)\tilde{\mathcal{Y}}+\frac{1}{p}\mathcal{A}m^3-\sqrt{\frac{p}{n}}m^2+o(p^{-1}),
			\end{split}
			\end{equation}
			which is exactly consistent with (\ref{y}). The corollary is then verified. 
		\end{proof}

	\subsection{Proof of Theorems \ref{null} and \ref{alternative}}
	\begin{proof}
		Theorem \ref{null} holds directly by Corollary \ref{cor2}. Hence, we only consider the alternative hypothesis and Theorem \ref{alternative}.
		Let $\tilde\Sbb^0$ be the renormalized separable sample covariance matrix corresponding to $\Ub=\tilde\bSigma_{\Ub}^{-1/2}$, i.e.,
		\[
		\tilde\Sbb^0=\sqrt{\frac{Tq}{p}}\bigg(\frac{1}{Tq}\sum_{t=1}^T\Xb_t\bSigma_{\Vb}\Xb_t^\prime-\Ib\bigg),
		\]
		and define $G_p^0(f)$  by (\ref{gpf}) using $\tilde\Sbb^0$. Therefore, under the alternative, with $f(x)=x^2$,
		\[
		\mathcal{T}=\frac{G_p(f)}{2\bar\lambda_{\bSigma_{\Vb}^2}}=\frac{1}{2\bar\lambda_{\bSigma_{\Vb}^2}}\bigg(\text{tr}\tilde\Sbb^2-\text{tr}(\tilde\Sbb^0)^2\bigg)+\frac{G_p^0(f)}{2\bar\lambda_{\bSigma_{\Vb}^2}},
		\]
		where
		\[
		\tilde\Sbb=\sqrt{\frac{Tq}{p}}\bigg(\frac{1}{Tq}\sum_{t=1}^T\tilde\bSigma_{\Ub}^{-1/2}\Ub\Xb_t\bSigma_{\Vb}\Xb_t^\prime\Ub^\prime\tilde\bSigma_{\Ub}^{-1/2}-\Ib\bigg).
		\]
		By Theorem \ref{null}, $G_p^0(f)=O_p(1)$. Therefore, it's sufficient to prove that
		\[
		\text{tr}\tilde\Sbb^2-\text{tr}(\tilde\Sbb^0)^2\overset{i.p.}{\rightarrow} \infty, \text{  as  } p\rightarrow \infty.
		\]
		Denote $\tilde\Ub=\tilde\bSigma_{\Ub}^{-1/2}\Ub$, $\Xb=(\Xb_1,\ldots,\Xb_T)$, $n=Tq$ and $\Bb_n=\Ib_T\otimes \bSigma_{\Vb}$, then by definition
		\[
		\begin{split}
		&\frac{1}{n}\bigg(\text{tr}\tilde\Sbb^2-\text{tr}(\tilde\Sbb^0)^2\bigg)\\
		=&\frac{1}{p}\text{tr}\bigg(\frac{1}{n^2}\big[\tilde\Ub\Xb\Bb_n\Xb^\prime\tilde\Ub^\prime-n\Ib_p\big]^2-\frac{1}{n^2}\big[\Xb\Bb_n\Xb^\prime-n\Ib_p\big]^2\bigg)\\
		=&\frac{1}{p}\text{tr}\bigg(\frac{1}{n^2}\bigg[(\tilde\Ub-\Ib)\Xb\Bb_n\Xb^\prime(\tilde\Ub^\prime-\Ib)+(\tilde\Ub-\Ib)\Xb\Bb_n\Xb^\prime+\Xb\Bb_n\Xb^\prime(\tilde\Ub^\prime-\Ib)\bigg]^2\bigg)\\
		&+\frac{2}{p}\text{tr}\bigg(\frac{1}{n^2}\bigg[(\tilde\Ub-\Ib)\Xb\Bb_n\Xb^\prime(\tilde\Ub^\prime-\Ib)+(\tilde\Ub-\Ib)\Xb\Bb_n\Xb^\prime+\Xb\Bb_n\Xb^\prime(\tilde\Ub^\prime-\Ib)\bigg]\\
		&\times\big[\Xb\Bb_n\Xb^\prime-n\Ib_p\big] \bigg)\\
		:=&L_1+L_2.
		\end{split}
		\]
		By the inequality $|p^{-1}\text{tr}\Ab_1\Ab_2|\le \|\Ab_1\|\|\Ab_2\|$ for any $p\times p$ matrices $\Ab_1$ and $\Ab_2$, we have
		\[
		\begin{split}
		|L_2|\le &C_0\bigg\|\frac{1}{n}\bigg[(\tilde\Ub-\Ib)\Xb\Bb_n\Xb^\prime(\tilde\Ub^\prime-\Ib)+(\tilde\Ub-\Ib)\Xb\Bb_n\Xb^\prime+\Xb\Bb_n\Xb^\prime(\tilde\Ub^\prime-\Ib)\bigg]\bigg\|\\
		&\times \bigg\|\frac{1}{n}\big[\Xb\Bb_n\Xb^\prime-n\Ib_p\big]\bigg\|\\
		\le &C_0\times \bigg\|\frac{1}{n}\Xb\Bb_n\Xb^\prime\bigg\|\times \bigg\|\frac{1}{n}\big[\Xb\Bb_n\Xb^\prime-n\Ib_p\big]\bigg\|\\=o_p(1),
		\end{split}
		\]
		where we use the facts $\|n^{-1}\Xb\Bb_n\Xb^\prime\|\le C_0$ for some constant $C_0$ and $\|n^{-1}\big[\Xb\Bb_n\Xb^\prime-n\Ib_p\big]\|=o_p(1)$ by Lemma \ref{lem1}. On the other hand, for $L_1$, denote
		\[
		\begin{split}
		J_1:=\frac{1}{n}\bigg[(\tilde\Ub-\Ib)\Xb\Bb_n\Xb^\prime(\tilde\Ub^\prime-\Ib)+(\tilde\Ub-\Ib)\Xb\Bb_n\Xb^\prime+\Xb\Bb_n\Xb^\prime(\tilde\Ub^\prime-\Ib)\bigg].
		\end{split}
		\]
		Then, by $\|n^{-1}\big[\Xb\Bb_n\Xb^\prime-n\Ib_p\big]\|=o_p(1)$, it's easy to see
		\[
		\begin{split}
		&\bigg\|J_1-\frac{1}{n}\bigg[n(\tilde\Ub-\Ib)^2+2n(\tilde\Ub-\Ib)\bigg]\bigg\|=o_p(1),
		\end{split}
		\]
		which further implies
		\[
		\bigg|L_1-\frac{1}{p}\text{tr}\bigg[(\tilde\Ub-\Ib)^2+2(\tilde\Ub-\Ib)\bigg]^2\bigg|=o_p(1).
		\]
		Note that
		\[
		\frac{1}{p}\text{tr}\bigg[(\tilde\Ub-\Ib)^2+2(\tilde\Ub-\Ib)\bigg]^2\bigg)=\frac{1}{p}\text{tr}[\tilde\bSigma_{\Ub}^{-1/2}\bSigma_{\Ub}\tilde\bSigma_{\Ub}^{-1/2}-\Ib]^2\ge c_0,
		\]
		for sufficiently large $p$. Therefore, with probability approaching  1,
		\[
		\frac{1}{n}\bigg|\text{tr}\tilde\Sbb^2-\text{tr}(\tilde\Sbb^0)^2\bigg|\ge c_0 \Rightarrow \bigg|\text{tr}\tilde\Sbb^2-\text{tr}(\tilde\Sbb^0)^2\bigg|\rightarrow \infty.
		\]
		That is, $\mathcal{T}\rightarrow \infty$, which concludes the theorem. 
	\end{proof}
	
	\subsection{Proof of Lemma \ref{consistcy: mu 1}}
	{	\begin{proof}
			Note that for any symmetric matrices $\Ab$ and $\Bb$,
			\[
			\|\Ab\|_F^2-\|\Bb\|_F^2=\text{tr}(\Ab-\Bb)^2+\text{tr}\Bb(\Ab-\Bb)+\text{tr}(\Ab-\Bb)\Bb.
			\]
			Then, it's sufficient to consider 
			\[
			\mathcal{J}_1:=\text{tr}\bigg(\frac{1}{Tp}\sum_{t=1}^T\tilde\Yb_t^\prime\tilde\Yb_t-\bSigma_{\Vb}\bigg)^2,\text{ and } \mathcal{J}_2:=\text{tr}\bigg(\frac{1}{Tp}\sum_{t=1}^T\tilde\Yb_t^\prime\tilde\Yb_t-\bSigma_{\Vb}\bigg)\bSigma_{\Vb}.
			\]
			For $\mathcal{J}_2$, under the null hypothesis, 
			\[
			\begin{split}
			\mathcal{J}_2=\text{tr}\bigg(\frac{1}{Tp}\sum_{t=1}^T\Xb_t^\prime\Xb_t-\Ib\bigg)\Vb^\prime\bSigma_{\Vb}\Vb.
			\end{split}
			\]
			Write $\tilde\Xb:=(\tilde\bX_{1\cdot},\ldots,\tilde\bX_{q\cdot})^\prime:=(\Xb_1^\prime,\ldots,\Xb_T^\prime)$, then
			\[
			\begin{split}
			\mathcal{J}_2=&\sum_{i=1}^q\bigg(\frac{1}{Tp}\tilde\bX_{i\cdot}^\prime\tilde\Xb^\prime-\be_i^\prime\bigg)\Vb^\prime\bSigma_{\Vb}\bV_i\\
			=&\sum_{i=1}^q\bigg(\frac{1}{Tp}\tilde\bX_{i\cdot}^\prime\tilde\bX_{i\cdot}-1\bigg)\bV_i^\prime\bSigma_{\Vb}\bV_i+\sum_{i=1}^q\sum_{j\ne i}^q\frac{1}{Tp}\tilde\bX_{i\cdot}^\prime\tilde\bX_{j\cdot}\bV_j^\prime\bSigma_{\Vb}\bV_i\\
			:=&\mathcal{J}_{21}+\mathcal{J}_{22},
			\end{split}
			\]
			where $\bV_i$ is the $i$-th column vector of $\Vb$, $\tilde\Xb_{(i)}$ and $\Vb_{(i)}$ are $\tilde\Xb$ and $\Vb$ after replacing the $i$-th rows with 0, respectively. By independence, 
			\[
			\frac{p^2}{q^2}\mathbb{E}|\mathcal{J}_{21}|^2\le O\bigg(\frac{p^2}{q^2}\times\frac{q}{Tp}\bigg)=O\bigg(\frac{p}{Tq}\bigg)\rightarrow 0.
			\]
			On the other hand,
			\[
			\begin{split}
			\mathbb{E}|\mathcal{J}_{22}|^2\le &C\sum_{i,j\ne i}^q\frac{1}{T^2p^2}\mathbb{E}(\tilde\bX_{i\cdot}^\prime\tilde\bX_{j\cdot})^2(\bV_j^\prime\bSigma_{\Vb}\bV_i)^2\le \frac{C}{Tp}\sum_{i,j}(\bV_j^\prime\bSigma_{\Vb}\bV_i)^2\\
			\le &\frac{C}{Tp}\text{tr}(\bSigma_{\Vb}^4)=O\bigg(\frac{q}{Tp}\bigg).
			\end{split}
			\]
			As a result, 
			\[
			\frac{p}{q}|\mathcal{J}_2|=O_p\bigg(\sqrt{\frac{p}{Tq}}\bigg).
			\]

			For $\mathcal{J}_1$, denote $\tilde\bX_{\cdot j}$ as the $j$-th column of $\tilde\Xb$. Then, under the null hypothesis,
			\[
			\begin{split}
			\mathcal{J}_1=\text{tr}\bigg(\sum_{i=1}^q\bV_i\big[\frac{1}{Tp}\tilde\bX_{i\cdot}^\prime\tilde\bX_{i\cdot}-1\big]\bV_i^\prime+\sum_{j\ne i}\frac{1}{Tp}\bV_i\tilde\bX_{i\cdot}^\prime\tilde\bX_{\cdot j}\bV_j^\prime\bigg)^2.
			\end{split}
			\]
			Let $\mathbb{E}_{k}(\cdot)=\mathbb{E}(\cdot\mid \tilde\bX_{1\cdot},\ldots,\tilde\bX_{k\cdot })$. Then,
			\[
			\begin{split}
			&\mathcal{J}_1-\mathbb{E}\mathcal{J}_1\\
			=&\sum_{k=1}^q(\mathbb{E}_k-\mathbb{E}_{k-1})\text{tr}\bigg(\sum_{i=1}^q\bV_i\big[\frac{1}{Tp}\tilde\bX_{i\cdot}^\prime\tilde\bX_{i\cdot}-1\big]\bV_i^\prime+\sum_{j\ne i}\frac{1}{Tp}\bV_i\tilde\bX_{i\cdot}^\prime\tilde\bX_{\cdot j}\bV_j^\prime\bigg)^2\\
			=&\sum_{k=1}^q(\mathbb{E}_k-\mathbb{E}_{k-1})\text{tr}\bigg(\bV_k\big[\frac{1}{Tp}\tilde\bX_{k\cdot}^\prime\tilde\bX_{k\cdot}-1\big]\bV_k^\prime+\sum_{i\ne k}\bV_i\big[\frac{1}{Tp}\tilde\bX_{i\cdot}^\prime\tilde\bX_{i\cdot}-1\big]\bV_i^\prime\\
			&+\frac{2}{Tp}\sum_{j\ne k}\bV_k\tilde\bX_{k\cdot}^\prime\tilde\bX_{j\cdot}\bV_j^\prime+\sum_{j\ne i,j\ne k,i\ne k}\frac{1}{Tp}\bV_i\tilde\bX_{i\cdot}^\prime\tilde\bX_{\cdot j}\bV_j^\prime\bigg)^2\\
			=&\sum_{k=1}^q(\mathbb{E}_k-\mathbb{E}_{k-1})\bigg(\big[\frac{1}{Tp}\tilde\bX_{k\cdot}^\prime\tilde\bX_{k\cdot}-1\big]\bV_k^\prime\bV_k\bigg)^2+\sum_{k=1}^q(\mathbb{E}_k-\mathbb{E}_{k-1})\text{tr}\bigg(\frac{2}{Tp}\sum_{j\ne k}\bV_k\tilde\bX_{k\cdot}^\prime\tilde\bX_{j\cdot}\bV_j^\prime\bigg)^2\\
			&+2\sum_{k=1}^q(\mathbb{E}_k-\mathbb{E}_{k-1})\text{tr}\bigg(\bV_k\big[\frac{1}{Tp}\tilde\bX_{k\cdot}^\prime\tilde\bX_{k\cdot}-1\big]\bV_k^\prime\frac{2}{Tp}\sum_{j\ne k}\bV_k\tilde\bX_{k\cdot}^\prime\tilde\bX_{j\cdot}\bV_j^\prime\bigg)\\
			&+2\sum_{k=1}^q(\mathbb{E}_k-\mathbb{E}_{k-1})\text{tr}\bigg(\bV_k\big[\frac{1}{Tp}\tilde\bX_{k\cdot}^\prime\tilde\bX_{k\cdot}-1\big]\bV_k^\prime+\frac{2}{Tp}\sum_{j\ne k}\bV_k\tilde\bX_{k\cdot}^\prime\tilde\bX_{j\cdot}\bV_j^\prime\bigg)\times\\
			&\bigg(\sum_{i\ne k}\bV_i\big[\frac{1}{Tp}\tilde\bX_{i\cdot}^\prime\tilde\bX_{i\cdot}-1\big]\bV_i^\prime+\sum_{j\ne i,j\ne k,i\ne k}\frac{1}{Tp}\bV_i\tilde\bX_{i\cdot}^\prime\tilde\bX_{\cdot j}\bV_j^\prime\bigg)\\
			:=&\mathcal{J}_{11}+\mathcal{J}_{12}+\mathcal{J}_{13}+\mathcal{J}_{14}.
			\end{split}
			\]
			We deal with them separately. Firstly, by Burkholder's inequality,
			\[
			\begin{split}
			\mathbb{E}|\mathcal{J}_{11}|^2\le& Cq\mathbb{E}\bigg|(\mathbb{E}_k-\mathbb{E}_{k-1})\text{tr}\bigg(\bV_k\big[\frac{1}{Tp}\tilde\bX_{k\cdot}^\prime\tilde\bX_{k\cdot}-1\big]\bV_k^\prime\bigg|^2\le \frac{Cq}{(Tp)^3}\mathbb{E}|x_{t,ij}|^8+\frac{Cq}{(Tp)^2}\\
			\le &\frac{Cq^2}{(Tp)^2}.
			\end{split} 
			\]
			For $\mathcal{J}_{12}$, we write
			\[
			\begin{split}
			\mathbb{E}|\mathcal{J}_{12}|^2\le& C\sum_k\mathbb{E}\bigg|\frac{1}{T^2p^2}(\mathbb{E}_k-\mathbb{E}_{k-1})\bigg(\sum_{j<k}\tilde\bX_{k\cdot}^\prime\tilde\bX_{j\cdot}\bV_j^\prime\bV_k+\sum_{j>k}\tilde\bX_{k\cdot}^\prime\tilde\bX_{j\cdot}\bV_j^\prime\bV_k\bigg)^2\bigg|^2\\
			\le &\frac{C}{T^4p^4}\sum_k\mathbb{E}\bigg|\tilde\bX_{k\cdot}^\prime(\sum_{j<k}\tilde\bX_{j\cdot}\bV_j^\prime\bV_k)(\sum_{j<k}\tilde\bX_{j\cdot}\bV_j^\prime\bV_k)^\prime\tilde\bX_{k\cdot}-\text{tr}(\sum_{j<k}\tilde\bX_{j\cdot}\bV_j^\prime\bV_k)(\sum_{j<k}\tilde\bX_{j\cdot}\bV_j^\prime\bV_k)^\prime\bigg|^2\\
			&+\frac{C}{T^4p^4}\sum_k\mathbb{E}\bigg|\sum_{j>k}(\bV_j^\prime\bV_k)^2(\|\tilde\bX_{k\cdot}\|^2-Tp)\bigg|^2\\
			\le &\frac{C}{T^4p^4}\sum_k\mathbb{E}\bigg\|(\sum_{j<k}\tilde\bX_{j\cdot}\bV_j^\prime\bV_k)(\sum_{j<k}\tilde\bX_{j\cdot}\bV_j^\prime\bV_k)^\prime\bigg\|_F^2+\frac{C}{T^3p^3}\sum_k\bigg(\sum_{j>k}\bSigma_{\Vb,kj}^2\bigg)^2\\
			\le &\frac{C}{T^2p^2}\sum_k\bigg(\sum_{j<k}(\bV_j^\prime\bV_k)^2\bigg)^2+\frac{Cq^2}{T^3p^3}\le \frac{Cq^2}{T^2p^2}.
			\end{split}
			\] 
			For $\mathcal{J}_{13}$, we write
			\[
			\begin{split}
			\mathbb{E}|\mathcal{J}_{13}|^2\le&C\sum_k\sum_{j<k}\mathbb{E}\bigg|\big[\frac{1}{Tp}\tilde\bX_{k\cdot}^\prime\tilde\bX_{k\cdot}-1\big]\bV_k^\prime\bV_k\frac{1}{Tp}\tilde\bX_{k\cdot}^\prime\tilde\bX_{j\cdot}\bV_j^\prime\bV_k\bigg|^2\\
			\le &C\sum_k\sum_{j<k}\bSigma_{Vb,kk}^2\bSigma_{\Vb,kj}^2\sqrt{\mathbb{E}\big[\frac{1}{Tp}\tilde\bX_{k\cdot}^\prime\tilde\bX_{k\cdot}-1\big]^4}\sqrt{\mathbb{E}\big[\frac{1}{Tp}\tilde\bX_{k\cdot}^\prime\tilde\bX_{j\cdot}\big]^4}\\
			\le &\frac{C}{T^2p^2}\sum_k\sum_{j<k}\bSigma_{\Vb,kj}^2\le \frac{Cq}{T^2p^2}.
			\end{split}
			\]
			Lastly, for $\mathcal{J}_{14}$, we write
			\[
			\begin{split}
			\mathbb{E}|\mathcal{J}_{141}|^2=&\sum_k\mathbb{E}\bigg|(\mathbb{E}_k-\mathbb{E}_{k-1})\text{tr}\bigg(\bV_k\big[\frac{1}{Tp}\tilde\bX_{k\cdot}^\prime\tilde\bX_{k\cdot}-1\big]\bV_k^\prime\sum_{i\ne k}\bV_i\big[\frac{1}{Tp}\tilde\bX_{i\cdot}^\prime\tilde\bX_{i\cdot}-1\big]\bV_i^\prime\bigg)\bigg|^2\\
			\le &\sum_k\sum_{i<k}\bSigma_{\Vb,ik}^4\sqrt{\mathbb{E}\big[\frac{1}{Tp}\tilde\bX_{k\cdot}^\prime\tilde\bX_{k\cdot}-1\big]^4}\sqrt{\mathbb{E}\big[\frac{1}{Tp}\tilde\bX_{i\cdot}^\prime\tilde\bX_{i\cdot}-1\big]^4}\\
			\le &\frac{C}{T^2p^2}\sum_k\sum_{j<k}\bSigma_{\Vb,kj}^2\le \frac{Cq}{T^2p^2}.
			\end{split}
			\]
			By tedious but elementary calculations, one can verify that
			\[
			\begin{split}
			\mathbb{E}|\mathcal{J}_{142}|^2=&\sum_k\mathbb{E}\bigg|(\mathbb{E}_k-\mathbb{E}_{k-1})\text{tr}\bigg(\bV_k\big[\frac{1}{Tp}\tilde\bX_{k\cdot}^\prime\tilde\bX_{k\cdot}-1\big]\bV_k^\prime\sum_{j\ne i,j\ne k,i\ne k}\frac{1}{Tp}\bV_i\tilde\bX_{i\cdot}^\prime\tilde\bX_{\cdot j}\bV_j^\prime\bigg)\bigg|^2\le \frac{Cq^2}{T^2p^2},\\
			\mathbb{E}|\mathcal{J}_{143}|^2=&\sum_k\mathbb{E}\bigg|(\mathbb{E}_k-\mathbb{E}_{k-1})\text{tr}\bigg(\frac{2}{Tp}\sum_{j\ne k}\bV_k\tilde\bX_{k\cdot}^\prime\tilde\bX_{j\cdot}\bV_j^\prime\sum_{i\ne k}\bV_i\big[\frac{1}{Tp}\tilde\bX_{i\cdot}^\prime\tilde\bX_{i\cdot}-1\big]\bV_i^\prime\bigg)\bigg|^2\le \frac{Cq^2}{T^2p^2},\\
			\mathbb{E}|\mathcal{J}_{144}|^2=&\sum_k\mathbb{E}\bigg|(\mathbb{E}_k-\mathbb{E}_{k-1})\text{tr}\bigg(\frac{2}{Tp}\sum_{j\ne k}\bV_k\tilde\bX_{k\cdot}^\prime\tilde\bX_{j\cdot}\bV_j^\prime\sum_{j\ne i,j\ne k,i\ne k}\frac{1}{Tp}\bV_i\tilde\bX_{i\cdot}^\prime\tilde\bX_{\cdot j}\bV_j^\prime\bigg)\bigg|^2\le \frac{Cq^2}{T^2p^2}.
			\end{split}
			\]
			Consequently, we conclude that
			\[
			\frac{p}{q}(\mathcal{J}_1-\mathbb{E}\mathcal{J}_1)=o_p(1).
			\]
			
			On the other hand,
			\[
			\begin{split}
			\mathbb{E}\mathcal{J}_1=&\mathbb{E}\text{tr}\bigg(\sum_{i=1}^q\bV_i\big[\frac{1}{Tp}\tilde\bX_{i\cdot}^\prime\tilde\bX_{i\cdot}-1\big]\bV_i^\prime\bigg)^2+\mathbb{E}\text{tr}\bigg(\sum_{i,j\ne i}\frac{1}{Tp}\bV_i\tilde\bX_{i\cdot}^\prime\tilde\bX_{\cdot j}\bV_j^\prime\bigg)^2\\
			=&\sum_i\bSigma_{\Vb,ii}^2\mathbb{E}\big[\frac{1}{Tp}\tilde\bX_{i\cdot}^\prime\tilde\bX_{i\cdot}-1\big]^2+\frac{1}{T^2p^2}\sum_{i\ne j}(\bSigma_{\Vb,ij}^2+\bSigma_{\Vb,ii}\bSigma_{\Vb,jj})\mathbb{E}(\tilde\bX_{i\cdot}^\prime\tilde\bX_{j\cdot})^2\\
			=&O\bigg(\frac{q}{Tp}\bigg)+\frac{1}{Tp}\text{tr}^2(\bSigma_{\Vb}).
			\end{split}
			\]
			Therefore,
			\[
			\frac{p}{q}\mathcal{J}_1=\frac{q}{T}+o_p(1),
			\]
			which concludes the lemma.
	\end{proof}}

	\subsection{Proof of Lemma \ref{consistency:nu}}
	{\begin{proof}
			We first prove the relationship (\ref{relation}) in the main paper. Define
			\[
			Z_t=\frac{1}{p}\text{tr}\tilde\Yb_t^\prime\tilde\Yb_t=\frac{1}{p}\sum_{i=1}^p\|\Vb\tilde\bY_{t,i\cdot}\|^2.
			\] 
			Clearly,
			
			\[
			\begin{split}
			\text{var}(Z_t)=\frac{1}{p^2}\sum_{i=1}^p\text{var}(\|\Vb\bX_{t,i\cdot}\|^2).
			\end{split}
			\]
			For each $i$,
			\[
			\begin{split}
			\text{var}(\|\Vb\bX_{t,i\cdot}\|^2)=&\mathbb{E}\bigg(\|\Vb\bX_{t,i\cdot}\|^2-\|\Vb\|_F^2\bigg)=\mathbb{E}\bigg(\sum_{j=1}^q\big[(\sum_{h=1}^qV_{jh}X_{t,ih})^2-\|\bV_j\|^2\big]\bigg)^2\\
			=&\mathbb{E}\bigg(\sum_{j=1}^q\big[\sum_{h=1}^qV_{jh}^2(X_{t,ih}^2-1)+\sum_{h_1\ne h_2}V_{jh_1}V_{jh_2}X_{t,ih_1}X_{t,ih_2}\big]\bigg)^2\\
			=&\mathbb{E}\bigg(\sum_{j=1}^q\sum_{h=1}^qV_{jh}^2(X_{t,ih}^2-1)\bigg)^2+\mathbb{E}\bigg(\sum_{j=1}^q\sum_{h_1\ne h_2}V_{jh_1}V_{jh_2}X_{t,ih_1}X_{t,ih_2}\bigg)^2\\
			=&\sum_{j_1,j_2}\sum_h(\nu_4-1)V_{j_1h}^2V_{j_2h}^2+2\sum_{j_1,j_2}\sum_{h_1\ne h_2}V_{j_1h_1}V_{j_2h_1}V_{j_1h_2}V_{j_2h_2}\\
			=&(\nu_4-3)\sum_h\sigma_h^4+2\|\bSigma_{\Vb}\|_F^2,
			\end{split}
			\]
			which concludes (\ref{relation}).
			
			Now we prove the consistency of $\nu_4$. 
			It's sufficient to prove that $q^{-1}\hat\zeta$, $q^{-1}\hat\omega$ and $q^{-1}\hat\tau$ are consistent estimators of $q^{-1}\zeta$, $q^{-1}\omega$ and $q^{-1}\tau$, respectively. The consistency of $q^{-1}\hat\tau$ has already been demonstrated by the proof of Lemma \ref{consistcy: mu 1}. Hence, we only focus on the other two terms.
			
			For $\hat\omega$, under the null hypothesis, $\tilde\Yb_t^\prime\tilde\Yb_t=\Vb\Xb_t^\prime\Xb_t\Vb$. Then, the model is exactly a vector-valued covariance model by regarding $q$ as dimension and $Tq$ as sample size. Then Theorem 1 in \cite{lopes2019bootstrapping} directly indicates that $q^{-1}\hat\omega$ is  consistent.
			
			For $\hat\zeta$,  we write
			\[
			\begin{split}
			q^{-1}\hat\zeta=&\frac{1}{Tpq}\sum_{t=1}^T\bigg[\text{tr}(\tilde Y_t^\prime\tilde Y_t)-\mathbb{E}\text{tr}(\tilde Y_t^\prime\tilde Y_t)-\frac{1}{T}\sum_{t=1}^T\{\text{tr}(\tilde Y_t^\prime\tilde Y_t)-\mathbb{E}\text{tr}(\tilde Y_t^\prime\tilde Y_t)\}\bigg]^2\\
			=&q^{-1}\zeta-\frac{1}{pq}\bigg[\frac{1}{T}\sum_{t=1}^T\{\text{tr}(\tilde Y_t^\prime\tilde Y_t)-\mathbb{E}\text{tr}(\tilde Y_t^\prime\tilde Y_t)\}\bigg]^2.
			\end{split}
			\]
			We already know that
			\[
			\mathbb{E}\frac{1}{pq}\bigg[\frac{1}{T}\sum_{t=1}^T\{\text{tr}(\tilde Y_t^\prime\tilde Y_t)-\mathbb{E}\text{tr}(\tilde Y_t^\prime\tilde Y_t)\}\bigg]^2=\frac{1}{T^2pq}\sum_t\text{var}(pZ_t)=\frac{p}{Tq}\text{var}(Z_1)\rightarrow 0.
			\]
			Therefore, $q^{-1}\hat\zeta-q^{-1}\zeta=o_p(1)$, which concludes the consistency of $\hat\nu_4$. For $\hat\mu_2$, the consistency holds directly because $q^{-1}\hat\zeta$ and $q^{-1}\hat\tau$ are consistent.
	\end{proof}}
	
	\subsection{Proof of Theorem \ref{bootstrap}}
	
	\begin{proof}
		
		{	Note that with probability tending to one $\hat\bSigma_{\Vb}$ satisfies the assumptions in Condition \ref{c2}. Therefore, by Theorem \ref{null}, 
			\[
			d_{LP}\bigg(\mathcal{N}(0,1),\frac{\mathcal{K}_1^*+\frac{p}{2\pi i}\int f(z)\hat{\mathcal{X}}_p(z)dz}{\hat\sigma}\bigg|\{\Xb_t\}_{t=1}^T\bigg)\overset{i.p.}{\longrightarrow } 0,
			\]
			where $\hat\sigma$ and $\hat{\mathcal{X}}_p(z)$ are calculated by replacing $\bSigma_{\Vb}$ and $\nu_4$ with $\hat\bSigma_{\Vb}$ and $\hat \nu_4$, respectively. 
			Therefore, it's sufficient to prove the respective convergence of $\hat\sigma$ and $p\int f(z)\hat{\mathcal{X}}_{p}(z)dz$ to $\sigma$ and $p\int f(z){\mathcal{X}}_{p}(z)dz$.  By Corollary \ref{cor1}, the consistency of $\hat\sigma$ is guaranteed as long as $\bar\lambda_{\hat\bSigma_{\Vb}^2}$ is consistent. This holds directly by the assumption $\|\hat\bSigma_{\Vb}-\bSigma_{\Vb}\|=o_p(1)$. Therefore, it suffices to consider $\hat{\mathcal{X}}_p(z)$. }
		
		Since the contour $\mathcal{C}$ encloses the interval $[-2\|\bSigma_{\Vb}\|+\epsilon_0,2\|\bSigma_{\Vb}\|+\epsilon_0]$, by (\ref{simplify}) ($\Bb_n=\Ib_T\otimes \bSigma_{\Vb}$), we conclude that with probability tending to 1
		\begin{equation}\label{ff1}
		\max_{z\in\mathcal{C}}| m_p(z)|\le C,\quad \max_{z\in\mathcal{C}}| m_p^\prime(z)|\le C,\quad\max_{z\in\mathcal{C}}|\hat m_p(z)|\le C,\quad \max_{z\in\mathcal{C}}|\hat m_p^\prime(z)|\le C.
		\end{equation}
		Then, with probability tending to 1, 
		\begin{equation}\label{ff2}
		\max_{z\in\mathcal{C}}|\mathcal{A}_p(z)|\le C.
		\end{equation}
		On the other hand, by (\ref{simplify}), we also have
		\begin{equation}\label{ff3}
		\begin{split}
		m_p(z)-\hat m_p(z)=&\frac{\bar\lambda_{\bSigma_{\Vb}^2}m_p(z)-\bar\lambda_{\hat\bSigma_{\Vb}^2}\hat m_p(z)}{(z+\bar\lambda_{\bSigma_{\Vb}^2}m_p(z))(z+\bar\lambda_{\hat\bSigma_{\Vb}^2}\hat m_p(z))}\\
		=&\frac{\bar\lambda_{\bSigma_{\Vb}^2}[m_p(z)-\hat m_p(z)]+[\bar\lambda_{\hat\bSigma_{\Vb}^2}-\bar\lambda_{\bSigma_{\Vb}^2}]\hat m_p(z)}{(z+\bar\lambda_{\bSigma_{\Vb}^2}m_p(z))(z+\bar\lambda_{\hat\bSigma_{\Vb}^2}\hat m_p(z))}.
		\end{split}
		\end{equation}
		For $z=\pm (2\|\bSigma_{\Vb}\|+\epsilon_0)$, 
		\[
		\begin{split}
		1-\frac{\bar\lambda_{\bSigma_{\Vb}^2}}{(z+\bar\lambda_{\bSigma_{\Vb}^2}m_p(z))(z+\bar\lambda_{\hat\bSigma_{\Vb}^2}\hat m_p(z))}=	1+\frac{\bar\lambda_{\bSigma_{\Vb}^2}\hat m_p(z)}{(z+\bar\lambda_{\bSigma_{\Vb}^2}m_p(z))}\ne 0
		\end{split}
		\]
		with probability tending to 1. Therefore,
		\[
		\max_{z\in\mathcal{C}}|	m_p(z)-\hat m_p(z)|=o_p(1).
		\]
		Similarly, we can prove that $	\max_{z\in\mathcal{C}}|	\mathcal{A}_p(z)-\hat {\mathcal{A}}_p(z)|=o_p(1).$ Therefore, by the dominated convergence theorem, 
		\[
		\int_{\mathcal{C}}f(z)[\mathcal{A}_p(z)m_p^3(z)-\hat{\mathcal{A}}_p(z)\hat m_p^3(z)]dz =o_p(1).
		\]
		
		Next, for the second part of $\mathcal{X}_p(z)$, we have
		\[
		\begin{split}
		&\frac{\mathcal{Y}_p(z)m_p^2(z)}{1+\mathcal{Y}_p(z)m_p(z)}-	\frac{\hat{\mathcal{Y}}_p(z)\hat m_p^2(z)}{1+\hat{\mathcal{Y}}_p(z)\hat m_p(z)}\\
		=&\frac{\mathcal{Y}_p(z)m_p^2(z)-\hat{\mathcal{Y}}_p(z)\hat m_p^2(z)}{[1+\mathcal{Y}_p(z)m_p(z)][1+\hat{\mathcal{Y}}_p(z)\hat m_p(z)]}+\frac{\mathcal{Y}_p(z)m_p(z)\hat{\mathcal{Y}}_p(z)\hat m_p(z)[m_p(z)-\hat m_p(z)]}{[1+\mathcal{Y}_p(z)m_p(z)][1+\hat{\mathcal{Y}}_p(z)\hat m_p(z)]}.
		\end{split}
		\]
		Similarly to (\ref{ff1}), the denominators are not equal to 0 with probability tending to 1. Moreover, by the definition of $\mathcal{Y}_p(z)$, 
		\[
		\begin{split}
		\mathcal{Y}_p(z)=&\frac{\bar\lambda_{\bSigma_{\Vb}^2}}{p}\mathcal{A}_p(z)m_p^3(z)+m_p(z)\bigg(\bar\lambda_{\bSigma_{\Vb}^2}-\frac{1}{q}\sum_{j=1}^q\frac{(\lambda_j^{\bSigma_{\Vb}})^2}{1+\lambda_j^{\bSigma_{\Vb}}\frac{1}{\sqrt{T}}\frac{m_p(z)}{1+\mathcal{Y}_p(z)m_p(z)}}\bigg)\\
		&+\frac{\mathcal{Y}_p(z)m_p^2(z)}{1+\mathcal{Y}_p(z)m_p(z)}\bigg(\frac{1}{q}\sum_{j=1}^q\frac{(\lambda_j^{\bSigma_{\Vb}})^2}{1+\lambda_j^{\bSigma_{\Vb}}\frac{1}{\sqrt{T}}\frac{m_p(z)}{1+\mathcal{Y}_p(z)m_p(z)}}\bigg)\\
		=&O(\frac{1}{\sqrt{T}})+\mathcal{Y}_p(z)m_p^2(z)\bar\lambda_{\bSigma_{\Vb}^2}[1+o(1)].
		\end{split}
		\]
		Then, similarly to (\ref{ff1}) and (\ref{ff2}), we conclude that $\max_{z\in\mathcal{C}}|\mathcal{Y}_p(z)|\le O(1/\sqrt{T})$. Similar result holds for $\hat{\mathcal{Y}}_p(z)$ with probability tending to 1. Moreover, 
		\[
		\begin{split}
		&\mathcal{Y}_p(z)-\hat{\mathcal{Y}}_p(z)\\
		=&\frac{\bar\lambda_{\bSigma_{\Vb}^2}}{p}\mathcal{A}_p(z)m_p^3(z)-\frac{\bar\lambda_{\hat\bSigma_{\Vb}^2}}{p}\hat{\mathcal{A}}_p(z)\hat m_p^3(z)\\
		&+m_p(z)\frac{1}{q}\sum_{j=1}^q\frac{\lambda_j^{\bSigma_{\Vb}}\frac{1}{\sqrt{T}}\frac{m_p(z)}{1+\mathcal{Y}_p(z)m_p(z)}}{1+\lambda_j^{\bSigma_{\Vb}}\frac{1}{\sqrt{T}}\frac{m_p(z)}{1+\mathcal{Y}_p(z)m_p(z)}}-\hat m_p(z)\frac{1}{q}\sum_{j=1}^q\frac{\lambda_j^{\hat\bSigma_{\Vb}}\frac{1}{\sqrt{T}}\frac{\hat m_p(z)}{1+\hat{\mathcal{Y}}_p(z)\hat m_p(z)}}{1+\lambda_j^{\hat\bSigma_{\Vb}}\frac{1}{\sqrt{T}}\frac{\hat m_p(z)}{1+\hat{\mathcal{Y}}_p(z)\hat m_p(z)}}\\
		&+\frac{\mathcal{Y}_p(z)m_p^2(z)}{1+\mathcal{Y}_p(z)m_p(z)}\bigg(\frac{1}{q}\sum_{j=1}^q\frac{(\lambda_j^{\bSigma_{\Vb}})^2}{1+\lambda_j^{\bSigma_{\Vb}}\frac{1}{\sqrt{T}}\frac{m_p(z)}{1+\mathcal{Y}_p(z)m_p(z)}}\bigg)\\
		&-\frac{\hat{\mathcal{Y}}_p(z)\hat m_p^2(z)}{1+\hat{\mathcal{Y}}_p(z)\hat m_p(z)}\bigg(\frac{1}{q}\sum_{j=1}^q\frac{(\lambda_j^{\hat\bSigma_{\Vb}})^2}{1+\lambda_j^{\hat\bSigma_{\Vb}}\frac{1}{\sqrt{T}}\frac{\hat m_p(z)}{1+\hat{\mathcal{Y}}_p(z)\hat m_p(z)}}\bigg).
		\end{split}
		\]
		Note that (\ref{ff3}) implies 
		\[
		\max_{z\in\mathcal{C}}|m_p(z)-\hat m_p(z)|\le \|\hat\bSigma_{\Vb}-\bSigma_{\Vb}\|=o\bigg(\frac{\sqrt{T}}{p}\bigg). 
		\]
		Then, after some tedious but elementary calculations, we conclude that
		\[
		\mathcal{Y}_p(z)-\hat{\mathcal{Y}}_p(z)=o_p\bigg(\frac{1}{p}\bigg)+[\mathcal{Y}_p(z)-\hat{\mathcal{Y}}_p(z)]\times o_p(1).
		\]
		That is, $\mathcal{Y}_p(z)-\hat{\mathcal{Y}}_p(z)=o_p(p^{-1})$. 
		As a result,
		\[
		\begin{split}
		&	\max_{z\in\mathcal{C}}|\mathcal{Y}_p(z)m_p^2(z)-\hat{\mathcal{Y}}_p(z)\hat m_p^2(z)|\\
		\le &	\max_{z\in\mathcal{C}}|[\mathcal{Y}_p(z)-\hat{\mathcal{Y}}_p(z)]m_p^2(z)|+\max_{z\in\mathcal{C}}|\hat{\mathcal{Y}}_p(z)[m_p^2(z)-\hat m_p^2(z)]|=o_p(p^{-1}).
		\end{split}
		\]
		On the other hand,
		\[
		\max_{z\in\mathcal{C}}|\mathcal{Y}_p(z)m_p(z)\hat{\mathcal{Y}}_p(z)\hat m_p(z)[m_p(z)-\hat m_p(z)]|=o_p(p^{-1}).
		\]
		Combining the above results, we conclude that
		\[
		\max_{z\in\mathcal{C}}\bigg|p\bigg\{\frac{\mathcal{Y}_p(z)m_p^2(z)}{1+\mathcal{Y}_p(z)m_p(z)}-	\frac{\hat{\mathcal{Y}}_p(z)\hat m_p^2(z)}{1+\hat{\mathcal{Y}}_p(z)\hat m_p(z)}\bigg\}\bigg|=o_p(1).
		\]
		{Then, the consistency of $p\int f(x)\hat{\mathcal{X}}(z)dz$ holds with the dominated convergence theorem, which further implies the consistency of bootstrap procedure.}
	\end{proof}

	\section{Proof of Results in Section \ref{sec:noise}: noised model}\label{sec: proof noise}
	{
		\subsection{Proof of Lemma \ref{sigma beta}: estimating $\sigma_\beta^2$}
		\begin{proof}
			By definition, 
			\begin{equation}\label{hat sigma beta}
			\hat\sigma_{\beta}^2=\mathbf{1}_p^\prime \text{Re}(\hat\Omega)\mathbf{1}_q/(pq)-\frac{\mathbf{1}_p^\prime \text{Re}(\hat\Omega)\hat v\times \hat u^\prime \text{Re}(\hat\Omega)\mathbf{1}_q}{pq\hat u^\prime \text{Re}(\hat\Omega)\hat v},
			\end{equation}
			where $\hat u$ and $\hat v$ are the leading left and right singular vectors of $(I-p^{-1}\mathbf{1}_p\mathbf{1}_p^\prime)\text{Re}(\hat\Omega)(I-q^{-1}\mathbf{1}_q\mathbf{1}_q^\prime)$ respectively, while $[\text{Re}(\hat\Omega)]_{ij}=T^{-1}\sum_{t=1}^T\hat Y_{ij}^2$. We first aim to find the limit of $\mathbf{1}_p^\prime \text{Re}(\hat\Omega)\mathbf{1}_q/(pq)$.  By definition,
			\[
			\mathbf{1}_p^\prime \text{Re}(\hat\Omega)\mathbf{1}_q/(pq)=\mathbf{1}_p^\prime (\vec{u}\vec{v}^\prime+\sigma_{\beta}^2\mathbf{1}_p\mathbf{1}_q^\prime)\mathbf{1}_q/(pq)+\mathbf{1}_p^\prime[\text{Re}(\hat\Omega)-\text{Re}(\Omega)]\mathbf{1}_q/(pq).
			\]
			Therefore, it suffices to consider the error term, or equivalently,
			\[
			\begin{split}
			\mathcal{E}_1:=&\frac{1}{pq}\sum_{i=1}^p\sum_{j=1}^q\bigg\{\frac{1}{T}\sum_{t=1}^T\bigg(\bU_{i\cdot}^\prime\Xb_t\bV_{j\cdot}+\sigma_{\beta}\phi_{t,ij}+\frac{\mathbf{1}_p^\prime(\Yb_t+\sigma_\beta\bPhi_t)\mathbf{1}_q}{pq}\bigg)^2-\bSigma_{\Ub,ii}\bSigma_{\Vb,jj}-\sigma_{\beta}^2\bigg\}\\
			=&\frac{1}{pq}\sum_{i=1}^p\sum_{j=1}^q\bigg\{\frac{1}{T}\sum_{t=1}^T\bigg(\bU_{i\cdot}^\prime\Xb_t\bV_{j\cdot}+\sigma_{\beta}\phi_{t,ij}\bigg)^2-\bSigma_{\Ub,ii}\bSigma_{\Vb,jj}-\sigma_{\beta}^2\bigg\}\\
			&+\frac{3}{T}\sum_{t=1}^T\bigg(\frac{\mathbf{1}_p^\prime(\Yb_t+\sigma_\beta\bPhi_t)\mathbf{1}_q}{pq}\bigg)^2.
			\end{split}
			\]
			One can easily verify that
			\[
			\begin{split}
			&\mathbb{E}\bigg(\frac{\mathbf{1}_p^\prime(\Yb_t+\sigma_\beta\bPhi_t)\mathbf{1}_q}{pq}\bigg)^2\le \frac{C}{pq},\quad\mathbb{E}\bigg(\frac{1}{Tpq}\sum_{t,i,j}\bU_{i\cdot}^\prime\Xb_t\bV_{j\cdot}\phi_{t,ij}\bigg)^2\le \frac{C}{Tpq}.
			\end{split}
			\]
			On the other hand,
			\[
			\begin{split}
			\mathbb{E}\bigg|\frac{1}{Tpq}\sum_{t,i,j}\{(\bU_{i\cdot}^\prime\Xb_t\bV_{j\cdot})^2-\bSigma_{\Ub,ii}\bSigma_{\Vb,jj}\}\bigg|^2=& \frac{1}{T^2p^2q^2}\sum_t\mathbb{E}\bigg|\text{tr}\Ub(\Xb_t\Vb^\prime\Vb\Xb_t^\prime-\text{tr}\bSigma_{\Vb})\Ub^\prime\bigg|^2\\
			\le &\frac{C}{Tpq}.
			\end{split}
			\]
			Therefore, we conclude that 
			\begin{equation}\label{mathcal E}
			|\mathcal{E}_1|=O_p\bigg(\frac{1}{pq}+\frac{1}{\sqrt{Tpq}}\bigg). 
			\end{equation}
			
			Next, we calculate the denominator of (\ref{hat sigma beta}). Recall that $\hat u$ and $\hat v$ are orthogonal of $\mathbf{1}_p$ and $\mathbf{1}_q$, respectively. Therefore,
			\[
			\frac{\hat u^\prime\text{Re}(\hat\Omega)\hat v}{\sqrt{pq}}=\frac{\hat u^\prime \vec{u}\hat v^\prime\vec{v}+\hat u^\prime[\text{Re}(\hat\Omega)-\text{Re}(\Omega)]\hat v}{\sqrt{pq}}.
			\]
			By similar technique in proving (\ref{mathcal E}), one can verify that
			\[
			\bigg\|\frac{1}{\sqrt{pq}}[\text{Re}(\hat\Omega)-\text{Re}(\Omega)]\bigg\|\le \bigg\|\frac{1}{\sqrt{pq}}[\text{Re}(\hat\Omega)-\text{Re}(\Omega)]\bigg\|_F\le O_p\bigg(\frac{1}{\sqrt{T}}+\frac{1}{pq}\bigg).
			\]
			Then, combining the famous Davis-Kahan's $\sin(\Theta)$ theorem, we conclude that
			\[
			\bigg\|\hat u-\frac{(\Ib-p^{-1}\mathbf{1}_p\mathbf{1}_p^\prime)\vec{u}}{\|(\Ib-p^{-1}\mathbf{1}_p\mathbf{1}_p^\prime)\vec{u}\|}\bigg\|=O_p\bigg(\frac{1}{\sqrt{T}}+\frac{1}{pq}\bigg), \quad \bigg\|\hat v-\frac{(\Ib-q^{-1}\mathbf{1}_q\mathbf{1}_q^\prime)\vec{v}}{\|(\Ib-q^{-1}\mathbf{1}_q\mathbf{1}_q^\prime)\vec{v}\|}\bigg\|=O_p\bigg(\frac{1}{\sqrt{T}}+\frac{1}{pq}\bigg).
			\]
			Therefore, 
			\[
			\begin{split}
			\frac{\hat u^\prime[\text{Re}(\hat\Omega)-\text{Re}(\Omega)]\hat v}{\sqrt{pq}}\le &\frac{1}{\sqrt{pq}}\bigg\|\hat u-\frac{(\Ib-p^{-1}\mathbf{1}_p\mathbf{1}_p^\prime)\vec{u}}{\|(\Ib-p^{-1}\mathbf{1}_p\mathbf{1}_p^\prime)\vec{u}\|}\bigg\|\|[\text{Re}(\hat\Omega)-\text{Re}(\Omega)]\hat v\|\\
			&+\frac{1}{\sqrt{pq}}\bigg\|\frac{\vec{u}^\prime(\Ib-p^{-1}\mathbf{1}_p\mathbf{1}_p^\prime)}{\|\vec{u}^\prime(\Ib-p^{-1}\mathbf{1}_p\mathbf{1}_p^\prime\|}[\text{Re}(\hat\Omega)-\text{Re}(\Omega)]\hat v\bigg\|.
			\end{split}
			\]
			On one hand, 
			\[
			\begin{split}
			\frac{1}{\sqrt{pq}}\|[\text{Re}(\hat\Omega)-\text{Re}(\Omega)]\hat v\|
			\le& \frac{1}{\sqrt{pq}}\bigg\|[\text{Re}(\hat\Omega)-\text{Re}(\Omega)]\frac{(\Ib-q^{-1}\mathbf{1}_q\mathbf{1}_q^\prime)\vec{v}}{\|(\Ib-q^{-1}\mathbf{1}_q\mathbf{1}_q^\prime)\vec{v}\|}\bigg\|\\
			&+\frac{1}{\sqrt{pq}}\|\text{Re}(\hat\Omega)-\text{Re}(\Omega)\|\bigg\|\hat v-\frac{(\Ib-q^{-1}\mathbf{1}_q\mathbf{1}_q^\prime)\vec{v}}{\|(\Ib-q^{-1}\mathbf{1}_q\mathbf{1}_q^\prime)\vec{v}\|}\bigg\|\\
			=&O_p\bigg(\frac{1}{pq}+\frac{1}{\sqrt{Tq}}\bigg)+O_p\bigg(\frac{1}{T}+\frac{1}{p^2q^2}\bigg)\\
			=&O_p\bigg(\frac{1}{pq}+\frac{1}{\sqrt{Tq}}+\frac{1}{T}\bigg),
			\end{split}
			\]
			where in the third line we use a similar technique in proving (\ref{mathcal E}). On the other hand, by parallel procedure, we can further conclude that
			\[
			\frac{1}{\sqrt{pq}}\bigg\|\frac{\vec{u}^\prime(\Ib-p^{-1}\mathbf{1}_p\mathbf{1}_p^\prime)}{\|\vec{u}^\prime(\Ib-p^{-1}\mathbf{1}_p\mathbf{1}_p^\prime\|}[\text{Re}(\hat\Omega)-\text{Re}(\Omega)]\hat v\bigg\|=O_p\bigg(\frac{1}{pq}+\frac{1}{\sqrt{Tpq}}+\frac{1}{T\sqrt{p}}\bigg).
			\]
			As a result, we conclude that
			\[
			\frac{\hat u^\prime[\text{Re}(\hat\Omega)-\text{Re}(\Omega)]\hat v}{\sqrt{pq}}=O_p\bigg(\frac{1}{pq}+\frac{1}{T\times\min\{\sqrt{T},\sqrt{p},\sqrt{q}\}}+\frac{1}{\sqrt{Tpq}}\bigg).
			\]
			
			Lastly, for the numerators of (\ref{hat sigma beta}), we write
			\[
			\begin{split}
			&\frac{1}{\sqrt{p}}\mathbf{1}_p^\prime\frac{1}{\sqrt{pq}}\text{Re}(\hat\Omega)\hat v= \frac{1}{\sqrt{p}}\mathbf{1}_p^\prime\frac{1}{\sqrt{pq}}\text{Re}(\Omega)\hat v+\frac{1}{\sqrt{p}}\mathbf{1}_p^\prime\frac{1}{\sqrt{pq}}[\text{Re}(\hat\Omega)-\text{Re}(\Omega)]\hat v\\
			=&\frac{1}{p}\mathbf{1}_p^\prime \vec{u}\times \frac{1}{\sqrt{q}}\hat v^\prime\vec{v}+O(1)\times\bigg\|\frac{1}{\sqrt{p}}\mathbf{1}_p^\prime\frac{1}{\sqrt{pq}}[\text{Re}(\hat\Omega)-\text{Re}(\Omega)]\bigg\|\bigg\|\hat v-\frac{(\Ib-q^{-1}\mathbf{1}_q\mathbf{1}_q^\prime)\vec{v}}{\|(\Ib-q^{-1}\mathbf{1}_q\mathbf{1}_q^\prime)\vec{v}\|}\bigg\|\\
			&+O(1)\times \bigg\|\frac{1}{\sqrt{p}}\mathbf{1}_p^\prime\frac{1}{\sqrt{pq}}[\text{Re}(\hat\Omega)-\text{Re}(\Omega)]\frac{(\Ib-q^{-1}\mathbf{1}_q\mathbf{1}_q^\prime)\vec{v}}{\|(\Ib-q^{-1}\mathbf{1}_q\mathbf{1}_q^\prime)\vec{v}\|}\bigg\|\\
			=&\frac{1}{p}\mathbf{1}_p^\prime \vec{u}\times \frac{1}{\sqrt{q}}\hat v^\prime\vec{v}+O_p\bigg(\frac{1}{pq}+\frac{1}{\sqrt{Tpq}}+\frac{1}{T\sqrt{p}}\bigg).
			\end{split}
			\]
			Similarly, we have
			\[
			\hat u^\prime\frac{1}{\sqrt{pq}}\text{Re}(\hat\Omega)\frac{1}{\sqrt{q}}\mathbf{1}_q=\frac{1}{\sqrt{p}}\hat u^\prime\vec{u}\times \frac{1}{q}\vec{v}^\prime\mathbf{1}_q+O_p\bigg(\frac{1}{pq}+\frac{1}{\sqrt{Tpq}}+\frac{1}{T\sqrt{q}}\bigg).
			\]
			Combining all the results and by the fact that $\hat u^\prime \vec{u}\hat v^\prime\vec{v}/\sqrt{pq}$ is of constant order, we conclude that
			\[
			\begin{split}
			\hat\sigma_{\beta}^2-\sigma_{\beta}^2=&\mathbf{1}_p^\prime (\vec{u}\vec{v}^\prime+\sigma_{\beta}^2\mathbf{1}_p\mathbf{1}_q^\prime)\mathbf{1}_q/(pq)-\frac{\frac{1}{\sqrt{p}}\hat u^\prime\vec{u}\times \frac{1}{q}\vec{v}^\prime\mathbf{1}_q\times \frac{1}{p}\mathbf{1}_p^\prime \vec{u}\times \frac{1}{\sqrt{q}}\hat v^\prime\vec{v}}{\frac{1}{\sqrt{pq}}\hat u^\prime\vec{u}\hat v^\prime\vec{v}}-\sigma_{\beta}^2\\
			&+O_p\bigg(\frac{1}{pq}+\frac{1}{T\times\min\{\sqrt{T},\sqrt{p},\sqrt{q}\}}+\frac{1}{\sqrt{Tpq}}\bigg)\\
			=&O_p\bigg(\frac{1}{pq}+\frac{1}{T\times\min\{\sqrt{T},\sqrt{p},\sqrt{q}\}}+\frac{1}{\sqrt{Tpq}}\bigg).
			\end{split}
			\]
			The lemma is then verified.
		\end{proof}\qed

		\subsection{Proof of Theorem \ref{noise clt}: asymptotic distribution under the null}
		\begin{proof}
			\noindent \textbf{Step 1: remove negligible errors.}\\
			We truncate $x_{t,ij}$ and $\phi_{t,ij}$ by $\delta_p\sqrt[4]{Tpq}$. Similarly to  Section \ref{secb}, this has minor effects on the results. By definition, 
			\[
			\begin{split}
			\bar{\mathcal{S}}=&\sqrt{\frac{1}{\bar Tpq}}\sum_{t=1}^{\bar T}\{\mathcal{Y}_{t}\mathcal{Y}_{t}^\prime-E(\mathcal{Y}_{t}\mathcal{Y}_{t}^\prime)\}.
			\end{split}
			\]
			Under the null hypothesis, without loss of generality we assume $\Ub=\bSigma_{\Ub}^{1/2}=\bSigma_{\Ub_0}^{1/2}$. Then,
			\begin{equation}\label{expansion}
			\mathcal{Y}_t\mathcal{Y}_t^\prime=\bar{\mathcal{S}}_{1,t}+\bar{\mathcal{S}}_{2,t}+\bar{\mathcal{S}}_{3,t}+\bar{\mathcal{S}}_{4,t}+\bar{\mathcal{S}}_{4,t}^\prime-\bar{\mathcal{S}}_{5,t}-\bar{\mathcal{S}}_{5,t}^\prime-\bar{\mathcal{S}}_{6,t}-\bar{\mathcal{S}}_{6,t}^\prime,
			\end{equation}
			where
			\[
			\begin{split}
			\bar{\mathcal{S}}_{1,t}=&\Xb_t\bSigma_\Vb\Xb_t^\prime,\quad 
			\bar{\mathcal{S}}_{2,t}=\sigma_{\beta}^2\Ub^{-1}\bPhi_t\bPhi_t^\prime \Ub^{-1},\quad 
			\bar{\mathcal{S}}_{3,t}=\bigg(\frac{\mathbf{1}_p^\prime(\Ub\Xb_t\Vb^\prime+\sigma_\beta \bPhi_t)\mathbf{1}_q}{pq}\bigg)^2q\Ub^{-1}\mathbf{1}_p\mathbf{1}_p^\prime \Ub^{-1},\\
			\bar{\mathcal{S}}_{4,t}=&\sigma_\beta \Xb_t\Vb^\prime\bPhi_t^\prime \Ub^{-1},\quad 
			\bar{\mathcal{S}}_{5,t}=\bigg(\frac{\mathbf{1}_p^\prime(\Ub\Xb_t\Vb^\prime+\sigma_\beta \bPhi_t)\mathbf{1}_q}{pq}\bigg)\Xb_t\Vb^\prime\mathbf{1}_q\mathbf{1}_p^\prime \Ub^{-1},\\
			\bar{\mathcal{S}}_{6,t}=&\bigg(\frac{\mathbf{1}_p^\prime(\Ub\Xb_t\Vb^\prime+\sigma_\beta \bPhi_t)\mathbf{1}_q}{pq}\bigg)\sigma_{\beta}\Ub^{-1}\bPhi_t\mathbf{1}_q\mathbf{1}_p^\prime\Ub^{-1}.
			\end{split}
			\]
			In the first step, we aim to remove some negligible errors from the expansion (\ref{expansion}). By definition and independence,
			\[
			\begin{split}
			&\mathbb{E}\bigg[\frac{1}{Tpq}\text{tr}\{\sum_{t=1}^T(\bar{\mathcal{S}}_{3,t}-\mathbb{E} \bar{\mathcal{S}}_{3,t})\}^2\bigg]\\
			=&\frac{q^2(\mathbf{1}_p^\prime\bSigma_\Ub^{-1}\mathbf{1}_p)^2}{Tpq}\sum_{t=1}^T\mathbb{E}\bigg[\bigg\{\frac{\mathbf{1}_p^\prime(\Ub\Xb_t\Vb^\prime+\sigma_\beta \bPhi_t)\mathbf{1}_q}{pq}\bigg\}^2-\mathbb{E}\bigg\{\frac{\mathbf{1}_p^\prime(\Ub\Xb_t\Vb^\prime+\sigma_\beta \bPhi_t)\mathbf{1}_q}{pq}\bigg\}^2\bigg]^2\\
			\le& \frac{C}{pq}\rightarrow 0.
			\end{split}
			\]
			On the other hand, by (\ref{mu and sigma2}), we conclude that 
			\[
			\frac{1}{Tpq}\text{tr}\{\sum_{t=1}^T(\bar{\mathcal{S}}_{1,t}-\mathbb{E} \bar{\mathcal{S}}_{1,t})\}^2=O_p (p).
			\]
			Therefore, by the Cauchy-Schwartz inequality, one can conclude that
			\[
			\frac{1}{Tpq}\text{tr}\{\sum_{t=1}^T(\bar{\mathcal{S}}_{1,t}-\mathbb{E} \bar{\mathcal{S}}_{1,t})\}^\prime\{\sum_{t=1}^T(\bar{\mathcal{S}}_{3,t}-\mathbb{E} \bar{\mathcal{S}}_{3,t})\}\le O_p\bigg(\frac{1}{\sqrt{q}}\bigg)\rightarrow 0. 
			\]
			Similar arguments will lead to 
			\[
			\begin{split}
			&\frac{1}{Tpq}\text{tr}\{\sum_{t=1}^T(\bar{\mathcal{S}}_{i,t}-\mathbb{E} \bar{\mathcal{S}}_{i,t})\}\{\sum_{t=1}^T(\bar{\mathcal{S}}_{3,t}-\mathbb{E} \bar{\mathcal{S}}_{3,t})\}=o_p(1),\\	&\frac{1}{Tpq}\text{tr}\{\sum_{t=1}^T(\bar{\mathcal{S}}_{i,t}-\mathbb{E} \bar{\mathcal{S}}_{i,t})\}^\prime\{\sum_{t=1}^T(\bar{\mathcal{S}}_{3,t}-\mathbb{E} \bar{\mathcal{S}}_{3,t})\}=o_p(1),\quad i\ne 3.
			\end{split}
			\]
			As a result, removing $\bar{\mathcal{S}}_{3,t}$ from $\bar{\mathcal{S}}_1$ will have asymptotically negligible effects on the corresponding limiting distribution. 
			
			Next, we consider the interaction term $\bar{\mathcal{S}}_{5,t}$. By independence and elementary calculation,
			\[
			\begin{split}
			&\mathbb{E}\bigg[\frac{1}{Tpq}\text{tr}\{\sum_{t=1}^T(\bar{\mathcal{S}}_{5,t}-\mathbb{E} \bar{\mathcal{S}}_{5,t})\}\{\sum_{t=1}^T(\bar{\mathcal{S}}_{5,t}-\mathbb{E} \bar{\mathcal{S}}_{5,t})\}^\prime\bigg]
			= \frac{1}{pq}\text{tr}\mathbb{E}(\bar{\mathcal{S}}_{5,1}-\mathbb{E} \bar{\mathcal{S}}_{5,1})(\bar{\mathcal{S}}_{5,1}-\mathbb{E} \bar{\mathcal{S}}_{5,1})^\prime\\
			\le &\frac{1}{pq}\text{tr}\mathbb{E}(\bar{\mathcal{S}}_{5,1}\bar{\mathcal{S}}_{5,1}^\prime)
			=\frac{1}{pq}\mathbb{E}\bigg(\frac{\mathbf{1}_p^\prime(\Ub\Xb_1\Vb^\prime+\sigma_\beta \bPhi_1)\mathbf{1}_q}{pq}\bigg)^2\mathbf{1}_q^\prime \Vb\Xb_1^\prime 	 \Xb_1\Vb^\prime\mathbf{1}_q\mathbf{1}_p^\prime \bSigma_\Ub^{-1}\mathbf{1}_p
			\le \frac{C}{q}\rightarrow 0.
			\end{split}
			\]
			By similar calculation and the Cauchy-Schwartz inequality, it's not hard to verify 
			\[
			\begin{split}
			&\frac{1}{Tpq}\text{tr}\{\sum_{t=1}^T(\bar{\mathcal{S}}_{5,t}-\mathbb{E} \bar{\mathcal{S}}_{5,t})\}^2=o_p(1),\\
			&\frac{1}{Tpq}\text{tr}\{\sum_{t=1}^T(\bar{\mathcal{S}}_{5,t}-\mathbb{E} \bar{\mathcal{S}}_{5,t})\}\{\sum_{t=1}^T(\bar{\mathcal{S}}_{6,t}-\mathbb{E} \bar{\mathcal{S}}_{6,t})\}=o_p(1),\\
			&\frac{1}{Tpq}\text{tr}\{\sum_{t=1}^T(\bar{\mathcal{S}}_{5,t}-\mathbb{E} \bar{\mathcal{S}}_{5,t})\}\{\sum_{t=1}^T(\bar{\mathcal{S}}_{6,t}-\mathbb{E} \bar{\mathcal{S}}_{6,t})\}^\prime=o_p(1).
			\end{split}
			\]
			It's more challenging to deal with the interaction between $\bar{\mathcal{S}}_{5,t}$ and $\bar{\mathcal{S}}_{i,t}$ for $i=1,2,4$. Firstly,
			\[
			\begin{split}
			&\frac{1}{Tpq}\text{tr}\{\sum_{t=1}^T(\bar{\mathcal{S}}_{5,t}-\mathbb{E} \bar{\mathcal{S}}_{5,t})\}\{\sum_{s=1}^T(\bar{\mathcal{S}}_{1,s}-\mathbb{E} \bar{\mathcal{S}}_{1,s})\}\\
			=&\frac{1}{Tpq}\text{tr}\{\sum_{t=1}^T(\bar{\mathcal{S}}_{5,t}-\mathbb{E} \bar{\mathcal{S}}_{5,t})(\bar{\mathcal{S}}_{1,t}-\mathbb{E} \bar{\mathcal{S}}_{1,t})\}+\frac{1}{Tpq}\text{tr}\{\sum_{t=1}^T(\bar{\mathcal{S}}_{5,t}-\mathbb{E} \bar{\mathcal{S}}_{5,t})\}\{\sum_{s\ne t}(\bar{\mathcal{S}}_{1,s}-\mathbb{E} \bar{\mathcal{S}}_{1,s})\}\\
			:=&\mathcal{L}_1+\mathcal{L}_2.
			\end{split}
			\]
			On one hand, for $\mathcal{L}_1$,  we have
			\[
			\begin{split}
			&\mathbb{E}\bigg\{\frac{1}{pq}\text{tr}(\bar{\mathcal{S}}_{5,t}-\mathbb{E} \bar{\mathcal{S}}_{5,t})(\bar{\mathcal{S}}_{1,t}-\mathbb{E} \bar{\mathcal{S}}_{1,t})\bigg\}^2\\
			\le& \frac{1}{p^2q^2}\mathbb{E}\text{tr}(\bar{\mathcal{S}}_{5,t}-\mathbb{E} \bar{\mathcal{S}}_{5,t})(\bar{\mathcal{S}}_{5,t}-\mathbb{E} \bar{\mathcal{S}}_{5,t})^\prime\times \text{tr}(\bar{\mathcal{S}}_{1,t}-\mathbb{E} \bar{\mathcal{S}}_{1,t})(\bar{\mathcal{S}}_{1,t}-\mathbb{E} \bar{\mathcal{S}}_{1,t})^\prime\\
			\le &\frac{1}{p^2q^2}\mathbb{E}\text{tr}(\bar{\mathcal{S}}_{1,t}^2)\times \text{tr}(\bar{\mathcal{S}}_{5,t}-\mathbb{E} \bar{\mathcal{S}}_{5,t})(\bar{\mathcal{S}}_{5,t}-\mathbb{E} \bar{\mathcal{S}}_{5 ,t})^\prime.
			\end{split}
			\]
			By random matrix theory,  there exists constant $C>0$ such that
			\begin{equation}\label{S1t}
			\mathbb{P}\bigg\{\max_t\frac{1}{p\vee q}\|\bar{\mathcal{S}}_{1,t}\|\le C\bigg\}\ge 1-(p\vee q)^{-d},
			\end{equation}
			for any $d>0$.  Consequently, we can truncate $\text{tr}(\bar{\mathcal{S}}_{1,t}^2)$ by $C(p\wedge q)(p\vee q)^2$, which further implies
			\[
			\mathbb{E}\bigg\{\frac{1}{pq}\text{tr}(\bar{\mathcal{S}}_{5,t}-\mathbb{E} \bar{\mathcal{S}}_{5,t})(\bar{\mathcal{S}}_{1,t}-\mathbb{E} \bar{\mathcal{S}}_{1,t})\bigg\}^2\le \frac{C}{pq^2}(p\wedge q)(p\vee q)^2.
			\]
			Then, by the law of large number and independence across $t$, we can conclude that
			\[
			\mathcal{L}_1\le O_p\bigg(\frac{1}{\sqrt{T}}\times \sqrt{\frac{(p\wedge q)(p\vee q)^2}{pq^2}}\bigg)=o_p(1),
			\]
			where the $o_p(1)$ is by $p/(Tq)\rightarrow 0$. On the other hand, for $\mathcal{L}_2$,
			\[
			\begin{split}
			\mathbb{E}\mathcal{L}_2^2=&\frac{1}{T^2p^2q^2}\mathbb{E}\bigg[\text{tr}\{\sum_{t=1}^T(\bar{\mathcal{S}}_{5,t}-\mathbb{E} \bar{\mathcal{S}}_{5,t})\}\{\sum_{s\ne t}(\bar{\mathcal{S}}_{1,s}-\mathbb{E} \bar{\mathcal{S}}_{1,s})\}\bigg]^2\\
			=&\frac{1}{T^2p^2q^2}\sum_{t=1}^T\sum_{s\ne t}\mathbb{E}\{\text{tr}(\bar{\mathcal{S}}_{5,t}-\mathbb{E} \bar{\mathcal{S}}_{5,t})(\bar{\mathcal{S}}_{1,s}-\mathbb{E} \bar{\mathcal{S}}_{1,s})\}^2\\
			&+\frac{1}{T^2p^2q^2}\sum_{t=1}^T\sum_{s\ne t}\mathbb{E}\{\text{tr}(\bar{\mathcal{S}}_{5,t}-\mathbb{E} \bar{\mathcal{S}}_{5,t})(\bar{\mathcal{S}}_{1,s}-\mathbb{E} \bar{\mathcal{S}}_{1,s})\}\{\text{tr}(\bar{\mathcal{S}}_{5,s}-\mathbb{E} \bar{\mathcal{S}}_{5,s})(\bar{\mathcal{S}}_{1,t}-\mathbb{E} \bar{\mathcal{S}}_{1,t})\}\\
			:=&\mathcal{L}_{21}+\mathcal{L}_{22}.	 
			\end{split}
			\]
			For $\mathcal{L}_{21}$, note that for any $s\ne t$,
			\[
			\begin{split}
			\mathbb{E}\{\text{tr}(\bar{\mathcal{S}}_{5,t}-\mathbb{E} \bar{\mathcal{S}}_{5,t})(\bar{\mathcal{S}}_{1,s}-\mathbb{E} \bar{\mathcal{S}}_{1,s})\}^2\le \mathbb{E}(\mathbf{1}_p^\prime\Ub^{-1}\Xb_t\bSigma_{\Vb}\Xb_t^\prime\Xb_t\Vb^\prime\mathbf{1}_q)^2\bigg(\frac{\mathbf{1}_p^\prime(\Ub\Xb_t\Vb^\prime+\sigma_\beta\bPhi_t)\mathbf{1}_q}{pq}\bigg)^2.
			\end{split}
			\]
			We can truncate $\Xb_t\Xb_t^\prime$ by $C(p\vee q)$ according to (\ref{S1t}). Then,
			\[
			\mathcal{L}_{11}\le \frac{C}{p^2q^2}\times (pq)^2(p\vee q)\times \frac{1}{pq}\rightarrow 0.
			\]
			On the other hand, for $\mathcal{L}_{22}$, we have
			\[
			\begin{split}
			\mathcal{L}_{22}\le &\frac{1}{T^2p^2q^2}\sum_{t=1}^T\sum_{s\ne t}\sqrt{\mathbb{E}\{\text{tr}(\bar{\mathcal{S}}_{5,t}-\mathbb{E} \bar{\mathcal{S}}_{5,t})(\bar{\mathcal{S}}_{1,s}-\mathbb{E} \bar{\mathcal{S}}_{1,s})\}^2\mathbb{E}\{\text{tr}(\bar{\mathcal{S}}_{5,s}-\mathbb{E} \bar{\mathcal{S}}_{5,s})(\bar{\mathcal{S}}_{1,t}-\mathbb{E} \bar{\mathcal{S}}_{1,t})\}^2}\\
			\le &C\mathcal{L}_{21}\rightarrow 0.
			\end{split}
			\]	
			Consequently, we conclude that
			\[
			\frac{1}{Tpq}\text{tr}\{\sum_{t=1}^T(\bar{\mathcal{S}}_{5,t}-\mathbb{E} \bar{\mathcal{S}}_{5,t})\}\{\sum_{s=1}^T(\bar{\mathcal{S}}_{1,s}-\mathbb{E} \bar{\mathcal{S}}_{1,s})\}=o_p(1).
			\]
			Further by similar calculations, one can verify that
			\[
			\begin{split}
			&\frac{1}{Tpq}\text{tr}\{\sum_{t=1}^T(\bar{\mathcal{S}}_{5,t}-\mathbb{E} \bar{\mathcal{S}}_{5,t})\}\{\sum_{s=1}^T(\bar{\mathcal{S}}_{i,s}-\mathbb{E} \bar{\mathcal{S}}_{i,s})\}=o_p(1),i=2,4,\\
			&\frac{1}{Tpq}\text{tr}\{\sum_{t=1}^T(\bar{\mathcal{S}}_{5,t}-\mathbb{E} \bar{\mathcal{S}}_{5,t})\}\{\sum_{s=1}^T(\bar{\mathcal{S}}_{4,s}-\mathbb{E} \bar{\mathcal{S}}_{4,s})\}^\prime=o_p(1).
			\end{split}
			\]
			That is, we can remove $\bar{\mathcal{S}}_{5,t}$ from the system. Similarly, we can also remove $\bar{\mathcal{S}}_{5,t}^\prime$, $\bar{\mathcal{S}}_{6,t}$ and $\bar{\mathcal{S}}_{6,t}^\prime$. Then, it remains to consider $\bar{\mathcal{S}}_{1,t}$, $\bar{\mathcal{S}}_{2,t}$ and $\bar{\mathcal{S}}_{4,t}$. 
			
			\vspace{1em}
			\noindent \textbf{Step 2: verify Lyapunov condition.}\\
			Now we remove the negligible terms and consider
			\[
			\begin{split}
			\bar{\mathcal{Q}}:=\frac{1}{Tpq}\text{tr}\bigg[\bigg\{\sum_{t=1}^T(\bar{\mathcal{S}}_{1,t}-\mathbb{E}\bar{\mathcal{S}}_{1,t}+\bar{\mathcal{S}}_{2,t}-\mathbb{E}\bar{\mathcal{S}}_{2,t}+\bar{\mathcal{S}}_{4,t}+\bar{\mathcal{S}}_{4,t}^\prime)\bigg\}\times \\
			\bigg\{\sum_{t=1}^T(\bar{\mathcal{S}}_{1,t}-\mathbb{E}\bar{\mathcal{S}}_{1,t}+\bar{\mathcal{S}}_{2,t}-\mathbb{E}\bar{\mathcal{S}}_{2,t}+\bar{\mathcal{S}}_{4,t}+\bar{\mathcal{S}}_{4,t}^\prime)\bigg\}^\prime\bigg].
			\end{split}
			\]
			Let  $\mathcal{F}_t$ be the $\sigma$-field by $\{\check\Yb_1,\ldots,\check\Yb_t\}$ and $\mathbb{E}_t(\cdot)=\mathbb{E}(\cdot\mid \mathcal{F}_t)$. Then,
			\[
			\bar{\mathcal{Q}}-\mathbb{E}\bar{\mathcal{Q}}=\sum_{t=1}^T(\mathbb{E}_t-\mathbb{E}_{t-1})\bar{\mathcal{Q}}:=\sum_{t=1}^T\bar{\mathcal{Q}}_t.
			\]
			We first focus on $\bar{\mathcal{S}}_{1,t}$. By independence across $t$, we have
			\begin{equation}\label{S1t first}
			\begin{split}
			&(\mathbb{E}_t-\mathbb{E}_{t-1})\frac{1}{Tpq}\text{tr}\bigg[\bar{\mathcal{S}}_{1,t}-\mathbb{E}\bar{\mathcal{S}}_{1,t}+\sum_{s\ne t}(\bar{\mathcal{S}}_{1,s}-\mathbb{E}\bar{\mathcal{S}}_{1,s})\bigg]^2\\
			=&\frac{1}{Tpq}\bigg\{\text{tr}(\bar{\mathcal{S}}_{1,t}-\mathbb{E}\bar{\mathcal{S}}_{1,t})^2-\mathbb{E}\text{tr}(\bar{\mathcal{S}}_{1,t}-\mathbb{E}\bar{\mathcal{S}}_{1,t})^2\bigg\}+\frac{2}{Tpq}\text{tr}(\bar{\mathcal{S}}_{1,t}-\mathbb{E}\bar{\mathcal{S}}_{1,t})\bigg\{\sum_{s<t}(\bar{\mathcal{S}}_{1,s}-\mathbb{E}\bar{\mathcal{S}}_{1,s})\bigg\}.
			\end{split}
			\end{equation}
			Now write $\bx_{i\cdot}=(\bx_{1,i\cdot}^\prime,\ldots,\bx_{T,i\cdot}^\prime)$ where $\bx_{t,i\cdot}^\prime$ is the $i$-th row of $\Xb_t$. Define $\mathbb{E}^i(\cdot)=\mathbb{E}(\cdot\mid \bx_{1\cdot},\ldots,\bx_{i\cdot})$. Then, by Burkholder's inequality,
			\[
			\begin{split}
			&\mathbb{E}\bigg|\frac{1}{Tpq}\bigg\{\text{tr}(\bar{\mathcal{S}}_{1,t}-\mathbb{E}\bar{\mathcal{S}}_{1,t})^2-\mathbb{E}\text{tr}(\bar{\mathcal{S}}_{1,t}-\mathbb{E}\bar{\mathcal{S}}_{1,t})^2\bigg\}\bigg|^4\\
			=&\mathbb{E}\bigg|\sum_{i=1}^p(\mathbb{E}^i-\mathbb{E}^{i-1})\frac{1}{Tpq}\text{tr}(\bar{\mathcal{S}}_{1,t}-\mathbb{E}\bar{\mathcal{S}}_{1,t})^2\bigg|^4
			\le C\bigg[\sum_{i=1}^p\mathbb{E}\bigg|(\mathbb{E}^i-\mathbb{E}^{i-1})\frac{1}{Tpq}\text{tr}(\bar{\mathcal{S}}_{1,t}-\mathbb{E}\bar{\mathcal{S}}_{1,t})^2\bigg|^2\bigg]^2.
			\end{split}
			\]
			For simplicity, write $\underline{\bx}_{t,i\cdot}=\bx_{t,i\cdot}^\prime\Vb^\prime$. Then,
			\[
			\begin{split}
			&	(\mathbb{E}^i-\mathbb{E}^{i-1})\frac{1}{Tpq}\text{tr}(\bar{\mathcal{S}}_{1,t}-\mathbb{E}\bar{\mathcal{S}}_{1,t})^2\\
			=&	(\mathbb{E}^i-\mathbb{E}^{i-1})\frac{1}{Tpq}\bigg[\sum_{j=1}^p(\|\underline{\bx}_{t,j\cdot}\|^2-q)^2+\sum_{j\ne k}(\underline{\bx}_{t,j\cdot}^\prime\underline{\bx}_{t,k\cdot})^2\bigg]\\
			=&(\mathbb{E}^i-\mathbb{E}^{i-1})\frac{1}{Tpq}\bigg[(\|\underline{\bx}_{t,i\cdot}\|^2-q)^2+2\underline{\bx}_{t,i\cdot}^\prime(\sum_{k\ne i}\underline{\bx}_{t,k\cdot}\underline{\bx}_{t,k\cdot}^\prime)\underline{\bx}_{t,i\cdot}\bigg]\\
			=&\frac{1}{Tpq}\bigg[(\|\underline{\bx}_{t,i\cdot}\|^2-q)^2-\mathbb{E}(\|\underline{\bx}_{t,i\cdot}\|^2-q)^2\bigg]+\frac{2}{Tpq}\bx_{t,i\cdot}^\prime (\sum_{k>i}\bSigma_{\Vb}^2-\text{tr}\sum_{k>i}\bSigma_{\Vb}^2)\bx_{t,i\cdot}\\
			&+\frac{2}{Tpq}\bx_{t,i\cdot}^\prime(\sum_{k< i}\Vb^\prime\underline{\bx}_{t,k\cdot}\underline{\bx}_{t,k\cdot}^\prime\Vb-\text{tr}\sum_{k< i}\Vb^\prime\underline{\bx}_{t,k\cdot}\underline{\bx}_{t,k\cdot}^\prime\Vb)\bx_{t,i\cdot}.
			\end{split}
			\]
			By large deviation bounds and the truncation $|x_{t,ij}|\le \delta_p\sqrt[4]{Tpq}$ for some $\delta_p\rightarrow 0$, we have
			\[
			\begin{split}
			&\mathbb{E}\bigg|(\mathbb{E}^i-\mathbb{E}^{i-1})\frac{1}{Tpq}\text{tr}(\bar{\mathcal{S}}_{1,t}-\mathbb{E}\bar{\mathcal{S}}_{1,t})^2\bigg|^2\\
			\le &\frac{C}{T^2p^2q^2}\mathbb{E}\{\bx_{t,i\cdot}^\prime(\bSigma_{\Vb}-\text{tr}\bSigma_{\Vb})\bx_{t,i\cdot}\}^4+\frac{C}{T^2p^2q^2}\bigg\|\sum_{k>i}\bSigma_{\Vb}^2\bigg\|_F^2+\frac{C}{T^2p^2q^2}\mathbb{E}\bigg\|\sum_{k<i}\bx_{t,k\cdot}\bx_{t,k\cdot}^\prime\bigg\|_F^2\\
			\le &\frac{C}{T^2p^2q^2}\mathbb{E}\bigg(\sum_{j=1}^q\bSigma_{\Vb,jj}(x_{t,ij}^2-1)+\sum_{j_1\ne j_2}\bSigma_{\Vb,j_1j_2}x_{t,ij_1}x_{t,ij_2}\bigg)^4+\frac{C}{T^2p^2q^2}\{(p-i)q+iq^2\}\\
			\le &\frac{C}{T^2p^2q^2}\bigg(q\times\delta_p^4Tpq+\|\bSigma_{\Vb}\|_F^4+(p-i)q+iq^2\bigg).
			\end{split}
			\]
			As a result, we conclude that
			\begin{equation}\label{squared}
			\sum_{i=1}^p\mathbb{E}\bigg|(\mathbb{E}^i-\mathbb{E}^{i-1})\frac{1}{Tpq}\text{tr}(\bar{\mathcal{S}}_{1,t}-\mathbb{E}\bar{\mathcal{S}}_{1,t})^2\bigg|^2\le \frac{\delta_p}{T}+\frac{1}{T^2}.
			\end{equation}
			Then, for the first term of (\ref{S1t first}) we conclude that
			\[
			\sum_{t=1}^T\mathbb{E}\bigg|\frac{1}{Tpq}\bigg\{\text{tr}(\bar{\mathcal{S}}_{1,t}-\mathbb{E}\bar{\mathcal{S}}_{1,t})^2-\mathbb{E}\text{tr}(\bar{\mathcal{S}}_{1,t}-\mathbb{E}\bar{\mathcal{S}}_{1,t})^2\bigg\}\bigg|^4=o(1).
			\]
			
			Now we consider the second term of (\ref{S1t first}). Similarly,
			\[
			\begin{split}
			&\mathbb{E}\bigg|\frac{1}{Tpq}\text{tr}(\bar{\mathcal{S}}_{1,t}-\mathbb{E}\bar{\mathcal{S}}_{1,t})\bigg\{\sum_{s<t}(\bar{\mathcal{S}}_{1,s}-\mathbb{E}\bar{\mathcal{S}}_{1,s})\bigg\}\bigg|^4\\
			\le &C\bigg[\sum_{i=1}^p\mathbb{E}\bigg|(\mathbb{E}^i-\mathbb{E}^{i-1})\frac{1}{Tpq}\text{tr}(\bar{\mathcal{S}}_{1,t}-\mathbb{E}\bar{\mathcal{S}}_{1,t})\bigg\{\sum_{s<t}(\bar{\mathcal{S}}_{1,s}-\mathbb{E}\bar{\mathcal{S}}_{1,s})\bigg\}\bigg|^2\bigg],
			\end{split}
			\]
			while by some elementary calculations,
			\[
			\begin{split}
			&(\mathbb{E}^i-\mathbb{E}^{i-1})\frac{1}{Tpq}\text{tr}(\bar{\mathcal{S}}_{1,t}-\mathbb{E}\bar{\mathcal{S}}_{1,t})\bigg\{\sum_{s<t}(\bar{\mathcal{S}}_{1,s}-\mathbb{E}\bar{\mathcal{S}}_{1,s})\bigg\}\\
			=&(\mathbb{E}^i-\mathbb{E}^{i-1})\frac{1}{Tpq}\sum_{s<t}\bigg[\sum_{j=1}^p(\|\underline{\bx}_{t,j\cdot}\|^2-q)(\|\underline{\bx}_{s,j\cdot}\|^2-q)+\sum_{k\ne j}\underline{\bx}_{t,j\cdot}^\prime\underline{\bx}_{t,k\cdot}\underline{\bx}_{s,j\cdot}^\prime\underline{\bx}_{s,k\cdot}\bigg]\\
			=&(\mathbb{E}^i-\mathbb{E}^{i-1})\frac{1}{Tpq}\sum_{s<t}\bigg[(\|\underline{\bx}_{t,i\cdot}\|^2-q)(\|\underline{\bx}_{s,i\cdot}\|^2-q)+2\sum_{k\ne i}\underline{\bx}_{t,i\cdot}^\prime\underline{\bx}_{t,k\cdot}\underline{\bx}_{s,i\cdot}^\prime\underline{\bx}_{s,k\cdot}\bigg]\\
			=&\frac{1}{Tpq}\sum_{s<t}\bigg[(\|\underline{\bx}_{t,i\cdot}\|^2-q)(\|\underline{\bx}_{s,i\cdot}\|^2-q)+2\sum_{k< i}\underline{\bx}_{t,i\cdot}^\prime\underline{\bx}_{t,k\cdot}\underline{\bx}_{s,i\cdot}^\prime\underline{\bx}_{s,k\cdot}\bigg].
			\end{split}
			\]
			Therefore,
			\[
			\begin{split}
			&\mathbb{E}\bigg|(\mathbb{E}^i-\mathbb{E}^{i-1})\frac{1}{Tpq}\text{tr}(\bar{\mathcal{S}}_{1,t}-\mathbb{E}\bar{\mathcal{S}}_{1,t})\bigg\{\sum_{s<t}(\bar{\mathcal{S}}_{1,s}-\mathbb{E}\bar{\mathcal{S}}_{1,s})\bigg\}\bigg|^2
			\le \frac{Ct}{T^2p^2q^2}\bigg(q^2+i\times q^2\bigg),
			\end{split}
			\]
			which further implies that
			\[
			\sum_{i=1}^p\mathbb{E}\bigg|(\mathbb{E}^i-\mathbb{E}^{i-1})\frac{1}{Tpq}\text{tr}(\bar{\mathcal{S}}_{1,t}-\mathbb{E}\bar{\mathcal{S}}_{1,t})\bigg\{\sum_{s<t}(\bar{\mathcal{S}}_{1,s}-\mathbb{E}\bar{\mathcal{S}}_{1,s})\bigg\}\bigg|^2\le\frac{Ct}{T^2},
			\]
			and 
			\[
			\sum_{t=1}^T\mathbb{E}\bigg|\frac{1}{Tpq}\text{tr}(\bar{\mathcal{S}}_{1,t}-\mathbb{E}\bar{\mathcal{S}}_{1,t})\bigg\{\sum_{s<t}(\bar{\mathcal{S}}_{1,s}-\mathbb{E}\bar{\mathcal{S}}_{1,s})\bigg\}\bigg|^4=o(1).
			\]
			Consequently, we can conclude that
			\begin{equation}\label{Lyapunov}
			\sum_{t=1}^T\mathbb{E}\bigg|(\mathbb{E}_t-\mathbb{E}_{t-1})\frac{1}{Tpq}\text{tr}\bigg[\bar{\mathcal{S}}_{1,t}-\mathbb{E}\bar{\mathcal{S}}_{1,t}+\sum_{s\ne t}(\bar{\mathcal{S}}_{1,s}-\mathbb{E}\bar{\mathcal{S}}_{1,s})\bigg]^2\bigg|^4=o(1).
			\end{equation}
			
			Indeed, by similar arguments, one can further verify that
			\begin{equation}\label{Lyapunov2}
			\begin{split}
			&\sum_{t=1}^T\mathbb{E}\bigg|(\mathbb{E}_t-\mathbb{E}_{t-1})\frac{1}{Tpq}\text{tr}\bigg[\sum_{s=1}^T(\bar{\mathcal{S}}_{2,s}-\mathbb{E}\bar{\mathcal{S}}_{2,s})\bigg]^2\bigg|^4=o(1),\\
			&	\sum_{t=1}^T\mathbb{E}\bigg|(\mathbb{E}_t-\mathbb{E}_{t-1})\frac{1}{Tpq}\text{tr}\bigg[\sum_{s=1}^T(\bar{\mathcal{S}}_{4,s}+\bar{\mathcal{S}}_{4,s}^\prime)\bigg]^2\bigg|^4=o(1).
			\end{split}
			\end{equation}
			With (\ref{Lyapunov}) and (\ref{Lyapunov2}), we have verified the Lyapunov condition. Therefore, it remains to calculate the limits of expectation and variance of $\text{tr}(\bar{\mathcal{S}}^2)$, respectively.
			
			\vspace{1em}
			\noindent\textbf{Step 3: calculate asymptotic mean and variance.}\\
			We start with mean.  Recall that  we have removed the negligible errors so it suffices to consider $\mathbb{E}\bar{\mathcal{Q}}$. 
			By independence between $\Xb_t$ and $\bPhi_t$,  we have
			\[
			\mathbb{E}(\bar{\mathcal{S}}_{i,t}-\mathbb{E}\bar{\mathcal{S}}_{i,t})\bar{\mathcal{S}}_{4,t}=0,i=1,2.
			\]
			Further by the independence across $t$, 
			\[
			\begin{split}
			\mathbb{E}\bar{\mathcal{Q}}=&\frac{1}{Tpq}\text{tr}\sum_{t=1}^T\mathbb{E}\bigg[(\bar{\mathcal{S}}_{1,t}-\mathbb{E}\bar{\mathcal{S}}_{1,t})^2+(\bar{\mathcal{S}}_{2,t}-\mathbb{E}\bar{\mathcal{S}}_{2,t})^2+(\bar{\mathcal{S}}_{4,t}+\bar{\mathcal{S}}_{4,t}^\prime)^2\bigg].
			\end{split}
			\]
			Firstly, for any $t$,
			\[
			\begin{split}
			\text{tr}\mathbb{E}(\bar{\mathcal{S}}_{1,t}-\mathbb{E}\bar{\mathcal{S}}_{1,t})^2=&\sum_{i=1}^p\mathbb{E}(\|\bx_{t,i\cdot}^\prime\bSigma_{\Vb}\bx_{t,i\cdot}-\text{tr}\bSigma_{\Vb})^2+\sum_{j\ne k}\mathbb{E}(\bx_{t,j\cdot}^\prime\bSigma_{\Vb}\bx_{t,k\cdot})^2\\
			=&p\bigg\{(\nu_4-3)\sum_{j=1}^q\bSigma_{\Vb,jj}^2+2\|\bSigma_{\Vb}\|_F^2\bigg\}+p(p-1)\|\bSigma_{\Vb}\|_F^2.
			\end{split}
			\]
			For $\bar{\mathcal{S}}_{2,t}$, we have
			\[
			\begin{split}
			&\sigma_{\beta}^{-4}\text{tr}\mathbb{E}(\bar{\mathcal{S}}_{2,t}-\mathbb{E}\bar{\mathcal{S}}_{2,t})^2\\
			=&\mathbb{E}\|\Ub^{-1}(\bPhi_t\bPhi_t^\prime-q \Ib)\Ub^{-1}\|_F^2=\mathbb{E}\text{tr}(\Ub^{-1}\bPhi_t\bPhi_t^\prime\Ub^{-2}\bPhi_t\bPhi_t^\prime\Ub^{-1})-q^2\text{tr}\bSigma_{\Ub}^{-2}\\
			=&\mathbb{E}\|\bPhi_t^\prime\Ub^{-2}\bPhi_t\|_F^2-q^2\text{tr}\bSigma_{\Ub}^{-2}\\
			=&q\bigg\{(\tilde\nu_4-3)\sum_{j=1}^p(\bSigma_{\Ub}^{-1})_{jj}^2+2\|\bSigma_{\Ub}^{-1}\|_F^2+(\text{tr}\bSigma_{\Ub}^{-1})^2\bigg\}+q(q-1)\|\bSigma_{\Ub}^{-1}\|_F^2-q^2\|\bSigma_{\Ub}^{-1}\|_F^2\\
			=&q\bigg\{(\tilde\nu_4-3)\sum_{j=1}^p(\bSigma_{\Ub}^{-1})_{jj}^2+\|\bSigma_{\Ub}^{-1}\|_F^2+(\text{tr}\bSigma_{\Ub}^{-1})^2\bigg\}.
			\end{split}
			\]
			For $\bar{\mathcal{S}}_{4,t}$, 
			\[
			\sigma_{\beta}^{-2}\text{tr}\mathbb{E}(\bar{\mathcal{S}}_{4,t}\bar{\mathcal{S}}_{4,t}^\prime)=p\text{tr}\bSigma_{\Vb}\text{tr}(\bSigma_{\Ub}^{-1}),
			\]
			while 
			\[
			\begin{split}
			\text{tr}\mathbb{E}(\bar{\mathcal{S}}_{4,t}^2)=&\text{tr}\mathbb{E}\bigg[\mathbb{E}(\bar{\mathcal{S}}_{4,t}^2\mid \bPhi_t)\bigg]=\mathbb{E}\text{tr}\Ub^{-2}\bPhi_t\Vb\Vb^\prime\bPhi_t^\prime=\text{tr}\bSigma_{\Vb}\text{tr}(\bSigma_{\Ub}^{-1}).
			\end{split}
			\]
			As a result,
			\[
			\begin{split}
			\mathbb{E}\bar{\mathcal{Q}}=&\frac{\nu_4-3}{q}\sum_{j=1}^q\bSigma_{\Vb,jj}^2+\frac{p+1}{q}\|\bSigma_{\Vb}\|_F^2\\
			&+\sigma_{\beta}^4\bigg(\frac{\tilde\nu_4-3}{p}\sum_{j=1}^p(\bSigma_{\Ub}^{-1})_{jj}^2+\frac{1}{p}\|\bSigma_{\Ub}^{-1}\|_F^2+\frac{1}{p}(\text{tr}\bSigma_{\Ub}^{-1})^2\bigg)+\frac{2(p+1)\sigma_{\beta}^2}{p}\text{tr}(\bSigma_{\Ub}^{-1}).
			\end{split}
			\]

			Next, we calculate the variance, i.e., $\text{Cov}(\bar{\mathcal{Q}})$.  We write
			\[
			\begin{split}
			\bar{\mathcal{Q}}_1=&\frac{1}{Tpq}\text{tr}\bigg\{\sum_{t=1}^T(\bar{\mathcal{S}}_{1,t}-\mathbb{E}\bar{\mathcal{S}}_{1,t})\bigg\}^2,\quad \bar{\mathcal{Q}}_2=\frac{1}{Tpq}\text{tr}\bigg\{\sum_{t=1}^T(\bar{\mathcal{S}}_{2,t}-\mathbb{E}\bar{\mathcal{S}}_{2,t})\bigg\}^2,\\
			\bar{\mathcal{Q}}_3=&\frac{1}{Tpq}\text{tr}\bigg\{\sum_{t=1}^T(\bar{\mathcal{S}}_{4,t}+\bar{\mathcal{S}}_{4,t}^\prime)\bigg\}^2,\quad
			\bar{\mathcal{Q}}_4=\frac{2}{Tpq}\text{tr}\bigg\{\sum_{t=1}^T(\bar{\mathcal{S}}_{1,t}-\mathbb{E}\bar{\mathcal{S}}_{1,t})\bigg\}\bigg\{\sum_{t=1}^T(\bar{\mathcal{S}}_{2,t}-\mathbb{E}\bar{\mathcal{S}}_{2,t})\bigg\},\\
			\bar{\mathcal{Q}}_5=&\frac{2}{Tpq}\text{tr}\bigg\{\sum_{t=1}^T(\bar{\mathcal{S}}_{1,t}-\mathbb{E}\bar{\mathcal{S}}_{1,t})\bigg\}\bigg\{\sum_{t=1}^T(\bar{\mathcal{S}}_{4,t}+\bar{\mathcal{S}}_{4,t}^\prime)\bigg\},\\
			\bar{\mathcal{Q}}_6=&\frac{2}{Tpq}\text{tr}\bigg\{\sum_{t=1}^T(\bar{\mathcal{S}}_{2,t}-\mathbb{E}\bar{\mathcal{S}}_{2,t})\bigg\}\bigg\{\sum_{t=1}^T(\bar{\mathcal{S}}_{4,t}+\bar{\mathcal{S}}_{4,t}^\prime)\bigg\}.
			\end{split}
			\]
			Corollary \ref{cor2} has already shown that 
			\[
			\text{Cov}\bar{\mathcal{Q}}_1=4\bar\lambda^2_{\bSigma_{\Vb}^2}+o(1).
			\]
			For $\bar{\mathcal{Q}}_2$, we write
			\[
			\bar{\mathcal{Q}}_2=\frac{1}{Tpq}\bigg\{\sum_{t=1}^T\text{tr}(\bar{\mathcal{S}}_{2,t}-\mathbb{E}\bar{\mathcal{S}}_{2,t})^2+\sum_{t,s\ne t}\text{tr}(\bar{\mathcal{S}}_{2,t}-\mathbb{E}\bar{\mathcal{S}}_{2,t})(\bar{\mathcal{S}}_{2,s}-\mathbb{E}\bar{\mathcal{S}}_{2,s})\bigg\}:=\bar{\mathcal{Q}}_{21}+\bar{\mathcal{Q}}_{22}.
			\]
			Then, $\text{Cov}(\bar{\mathcal{Q}}_{21},\bar{\mathcal{Q}}_{22})=0$.  Furthermore, using similar technique in (\ref{squared}), we have
			\[
			\begin{split}
			\text{Cov}(\bar{\mathcal{Q}}_{21})=&\frac{1}{Tp^2q^2}\text{Cov}\{\text{tr}(\bar{\mathcal{S}}_{2,t}-\mathbb{E}\bar{\mathcal{S}}_{2,t})^2\}\rightarrow 0,
			\end{split}
			\]
			while
			\[
			\begin{split}
			&\text{Cov}(\bar{\mathcal{Q}}_{22})=\frac{2}{T^2p^2q^2}\sum_{t,s\ne t}\mathbb{E}\text{tr}^2(\bar{\mathcal{S}}_{2,t}-\mathbb{E}\bar{\mathcal{S}}_{2,t})(\bar{\mathcal{S}}_{2,s}-\mathbb{E}\bar{\mathcal{S}}_{2,s})\\
			=&\frac{2\sigma_{\beta}^8}{T^2p^2q^2}\sum_{t,s\ne t}\mathbb{E}\bigg[\text{tr}^2\bigg\{\Ub^{-1}(\bPhi_t\bPhi_t^\prime-q\Ib)\Ub^{-2}(\bPhi_s\bPhi_s^\prime-q\Ib)\Ub^{-1}\bigg\}\bigg]\\
			=&\frac{2\sigma_{\beta}^8}{T^2p^2q^2}\sum_{t,s\ne t}\mathbb{E}\text{tr}\bigg\{\Ub^{-1}\bPhi_t\bPhi_t^\prime\bSigma_{\Ub}^{-1}\bPhi_s\bPhi_s^\prime\Ub^{-1}\bigg\}\times\text{tr}\bigg\{\Ub^{-1}(\bPhi_t\bPhi_t^\prime-q\Ib)\bSigma_{\Ub}^{-1}(\bPhi_s\bPhi_s^\prime-q\Ib)\Ub^{-1}\bigg\}.
			\end{split}
			\]
			Given $s\ne t$, by elementary calculations,
			\begin{equation}\label{Q221}
			\begin{split}
			&\mathbb{E}\text{tr}^2\{\bPhi_s^\prime\bSigma_{\Ub}^{-1}\bPhi_t\bPhi_t^\prime\bSigma_{\Ub}^{-1}\bPhi_s\}= \mathbb{E}\{\sum_{j=1}^q\bPhi_{s,\cdot j}^\prime\bSigma_{\Ub}^{-1}\bPhi_t\bPhi_t^\prime\bSigma_{\Ub}^{-1}\bPhi_{s,\cdot j}\}^2\\
			=&\mathbb{E}\{\sum_{j=1}^q(\bPhi_{s,\cdot j}^\prime\bSigma_{\Ub}^{-1}\bPhi_t\bPhi_t^\prime\bSigma_{\Ub}^{-1}\bPhi_{s,\cdot j}-\text{tr}(\bSigma_{\Ub}^{-1}\bPhi_t\bPhi_t^\prime\bSigma_{\Ub}^{-1})\}^2+q^2\mathbb{E}\text{tr}^2(\bSigma_{\Ub}^{-1}\bPhi_t\bPhi_t^\prime\bSigma_{\Ub}^{-1})\\
			=&q\mathbb{E}(\bPhi_{s,\cdot 1}^\prime\bSigma_{\Ub}^{-1}\bPhi_t\bPhi_t^\prime\bSigma_{\Ub}^{-1}\bPhi_{s,\cdot 1}-\text{tr}\bSigma_{\Ub}^{-1}\bPhi_t\bPhi_t^\prime\bSigma_{\Ub}^{-1})\}^2+q^2\mathbb{E}\text{tr}^2(\bSigma_{\Ub}^{-1}\bPhi_t\bPhi_t^\prime\bSigma_{\Ub}^{-1}).
			\end{split}
			\end{equation}
			Note that for symmetric matrix $\Ab$ independent of $\bPhi_s$, we have
			\begin{equation}\label{quatratic}
			\begin{split}
			\mathbb{E}(\bPhi_{s,\cdot 1}^\prime\Ab\bPhi_{s,\cdot 1}-\text{tr}\Ab)^2=(\tilde\nu_4-3)\sum_{j=1}^p\Ab_{jj}^2+2\text{tr}\Ab^2.
			\end{split}
			\end{equation}
			Then, (\ref{Q221}) can be written as 
			\[
			q(\tilde\nu_4-3)\sum_{j=1}^p\mathbb{E}(\bSigma_{\Ub}^{-1}\bPhi_t\bPhi_t^\prime\bSigma_{\Ub}^{-1})_{jj}^2+2q\text{tr}(\bSigma_{\Ub}^{-1}\bPhi_t\bPhi_t^\prime\bSigma_{\Ub}^{-1})^2+q^2\mathbb{E}\text{tr}^2(\bSigma_{\Ub}^{-1}\bPhi_t\bPhi_t^\prime\bSigma_{\Ub}^{-1}).
			\]
			Further,
			\[
			\begin{split}
			\mathbb{E}\text{tr}\bigg\{\Ub^{-1}\bPhi_t\bPhi_t^\prime\bSigma_{\Ub}^{-1}\bPhi_s\bPhi_s^\prime\Ub^{-1}\bigg\}\times \text{tr}\bigg\{q\Ub^{-1}\bPhi_t\bPhi_t^\prime\bSigma_{\Ub}^{-1}\Ub^{-1}\bigg\}=q^2\mathbb{E} \text{tr}^2\bigg\{\bSigma_{\Ub}^{-1}\bPhi_t\bPhi_t^\prime\bSigma_{\Ub}^{-1}\bigg\}.
			\end{split}
			\]
			Therefore, 
			\[
			\begin{split}
			\sigma_{\beta}^{-8}\text{Cov}(\bar{\mathcal{Q}}_{22})
			=&\frac{2(T-1)}{Tp^2q}\bigg[ (\tilde\nu_4-3)\sum_{j=1}^p\mathbb{E}(\bSigma_{\Ub}^{-1}\bPhi_t\bPhi_t^\prime\bSigma_{\Ub}^{-1})_{jj}^2+2\text{tr}(\bSigma_{\Ub}^{-1}\bPhi_t\bPhi_t^\prime\bSigma_{\Ub}^{-1})^2\\
			&-q\mathbb{E}\text{tr}^2(\bPhi_t^\prime\bSigma_{\Ub}^{-2}\bPhi_t)+q^3\text{tr}^2(\bSigma_{\Ub}^{-2})\bigg]\\
			=&\frac{2(T-1)}{Tp^2q}\bigg[ (\tilde\nu_4-3)\sum_{j=1}^p\mathbb{E}(\bSigma_{\Ub}^{-1}\bPhi_t\bPhi_t^\prime\bSigma_{\Ub}^{-1})_{jj}^2+2\text{tr}(\bSigma_{\Ub}^{-1}\bPhi_t\bPhi_t^\prime\bSigma_{\Ub}^{-1})^2\\
			&-q^2\{(\tilde\nu_4-3)\sum_{j=1}^p(\bSigma_{\Ub}^{-2})_{jj}^2+2\text{tr}(\bSigma_{\Ub}^{-4})\}\bigg],
			\end{split}
			\]
			where the last line is by (\ref{quatratic}). Moreover, for deterministic vector $\ba$, 
			\[
			\begin{split}
			\mathbb{E}(\ba^\prime\bPhi_t\bPhi_t^\prime\ba)^2=&\mathbb{E}\{\sum_{i=1}^q(\ba^\prime\bPhi_{t,\cdot j})^2\}^2=\mathbb{E}[\sum_{i=1}^q\{(\ba^\prime\bPhi_{t,\cdot j})^2-\mathbb{E}(\ba^\prime\bPhi_{t,\cdot j})^2\}]^2+q^2[\mathbb{E}(\ba^\prime\bPhi_{t,\cdot j})^2]^2\\
			=&q\mathbb{E}\{(\ba^\prime\bPhi_{t,\cdot j})^2-\mathbb{E}(\ba^\prime\bPhi_{t,\cdot j})^2\}^2+q^2[\|\ba\|^2]^2\\
			=&q[(\tilde\nu_4-3)\sum_{j=1}^p(\ba\ba^\prime)_{jj}^2+2\text{tr}(\ba\ba^\prime)^2]+q^2\|\ba\|^4.
			\end{split}
			\]
			Then,
			\[
			\sum_{j=1}^p\mathbb{E}(\bSigma_{\Ub}^{-1}\bPhi_t\bPhi_t^\prime\bSigma_{\Ub}^{-1})_{jj}^2=O(pq)+q^2\sum_{j=1}^p(\bSigma_{\Ub}^{-2})_{jj}^2.
			\]
			On the other hand,
			\[
			\mathbb{E}\text{tr}(\bSigma_{\Ub}^{-1}\bPhi_t\bPhi_t^\prime\bSigma_{\Ub}^{-1})^2=O(pq)+q\text{tr}^2(\bSigma_{\Ub}^{-2})+q^2\text{tr}(\bSigma_{\Ub}^{-4}).
			\]
			Consequently,
			\[
			\text{Cov}(\bar{\mathcal{Q}}_2)=\text{Cov}(\bar{\mathcal{Q}}_{22})+o(1)=4\sigma_{\beta}^8\bar\lambda_{\bSigma_{\Ub}^{-2}}^2+o(1).
			\] 
			
			Next, we move to $\bar{\mathcal{Q}}_3$. Similarly to $\bar{\mathcal{Q}}_2$, we write
			\[
			\bar{\mathcal{Q}}_3= \frac{1}{Tpq}\text{tr}\{\sum_{t=1}^T(\bar{\mathcal{S}}_{4,t}+\bar{\mathcal{S}}_{4,t}^\prime)^2+\sum_{t,s\ne t}(\bar{\mathcal{S}}_{4,t}+\bar{\mathcal{S}}_{4,t}^\prime)(\bar{\mathcal{S}}_{4,s}+\bar{\mathcal{S}}_{4,s}^\prime)\}:=\bar{\mathcal{Q}}_{31}+\bar{\mathcal{Q}}_{32},
			\]
			while $\text{Cov}(\bar{\mathcal{Q}}_{31})=o(1)$ and $\text{Cov}(\bar{\mathcal{Q}}_{31},\bar{\mathcal{Q}}_{32})=0$. Moreover,
			\[
			\begin{split}
			\text{Cov}(\bar{\mathcal{Q}}_{32})=&\frac{2}{T^2p^2q^2}\sum_{t,s\ne t}\mathbb{E}\text{tr}^2(\bar{\mathcal{S}}_{4,t}+\bar{\mathcal{S}}_{4,t}^\prime)(\bar{\mathcal{S}}_{4,s}+\bar{\mathcal{S}}_{4,s}^\prime)=\frac{8}{T^2p^2q^2}\sum_{t,s\ne t}\mathbb{E}\text{tr}^2(\bar{\mathcal{S}}_{4,t}\bar{\mathcal{S}}_{4,s}+\bar{\mathcal{S}}_{4,t}\bar{\mathcal{S}}_{4,s}^\prime).
			\end{split}
			\]
			Given any $t$ and $s\ne t$, we have
			\[
			\begin{split}
			\sigma_{\beta}^{-4}\mathbb{E}\text{tr}^2(\bar{\mathcal{S}}_{4,t}\bar{\mathcal{S}}_{4,s})=&\mathbb{E}\text{tr}^2(\Xb_t\Vb^\prime\bPhi_t^\prime\Ub^{-1}\Xb_s\Vb^\prime\bPhi_s^\prime\Ub^{-1})=\mathbb{E}\|\Vb^\prime\bPhi_t^\prime\Ub^{-1}\Xb_s\Vb^\prime\bPhi_s^\prime\Ub^{-1}\|_F^2 \\
			=&\mathbb{E}\text{tr}(\Vb\Xb_s^\prime\Ub^{-1}\bPhi_t\Vb\Vb^\prime\bPhi_t^\prime\Ub^{-1}\Xb_s\Vb^\prime\bPhi_s^\prime\Ub^{-2}\bPhi_s)\\
			=&\mathbb{E}\text{tr}(\bSigma_{\Vb})\text{tr}(\bSigma_{\Ub}^{-1})\mathbb{E}\text{tr}(\Vb\Xb_s^\prime\bSigma_{\Ub}^{-1}\Xb_s\Vb^\prime)=\text{tr}^2(\bSigma_{\Vb})\text{tr}^2(\bSigma_{\Ub}^{-1}).
			\end{split}
			\]
			On the other hand,
			\[
			\begin{split}
			\sigma_{\beta}^{-4}\mathbb{E}\text{tr}^2(\bar{\mathcal{S}}_{4,t}\bar{\mathcal{S}}_{4,s}^\prime)=&\mathbb{E}\text{tr}^2(\Xb_t\Vb^\prime\bPhi_t^\prime\bSigma_{\Ub}^{-1}\bPhi_s\Vb\Xb_s^\prime)=\mathbb{E}\|\Vb^\prime\bPhi_t^\prime\bSigma_{\Ub}^{-1}\bPhi_s\Vb\Xb_s^\prime\|_F^2=p \text{tr}(\bSigma_{\Ub}^{-2})\text{tr}^2(\bSigma_{\Vb}),
			\end{split}
			\]
			while
			\[
			\begin{split}
			\sigma_{\beta}^{-4}\mathbb{E}\text{tr}(\bar{\mathcal{S}}_{4,t}\bar{\mathcal{S}}_{4,s})\text{tr}(\bar{\mathcal{S}}_{4,t}\bar{\mathcal{S}}_{4,s}^\prime)=&\mathbb{E}\text{tr}(\Xb_t\Vb^\prime\bPhi_t^\prime\Ub^{-1}\Xb_s\Vb^\prime\bPhi_s^\prime\Ub^{-1})\text{tr}(\Xb_t\Vb^\prime\bPhi_t^\prime\bSigma_{\Ub}^{-1}\bPhi_s\Vb\Xb_s^\prime)\\
			=&\text{tr}^2(\bSigma_{\Vb})\text{tr}(\bSigma_{\Ub}^{-2})=O(pq^2).
			\end{split}
			\]
			Therefore,
			\[
			\text{Cov}(\bar{\mathcal{Q}}_3)=\frac{8\sigma_{\beta}^4}{p^2}\bigg[\text{tr}^2(\bSigma_{\Ub}^{-1})+p\text{tr}(\bSigma_{\Ub}^{-2})\bigg]+o(1). 
			\]
			
			Next, for $\bar{\mathcal{Q}}_4$, 
			\[
			\begin{split}
			\text{Cov}(\bar{\mathcal{Q}}_4)=&\frac{4}{T^2p^2q^2}\sum_{t,s}\mathbb{E}\text{tr}^2[(\bar{\mathcal{S}}_{1,t}-\mathbb{E}\bar{\mathcal{S}}_{1,t})(\bar{\mathcal{S}}_{2,s}-\mathbb{E}\bar{\mathcal{S}}_{2,s})].
			\end{split}
			\]
			Given any $t,s$,
			\[
			\begin{split}
			& \sigma_{\beta}^{-4}\mathbb{E}\text{tr}^2[(\bar{\mathcal{S}}_{1,t}-\mathbb{E}\bar{\mathcal{S}}_{1,t})(\bar{\mathcal{S}}_{2,s}-\mathbb{E}\bar{\mathcal{S}}_{2,s})]= \mathbb{E}\text{tr}^2[(\Xb_t\bSigma_{\Vb}\Xb_t^\prime-q\Ib)\Ub^{-1}(\bPhi_s\bPhi_s^\prime-q\Ib)\Ub^{-1}]\\
			=&\mathbb{E}\text{tr}^2[\bPhi_s^\prime\Ub^{-1}(\Xb_t\bSigma_{\Vb}\Xb_t^\prime-q\Ib)\Ub^{-1}\bPhi_s-q\Ub^{-1}(\Xb_t\bSigma_{\Vb}\Xb_t^\prime-q\Ib)\Ub^{-1}]\\
			=&q(\tilde\nu_4-3)\sum_{j=1}^p\mathbb{E}[\Ub^{-1}(\Xb_t\bSigma_{\Vb}\Xb_t^\prime-q\Ib)\Ub^{-1}]_{jj}^2+2q\text{tr}\mathbb{E}[\Ub^{-1}(\Xb_t\bSigma_{\Vb}\Xb_t^\prime-q\Ib)\Ub^{-1}]^2\\
			=&o(p^2q^2)+2q\text{tr}^2(\bSigma_{\Ub}^{-1})\text{tr}(\bSigma_{\Vb}^2),
			\end{split}
			\]
			where the last line is by some tedious but elementary calculations. Therefore,
			\[
			\text{Cov}(\bar{\mathcal{Q}}_4)=\frac{8\sigma_{\beta}^4}{p^2q}\text{tr}^2(\bSigma_{\Ub}^{-1})\text{tr}(\bSigma_{\Vb}^2)+o(1).
			\]
			
			For $\bar{\mathcal{Q}}_5$, note that
			\[
			\begin{split}
			\sigma_{\beta}^{-2}\text{Cov}(\bar{\mathcal{Q}}_5)=&\frac{4}{T^2p^2q^2}\sum_{t,s}\mathbb{E}\text{tr}^2[(\Xb_t\bSigma_{\Vb}\Xb_t^\prime-\text{tr}\bSigma_{\Vb})(\Xb_s\Vb^\prime\bPhi_s^\prime\Ub^{-1}+\Ub^{-1}\bPhi_s\Vb\Xb_s^\prime)]\\
			=&\frac{16}{T^2p^2q^2}\sum_{t,s}\mathbb{E}\text{tr}^2[(\Xb_t\bSigma_{\Vb}\Xb_t^\prime-\text{tr}\bSigma_{\Vb})(\Xb_s\Vb^\prime\bPhi_s^\prime\Ub^{-1})].
			\end{split}
			\]
			Given any $t\ne s$,
			\[
			\begin{split}
			&\mathbb{E}\text{tr}^2[(\Xb_t\bSigma_{\Vb}\Xb_t^\prime-\text{tr}\bSigma_{\Vb})\Xb_s\Vb^\prime\bPhi_s^\prime\Ub^{-1}]=\mathbb{E}\text{tr}^2[\bPhi_s^\prime\Ub^{-1}(\Xb_t\bSigma_{\Vb}\Xb_t^\prime-\text{tr}\bSigma_{\Vb})\Xb_s\Vb^\prime]\\
			=&\mathbb{E}\|\Ub^{-1}(\Xb_t\bSigma_{\Vb}\Xb_t^\prime-\text{tr}\bSigma_{\Vb})\Xb_s\Vb^\prime\|_F^2=\text{tr}(\bSigma_{\Vb})\mathbb{E}\|\Ub^{-1}(\Xb_t\bSigma_{\Vb}\Xb_t^\prime-\text{tr}\bSigma_{\Vb})\|_F^2\\
			=&O(pq^2)+pq\text{tr}(\bSigma_{\Ub}^{-1})\text{tr}(\bSigma_{\Vb}^2),
			\end{split}
			\]
			while the case $t=s$ is asymptotically negligible. Therefore, 
			\[
			\text{Cov}(\bar{\mathcal{Q}}_5)=\frac{16\sigma_{\beta}^2}{pq}\text{tr}(\bSigma_{\Ub}^{-1})\text{tr}(\bSigma_{\Vb}^2)+o(1).
			\]
			Similarly, for $\bar{\mathcal{Q}}_6$, we conclude that
			\[
			\text{Cov}(\bar{\mathcal{Q}}_6)=\frac{16\sigma_\beta^6}{p^2}\text{tr}(\bSigma_{\Ub}^{-1})\text{tr}(\bSigma_{\Ub}^{-2})+o(1).
			\]
			Furthermore, one can verify 
			\[
			\text{Cov}(\bar{\mathcal{Q}}_1,\bar{\mathcal{Q}}_i)=0,\quad i=2,4;\quad \text{Cov}(\bar{\mathcal{Q}}_1,\bar{\mathcal{Q}}_j)\rightarrow 0, \quad i=3,5,6,
			\]
			and further $\text{Cov}(\bar{\mathcal{Q}}_i,\bar{\mathcal{Q}}_j)\rightarrow 0$ for any $i\ne j$. Therefore,
			\[
			\begin{split}
			&\text{Cov}(\text{tr}\bar{\mathcal{S}}^2)
			=4\bar\lambda^2_{\bSigma_{\Vb}^2}+4\sigma_{\beta}^8\bar \lambda^2_{\bSigma_{\Ub}^{-2}}+8\sigma_{\beta}^4(\bar\lambda^2_{\bSigma_{\Ub}^{-1}}+\bar\lambda_{\bSigma_{\Ub}^{-2}}+\bar\lambda^2_{\bSigma_{\Ub}^{-1}}\bar\lambda_{\bSigma_{\Vb}^2})+16\sigma_{\beta}^2(\bar\lambda_{\bSigma_{\Ub}^{-1}}\bar\lambda_{\bSigma_{\Vb}^2}+\sigma_{\beta}^4\bar\lambda^2_{\bSigma_{\Ub}^{-1}}).
			\end{split}
			\]
			The asymptotic distribution under the null hypothesis then follows.
		\end{proof}

		\subsection{Proof of Theorem \ref{noise power}: asymptotic power}
		\begin{proof}
			Under the alternative hypothesis, $\bSigma_{\Ub_0}$ is not necessarily equal to $\bSigma_{\Ub}$. Then,
			\[
			\bar{\mathcal{S}}=\{( Tq)/p\}^{1/2}\{(Tq)^{-1}\sum_{t=1}^{\bar T}\mathcal{Y}_{t}\mathcal{Y}_{t}^\prime-{\mathcal{E}}_0\}=\frac{1}{\sqrt{Tpq}}\sum_{t=1}^T(\mathcal{Y}_t\mathcal{Y}_t^\prime-q\tilde{\mathcal{E}}_0+q\tilde{\mathcal{E}}_0-q\mathcal{E}_0),
			\]
			where $\tilde{\mathcal{Y}}_t$ and $\tilde{\mathcal{E}}_0$ are given by
			\[
			\begin{split}
			{\mathcal{Y}}_t=&\bSigma_{\Ub_0}^{-1/2}\Ub\Xb_t\Vb^\prime+\sigma_\beta\bSigma_{\Ub_0}^{-1/2}\bPhi_t,\\
			\tilde{\mathcal{E}}_0=&\bSigma_{\Ub_0}^{-1/2}\bSigma_{\Ub}\bSigma_{\Ub_0}^{-1/2}-\sigma_\beta^2\bSigma_{\Ub_0}^{-1}+\frac{p^{-1}\mathbf{1}_p\prime\bSigma_{\Ub}\mathbf{1}_p+\sigma_\beta^2}{p}\bSigma_{\Ub_0}^{-1/2}\mathbf{1}_p\mathbf{1}_p^\prime\bSigma_{\Ub_0}^{-1/2}\\
			&-\frac{2\mathbf{1}_q^\prime\bSigma_{\Vb}\mathbf{1}_q}{pq}\bSigma_{\Ub_0}^{-1/2}\bSigma_{\Ub}\mathbf{1}_p\mathbf{1}_p^\prime\bSigma_{\Ub_0}^{-1/2}.
			\end{split}
			\]
			In the following, we aim to prove that $\text{tr}\bar{\mathcal{S}}^2-\bar\mu\rightarrow \infty$, where $\bar\mu$ is given in Theorem \ref{noise clt}. 
			
			Define $\bar{\mathcal{S}}_A$ as 
			\[
			\bar{\mathcal{S}}_A=\frac{1}{\sqrt{Tpq}}\sum_{t=1}^T(\mathcal{Y}_t\mathcal{Y}_t^\prime-q\tilde{\mathcal{E}}_0).
			\]
			Then, following the same technique in proving Theorem \ref{noise clt}, one can verify that the expectation of $\text{tr}\bar{\mathcal{S}}_A^2$ diverges with rate $O(p)$ while its variance is bounded by some constant.  That is to say, $\text{tr}\bar{\mathcal{S}}_A^2=O_p(p)$. On the other hand, 
			\[
			\begin{split}
			\text{tr}\bigg(\frac{Tq}{\sqrt{Tpq}}(\tilde{\mathcal{E}}_0-\mathcal{E}_0)\bigg)^2=&\frac{Tq}{p}\text{tr}\bigg(\bSigma_{\Ub_0}^{-1/2}\bSigma_{\Ub}\bSigma_{\Ub_0}^{-1/2}-\Ib+\frac{\bSigma_{\Ub_0}^{-1/2}\mathbf{1}_p\mathbf{1}_p^\prime\bSigma_{\Ub_0}^{-1/2}}{p}\frac{\mathbf{1}_p^\prime(\bSigma_{\Ub}-\bSigma_{\Ub_0})\mathbf{1}_p}{p}\\
			&-\frac{2\mathbf{1}_q^\prime\bSigma_{\Vb}\mathbf{1}_q}{pq}\bSigma_{\Ub_0}^{-1/2}(\bSigma_{\Ub}-\bSigma_{\Ub_0})\mathbf{1}_p\mathbf{1}_p^\prime\bSigma_{\Ub_0}^{-1/2}\bigg)^2.
			\end{split}
			\]
			By elementary calculations, 
			\[
			\text{tr}\bigg(\frac{Tq}{\sqrt{Tpq}}(\tilde{\mathcal{E}}_0-\mathcal{E}_0)\bigg)^2=\frac{Tq}{p}\text{tr}\bigg(\bSigma_{\Ub_0}^{-1/2}\bSigma_{\Ub}\bSigma_{\Ub_0}^{-1/2}-\Ib\bigg)^2+O_p\bigg(\frac{Tq}{p}\bigg)\ge cTq[1+o_p(1)].
			\]
			Therefore, it suffices to calculate the interaction term, given by
			\[
			\text{tr}\bigg(\bar{\mathcal{S}}_A\times\frac{Tq}{\sqrt{Tpq}}(\tilde{\mathcal{E}}_0-\mathcal{E}_0)\bigg).
			\]
			Expanding $\bar{\mathcal{S}}_A$ and $(\tilde{\mathcal{E}}_0-\mathcal{E}_0)$, here we only show to bound
			\[
			\mathcal{L}_1:=\text{tr}\bigg(\frac{1}{p}\sum_{t=1}^T(\bSigma_{\Ub_0}^{-1/2}\Ub\Xb_t\Vb^\prime\Vb\Xb_t^\prime\Ub^\prime\bSigma_{\Ub_0}^{-1/2}-q\bSigma_{\Ub_0}^{-1/2}\bSigma_{\Ub}\bSigma_{\Ub_0}^{-1/2})\times(\bSigma_{\Ub_0}^{-1/2}\bSigma_{\Ub}\bSigma_{\Ub_0}^{-1/2}-\Ib)\bigg).
			\]
			By the independence across $t$ and the Cauchy-Schwartz inequality,
			\[
			\begin{split}
			\mathbb{E}\mathcal{L}_1^2=&\frac{T}{p^2}\bigg[\mathbb{E}\text{tr}^2\bSigma_{\Ub_0}^{-1/2}\Ub(\Xb_t\Vb^\prime\Vb\Xb_t^\prime-q\Ib)\Ub^\prime\bSigma_{\Ub_0}^{-1/2}\bigg]\text{tr}^2(\bSigma_{\Ub_0}^{-1/2}\bSigma_{\Ub}\bSigma_{\Ub_0}^{-1/2}-\Ib)\\
			\le &CTq\times \frac{1}{p}\text{tr}^2(\bSigma_{\Ub_0}^{-1/2}\bSigma_{\Ub}\bSigma_{\Ub_0}^{-1/2}-\Ib).
			\end{split}
			\]
			As a result,
			\[
			\mathcal{L}_1=O_p\bigg(\sqrt{Tq\times \frac{1}{p}\text{tr}^2(\bSigma_{\Ub_0}^{-1/2}\bSigma_{\Ub}\bSigma_{\Ub_0}^{-1/2}-\Ib)}\bigg)=o_p\bigg(Tq\times \frac{1}{p}\text{tr}^2(\bSigma_{\Ub_0}^{-1/2}\bSigma_{\Ub}\bSigma_{\Ub_0}^{-1/2}-\Ib)\bigg).
			\]
			The other interaction terms after expanding can be handled similarly, which leads to 
			\[
			\text{tr}\bigg(\bar{\mathcal{S}}_A\times\frac{Tq}{\sqrt{Tpq}}(\tilde{\mathcal{E}}_0-\mathcal{E}_0)\bigg)=o_p\bigg(Tq\times \frac{1}{p}\text{tr}^2(\bSigma_{\Ub_0}^{-1/2}\bSigma_{\Ub}\bSigma_{\Ub_0}^{-1/2}-\Ib)\bigg).
			\]
			Note that $(Tq)\gg p$. Then, under the alternative hypothesis,
			\[
			\text{tr}\bar{\mathcal{S}}^2-\bar\mu\ge cTq[1+o(1)]-O(p)\rightarrow \infty,
			\]
			with probability tending to one, which concludes the theorem.
		\end{proof}
		
		\subsection{Proof of lemmas in Section \ref{sec: estimate}: estimating unknown parameters}\label{sec: proof estimate}
		We first prove Lemma \ref{lem: plug-in 1}. 
		\begin{proof}
			We start with the estimation error of $q^{-1}\mathbf{1}_q^\prime\bSigma_{\Vb}\mathbf{1}_q$. Under the null hypothesis,
			\[
			\begin{split}
			&\bigg(\ref{plug-in 1}-\frac{1}{q}\mathbf{1}_q^\prime\bSigma_{\Vb}\mathbf{1}_q\bigg)\times \bigg(1+\frac{\mathbf{1}_p^\prime\bSigma_{\Ub_0}^{-1}\mathbf{1}_p}{p^2}\frac{\mathbf{1}_p^\prime\bSigma_{\Ub_0}\mathbf{1}_p}{p}-\frac{2}{p}\bigg)\\
			=&\frac{1}{Tpq}\sum_t\mathbf{1}_q^\prime\Vb(\Xb_t^\prime\Xb_t-p\Ib)\Vb^\prime\mathbf{1}_q+\frac{1}{Tpq}\sum_t\sigma_{\beta}^2\mathbf{1}_q^\prime(\bPhi_t^\prime\bSigma_{\Ub_0}^{-1}\bPhi_t-\text{tr}\bSigma_{\Ub_0}^{-1}\Ib)\mathbf{1}_q\\
			&+O_p(|\sigma_\beta^2-\hat\sigma_\beta^2|)+\frac{2\sigma_\beta}{Tpq}\sum_t\mathbf{1}_q^\prime\Vb\Xb_t^\prime\bSigma_{\Ub_0}^{-1/2}\bPhi_t\mathbf{1}_q\\
			&+\frac{q\mathbf{1}_p^\prime\bSigma_{\Ub_0}^{-1}\mathbf{1}_p}{Tp}\sum_t\bigg\{\bigg(\frac{\mathbf{1}_p^\prime(\Yb_t+\sigma_\beta\bPhi_t)\mathbf{1}_q}{pq}\bigg)^2-\mathbb{E}\bigg(\frac{\mathbf{1}_p^\prime(\Yb_t+\sigma_\beta\bPhi_t)\mathbf{1}_q}{pq}\bigg)^2\bigg\}\\
			&-\frac{2q}{Tpq}\sum_t\mathbf{1}_q^\prime(\Vb\Xb_t^\prime\Ub^\prime+\sigma_\beta\bPhi_t^\prime)\bSigma_{\Ub_0}^{-1}\mathbf{1}_p\times \frac{\mathbf{1}_p^\prime(\Yb_t+\sigma_\beta\bPhi_t)\mathbf{1}_q}{pq}\\
			&+\frac{2q}{Tpq}\sum_t\mathbb{E}\bigg\{\mathbf{1}_q^\prime(\Vb\Xb_t^\prime\Ub^\prime+\sigma_\beta\bPhi_t^\prime)\bSigma_{\Ub_0}^{-1}\mathbf{1}_p\times \frac{\mathbf{1}_p^\prime(\Yb_t+\sigma_\beta\bPhi_t)\mathbf{1}_q}{pq}\bigg\}\\
			=&O_p\bigg(\frac{1}{\sqrt{Tp}}+|\sigma_\beta^2-\hat\sigma_\beta^2|\bigg).
			\end{split}
			\]
			Therefore, 
			\begin{equation}\label{difference E0}
			\begin{split}
			\hat{\mathcal{E}}_0-\mathcal{E}_0
			=&(\hat\sigma_\beta^2-\sigma_\beta^2)\bigg(\bSigma_{\Ub_0}^{-1}-\frac{1}{p}\bSigma_{\Ub_0}^{-1/2}\mathbf{1}_p\mathbf{1}_p^\prime\bSigma_{\Ub_0}^{-1/2}\bigg)\\
			&+O_p\bigg(\frac{1}{\sqrt{Tp}}+|\hat\sigma_\beta^2-\sigma_\beta^2|\bigg)\times \frac{1}{p}\bSigma_{\Ub_0}^{-1/2}\mathbf{1}_p\mathbf{1}_p^\prime\bSigma_{\Ub_0}^{-1/2}.
			\end{split}
			\end{equation}
			Recall that $\hat{\mathcal{S}}=\bar{\mathcal{S}}+\sqrt{(Tq)/p}(\hat{\mathcal{E}}_0-\mathcal{E}_0)$. On one hand,
			\[
			\frac{Tq}{p}\text{tr}(	\hat{\mathcal{E}}_0-\mathcal{E}_0)^2\le(Tq)(\hat\sigma_\beta^2-\sigma_\beta^2)^2+\frac{Tq}{p}\times O_p\bigg(\frac{1}{Tp}\bigg)=o_p(1).
			\]
			On the other hand,
			\[
			\text{tr}\bigg\{\bar{\mathcal{S}}\times \sqrt{\frac{Tq}{p}}(	\hat{\mathcal{E}}_0-\mathcal{E}_0)\bigg\}=\text{tr}\frac{1}{p}\sum_{t=1}^T(\mathcal{Y}_t\mathcal{Y}_t^\prime-\mathbb{E}\mathcal{Y}_t\mathcal{Y}_t^\prime)(\hat{\mathcal{E}}_0-\mathcal{E}_0).
			\]
			Then, expand $\mathcal{Y}_t\mathcal{Y}_t^\prime-\mathbb{E}\mathcal{Y}_t\mathcal{Y}_t^\prime$ and replace $(\hat{\mathcal{E}}_0-\mathcal{E}_0)$ with (\ref{difference E0}). We only show the calculation of the first term, i.e.,
			\[
			\begin{split}
			\frac{\hat\sigma_\beta^2-\sigma_\beta^2}{p}\text{tr}\sum_{t=1}^T(\Xb_t\Vb^\prime\Vb\Xb_t^\prime-q\text{tr}\bSigma_{\Vb}\Ib)\bSigma_{\Ub_0}^{-1}=O_p\bigg(|\hat\sigma_\beta^2-\sigma_\beta^2|\times \frac{\sqrt{Tpq}}{p}\bigg)=o_p(1).
			\end{split}
			\]
			The remaining terms can be proved similarly  to be negligible, which concludes the lemma.
		\end{proof}
		
		The next step is to prove Lemma \ref{tilde mu}.
		\begin{proof}
			Similarly to the proof of Theorem \ref{noise clt},  there exists constant $C>0$ such that
			\[
			\text{Cov}\bigg(\frac{Tp}{q}	\bigg\|\frac{1}{Tp}\sum_{t=1}^T(\mathcal{Y}_t^\prime\mathcal{Y}_t-\mathbb{E}\mathcal{Y}_t^\prime\mathcal{Y}_t)\bigg\|_F^2\bigg)\le C.
			\]
			Therefore, we conclude that
			\begin{equation}\label{estimate tilde mu}
			\frac{p}{q}	\bigg\|\frac{1}{Tp}\sum_{t=1}^T(\mathcal{Y}_t^\prime\mathcal{Y}_t-\mathbb{E}\mathcal{Y}_t^\prime\mathcal{Y}_t)\bigg\|_F^2=\mathbb{E}\bigg(\frac{p}{q}	\bigg\|\frac{1}{Tp}\sum_{t=1}^T(\mathcal{Y}_t^\prime\mathcal{Y}_t-\mathbb{E}\mathcal{Y}_t^\prime\mathcal{Y}_t)\bigg\|_F^2\bigg)+o_p(1).
			\end{equation}
			We first deal with the left hand side, which can be written as 
			\[
			\frac{p}{q}\bigg\|\frac{1}{Tp}\sum_t\mathcal{Y}_t^\prime\mathcal{Y}_t\bigg\|_F^2-\frac{p}{q}\bigg\|\frac{1}{p}\mathbb{E}\mathcal{Y}_1^\prime\mathcal{Y}_1\bigg\|_F^2-\frac{2p}{q}\frac{1}{Tp^2}\text{tr}\sum_t(\mathcal{Y}_t^\prime\mathcal{Y}_t-\mathbb{E}\mathcal{Y}_t^\prime\mathcal{Y}_t)\mathbb{E}\mathcal{Y}_1^\prime\mathcal{Y}_1.
			\]
			Recall that 
			\[
			\mathcal{Y}_t^\prime=\Vb\Xb_t^\prime+\sigma_\beta\bPhi_t^\prime\bSigma_{\Ub_0}^{-1/2}-\frac{\mathbf{1}_p^\prime(\Yb_t+\sigma_\beta\bPhi_t)\mathbf{1}_q}{pq}\mathbf{1}_q\mathbf{1}_p^\prime\bSigma_{\Ub_0}^{-1/2}.
			\]
			Then, elementary calculations lead to 
			\begin{equation}\label{left}
			\begin{split}
			\frac{p}{q}\bigg\|\frac{1}{p}\mathbb{E}\mathcal{Y}_1^\prime\mathcal{Y}_1\bigg\|_F^2=&\frac{p}{q}\bigg\|\bSigma_{\Vb}+\frac{\sigma_\beta^2\text{tr}\bSigma_{\Ub_0}^{-1}}{p}\Ib+\frac{\mathbf{1}_p^\prime\bSigma_{\Ub_0}^{-1}\mathbf{1}_p\mathbf{1}_p^\prime\bSigma_{\Ub_0}\mathbf{1}_p\mathbf{1}_q^\prime\bSigma_{\Vb}\mathbf{1}_q}{p^3q^2}\mathbf{1}_q\mathbf{1}_q^\prime\\
			&-\frac{1}{pq}(\bSigma_{\Vb}\mathbf{1}_q\mathbf{1}_q^\prime+\mathbf{1}_q\mathbf{1}_q^\prime\bSigma_{\Vb})-\frac{\sigma_\beta^2\mathbf{1}_p^\prime\bSigma_{\Ub_0}^{-1}\mathbf{1}_p}{p^2q}\mathbf{1}_q\mathbf{1}_q^\prime\bigg\|_F^2\\
			=&p\bar\lambda_{\bSigma_{\Vb}^2}+p\sigma_\beta^4\bar\lambda^2_{\bSigma_{\Ub_0}^{-1}}+2p\sigma_\beta^2\bar\lambda_{\bSigma_{\Ub_0}^{-1}}+o(1).
			\end{split}
			\end{equation}
			On the other hand, for the interaction term, we have
			\[
			\frac{p}{q}\frac{1}{Tp^2}\text{tr}\sum_t(\mathcal{Y}_t^\prime\mathcal{Y}_t-\mathbb{E}\mathcal{Y}_t^\prime\mathcal{Y}_t)\mathbb{E}\mathcal{Y}_1^\prime\mathcal{Y}_1=O_p\bigg(\frac{p}{\sqrt{Tpq}}\bigg)\rightarrow 0.
			\]
			
			Next, we calculate the right hand side of (\ref{estimate tilde mu}). Indeed, under the null hypothesis, this is totally parallel to the calculation of $\tilde\mu$. Therefore, we easily conclude from Theorem \ref{noise clt} that 
			\begin{equation}\label{right}
			\mathbb{E}\bigg(\frac{p}{q}	\bigg\|\frac{1}{Tp}\sum_{t=1}^T(\mathcal{Y}_t^\prime\mathcal{Y}_t-\mathbb{E}\mathcal{Y}_t^\prime\mathcal{Y}_t)\bigg\|_F^2\bigg)=\frac{q}{T}\bigg(\sigma_\beta^4\bar\lambda_{\bSigma_{\Ub_0}^{-2}}+1+2\sigma_\beta^2\bar\lambda_{\bSigma_{\Ub_0}^{-1}}\bigg)+o(1).
			\end{equation}
			The lemma follows from (\ref{estimate tilde mu}), (\ref{left}) and (\ref{right}).
		\end{proof}
	}
\end{appendices}

\bibliographystyle{model2-names}
\bibliography{Ref}

\end{document}